\documentclass[11pt]{amsart}

\usepackage[top=27 mm, bottom=20 mm, left=25 mm, right=25mm]{geometry}

\usepackage[all]{xy}
\usepackage{mathrsfs}
\usepackage{amsmath,amsthm, amssymb, amsgen,amsxtra,amsfonts, amsbsy} 
\usepackage[all]{xy}
\usepackage{verbatim}
\usepackage{arydshln}

\usepackage{graphicx}
\usepackage{caption,rotating}
\usepackage{color}
\usepackage{subcaption}

\usepackage{tikz}
\usetikzlibrary{trees}
\usetikzlibrary{decorations.pathmorphing}
\usetikzlibrary{decorations.markings}

\usepackage{times}

\usepackage[font=small,labelfont=small]{caption}

\usepackage{enumitem}

\usepackage{mathtools,stackengine}

\usepackage{float}
\restylefloat{table}

\theoremstyle{definition}
\newtheorem{Definition}{Definition}[section]

\newtheorem{Remark}[Definition]{Remark}
\newtheorem{Remark/Notation}[Definition]{Remark/Notation}
\newtheorem{Notation}[Definition]{Notation}
\newtheorem{Procedure}[Definition]{Procedure}

\newtheorem*{Acknowledgements}{Acknowledgements}

\theoremstyle{plain}

\newtheorem*{Main Theoremx}{Main Theorem}

\newtheorem{Lemma}[Definition]{Lemma}
\newtheorem{Corollary}[Definition]{Corollary}

\newtheoremstyle{voiditstyle}{3pt}{3pt}{\itshape}{\parindent}%
{\bfseries}{.}{ }{\thmnote{#3}}%
\theoremstyle{voiditstyle}

\newtheoremstyle{voidromstyle}{3pt}{3pt}{\rm}{\parindent}%
{\bfseries}{.}{ }{\thmnote{#3}}%
\theoremstyle{voidromstyle}

\newcommand{\prf}{\par\noindent{\sc Proof.}\quad}

\newcommand{\cal}{\mathcal}

\numberwithin{equation}{section}

\newcommand{\bbP}{{\mathbb{P}}}
\newcommand{\bbG}{{\mathbb{G}}}

\newcommand{\bbF}{{\mathbb{F}}}

\newcommand{\GL}{\operatorname{GL}}

\newcommand{\PGL}{\operatorname{PGL}}
\newcommand{\Aut}{\operatorname{Aut}}

\newcommand{\Spec}{\operatorname{Spec}}

\newcommand{\id}{\operatorname{id}}
\newcommand{\Pic}{\operatorname{Pic}}
\newcommand{\red}{\operatorname{red}}

\newcommand{\Char}{\operatorname{char}}
\newcommand{\Stab}{\operatorname{Stab}}

\newcommand{\bsm}{\left(\begin{smallmatrix}}
\newcommand{\esm}{\end{smallmatrix}\right)}

\newcommand{\calC}{\mathcal{C}}

\newcommand{\beq}{\begin{equation}}
\newcommand{\eeq}{\end{equation}}
\usepackage[
            pdftex,
           colorlinks=true,
            urlcolor=blue,       
            filecolor=blue,     
            linkcolor=blue,       
            citecolor=blue,         
            pdftitle= {Title},
            pdfauthor={Author},
            pdfsubject={Subject},
            pdfkeywords={Key}
            pagebackref,
            bookmarksopen=true]
            {hyperref}
\usepackage{hyperref}

\begin{document}

\title[Weak del Pezzo surfaces with global vector fields]{Weak del Pezzo surfaces with global vector fields}

\author{Gebhard Martin}
\address{Mathematisches Institut der Universit\"at Bonn, Endenicher Allee 60, 53115 Bonn, Germany}
\curraddr{}
\email{gmartin@math.uni-bonn.de}

\author{Claudia Stadlmayr}
\address{TU M\"unchen, Zentrum Mathematik - M11, Boltzmannstra{\ss}e 3, 85748 Garching bei M\"unchen, Germany}
\curraddr{}
\email{claudia.stadlmayr@ma.tum.de}

\date{\today}

\maketitle

\begin{abstract}
We classify smooth weak del Pezzo surfaces with global vector fields over an arbitrary algebraically closed field $k$ of arbitrary characteristic $p \geq 0$. We give a complete description of the configuration of $(-1)$- and $(-2)$-curves on these surfaces and calculate the identity component of their automorphism schemes. It turns out that there are $53$ distinct families of such surfaces if $p \neq 2,3$, while there are $61$ such families if $p = 3$, and $75$ such families if $p = 2$. Each of these families has at most one moduli. As a byproduct of our classification, it follows that weak del Pezzo surfaces with non-reduced automorphism schemes exist over $k$ if and only if $p \in \{2,3\}$.
\end{abstract}

 \section{Introduction}
 Recall that a weak del Pezzo surface over an algebraically closed field $k$ is a smooth projective surface $X$ with anticanonical divisor class $-K_X$ big and nef, or, equivalently, $X$ is $\bbP^1 \times \bbP^1$, the second Hirzebruch surface $\bbF_2$, or the blow-up of at most $8$ points in $\bbP^2$ in almost general position. More classically, weak del Pezzo surfaces appear as the minimal resolution of surfaces of degree $d$ in $\bbP^d$ which are neither cones nor projections of surfaces of minimal degree $d$ in $\bbP^{d+1}$ \cite[Definition 8.1.5]{Dolgachev}.

 By a result of Matsumura and Oort \cite{MatsumuraOort}, the automorphism functor $\Aut_X$ of a proper variety $X$ over $k$ is representable by a group scheme locally of finite type over $k$. 
 Since $\Aut_X$ is well-known for surfaces of minimal degree (that is for quadric surfaces, the Veronese surface, and rational normal scrolls \cite[Corollary 8.1.2]{Dolgachev}), weak del Pezzo surfaces form the first class of smooth projective surfaces for which the study of $\Aut_X$ is interesting. In this paper, we are concerned with the identity component $\Aut_X^0$ of $\Aut_X$, which can be non-reduced in positive characteristic.

 While this non-reducedness phenomenon does not occur for smooth projective curves, we will see that it appears for one of the first non-trivial classes of smooth projective surfaces, namely for weak del Pezzo surfaces (see also \cite{Neuman}), at least in characteristic $2$ and $3$. This means that for a weak del Pezzo surface $X$ in characteristic $2$ and $3$ we may have $h^0(X,T_X) > {\rm dim} \Aut_X^0$, that is, $X$ may have more global vector fields than expected.

 More classically, automorphisms of (weak) del Pezzo surfaces are being studied in the context of the plane Cremona group, i.e. the group of birational automorphisms of $\bbP^2$. The main reason for this is that automorphisms of (weak) del Pezzo surfaces yield birational automorphisms of $\bbP^2$ that do not necessarily extend to biregular automorphisms. For the action of $\Aut_X^0$ on a weak del Pezzo surface $X$, the situation is very different, since this action always descends to an action on the whole minimal model of $X$ by Blanchard's Lemma \ref{L BlanchardLemma}.
 
 This special feature of the connected automorphism scheme $\Aut_X^0$ will enable us to calculate it explicitly for all weak del Pezzo surfaces that are blow-ups of $\bbP^2$ in terms of stabilizers as a subgroup scheme of $\PGL_3$. Using this, we will classify all weak del Pezzo surfaces $X$ with non-trivial $\Aut_X^0$ and determine their configurations of $(-2)$- and $(-1)$-curves, as well as their number of moduli, which is the content of the following \hyperref[MainTheorem]{Main Theorem}:

 \begin{Main Theoremx} \label{MainTheorem} 
 Let $X$ be a weak del Pezzo surface over an algebraically closed field. If $h^0(X,T_X) \neq 0$, then $X$ is one of the surfaces in Table \ref{Table98765}, \ref{TableRest}, \ref{Table4}, \ref{Table3}, \ref{Table2}, or Table \ref{Table1}. All cases exist and have an irreducible moduli space of the stated dimension.
 \end{Main Theoremx} 
 
 In Tables \ref{Table98765}, \ref{Table4}, \ref{Table3}, \ref{Table2}, \ref{Table1}, the configuration of $(-2)$-curves and $(-1)$-curves (``lines") on these surfaces is given in columns $2$-$4$. In the corresponding figures, a ``thick'' curve denotes a $(-2)$-curve, while a ``thin'' curve denotes a $(-1)$-curve. The intersection multiplicity of two such curves is no more than $3$ at every point; intersection multiplicities $1$ and $2$ will be clear from the picture, whereas we write a small $3$ next to the point of intersection if the intersection multiplicity is $3$. In column $5$ of the tables, we describe a general $S$-valued point of $\Aut_X^0$, where $S$ is a $k$-scheme. In particular, the dimension of $H^0(X,T_X) = \Aut_X^0(k[\epsilon]/(\epsilon^2))$ can be read off from this description and is listed in column $6$ for the convenience of the reader.
 Comparing this with the dimension of $\Aut_X^0$, it can be checked whether $\Aut_X^0$ is smooth or not. This is done in column $7$. If there is more than one weak del Pezzo surface with the configuration of curves and with the automorphism scheme as in the previous columns, we give the dimension of a modular family of such surfaces in column $8$. If, instead, there is a unique surface of this type, we write ``\{pt\}'' in column $8$ in order to emphasize that the surface is unique. Finally, in column $9$, we give the characteristic(s) in which the respective surface(s) exist(s).

 In particular, our classification also gives a complete list of weak del Pezzo surfaces with non-reduced automorphism schemes. In the following corollary, we list the characteristics $p$ and degrees $d$ for which every weak del Pezzo surface of degree $d$ in characteristic $p$ has reduced automorphism scheme.

 \begin{Corollary}
 Let $k$ be an algebraically closed field of characteristic $p \geq 0$. Then, every weak del Pezzo surface of degree $d$ over $k$ has reduced automorphism scheme if and only if one of the following three conditions holds:
 \begin{enumerate}
     \item $p \neq 2,3$,
     \item $p = 3$ and $d \geq 4$,
     \item $p = 2$ and $d \geq 5$.
 \end{enumerate}
 \end{Corollary}
 
 \begin{Remark}
 Since every Jacobian rational (quasi-)elliptic surface $X'$ is the blow-up of a weak del Pezzo surface $X$ of degree $1$ in the unique basepoint of its anticanonical linear system, Lemma \ref{L NeumanLemma} yields an isomorphism $\Aut_{X'}^0 \cong \Aut_X^0$. In particular, our \hyperref[MainTheorem]{Main Theorem} gives a complete classification of Jacobian rational (quasi-)elliptic surfaces with global vector fields. The non-Jacobian case is more involved and will be treated by the second named author in an upcoming article.
 \end{Remark}
 
\begin{Acknowledgements}
We thank Christian Liedtke for comments on a first version of this article.
The first named author would like to thank the Department of Mathematics at the
University of Utah for its hospitality while this article was written. \\
Research of the first named author is funded by the DFG Research Grant MA 8510/1-1 ``Infinitesimal Automorphisms of Algebraic Varieties''.
Research of the second named author is funded by the DFG Sachbeihilfe LI 1906/5-1 ``Geometrie von rationalen Doppelpunkten''. The second named author also gratefully acknowledges support by the doctoral program TopMath and the TUM Graduate School.
\end{Acknowledgements}
 
 \newpage
 
 \begin{table}[H] 
 $$
 \begin{array}{|c|c|c|c|c|c|c|c|c|} 
 \hline
 \text{Case}
& \text{Figure}
&(-2) \text{-curves} 
& \# \{\text{lines}\}
& \Aut_X^0 \subseteq \PGL_3
&h^0(X, T_X)
&
\addstackgap[1.8pt]{\begin{tabular}{c}
$\Aut_X^0$ \\
\text{smooth?}
\end{tabular}}
& \text{Moduli}
& \Char(k)
\\

\hline  \hline
\multicolumn{9}{|c|}{\textbf{degree }\bf{9}}
\\
\hline

9A \label{Tab9A} & \text{} & \emptyset & 0 & \PGL_3 &  8 & \checkmark & \{ \text{pt} \} & \text{any}\\ \hline

\hline \hline 

\multicolumn{9}{|c|}{\textbf{degree }\bf{8}}

\\

\hline

8A \label{Tab8A}
 & \text{Fig. \ref{Conf8A}} & \emptyset & 1 & 
\addstackgap[1.8pt]{\resizebox{!}{0.18in}{ $   \left( \begin{smallmatrix}
1 & b & c \\
 & e & f \\
 & h & i
\end{smallmatrix} \right)  $}}
&  6 & \checkmark & \{ \text{pt} \}& \text{any}\\  
\hline \hline
\multicolumn{9}{|c|}{\textbf{degree }\bf{7}}
\\
\hline

7A \label{Tab7A}
& \text{Fig. \ref{Conf7A}} & \emptyset & 3 &
\addstackgap[1.8pt]{\resizebox{!}{0.18in}{
$
  \left( \begin{smallmatrix}
1 &  & c \\
 & e & f \\
 &  & i
\end{smallmatrix}  
\right)
$
}}
&  4 & \checkmark & \{ \text{pt} \}& \text{any}\\ \hline

7B \label{Tab7B}
& \text{Fig. \ref{Conf7B}} & A_1 & 2 &
\addstackgap[1.8pt]{\resizebox{!}{0.18in}{$   \left( \begin{smallmatrix}
1 & b & c \\
 & e & f \\
 &  & i
\end{smallmatrix}\right)  $}}
&  5 & \checkmark & \{ \text{pt} \}& \text{any}\\ \hline

\hline \hline
\multicolumn{9}{|c|}{\textbf{degree }\bf{6}}
\\
\hline

6A \label{Tab6A}
& \text{Fig. \ref{Conf6A}} & \emptyset & 6 &
\addstackgap[1.8pt]{\resizebox{!}{0.18in}{$   \left( \begin{smallmatrix}
1 &  &  \\
 & e &  \\
 &  & i
\end{smallmatrix} \right)  $}}
&  2 & \checkmark & \{ \text{pt} \}& \text{any}\\ \hline

6B \label{Tab6B}
& \text{Fig. \ref{Conf6B}} & A_1 & 4 &
\addstackgap[1.8pt]{\resizebox{!}{0.18in}{$   \left( \begin{smallmatrix}
1 &  & c \\
 & e &  \\
 &  & i
\end{smallmatrix} \right)  $}}
&  3 &  \checkmark & \{ \text{pt} \}& \text{any}\\ \hline

6C \label{Tab6C}
& \text{Fig. \ref{Conf6C}} & A_1 & 3 &
\addstackgap[1.8pt]{\resizebox{!}{0.18in}{$   \left( \begin{smallmatrix}
1 &  & c \\
 & 1 & f \\
 &  & i
\end{smallmatrix} \right)  $}}
&  3 &  \checkmark & \{ \text{pt} \}& \text{any}\\ \hline

6D \label{Tab6D}
& \text{Fig. \ref{Conf6D}} & 2A_1 & 2 &
\addstackgap[1.8pt]{\resizebox{!}{0.18in}{$   \left( \begin{smallmatrix}
1 &  & c \\
 & e & f \\
 &  & i
\end{smallmatrix} \right)  $}}
&  4 & \checkmark & \{ \text{pt} \}& \text{any}\\ \hline

6E \label{Tab6E}
& \text{Fig. \ref{Conf6E}} & A_2 & 2 &
\addstackgap[1.8pt]{\resizebox{!}{0.18in}{$   \left( \begin{smallmatrix}
1 & b & c \\
 & e & f \\
 &  & e^2
\end{smallmatrix} \right)  $}}
&  4 &  \checkmark & \{ \text{pt} \}& \text{any}\\ \hline

6F \label{Tab6F}
& \text{Fig. \ref{Conf6F}} & A_2 + A_1 & 1 &
\addstackgap[1.8pt]{\resizebox{!}{0.18in}{$   \left( \begin{smallmatrix}
1 & b & c \\
 & e & f \\
 &  & i
\end{smallmatrix} \right)  $}}
&  5& \checkmark & \{ \text{pt} \}& \text{any}\\ \hline

\hline \hline
\multicolumn{9}{|c|}{\textbf{degree }\bf{5}}
\\
\hline

5A \label{Tab5A}
& \text{Fig. \ref{Conf5A}} & A_1 & 7 &
\addstackgap[1.8pt]{\resizebox{!}{0.18in}{$   \left( \begin{smallmatrix}
1 &  &  \\
 & 1 &  \\
 &  & i
\end{smallmatrix} \right)  $}}
&  1 & \checkmark & \{ \text{pt} \}& \text{any}\\ \hline

5B \label{Tab5B}
& \text{Fig. \ref{Conf5B}} & 2A_1 & 5 &
\addstackgap[1.8pt]{\resizebox{!}{0.18in}{$   \left( \begin{smallmatrix}
1 &  &  \\
 & e &  \\
 &  & i
\end{smallmatrix} \right)  $}}
&  2 & \checkmark & \{ \text{pt} \}& \text{any}\\ \hline

5C \label{Tab5C}
& \text{Fig. \ref{Conf5C}} & A_2 & 4 &
\addstackgap[1.8pt]{\resizebox{!}{0.18in}{$   \left( \begin{smallmatrix}
1 &  & c \\
 & 1 &  \\
 &  & i
\end{smallmatrix} \right)  $}}
&  2 & \checkmark & \{ \text{pt} \}& \text{any}\\ \hline

5D \label{Tab5D}
& \text{Fig. \ref{Conf5D}} & A_2+A_1 & 3 &
\addstackgap[1.8pt]{\resizebox{!}{0.18in}{$   \left( \begin{smallmatrix}
1 &  &  \\
 & e & f \\
 &  & i
\end{smallmatrix} \right)  $}}
&  3 & \checkmark & \{ \text{pt} \}& \text{any}\\ \hline

5E \label{Tab5E}
& \text{Fig. \ref{Conf5E}} & A_3 & 2 &
\addstackgap[1.8pt]{\resizebox{!}{0.18in}{$   \left( \begin{smallmatrix}
1 &  & c \\
 & e & f \\
 &  & e^2
\end{smallmatrix} \right)  $}}
&  3 & \checkmark & \{ \text{pt} \}& \text{any}\\ \hline

5F \label{Tab5F}
& \text{Fig. \ref{Conf5F}} & A_4 & 1 &
\addstackgap[1.8pt]{\resizebox{!}{0.18in}{$   \left( \begin{smallmatrix}
1 & b & c \\
 & e & f \\
 &  & e^3
\end{smallmatrix} \right)  $}}
&  4 & \checkmark & \{ \text{pt} \}& \text{any}\\ \hline

\end{array}$$

\caption{Weak del Pezzo surfaces of degree $\geq 5$ that are blow-ups of $\bbP^2$} \label{Table98765}

\end{table}

\begin{table}[H]
 $$
 \begin{array}{|c|c|c|c|c|c|c|c|}
 \hline
 \text{Case}
&(-2) \text{-curves} 
& \# \{\text{lines}\}
& \Aut_X^0
&h^0(X, T_X)
&
\addstackgap[1.8pt]{\begin{tabular}{c}
$\Aut_X^0$ \\
\text{smooth?}
\end{tabular}}
& \text{Moduli}
& \Char(k)
\\

\hline \hline 
\bbP^1 \times \bbP^1  & \emptyset & 0 &
\PGL_2 \times \PGL_2
&  6 & \checkmark & \{ \text{pt} \}& \text{any}\\ \hline
\bbF_2 & A_1 & 0 &
 \addstackgap[1.8pt]{\begin{tabular}{c}
 $(\Aut_{\bbP(1,1,2)})_{\rm red} $\\
  $=(\bbG_a^3 \rtimes \GL_2)/\mu_2$
 \end{tabular}}
&  7 & \checkmark & \{ \text{pt} \}& \text{any}\\ \hline
\end{array}$$
\caption{Weak del Pezzo surfaces of degree $8$ that are not blow-ups of $\bbP^2$} \label{TableRest}
\end{table}

\begin{table}[H]
 $$
 \begin{array}{|c|c|c|c|c|c|c|c|c|}
 \hline
 \text{Case}
& \text{Figure}
&(-2) \text{-curves} 
& \# \{\text{lines}\}
& \Aut_X^0 \subseteq \PGL_3
&h^0(X, T_X)
&
\addstackgap[1.8pt]{\begin{tabular}{c}
$\Aut_X^0$ \\
\text{smooth?}
\end{tabular}}
& \text{Moduli}
& \Char(k)
\\

\hline \hline
4A \label{Tab4A}
& \text{Fig. \ref{Conf4A}} & 2A_1 & 8 &
\addstackgap[1.8pt]{\resizebox{!}{0.18in}{$   \left( \begin{smallmatrix}
1 &  &  \\
 & 1 &  \\
 &  & i
\end{smallmatrix} \right)  $}}
&  1 & \checkmark & 1 \text{ dim}& \text{any}\\ \hline

4B \label{Tab4B}
& \text{Fig. \ref{Conf4B}} & 3A_1 & 6 &
\addstackgap[1.8pt]{\resizebox{!}{0.18in}{$   \left( \begin{smallmatrix}
1 &  &  \\
 & 1 &  \\
 &  & i
\end{smallmatrix} \right)  $}}
&  1 & \checkmark & \{ \text{pt} \}& \text{any}\\ \hline

4C \label{Tab4C}
& \text{Fig. \ref{Conf4C}} & A_2+A_1 & 6 &
\addstackgap[1.8pt]{\resizebox{!}{0.18in}{$   \left( \begin{smallmatrix}
1 &  &  \\
 & 1 &  \\
 &  & i
\end{smallmatrix} \right)  $}}
&  1 & \checkmark & \{ \text{pt} \}& \text{any}\\ \hline

4D \label{Tab4D}
& \text{Fig. \ref{Conf4D}} & A_3 & 5 &
\addstackgap[1.8pt]{\resizebox{!}{0.18in}{$   \left( \begin{smallmatrix}
1 &  &  \\
 & 1 &  \\
 &  & i
\end{smallmatrix} \right)  $}}
&  1 & \checkmark & \{ \text{pt} \}& \text{any}\\ \hline

4E \label{Tab4E}
& \text{Fig. \ref{Conf4E4M}} & A_3 & 4 &
\addstackgap[1.8pt]{\resizebox{!}{0.18in}{$   \left( \begin{smallmatrix}
1 &  & c \\
 & 1 &  \\
 &  & 1
\end{smallmatrix} \right) $}}
&  1 & \checkmark & \{ \text{pt} \}& \neq 2 \\ \hline

4F \label{Tab4F}
& \text{Fig. \ref{Conf4F}} & 4A_1 & 4 &
\addstackgap[1.8pt]{\resizebox{!}{0.18in}{$   \left( \begin{smallmatrix}
1 &  &  \\
 & e &  \\
 &  & i
\end{smallmatrix} \right) $}}
&  2 & \checkmark & \{ \text{pt} \}& \text{any} \\ \hline

4G \label{Tab4G}
& \text{Fig. \ref{Conf4G}} & A_2+2A_1 & 4 &
\addstackgap[1.8pt]{\resizebox{!}{0.18in}{$   \left( \begin{smallmatrix}
1 &  &  \\
 & e &  \\
 &  & i
\end{smallmatrix} \right)  $}}
&  2 & \checkmark & \{ \text{pt} \}& \text{any}\\ \hline

4H \label{Tab4H}
& \text{Fig. \ref{Conf4H}} & A_3+A_1 & 3 &
\addstackgap[1.8pt]{\resizebox{!}{0.18in}{$   \left( \begin{smallmatrix}
1 &  & c \\
 & 1 &  \\
 &  & i
\end{smallmatrix} \right)  $}}
&  2 & \checkmark & \{ \text{pt} \}& \text{any} \\ \hline

4I \label{Tab4I}
& \text{Fig. \ref{Conf4I}} & A_4 & 3 &
\addstackgap[1.8pt]{\resizebox{!}{0.18in}{$   \left( \begin{smallmatrix}
1 &  &  \\
 & e & f \\
 &  & e^2
\end{smallmatrix} \right) $}}
&  2 & \checkmark & \{ \text{pt} \}& \text{any} \\ \hline

4J \label{Tab4J}
& \text{Fig. \ref{Conf4N4J4O}} & D_4 & 2 &
\addstackgap[1.8pt]{\resizebox{!}{0.18in}{$   \left( \begin{smallmatrix}
1 &  & c \\
 & e &  \\
 &  & e^2
\end{smallmatrix} \right)  $}}
&  2 & \checkmark & \{ \text{pt} \}& \neq 2 \\ \hline

4K \label{Tab4K}
& \text{Fig. \ref{Conf4K}} & A_3+2A_1 & 2 &
\addstackgap[1.8pt]{\resizebox{!}{0.18in}{$   \left( \begin{smallmatrix}
1 &  &  \\
 & e & f \\
 &  & i
\end{smallmatrix} \right) $}}
&  3 & \checkmark & \{ \text{pt} \}& \text{any} \\ \hline

4L \label{Tab4L}
& \text{Fig. \ref{Conf4P4Q4L}} & D_5 & 1 &
\addstackgap[1.8pt]{\resizebox{!}{0.18in}{$   \left( \begin{smallmatrix}
1 &  & c \\
 & e & f \\
 &  & e^3
\end{smallmatrix} \right) $}}
&  3 & \checkmark & \{ \text{pt} \}& \neq 2 \\ \hline \hline

4M \label{Tab4M}
& \text{Fig. \ref{Conf4E4M}} & A_3 & 4 &
\addstackgap[1.8pt]{\resizebox{!}{0.18in}{$   \left( \begin{smallmatrix}
1 &  & c \\
 & 1 &  \\
 &  & i
\end{smallmatrix} \right) , i^2=1  $ }}
&  2 & \times & \{ \text{pt} \}& = 2 \\ \hline

4N \label{Tab4N}
& \text{Fig. \ref{Conf4N4J4O}} & D_4 & 2 &
\addstackgap[1.8pt]{\resizebox{!}{0.18in}{$   \left( \begin{smallmatrix}
1 &  & c \\
 & 1 & f \\
 &  & 1
\end{smallmatrix} \right) $}}
&  2 & \checkmark & \{ \text{pt} \}& = 2 \\ \hline

4O \label{Tab4O}
& \text{Fig. \ref{Conf4N4J4O}} & D_4 & 2 &
\addstackgap[1.8pt]{\resizebox{!}{0.18in}{$   \left( \begin{smallmatrix}
1 &  & c \\
 & e & f \\
 &  & e^2
\end{smallmatrix} \right)  $}}
&  3 & \checkmark & \{ \text{pt} \}& = 2 \\ \hline

4P \label{Tab4P}
& \text{Fig. \ref{Conf4P4Q4L}} & D_5 & 1 &
\addstackgap[1.8pt]{\resizebox{!}{0.18in}{$   \left( \begin{smallmatrix}
1 & b & c \\
 & 1 & f \\
 &  & 1
\end{smallmatrix} \right) $}}
&  3 & \checkmark & \{ \text{pt} \}& = 2 \\ \hline

4Q \label{Tab4Q}
& \text{Fig. \ref{Conf4P4Q4L}} & D_5 & 1 &
\addstackgap[1.8pt]{\resizebox{!}{0.18in}{$    \left( \begin{smallmatrix}
1 & b & c \\
 & e & f \\
 &  & e^3
\end{smallmatrix} \right) $}}
&  4 & \checkmark & \{ \text{pt} \}& = 2 \\ \hline

\end{array}$$
\caption{Weak del Pezzo surfaces of degree $4$}
\label{Table4}
\end{table}

\begin{table}[H]
 $$
 \begin{array}{|c|c|c|c|c|c|c|c|c|}
 \hline
 \text{Case}
& \text{Figure}
&(-2) \text{-curves} 
& \# \{\text{lines}\}
& \Aut_X^0 \subseteq \PGL_3
&h^0(X, T_X)
&
\addstackgap[1.8pt]{\begin{tabular}{c}
$\Aut_X^0$ \\
\text{smooth?}
\end{tabular}}
& \text{Moduli}
& \Char(k)
\\

\hline \hline

3A \label{Tab3A}
& \text{Fig. \ref{Conf3A}} & 2A_2 & 7 &
\addstackgap[1.8pt]{\resizebox{!}{0.18in}{$   \left( \begin{smallmatrix}
1 &  &  \\
 & 1 &  \\
 &  & i
\end{smallmatrix} \right) $}}
&  1 & \checkmark & 1 \text{ dim}& \text{any}\\ \hline

3B \label{Tab3B}
& \text{Fig. \ref{Conf3B}} & D_4 & 6 &
\addstackgap[1.8pt]{\resizebox{!}{0.18in}{$   \left( \begin{smallmatrix}
1 &  &  \\
 & 1 &  \\
 &  & i
\end{smallmatrix} \right)  $}}
&  1 & \checkmark & \{ \text{pt} \}& \text{any}\\ \hline

3C \label{Tab3C}
& \text{Fig. \ref{Conf3C}} & 2A_2+A_1 & 5 &
\addstackgap[1.8pt]{\resizebox{!}{0.18in}{ $  \left( \begin{smallmatrix}
1 &  &  \\
 & 1 &  \\
 &  & i
\end{smallmatrix} \right) $}}
&  1 & \checkmark & \{ \text{pt} \}& \text{any}\\ \hline

3D \label{Tab3D}
& \text{Fig. \ref{Conf3D}} & A_3+2A_1 & 5 &
\addstackgap[1.8pt]{\resizebox{!}{0.18in}{ $  \left( \begin{smallmatrix}
1 &  &  \\
 & 1 &  \\
 &  & i
\end{smallmatrix} \right) $}}
&  1 & \checkmark & \{ \text{pt} \}& \text{any}\\ \hline

3E \label{Tab3E}
& \text{Fig. \ref{Conf3E}} & A_4+A_1 & 4 &
\addstackgap[1.8pt]{\resizebox{!}{0.18in}{$   \left( \begin{smallmatrix}
1 &  &  \\
 & 1 &  \\
 &  & i
\end{smallmatrix} \right) $}}
&  1 & \checkmark & \{ \text{pt} \}& \text{any}\\ \hline

3F \label{Tab3F}
& \text{Fig. \ref{Conf3F3K}} & A_5 & 3 &
\addstackgap[1.8pt]{\resizebox{!}{0.18in}{$   \left( \begin{smallmatrix}
1 &  &  \\
 & 1 & f \\
 &  & 1
\end{smallmatrix} \right) $}}
&  1 & \checkmark & \{ \text{pt} \}& \neq 3 \\ \hline

3G \label{Tab3G}
& \text{Fig. \ref{Conf3G3O3P}} & D_5 & 3 &
\addstackgap[1.8pt]{\resizebox{!}{0.18in}{$   \left( \begin{smallmatrix}
1 &  &  \\
 & e &  \\
 &  & e^2
\end{smallmatrix} \right) $}}
&  1 & \checkmark & \{ \text{pt} \}& \neq 2 \\ \hline

3H \label{Tab3H}
& \text{Fig. \ref{Conf3H}} & 3A_2 & 3 &
\addstackgap[1.8pt]{\resizebox{!}{0.18in}{$   \left( \begin{smallmatrix}
1 &  &  \\
 & e &  \\
 &  & i
\end{smallmatrix} \right) $}}
&  2 & \checkmark & \{ \text{pt} \}& \text{any} \\ \hline

3I \label{Tab3I}
& \text{Fig. \ref{Conf3I}} & A_5+A_1 & 2 &
\addstackgap[1.8pt]{\resizebox{!}{0.18in}{$   \left( \begin{smallmatrix}
1 &  &  \\
 & e & f \\
 &  & e^2
\end{smallmatrix} \right)  $}}
&  2 & \checkmark & \{ \text{pt} \}& \text{any} \\ \hline

3J \label{Tab3J}
& \text{Fig. \ref{Conf3J3L3M3R3Q}} & E_6 & 1 &
\addstackgap[1.8pt]{\resizebox{!}{0.18in}{$   \left( \begin{smallmatrix}
1 &  & c \\
 & e &  \\
 &  & e^3
\end{smallmatrix} \right) $}}
&  2 & \checkmark & \{ \text{pt} \}& \neq 2,3 \\ \hline \hline

3K \label{Tab3K}
& \text{Fig. \ref{Conf3F3K}} & A_5 & 3 &
\addstackgap[1.8pt]{\resizebox{!}{0.18in}{$   \left( \begin{smallmatrix}
1 &  &  \\
 & e & f \\
 &  & e^2
\end{smallmatrix} \right) , e^3=1  $}}
&  2 & \times & \{ \text{pt} \}& = 3 \\ \hline

3L \label{Tab3L}
& \text{Fig. \ref{Conf3J3L3M3R3Q}} & E_6 & 1 &
\addstackgap[1.8pt]{\resizebox{!}{0.18in}{$   \left( \begin{smallmatrix}
1 &  & c \\
 & 1 & f \\
 &  & 1
\end{smallmatrix} \right)  $}}
&  2 & \checkmark & \{ \text{pt} \}& = 3 \\ \hline

3M \label{Tab3M}
& \text{Fig. \ref{Conf3J3L3M3R3Q}} & E_6 & 1 &
\addstackgap[1.8pt]{\resizebox{!}{0.18in}{$   \left( \begin{smallmatrix}
1 &  & c \\
 & e & f \\
 &  & e^3
\end{smallmatrix} \right) $}}
&  3 & \checkmark & \{ \text{pt} \}& = 3 \\ \hline \hline

3N \label{Tab3N}
& \text{Fig. \ref{Conf3N}} & A_4 & 6 &
\addstackgap[1.8pt]{\resizebox{!}{0.18in}{$   \left( \begin{smallmatrix}
1 &  &  \\
 & 1 &  \\
 &  & i
\end{smallmatrix} \right) , i^2=1  $}}
&  1  & \times & \{ \text{pt} \}& = 2 \\ \hline

3O \label{Tab3O}
& \text{Fig. \ref{Conf3G3O3P}} & D_5 & 3 &
\addstackgap[1.8pt]{\resizebox{!}{0.18in}{$   \left( \begin{smallmatrix}
1 &  &  \\
 & 1 & f \\
 &  & 1
\end{smallmatrix} \right)  $}}
&  1 & \checkmark & \{ \text{pt} \}& = 2 \\ \hline

3P \label{Tab3P}
& \text{Fig. \ref{Conf3G3O3P}} & D_5 & 3 &
\addstackgap[1.8pt]{\resizebox{!}{0.18in}{$   \left( \begin{smallmatrix}
1 &  &  \\
 & e & f \\
 &  & e^2
\end{smallmatrix} \right) $}}
&  2 & \checkmark & \{ \text{pt} \}& = 2 \\ \hline

3Q \label{Tab3Q}
& \text{Fig. \ref{Conf3J3L3M3R3Q}} & E_6 & 1 &
\addstackgap[1.8pt]{\resizebox{!}{0.18in}{$    \left( \begin{smallmatrix}
1 & b & c \\
 & 1 & b^2+b \\
 &  & 1
\end{smallmatrix} \right) $}}
&  2 & \checkmark & \{ \text{pt} \}& = 2 \\ \hline

3R \label{Tab3R}
& \text{Fig. \ref{Conf3J3L3M3R3Q}} & E_6 & 1 &
\addstackgap[1.8pt]{\resizebox{!}{0.18in}{$    \left( \begin{smallmatrix}
1 & b & c \\
 & e & b^2e \\
 &  & e^3
\end{smallmatrix} \right) $}}
&  3 & \checkmark & \{ \text{pt} \}& = 2 \\ \hline

\end{array}$$
\caption{Weak del Pezzo surfaces of degree $3$}
\label{Table3}
\end{table}

\begin{table}[H]
 $$
 \begin{array}{|c|c|c|c|c|c|c|c|c|}
 \hline
 \text{Case}
& \text{Figure}
&(-2) \text{-curves} 
& \# \{\text{lines}\}
& \Aut_X^0 \subseteq \PGL_3
&h^0(X, T_X)
&
\addstackgap[1.8pt]{\begin{tabular}{c}
$\Aut_X^0$ \\
\text{smooth?}
\end{tabular}}
& \text{Moduli}
& \Char(k)
\\

\hline \hline

2A \label{Tab2A}
& \text{Fig. \ref{Conf2A}} & 2A_3 & 6 &
\addstackgap[1.8pt]{\resizebox{!}{0.18in}{$   \left( \begin{smallmatrix}
1 &  &  \\
 & 1 &  \\
 &  & i
\end{smallmatrix} \right)  $}}
&  1 & \checkmark & 1 \text{ dim}& \text{any}\\ \hline

2B \label{Tab2B}
& \text{Fig. \ref{Conf2B}} & D_5+A_1 & 5 &
\addstackgap[1.8pt]{\resizebox{!}{0.18in}{$   \left( \begin{smallmatrix}
1 &  &  \\
 & 1 &  \\
 &  & i
\end{smallmatrix} \right) $}}
&  1 & \checkmark & \{ \text{pt} \}& \text{any} \\ \hline

2C \label{Tab2C}
& \text{Fig. \ref{Conf2C2S}} & E_6 & 4 &
\addstackgap[1.8pt]{\resizebox{!}{0.18in}{$    \left( \begin{smallmatrix}
1 &  &  \\
 & e &  \\
 &  & e^2
\end{smallmatrix} \right) $}}
&  1 & \checkmark & \{ \text{pt} \}& \neq 2 \\ \hline

2D \label{Tab2D}
& \text{Fig. \ref{Conf2D}} & 2A_3+A_1 & 4 &
\addstackgap[1.8pt]{\resizebox{!}{0.18in}{$   \left( \begin{smallmatrix}
1 &  &  \\
 & 1 &  \\
 &  & i
\end{smallmatrix} \right) $}}
&  1 & \checkmark & \{ \text{pt} \}& \text{any} \\ \hline

2E \label{Tab2E}
& \text{Fig. \ref{Conf2E}} & D_4+3A_1 & 4 &
\addstackgap[1.8pt]{\resizebox{!}{0.18in}{$   \left( \begin{smallmatrix}
1 &  &  \\
 & 1 &  \\
 &  & i
\end{smallmatrix} \right) $}}
&  1 & \checkmark & \{ \text{pt} \}& \text{any}\\ \hline

2F \label{Tab2F}
& \text{Fig. \ref{Conf2F}} & A_5+A_2 & 3 &
\addstackgap[1.8pt]{\resizebox{!}{0.18in}{$   \left( \begin{smallmatrix}
1 &  &  \\
 & 1 &  \\
 &  & i
\end{smallmatrix} \right) $ }}
&  1& \checkmark & \{ \text{pt} \}& \text{any} \\ \hline

2G \label{Tab2G}
& \text{Fig. \ref{Conf2G2T2U}} & D_6+A_1 & 2 &
\addstackgap[1.8pt]{\resizebox{!}{0.18in}{$   \left( \begin{smallmatrix}
1 &  &  \\
 & e &  \\
 &  & e^2
\end{smallmatrix} \right)  $ }}
&  1 & \checkmark& \{ \text{pt} \}& \neq 2 \\ \hline

2H \label{Tab2H}
& \text{Fig. \ref{Conf2H2V}} & A_7 & 2 &
\addstackgap[1.8pt]{\resizebox{!}{0.18in}{$   \left( \begin{smallmatrix}
1 &  &  \\
 & 1 & f \\
 &  & 1
\end{smallmatrix} \right) $ }}
&  1 & \checkmark& \{ \text{pt} \}& \neq 2 \\ \hline

2I \label{Tab2I}
& \text{Fig. \ref{Conf2I2L2M2W2X2Y}} & E_7 & 1 &
\addstackgap[1.8pt]{\resizebox{!}{0.18in}{$    \left( \begin{smallmatrix}
1 &  &  \\
 & e &  \\
 &  & e^3
\end{smallmatrix} \right)  $}}
&  1 & \checkmark & \{ \text{pt} \}& \neq 2,3 \\ \hline \hline

2J \label{Tab2J}
& \text{Fig. \ref{Conf2J}} & A_6 & 4 &
\addstackgap[1.8pt]{\resizebox{!}{0.18in}{$   \left( \begin{smallmatrix}
1 &  &  \\
 & e &  \\
 &  & e^2
\end{smallmatrix} \right) , e^3=1  $}}
&  1 & \times & \{ \text{pt} \}& = 3 \\ \hline

2K \label{Tab2K}
& \text{Fig. \ref{Conf2K2R}} & D_6 & 3 &
\addstackgap[1.8pt]{\resizebox{!}{0.18in}{$   \left( \begin{smallmatrix}
1 &  &  \\
 & e &  \\
 &  & e^2
\end{smallmatrix} \right) , e^3=1  $}}
&  1 & \times & \{ \text{pt} \}& = 3 \\ \hline

2L \label{Tab2L}
& \text{Fig. \ref{Conf2I2L2M2W2X2Y}} & E_7 & 1 &
\addstackgap[1.8pt]{\resizebox{!}{0.18in}{$    \left( \begin{smallmatrix}
1 &  &  \\
 & 1 & f \\
 &  & 1
\end{smallmatrix} \right) $}}
&  1 & \checkmark & \{ \text{pt} \}& = 3 \\ \hline

2M \label{Tab2M}
& \text{Fig. \ref{Conf2I2L2M2W2X2Y}} & E_7 & 1 &
\addstackgap[1.8pt]{\resizebox{!}{0.18in}{$    \left( \begin{smallmatrix}
1 &  &  \\
 & e & f \\
 &  & e^3
\end{smallmatrix} \right)  $}}
&  2 & \checkmark & \{ \text{pt} \}& = 3 \\ \hline \hline

2N \label{Tab2N}
& \text{Fig. \ref{Conf2N}} & A_5 & 7 &
\addstackgap[1.8pt]{\resizebox{!}{0.18in}{$   \left( \begin{smallmatrix}
1 &  &  \\
 & 1 &  \\
 &  & i
\end{smallmatrix} \right) , i^2=1  $}}
&  1 & \times & 1 \text{ dim}& =2 \\ \hline

2O \label{Tab2O}
& \text{Fig. \ref{Conf2O}} & D_5 & 8 &
\addstackgap[1.8pt]{\resizebox{!}{0.18in}{$   \left( \begin{smallmatrix}
1 &  &  \\
 & 1 &  \\
 &  & i
\end{smallmatrix} \right) , i^2=1  $}}
&  1 & \times& \{ \text{pt} \}& =2 \\ \hline

2P \label{Tab2P}
& \text{Fig. \ref{Conf2P}} & A_5+A_1 & 6 &
\addstackgap[1.8pt]{\resizebox{!}{0.18in}{$   \left( \begin{smallmatrix}
1 &  &  \\
 & 1 &  \\
 &  & i
\end{smallmatrix} \right) , i^2=1  $}}
&  1 & \times & \{ \text{pt} \}& = 2 \\ \hline

2Q \label{Tab2Q}
& \text{Fig. \ref{Conf2Q}} & A_5+A_1 & 5 &
\addstackgap[1.8pt]{\resizebox{!}{0.18in}{$   \left( \begin{smallmatrix}
1 &  &  \\
 & 1 &  \\
 &  & i
\end{smallmatrix} \right) , i^2=1  $ }}
&  1 & \times & \{ \text{pt} \}& =2 \\ \hline

2R \label{Tab2R}
& \text{Fig. \ref{Conf2K2R}} & D_6 & 3 &
\addstackgap[1.8pt]{\resizebox{!}{0.18in}{$   \left( \begin{smallmatrix}
1 &  &  \\
 & 1 & f \\
 &  & 1
\end{smallmatrix} \right) $ }}
&  1 & \checkmark & 1 \text{ dim}& =2 \\ \hline

2S \label{Tab2S}
& \text{Fig. \ref{Conf2C2S}} & E_6 & 4 &
\addstackgap[1.8pt]{\resizebox{!}{0.18in}{$    \left( \begin{smallmatrix}
1 &  &  \\
 & e & f \\
 &  & e^2
\end{smallmatrix} \right) , f^2=0 $}}
&  2 & \times & \{ \text{pt} \}& = 2 \\ \hline

2T \label{Tab2T}
& \text{Fig. \ref{Conf2G2T2U}} & D_6+A_1 & 2 &
\addstackgap[1.8pt]{\resizebox{!}{0.18in}{$   \left( \begin{smallmatrix}
1 &  &  \\
 & 1 & f \\
 &  & 1
\end{smallmatrix} \right) $}}
&  1 & \checkmark & \{ \text{pt} \}& = 2 \\ \hline

2U \label{Tab2U}
& \text{Fig. \ref{Conf2G2T2U}} & D_6+A_1 & 2 &
\addstackgap[1.8pt]{\resizebox{!}{0.18in}{$    \left( \begin{smallmatrix}
1 &  &  \\
 & e & f \\
 &  & e^2
\end{smallmatrix} \right) $}}
&  2 & \checkmark& \{ \text{pt} \}& = 2 \\ \hline

2V \label{Tab2V}
& \text{Fig. \ref{Conf2H2V}} & A_7 & 2 &
\addstackgap[1.8pt]{\resizebox{!}{0.18in}{$   \left( \begin{smallmatrix}
1 &  &  \\
 & e & f \\
 &  & e^2
\end{smallmatrix} \right) , e^4=1 $}}
&  2 & \times& \{ \text{pt} \}& =2 \\ \hline

2W \label{Tab2W}
& \text{Fig. \ref{Conf2I2L2M2W2X2Y}} & E_7 & 1 &
\addstackgap[1.8pt]{\resizebox{!}{0.18in}{$    \left( \begin{smallmatrix}
1 &  & c \\
 & 1 &  \\
 &  & 1
\end{smallmatrix} \right) $}}
&  1 & \checkmark & \{ \text{pt} \}& =2 \\ \hline

2X \label{Tab2X}
& \text{Fig. \ref{Conf2I2L2M2W2X2Y}} & E_7 & 1 &
\addstackgap[1.8pt]{\resizebox{!}{0.18in}{$    \left( \begin{smallmatrix}
1 & b & c \\
 & 1 & b^2 \\
 &  & 1
\end{smallmatrix} \right) $}}
&  2 & \checkmark & \{ \text{pt} \}& = 2 \\ \hline

2Y \label{Tab2Y}
& \text{Fig. \ref{Conf2I2L2M2W2X2Y}} & E_7 & 1 &
\addstackgap[1.8pt]{\resizebox{!}{0.18in}{$    \left( \begin{smallmatrix}
1 & b & c \\
 & e & b^2e \\
 &  & e^3
\end{smallmatrix} \right) $}}
&  3 & \checkmark & \{ \text{pt} \}& = 2 \\ \hline

\end{array}$$
\caption{Weak del Pezzo surfaces of degree $2$}
\label{Table2}
\end{table}

\begin{table}[H]
 $$
 \begin{array}{|c|c|c|c|c|c|c|c|c|}
 \hline
 \text{Case}
& \text{Figure}
&(-2) \text{-curves} 
& \# \{\text{lines}\}
& \Aut_X^0 \subseteq \PGL_3
&h^0(X, T_X)
&
\addstackgap[1.8pt]{\begin{tabular}{c}
$\Aut_X^0$ \\
\text{smooth?}
\end{tabular}}
& \text{Moduli}
& \Char(k)
\\

\hline \hline
1A \label{Tab1A}
& \text{Fig. \ref{Conf1A}} & 2D_4 & 5 &
\addstackgap[1.8pt]{\resizebox{!}{0.18in}{$   \left( \begin{smallmatrix}
1 &  &  \\
 & 1 &  \\
 &  & i
\end{smallmatrix} \right)  $}}
&  1 & \checkmark & 1 \text{ dim} & \text{any} \\ \hline

1B \label{Tab1B}
& \text{Fig. \ref{Conf1B}} & E_6+A_2 & 4 &
\addstackgap[1.8pt]{\resizebox{!}{0.18in}{$   \left( \begin{smallmatrix}
1 &  &  \\
 & 1 &  \\
 &  & i
\end{smallmatrix} \right)  $}}
&  1 & \checkmark & \{ \text{pt} \} & \text{any} \\ \hline

1C \label{Tab1C}
& \text{Fig. \ref{Conf1C1P}} & E_7+A_1 & 3 &
\addstackgap[1.8pt]{\resizebox{!}{0.18in}{$   \left( \begin{smallmatrix}
1 &  &  \\
 & e &  \\
 &  & e^2
\end{smallmatrix} \right)  $}}
&  1 & \checkmark & \{ \text{pt} \} & \neq 2 \\ \hline

1D \label{Tab1D}
& \text{Fig. \ref{Conf1D1H1I1S1T}} & E_8 & 1 &
\addstackgap[1.8pt]{\resizebox{!}{0.18in}{$   \left( \begin{smallmatrix}
1 &  &  \\
 & e &  \\
 &  & e^3
\end{smallmatrix} \right) $}}
&  1 & \checkmark & \{ \text{pt} \} & \neq 2,3 \\ \hline \hline

1E \label{Tab1E}
& \text{Fig. \ref{Conf1E}} & D_7 & 5 &
\addstackgap[1.8pt]{\resizebox{!}{0.18in}{$   \left( \begin{smallmatrix}
1 &  &  \\
 & 1 &  \\
 &  & i
\end{smallmatrix} \right) , i^3=1 $}}
&  1 & \times & \{ \text{pt} \} & =3 \\ \hline

1F \label{Tab1F}
& \text{Fig. \ref{Conf1F}} & E_7 & 5 &
\addstackgap[1.8pt]{\resizebox{!}{0.18in}{$   \left( \begin{smallmatrix}
1 &  &  \\
 & e &  \\
 &  & e^2
\end{smallmatrix} \right) , e^3=1  $}}
&  1 & \times & \{ \text{pt} \} & = 3 \\ \hline

1G \label{Tab1G}
& \text{Fig. \ref{Conf1G}} & A_8 & 3 &
\addstackgap[1.8pt]{\resizebox{!}{0.18in}{$   \left( \begin{smallmatrix}
1 &  &  \\
 & e &  \\
 &  & e^2
\end{smallmatrix} \right) , e^3=1 $}}
&  1 & \times & \{ \text{pt} \} & =3 \\ \hline

1H \label{Tab1H}
& \text{Fig. \ref{Conf1D1H1I1S1T}} & E_8 & 1 &
\addstackgap[1.8pt]{\resizebox{!}{0.18in}{$   \left( \begin{smallmatrix}
1 &  &  \\
 & 1 & f \\
 &  & 1
\end{smallmatrix} \right)  $}}
&  1 & \checkmark & \{ \text{pt} \} & =3 \\ \hline

1I \label{Tab1I}
& \text{Fig. \ref{Conf1D1H1I1S1T}} & E_8 & 1 &
\addstackgap[1.8pt]{\resizebox{!}{0.18in}{$   \left( \begin{smallmatrix}
1 &  &  \\
 & e & f \\
 &  & e^3
\end{smallmatrix} \right) $}}
&  2 & \checkmark & \{ \text{pt} \} & =3 \\ \hline \hline

1J \label{Tab1J}
& \text{Fig. \ref{Conf1J}} & E_6 & 13 &
\addstackgap[1.8pt]{\resizebox{!}{0.18in}{$   \left( \begin{smallmatrix}
1 &  &  \\
 & 1 &  \\
 &  & i
\end{smallmatrix} \right) , i^2=1  $}}
&  1 & \times & 1 \text{ dim}& =2 \\ \hline

1K \label{Tab1K}
& \text{Fig. \ref{Conf1K}} & E_6+A_1 & 8 &
\addstackgap[1.8pt]{\resizebox{!}{0.18in}{$   \left( \begin{smallmatrix}
1 &  &  \\
 & 1 &  \\
 &  & i
\end{smallmatrix} \right) , i^2=1 $}}
&  1 & \times & \{ \text{pt} \} & =2 \\ \hline

1L \label{Tab1L}
& \text{Fig. \ref{Conf1L}} & A_7 & 8 &
\addstackgap[1.8pt]{\resizebox{!}{0.18in}{$   \left( \begin{smallmatrix}
1 &  &  \\
 & 1 &  \\
 &  & i
\end{smallmatrix} \right) , i^2=1  $}}
&  1 & \times & 1 \text{ dim}& =2 \\ \hline

1M \label{Tab1M}
& \text{Fig. \ref{Conf1M}} & E_7 & 5 &
\addstackgap[1.8pt]{\resizebox{!}{0.18in}{$   \left( \begin{smallmatrix}
1 &  &  \\
 & 1 & f \\
 &  & 1
\end{smallmatrix} \right) , f^2=0  $ }}
&  1 & \times & \{ \text{pt} \} & =2 \\ \hline

1N \label{Tab1N}
& \text{Fig. \ref{Conf1N}} & D_6+2A_1 & 6 &
\addstackgap[1.8pt]{\resizebox{!}{0.18in}{$   \left( \begin{smallmatrix}
1 &  &  \\
 & 1 &  \\
 &  & i
\end{smallmatrix} \right) , i^2=1  $}}
&  1 & \times& \{ \text{pt} \} & =2 \\ \hline

1O \label{Tab1O}
& \text{Fig. \ref{Conf1O}} & A_7+A_1 & 5 &
\addstackgap[1.8pt]{\resizebox{!}{0.18in}{$   \left( \begin{smallmatrix}
1 &  &  \\
 & 1 &  \\
 &  & i
\end{smallmatrix} \right) , i^2=1  $ }}
&  1 & \times & \{ \text{pt} \} & =2 \\ \hline

1P \label{Tab1P}
& \text{Fig. \ref{Conf1C1P}} & E_7+A_1 & 3 &
\addstackgap[1.8pt]{\resizebox{!}{0.18in}{$   \left( \begin{smallmatrix}
1 &  &  \\
 & e & f \\
 &  & e^2
\end{smallmatrix} \right) , f^2=0  $}}
&  2 & \times & \{ \text{pt} \} & =2 \\ \hline

1Q \label{Tab1Q}
& \text{Fig. \ref{Conf1Q1R}} & D_8 & 2 &
\addstackgap[1.8pt]{\resizebox{!}{0.18in}{$   \left( \begin{smallmatrix}
1 &  &  \\
 & 1 & f \\
 &  & 1
\end{smallmatrix}  \right) $}}
&  1 & \checkmark & 1 \text{ dim} & =2 \\ \hline

1R \label{Tab1R}
& \text{Fig. \ref{Conf1Q1R}} & D_8 & 2 &
\addstackgap[1.8pt]{\resizebox{!}{0.18in}{$   \left( \begin{smallmatrix}
1 &  &  \\
 & e & f \\
 &  & e^2
\end{smallmatrix} \right), e^4=1  $}}
&  2 & \times & \{ \text{pt} \} & =2 \\ \hline

1S \label{Tab1S}
& \text{Fig. \ref{Conf1D1H1I1S1T}} & E_8 & 1 &
\addstackgap[1.8pt]{\resizebox{!}{0.18in}{$   \left( \begin{smallmatrix}
1 &  & c \\
 & 1 &  \\
 &  & 1
\end{smallmatrix} \right) $}}
&  1 & \checkmark & \{ \text{pt} \} & =2 \\ \hline

1T \label{Tab1T}
& \text{Fig. \ref{Conf1D1H1I1S1T}} & E_8 & 1 &
\addstackgap[1.8pt]{\resizebox{!}{0.18in}{$   \left( \begin{smallmatrix}
1 & b & c \\
 & e & b^2e \\
 &  & e^3
\end{smallmatrix} \right),b^4=0  $}}
&  3 & \times & \{ \text{pt} \} & = 2 \\ \hline

\end{array}$$
\caption{Weak del Pezzo surfaces of degree $1$}
\label{Table1}
\end{table}

\newpage

 \section{Generalities}
This section provides the necessary background on the two main topics of this paper: weak del Pezzo surfaces and automorphism schemes. Throughout, we will be working over an algebraically closed field $k$.

 \subsection{Geometry of weak del Pezzo surfaces and their ``height''}
 We recall that every weak del Pezzo surface $X$ (except $X = \bbP^1 \times \bbP^1$ and the second Hirzebruch surface $X = \bbF_2$) is a successive blow-up of $\bbP^2$ satisfying certain properties, and we define the notion of ``height'', which is a measure for the complexity of $X$. We describe the set of all $(-2)$- and $(-1)$-curves on $X$ in terms of a realization of $X$ as a blow-up of $\bbP^2$.

 \begin{Definition}
 A \emph{weak del Pezzo surface} is a smooth projective surface $X$ with nef and big anticanonical class $-K_X$. The number $\deg(X) =K_X^2$ is called the \emph{degree of} $X$.
 \end{Definition}
 
 Recall that every birational morphism $\pi: X' \to X$ of smooth projective surfaces can be factored as
 $$
 \pi: X' \overset{\varphi}{\longrightarrow} X^{(n)} \overset{\pi^{(n-1)}}{\longrightarrow} X^{(n-1)} \overset{\pi^{(n-2)}}{\longrightarrow} \hdots \overset{\pi^{(1)}}{\longrightarrow} X^{(1)} \overset{\pi^{(0)}}{\longrightarrow} X^{(0)} = X \text{,}
 $$
 where $\varphi$ is an isomorphism and each $\pi^{(i)}: X^{(i+1)} \to X^{(i)}$ is the blow-up of a number of distinct closed points on $X^{(i)}$. The isomorphism $\varphi$ can be neglected by identifying $X'$ with $X^{(n)}$ via $\varphi$. Then, the above factorization becomes unique if in each step $\pi^{(i)}$ the maximal number of distinct closed points of $X^{(i)}$ is blown up.
 In this case, we call the above factorization of $\pi$ \emph{minimal}.
 
\begin{Definition}
Let $X$ and $X'$ be two smooth projective surfaces.
 \begin{itemize}[leftmargin=20pt]
 \item
For every birational morphism $\pi: X' \to X$, let $\pi = \pi^{(0)} \circ \hdots \circ \hdots \pi^{(n-1)}$ be its minimal factorization.  The \emph{height of} $\pi$ is defined as
 $$
 {\rm ht}(\pi) :=n.
 $$
 \item
 If $X'$ admits some birational morphism to $X$, we define the \emph{height of} $X'$ \emph{ over } $X$ as
 $$
 {\rm ht}(X'/X) := {\rm min}_{\pi: X' \to X} \{ {\rm ht}(\pi) \},
 $$
 where the minimum is taken over all birational morphisms $\pi: X' \to X$.
 \item If $X$ is a weak del Pezzo surface which is a successive blow-up of $\bbP^2$, then we define
 $$
 {\rm ht}(X) := {\rm ht}(X/\bbP^2)
 $$
 and if $X$ is not a blow-up of $\bbP^2$, we set ${\rm ht}(X) = 0$.
 \end{itemize}
\end{Definition}

 \begin{Remark}
 The reader should compare our notion of height with the height function on the bubble space of $X$ considered in \cite[Section 7.3.2]{Dolgachev}.
 \end{Remark}
 
 \begin{Notation} \label{N Notation}
 Let $\pi: X \to \bbP^2$ be a birational morphism of height $n$, and let $\pi = \pi^{(0)} \circ \hdots \circ \pi^{(n-1)}$ be its minimal factorization. Then, we fix the following notation:
 \begin{itemize}[leftmargin=20pt]
     \item 
     For each $0 \leq i < n$, we let $p_{1,i},\hdots,p_{n_i,i} \in X^{(i)}$ be the points blown up under $\pi^{(i)}$.       
     \item
     The exceptional divisor $(\pi^{(i)})^{-1}(p_{j,i}) \subseteq X^{(i+1)}$ over a closed point $p_{j,i} \in X^{(i)}$ will be denoted by $E_{j,i}$ for $j=1, \hdots, n_i$.
     \item
     For every $0 \leq i \leq k \leq n$, the strict transform of a curve $C \subseteq X^{(i)}$ along $\pi^{(i)} \circ \hdots \circ \pi^{(k-1)} $ is denoted by $C^{(k)}$.
 \end{itemize}
 \end{Notation}

Using this notation, we can now state a necessary and sufficient criterion for a successive blow-up of $\bbP^2$ to be a weak del Pezzo surface.
 
 \begin{Lemma}\label{L almost general position} \cite[Section 8.1.3]{Dolgachev}
 With Notation \ref{N Notation}, let $\pi:X \to \bbP^2$ be a birational morphism of height $n$. Then, $X$ is a weak del Pezzo surface if and only if the following three conditions hold.
 \begin{itemize}[leftmargin=20pt]
    \item
     On each $E_{j,i}$ there is at most one $p_{k,i+1}$.
     \item 
     For every line $\ell \subseteq \bbP^2$, there are at most three $p_{j,i}$ with $p_{j,i} \in \ell^{(i)}$, where $i$ ranges over $0, \hdots,  n-1$.
     \item
     For every irreducible conic $Q \subseteq \bbP^2$, there are at most six $p_{j,i}$ with $p_{j,i} \in Q^{(i)}$, where $i$ ranges over $0, \hdots,  n-1$.
 \end{itemize}
 \end{Lemma}
 
 \begin{Notation}
By Lemma \ref{L almost general position}, there is at most one $p_{k,i+1}$ on each $E_{j,i}$. Therefore, it makes sense to rename the $p_{k,i+1}$ so that $p_{k,i+1}$ lies on $E_{k,i}$. We will adopt this convention from now on.
 \end{Notation}

 If the above three conditions of Lemma \ref{L almost general position} are satisfied, we say that the points $p_{j,i}$ are in \emph{almost general position}. Using this terminology, there is the following well-known characterization of weak del Pezzo surfaces.

 \begin{Lemma}  \cite[Section 8.1.3]{Dolgachev} \label{L characterization}
 If $X$ is a weak del Pezzo surface, then
 \begin{enumerate}
     \item[(i)]
     $X \cong \bbP^1 \times \bbP^1$, or
     \item[(ii)]
     $X \cong \bbF_2$, the \emph{second Hirzebruch surface}, or
     \item[(iii)]
     $X$ is the successive blow-up of $\bbP^2$ in $n \leq 8$ points in \emph{almost general position}.
 \end{enumerate}
 In particular, we have $1 \leq \deg(X) \leq 9$, and ${\rm ht}(X)=0$ if and only if 
$X \in \{\bbP^2,\bbP^1 \times \bbP^1, \bbF_2\}$.
 \end{Lemma}

All the possible classes of $(-2)$- and $(-1)$-curves in the odd unimodular lattice $\Pic(X) = {\rm I}_{1,9-\deg(X)}$ of signature $(1,9-\deg(X))$ are well-known and described in \cite[Definition 23.7., Proposition 26.1.]{Manin} and \cite[Proposition 8.2.7]{Dolgachev}. This lattice-theoretic description can be translated into geometry (see \cite[Theorem 26.2. (ii)]{Manin} for the case of del Pezzo surfaces). A straightforward adaption of Manin's approach to our situation of \emph{weak} del Pezzo surfaces yields the following description of $(-2)$- and $(-1)$-curves on $X$.  

\begin{Lemma} \label{L criterion}
Let $X$ be a weak del Pezzo surface and let $\pi: X=X^{(n)} \to \bbP^2$ be a birational morphism of height $n$.
\begin{enumerate}[leftmargin=*]     
    \item[(i)] \label{L (-2)curves}
A curve on $X$ is a $(-2)$-curve if and only if it is of one of the following four types:
    \begin{itemize}[leftmargin=12pt]
        \item 
        the strict transform $E_{j,i}^{(n)}$ of an exceptional curve such that there is exactly one $p_{j,i+1}$ on $E_{j,i}$,
        \item 
        the strict transform $\ell^{(n)}$ of a line $\ell \subseteq \bbP^2$ such that there are exactly three $p_{j,i}$ with $p_{j,i} \in \ell^{(i)}$,
        \item
        the strict transform $C^{(n)}$ of an irreducible conic $C \subseteq \bbP^2$ such that there are exactly six $p_{j,i}$ with $p_{j,i} \in C^{(i)}$, or
        \item
        the strict transform $C^{(n)}$ of an irreducible singular cubic $C \subseteq \bbP^2$ such that there are exactly eight $p_{j,i}$ with $p_{j,i} \in C^{(i)}$, and such that one of the $p_{j,0}$ is the singular point of $C$.
    \end{itemize}
    \item[(ii)] \label{L (-1)curves}
 A curve on $X$ is a $(-1)$-curve if and only if it is of one of the following seven types:
    \begin{itemize}[leftmargin=12pt]
        \item 
        the strict transform $E_{j,i}^{(n)}$ of an exceptional curve such that there is no $p_{k,i+1}$ on $E_{j,i}$,
        \item 
        the strict transform $\ell^{(n)}$ of a line $\ell \subseteq \bbP^2$ such that there are exactly two $p_{j,i}$ with $p_{j,i} \in \ell^{(i)}$,
        \item
        the strict transform $C^{(n)}$ of an irreducible conic $C \subseteq \bbP^2$ such that there are exactly five $p_{j,i}$ with $p_{j,i} \in C^{(i)}$,
        \item
        the strict transform $C^{(n)}$ of an irreducible singular cubic $C \subseteq \bbP^2$ such that there are exactly seven $p_{j,i}$ with $p_{j,i} \in C^{(i)}$, and such that one of the $p_{j,0}$ is the singular point of $C$,
        \item
        the strict transform $C^{(n)}$ of an irreducible singular quartic $C \subseteq \bbP^2$ such that there are exactly eight $p_{j,i}$ with $p_{j,i} \in C^{(i)}$, and such that exactly three of the $p_{j,i}$ are double points of $C^{(i)}$,
        \item
        the strict transform $C^{(n)}$ of an irreducible singular quintic $C \subseteq \bbP^2$ such that there are exactly eight $p_{j,i}$ with $p_{j,i} \in C^{(i)}$, and such that exactly six of the $p_{j,i}$ are double points of $C^{(i)}$, or
        \item
        the strict transform $C^{(n)}$ of an irreducible singular sextic $C \subseteq \bbP^2$ such that there are exactly eight $p_{j,i}$ with $p_{j,i} \in C^{(i)}$, and such that exactly seven of the $p_{j,i}$ are double points of $C^{(i)}$ and exactly one of the $p_{j,0}$ is a triple point of $C$.
    \end{itemize}
\end{enumerate}
\end{Lemma}

\begin{Remark}
In particular, it can be seen that the criterion given in Lemma \ref{L almost general position} simply tells us that a successive blow-up of $\bbP^2$ in at most $8$ points is a weak del Pezzo surface if and only if we have never blown up a point on a $(-2)$-curve.
\end{Remark}

 \subsection{Automorphism schemes of blow-ups of smooth surfaces}
 
By a result of Matsumura and Oort \cite{MatsumuraOort}, the automorphism functor $\Aut_X^0$ of a proper variety over $k$ is representable and it is well-known that the tangent space of $\Aut_X^0$ can be identified naturally with $H^0(X,T_X)$. The main tool in our study of automorphism schemes of weak del Pezzo surfaces is the following lemma of Blanchard (see \cite[Theorem 7.2.1]{Brion1}).
 
 \begin{Lemma}\emph{(Blanchard's Lemma)} \label{L BlanchardLemma}
Let $f: Y \to X$ be a morphism of proper schemes over $k$ with $f_* \cal{O}_Y = \cal{O}_X$. Then, $f$ induces a homomorphism of group schemes $f_*: \Aut_Y^0 \to \Aut_X^0$. If $f$ is birational, then $f_*$ is a closed immersion.
\end{Lemma}

Thus, if $f$ is birational, we can and will identify $\Aut_Y^0$ with its image under $f_*$ in the following.
If $f$ is the blow-up of a smooth surface $X$ in a closed point $p$, it is possible to describe the image of $f_*$ (see \cite[Lemma 1.1]{Neuman} and \cite[Proposition 2.7]{Martin}). 

\begin{Lemma}\label{L NeumanLemma}
Let $f: Y \to X$ be the blow-up of a smooth projective surface $X$ in $n$ distinct points $p_1,\hdots,p_n \in X$. Then, we have $\Aut_Y^0 = (\bigcap_{i=1}^n{\rm Stab}_{p_i}^0)^0$.
\end{Lemma}

\prf
We prove the claim by induction on $n$ with the case $n=0$ being trivial.
For the inductive step, let $Y'$ be the blow-up of $X$ in $p_1,\hdots,p_{n-1}$. Then, $f': Y \to Y'$ is the blow-up in $p_n$ and we have $\Aut_{Y'}^0 = (\bigcap_{i=1}^{n-1}{\rm Stab}_{p_i}^0)^0$ by the induction hypothesis. Note that the identity component of the stabilizer of $p_n \in Y'$, with respect to the action of $\Aut_{Y'}^0$, is precisely $(\bigcap_{i=1}^{n}{\rm Stab}_{p_i}^0)^0$. By \cite[Remark 2.8]{Martin}, the $\Aut_Y^0$-action on $Y$ preserves the exceptional divisor of $f'$, hence $\Aut_Y^0$, being connected, is contained in $(\bigcap_{i=1}^{n}{\rm Stab}_{p_i}^0)^0$. Conversely, by \cite[Proposition 2.7]{Martin}, the $(\bigcap_{i=1}^{n}{\rm Stab}_{p_i}^0)^0$-action on $Y'$ lifts to $Y$ and since $(\bigcap_{i=1}^{n}{\rm Stab}_{p_i}^0)^0$ is connected, it actually lifts to a subgroup scheme of $\Aut_Y^0$. This finishes the proof.
\qed
\vspace{5mm}

 Let $\pi: X^{(n)} \to X$ be a birational morphism of smooth projective surfaces $X$ and $X^{(n)}$. Let $E \subseteq X^{(n)}$ be a $\pi$-exceptional irreducible curve. 
 Recall that the left-action of $\Aut_X^0$ on ${\rm Hilb}_X$ is given on $S$-valued points by
\begin{eqnarray*}
\Aut_X^0(S) \times {\rm Hilb}_X(S) &\overset{\rho(S)}{\longrightarrow} &{\rm Hilb}_X(S) \\
(g: X_S \to X_S, \iota: Z \hookrightarrow X_S) &\longmapsto &(Z \times_{\iota,X_S,g^{-1}} X_S \hookrightarrow X_S),
\end{eqnarray*}
where $X_S:= X \times S$, and this induces a natural action $\rho$ of $\Aut_{X^{(n)}}^0 \subseteq \Aut_X^0$ on ${\rm Hilb}_X$.
For a \emph{pencil} (that is, a $1$-dimensional linear system) $f: \calC \to \bbP^1 \subseteq {\rm Hilb}_X$ of curves on $X$  we will identify a point $p \in \bbP^1(S)$ with its fiber $\cal{C}_p$ under $f$. 
Let $V \subseteq \bbP^1$ be an open subset such that any two fibers $\cal{C}_p$ and $\cal{C}_{q}$ with $p,q \in V$ (as well as their strict transforms to all the $X^{(i)}$) have the same multiplicity at the $p_{j,i}$. Then, the rational map
\begin{eqnarray}\label{eq adaptedequation}
  \bbP^1 \supseteq & V  &{\longrightarrow} {\rm Hilb}_E \\
    & p & \longmapsto \mathcal{C}_p^{(n)} \cap E, \nonumber
\end{eqnarray}
can be extended to a morphism $\varphi$ from $\bbP^1$, since every irreducible component of ${\rm Hilb}_E$ is proper. In fact, the morphism $\varphi$ extends naturally to an $\Aut_{X^{(n)}}^0$-equivariant morphism
$$
B := (\Aut_{X^{(n)}}^0 \times {\rm Hilb}_X) \times_{\rho, {\rm Hilb}_X} \bbP^1 \longrightarrow \bbP^1 \overset{\varphi}{\longrightarrow} {\rm Hilb}_E,
$$
which we also call $\varphi$ by abuse of notation. The pullback of the family $\calC \to \bbP^1$ to $B$ yields a family of curves $\calC' \to B$, whose fiber over an $S$-valued point $p: S \to B \to \bbP^1$ is $\calC_p \times_{\iota,X_S,g^{-1}} X_S$, where $g \in \Aut_{X^{(n)}}^0(S)$. 

\begin{Definition}
Let $\pi: X^{(n)} \to X$ be a birational morphism of smooth projective surfaces $X$ and $X^{(n)}$. Let $E \subseteq X^{(n)}$ be a $\pi$-exceptional irreducible curve. A pencil of curves $f: \mathcal{C} \to \bbP^1$ is called \emph{adapted} (to $E$ and $\pi$) (or $E$\emph{-adapted}), if the morphism $\varphi$ of \eqref{eq adaptedequation} factors through an isomorphism $\bbP^1 \overset{\cong}{\to} E \subseteq {\rm Hilb}_E$.
\end{Definition}

\begin{Remark} \label{R obwohlsnichtadaptedist}
In most of the cases occurring in our classification we can choose the adapted pencil $\calC \to \bbP^1$ to be stable under the natural action of $\Aut_{X^{(n)}}^0$ on ${\rm Hilb}_X$. But even if this is not possible (i.e. if there exists $\calC_p \in \bbP^1(S)$ with $p:S \to \bbP^1$ such that $(\calC_p \times_{\iota,X_S,g^{-1}} X_S) \not \in \bbP^1(S)$), there is a unique action of $\Aut_{X^{(n)}}^0$ on $\bbP^1$ which makes the map $\varphi$ equivariant. This action is given on $S$-valued points as follows: the element $g \in \Aut_{X^{(n)}}^0(S)$ sends $\calC_p \in \bbP^1(S)$ with embedding $\iota:\calC_p \to X_S$ to the unique curve $\calC_{g(p)} \in \bbP^1(S)$ such that $\varphi( \calC_{g(p)}) = \varphi(\calC_p \times_{\iota,X_S,g^{-1}} X_S)$. In particular, orbits and stabilizers of the $\Aut_{X^{(n)}}^0$-action on $E$ can be calculated on $\bbP^1$, which we are going to exploit throughout.
\end{Remark}

 \begin{Remark/Notation} If $X = \bbP^2$ and $f_1,f_2$ are homogeneous equations of the same degree, we say that $\lambda f_1 + \mu f_2$ is \emph{adapted} (to $\pi$ and $E$)
 if the pencil spanned by $C_1= \cal{V}(f_1)$ and $C_2= \cal{V}(f_2)$ is adapted to $\pi$ and $E$ and if, in addition, we identified $C_1$ and $C_2$ with $[1:0]$ and $[0:1]$ in $\bbP^1$, respectively.
We will use this choice of coordinates to determine the orbits and stabilizers of the $\Aut_{X^{(n)}}^0$-action on $E$ explicitly by reducing it to a calculation on the pencil $[\lambda: \mu]$.
\end{Remark/Notation}

 \section{Strategy of proof}

 For the proof of our \hyperref[MainTheorem]{Main Theorem} we are going to argue inductively by going through all possible weak del Pezzo surfaces with non-trivial connected automorphism scheme in the order given by their height, i.e., we start with del Pezzo surfaces of height $0$, which are $\bbP^2, \bbP^1 \times \bbP^1$ and $\bbF_2$. Then, by Lemma \ref{L characterization}, to study del Pezzo surfaces of height $1$ we have to study blow-ups of $\bbP^2$ in a number of distinct ``honest'' points. After that, for height $2$, we have to consider del Pezzos that arise as blow-ups of points on exceptional divisors of blow-ups of points in $\bbP^2$ (sometimes we will also refer to such points as \emph{infinitely near} points of the first order, as was introduced in \cite[Section 7.3.2, p. 307]{Dolgachev}). 
Continuing this pattern, increasing the height by $1$ means that we have to study those surfaces that arise as blow-ups of points on the ``latest exceptional divisor".

In this subsection, we are going to further specify our strategy of proof and explain why the classification of weak del Pezzo surfaces with non-trivial vector fields obtained via our inductive procedure is indeed complete.

\subsection{Inductive strategy}
\hfill

Assume we have a complete set $\cal{L}_{i} = \{X_1,\hdots,X_{n_i}\}$ of representatives of weak del Pezzo surfaces of height $i$ that are blow-ups of $\bbP^2$, where for every $X_k$ we have fixed a birational morphism $\psi_k: X_k \to \bbP^1$ of height $i$. Further assume that we have calculated $(\psi_{k})_*(\Aut_{X_k}^0) \subseteq \PGL_3$ (see Lemma \ref{L NeumanLemma}) for every $k$.
If $i = 0$, such a list is given by $\cal{L}_0 = \{\bbP^2\}$ with $\Aut_{\bbP^2}^0 = \PGL_3$. 
Using the list $\cal{L}_i$, we produce a list $\cal{L}_{i+1}$ as follows:

\begin{Procedure} \label{Procedure} \hfill
\vspace{-2mm}
\begin{enumerate}
    \item[Step $1$:] Choose $X \in \cal{L}_i$ with $\psi:X \to \bbP^2$ and let
    $
    \psi: X \overset{\psi^{(i-1)}}{\longrightarrow} X^{(i-1)} \overset{\psi^{(i-2)}}{\longrightarrow} \hdots \overset{\psi^{(0)}}{\longrightarrow} X^{(0)}=\bbP^2
    $
    be the minimal factorization of $\psi$.
    \item[Step $2$:] If $i = 0$, let $E := X = \bbP^2$. 
    Otherwise, 
    let 
    $$
    E := \left( {\rm Exc}(\psi^{(i-1)}) - \bigcup_{j=0}^{i-2}{\rm Exc}(\psi^{(j)}) \right) - D, 
    $$
    where $D$ is the union of all $(-2)$-curves on $X$.
    Note that, if $i > 0$, then $E$ is the set of points on the ``latest" exceptional divisors that do not lie on $(-2)$-curves.
    Using the description of $\Aut_X^0$ as a subgroup scheme of $\PGL_3$, we calculate the orbits and stabilizers of the action of $\Aut_X^0$ on $E$ using $E_{j,i-1}$-adapted pencils.
    \item[Step $3$:] Choose a set of points $\{p_{j,i}\}_{j \in J} \subseteq E$ such that $(\bigcap_{j \in J} {\rm Stab}_{p_{j,i}}^0)^0$ is non-trivial and such that the blow-up $\psi': X' \to X$ in these points is still a weak del Pezzo surface (see the criterion given in Lemma \ref{L criterion}). In particular, since there is at most one of the $p_{j,i}$ on every exceptional curve, we may assume that $p_{j,i} \in E_{j,i-1}$.
     Note that we obtain isomorphic surfaces if we replace a point $p_{j,i}$ by a point in the same orbit under the action of $\bigcap_{k \neq j} \Stab_{p_{k,i}} \subseteq \Aut_X$.
    \item[Step $4$:] If $X'$ is isomorphic to a surface already contained in $\cal{L}_{j}$ for some $j \leq i+1$, discard this case. Otherwise, add $X'$ to $\cal{L}_{i+1}$, choose the blow-up realization $\psi \circ \psi': X' \to \bbP^2$, and calculate
    $$
     (\psi \circ \psi')_*(\Aut_{X'}^0) = (\psi_*) (\bigcap_{j=1}^{n_i} {\rm Stab}_{p_{j,i}}^0)^0 \subseteq \PGL_3.
    $$
    We do this by describing the group $\Aut_{X'}^0(R)$ for an arbitrary local $k$-algebra $R$ (see Subsection \ref{S Stabilizer}).
    \item[Step $5$:] Repeat Steps $3$ and $4$ until all possible point combinations $\{p_{j,i}\}_{j \in J}$ are exhausted. 
    \item[Step $6$:] Then, repeat Steps $1$-$5$ until all possible $X \in \cal{L}_i$ are exhausted.
\end{enumerate}
\end{Procedure}

\begin{Lemma}
For every $i$, the above Procedure \ref{Procedure} yields a complete set $\cal{L}_{i+1} = \{X_1,\hdots,X_{n_{i+1}}\}$ of representatives of isomorphism classes of weak del Pezzo surfaces of height $(i+1)$ with non-trivial global vector fields, that are blow-ups of $\bbP^2$.
\end{Lemma}

\prf
We prove the claim by induction on the height $i$. The case $i = 0$ with $\cal{L}_0=\{\bbP^2\}$ is clear. Therefore, assume that the claim holds for $i>0$ and that we have a list $\cal{L}_i$.

Let $X'$ be a weak del Pezzo surface of height $(i+1)$ with $h^0(X',T_{X'}) \neq 0$. Choose a birational morphism $\pi: X' \to \bbP^2$ with minimal factorization
$$
\pi: X' = X'^{(i+1)} \overset{\pi^{(i)}}{\longrightarrow} X'^{(i)} \overset{\pi^{(i-1)}}{\longrightarrow} \hdots \overset{\pi^{(0)}}{\longrightarrow} X'^{(0)}= \bbP^2
$$
such that for every birational morphism $\pi': X' \to \bbP^2$, the number of exceptional curves for $\pi'^{(i)}$ is at least as great as the number of exceptional curves for $\pi^{(i)}$, i.e. such that the number of points blown up by the last step $\pi^{(i)}$ is minimal.
By Lemma \ref{L BlanchardLemma}, there is an inclusion
$$
(\pi^{(i)})_* (\Aut_{X'}^0) \subseteq \Aut_{X'^{(i)}}^0.
$$
In particular, we have $h^0(X'^{(i)},T_{X'^{(i)}}) \neq 0$ since $\Aut_{X'}^0 \neq \{\id\}$ and $\pi^{(i)}_*$ is a closed immersion. Hence, by the induction hypothesis, there is $X \in \cal{L}_i$ such that there exists an isomorphism $\phi: X'^{(i)} \to X$ and $X$ comes with a birational morphism $\psi: X \to \bbP^2$.

To prove the claim, it suffices to show that $\phi \circ \pi^{(i)}$ is the blow-up of $X$ in a set of points $p_{1,i},\hdots,p_{n_i,i}$ on $E$ defined as in Procedure \ref{Procedure}.
Indeed, once we prove this, it will follow from Lemma \ref{L NeumanLemma} and the assumption $h^0(X',T_{X'}) \neq 0$ that $\Aut_{X'}^0 = (\bigcap_{j=1}^{n_i} \Stab_{p_{j,i}}^0)^0$ is non-trivial.

Now, note that the condition that the $p_{j,i}$ lie on $E$ is trivially satisfied if $i = 1$, and equivalent to $\phi \circ \pi^{(i)}$ being the first step in the minimal factorization of 
$$
\psi' := \psi \circ \phi \circ \pi^{(i)}: X' \to X^{(i)} \to X \to \bbP^2
$$
if $i > 1$.
Thus, to prove the case $i > 1$, we take the minimal factorization of $\psi'$ and let $\psi'^{(i)}: X' \to X''$ be the first morphism in the minimal factorization of $\psi'$. Since $X$ has height $i$, the morphism $\phi \circ \pi^{(i)}: X' \to X$ factors through $\psi'^{(i)}$, which means there is a morphism $f: X'' \to X$ such that $f \circ \psi'^{(i)} = \phi \circ \pi^{(i)}$. In particular, the number of points blown up under $\psi'^{(i)}$ is at most the number of points blown up under $\pi^{(i)}$. As we chose the number of points blown up under $\pi^{(i)}$ to be minimal, this shows that $f$ is an isomorphism. In fact, since $f$ is an isomorphism over $\bbP^2$, this isomorphism is unique and we can identify $X''$ with $X$. 
\qed

\begin{Remark}\label{R IsomorphismCheck}
In order for $X'$ to be isomorphic to a weak del Pezzo surface $X$ in our lists, a necessary condition is that the configurations of $(-1)$- and $(-2)$-curves on $X'$ and $X$ are the same, and that $\Aut_X^0 \cong \Aut_{X'}^0$. Then, one can blow down $(-1)$-curves on $X'$ to get a realization $\pi: X' \to \bbP^2$ as a blow-up of $\bbP^2$ similar to the realization $\psi: X \to \bbP^2$ we fixed for $X$. Finally, it has to be checked that the points blown up by $\pi$ and $\psi$ are the same up to automorphisms of $\bbP^2$. This last part is straightforward but tedious, and we leave the details to the reader in each case.
\end{Remark}

\subsection{On the calculation of stabilizers} \label{S Stabilizer}
Before starting our classification, let us explain how to calculate the scheme-theoretic stabilizers of the points $p_{j,i} \in E_{j,i-1}$ occurring in Step 4 of Procedure \ref{Procedure}. First, recall the definition of the scheme-theoretic stabilizer.

\begin{Definition}
Let $\rho: G \times X \to X$ be an action of a group scheme $G$ on a scheme $X$ over $k$. Let $p: \Spec k \to X$ be a $k$-valued point.
The stabilizer ${\rm Stab}_p \subseteq G$ of $p$ with respect to $\rho$ is defined as
\begin{eqnarray*}
{\rm Stab}_p: & ( Sch/k) &\to ( Sets) \\
&S &\mapsto \{ g \in G(S) \mid g(p_S) = p_S \}
\end{eqnarray*}
where $p_S: S \to \Spec k \to X$. 
\end{Definition}

The stabilizer ${\rm Stab}_p \subseteq G$ is a closed subgroup scheme of $G$. As mentioned in Step 4 of Procedure \ref{Procedure}, we will describe only the $R$-valued points of the stabilizers occuring in our classification, where $R$ is a local $k$-algebra. This is sufficient, since in each case -- all the conditions on the matrices in $\PGL_3(R)$ of Tables \ref{Table98765}, \ref{Table4}-\ref{Table1} being given by polynomial equations which respect the group structure on $\PGL_3$ -- there will be an obvious closed subgroup scheme $G$ of $\PGL_3$ that admits the same $R$-valued points as the given stabilizer. The group scheme $G$ will then be equal to the stabilizer because of the following well-known lemma.

\begin{Lemma} \label{L well-known}
Let $Z_1,Z_2 \subseteq X$ be two closed subschemes of a scheme $X$ over a field $k$. If $Z_1(R) = Z_2(R) \subseteq X(R)$ for all local $k$-algebras $R$, then $Z_1 = Z_2$ as closed subschemes of $X$.
\end{Lemma}

The advantage of only considering $R$-valued points of $\PGL_n$ lies in the fact that $R$-valued points $\bbP^n$ are simply given by $(n+1)$-tuples of elements in $R$ up to units in $R$ such that at least one of the elements in the $(n+1)$-tuple is a unit. This allows us to describe the action of $\Aut_X^0(R)$ on $E_{j,i-1}(R) \cong \bbP^1(R)$ explicitly using adapted pencils, so that the calculation of the scheme-theoretic stabilizer of a $k$-valued point $p_{j,i} \in E_{j,i-1}$ becomes straightforward (by Lemma \ref{L well-known}). Thus, $R$ will denote a local $k$-algebra from now on.

\section{Proof of \hyperref[MainTheorem]{Main Theorem}: Classification}
 
 In this section, we will carry out Procedure \ref{Procedure} to obtain the classification of weak del Pezzo surfaces with global regular vector fields and prove our \hyperref[MainTheorem]{Main Theorem}.
 
 Firstly, note that there are two weak del Pezzo surfaces which do not fit into the framework of Procedure \ref{Procedure}, namely those which are not blow-ups of $\bbP^2$. By Lemma \ref{L characterization}, these are $\bbP^1 \times \bbP^1$ and $\bbF_2$. As is well-known, we have $\Aut_{\bbP^1 \times \bbP^1} = \PGL_2 \times \PGL_2$. As for $\Aut_{\bbF_2}$, we make use of the fact that this group scheme is smooth and connected by \cite[Theorem 1 and Lemma 10]{Maruyama}. An explicit description of this group scheme is given in \cite{Maruyama}. Alternatively, one can blow-down the unique $(-2)$-curve on $\bbF_2$ to obtain the weighted projective plane $\bbP(1,1,2)$ and use the fact that $(\Aut_{\bbP(1,1,2)})_{\red}$ fixes the unique singular point on $\bbP(1,1,2)$. Hence, this action lifts to $\bbF_2$ and we get $\Aut_{\bbF_2} = (\Aut_{\bbP(1,1,2)})_{\red}$. These results are summarized in Table \ref{TableRest}.
 
 For the remaining cases, we can apply Procedure \ref{Procedure} and we will subdivide the proof into subsections according to the height of our weak del Pezzo surfaces. Throughout, we write $\ell_{f}:= \cal{V}(f) $ for the line given by $f=0$ in $\bbP^2$. Recall that in the following figures a ``thick'' curve denotes a $(-2)$-curve, while a ``thin'' curve denotes a $(-1)$-curve. The intersection multiplicity of two such curves is at most $3$ at every point; intersection multiplicities $1$ and $2$ will be clear from the picture, whereas we write a small $3$ next to the point of intersection if the intersection multiplicity is $3$.

\subsection{Height 0}
We have $\cal{L}_0 = \{X_{9A}\}$, where $X_{9A} := \bbP^2$ with $\Aut_{\bbP^2} = \PGL_3$.

\subsection{Height 1}

\subsubsection*{\underline{Case \hyperref[Tab9A]{$9A$}}}
    
    In this case, $X = \bbP^2$ and $\psi = \id$. We have $E = \bbP^2$ and the action of $\Aut_X^0 = \PGL_3$ on $E$ is transitive. Now, note that if $p_{1,0},\hdots,p_{n_0,0} \in \bbP^2$ are points such that at least four of them are in general position, then
    $$
    \Aut_{X'}^0 = (\bigcap_{j=1}^{n_0} \Stab_{p_{j,0}}^0) = \{ \ast \}.
    $$
    On the other hand, according to Lemma \ref{L almost general position}, to guarantee that $X'$ is a weak del Pezzo surface, no more than three of the $p_{j,0}$ may be on a line. Up to isomorphism, this leaves the following five possibilities for $p_{1,0},\hdots,p_{n,0}$:

    \begin{enumerate}[leftmargin=*]
    \noindent     \begin{minipage}{0.65\textwidth}
        \item \vspace{2mm}
        $n = 4$ and $ p_{1,0},p_{2,0},p_{4,0}$ on a line $\ell$, $p_{3,0} \not\in \ell$:
        Using the action of $\PGL_3$, we may assume that $p_{1,0} = [1:0:0], p_{2,0} = [0:1:0], \\p_{3,0} = [0:0:1], p_{4,0} = [1:1:0]$ and $\ell = \ell_{z}$.

    \begin{itemize}[leftmargin=20pt]    
        \item
        $
        \Aut_{X'}^0(R) =    
        \left\{  \left( \begin{smallmatrix}
1 &  &  \\
 & 1 &  \\
 &  & i
\end{smallmatrix} \right) 
\in \PGL_3(R) \right\}
        $
        \item $(-2)$-curves: $\ell_z^{(1)}$
        \item $(-1)$-curves: $E_{1,0}^{}, E_{2,0}^{}, E_{3,0}^{},E_{4,0}^{}, \ell_x^{(1)}, \ell_y^{(1)}, \ell_{x-y}^{(1)}$
        \item
        with configuration as in Figure \ref{Conf5A}.
    \end{itemize} 
    This is case \hyperref[Tab5A]{$5A$}.
    \end{minipage} \hspace{2mm} \begin{minipage}{0.3\textwidth}
\begin{center}
\includegraphics[width=0.65\textwidth]{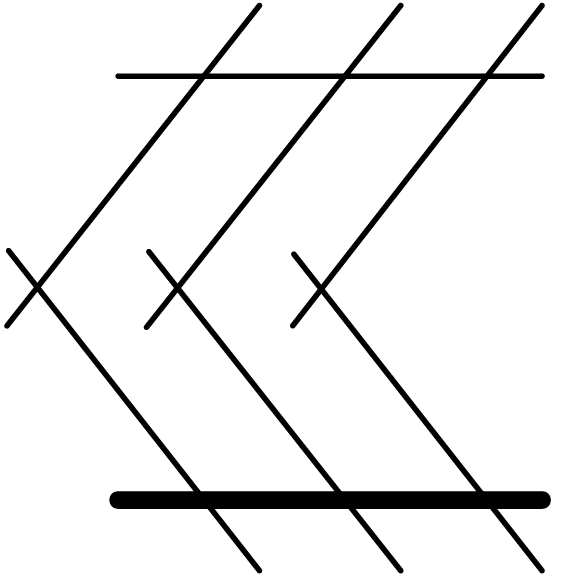}
\vspace{-3.5mm}\captionof{figure}{} \label{Conf5A}
\end{center}
\end{minipage}

        \noindent     \begin{minipage}{0.65\textwidth}
        \item \vspace{2mm}
        $n = 3$, all points on a line $\ell$: 
        We may assume that $p_{1,0} = [1:0:0],\\ p_{2,0} = [0:1:0], p_{3,0} = [1:1:0]$ and $\ell = \ell_z$.

    \begin{itemize}[leftmargin=20pt]
        \item
        $
        \Aut_{X'}^0(R) =    
        \left\{  \left( \begin{smallmatrix}
1 &  & c \\
 & 1 & f \\
 &  & i
\end{smallmatrix} \right) 
\in \PGL_3(R) \right\}
        $
        \item $(-2)$-curves: $\ell_z^{(1)}$
        \item $(-1)$-curves: $E_{1,0}^{}, E_{2,0}^{}, E_{3,0}^{}$
        \item
        with configuration as in Figure \ref{Conf6C}.
    \end{itemize} 
    This is case \hyperref[Tab6C]{$6C$}.
    \end{minipage} \hspace{2mm} \begin{minipage}{0.3\textwidth} \begin{center}
\includegraphics[width=0.9\textwidth]{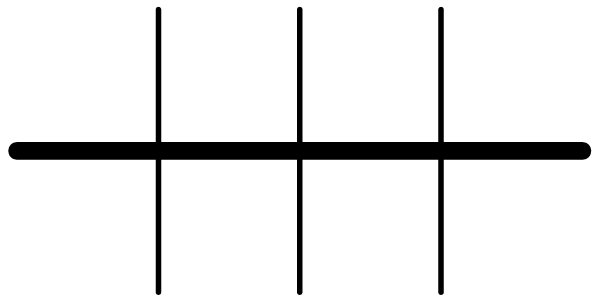}
\vspace{-3.5mm}\captionof{figure}{} \label{Conf6C}
\end{center} \end{minipage}    
        
        \noindent     \begin{minipage}{0.65\textwidth}
        \item \vspace{2mm}
        $n = 3$, not all points on a line: We may assume that $p_{1,0} = [1:0:0],\\ p_{2,0} = [0:1:0], p_{3,0} = [0:0:1]$.

    \begin{itemize}[leftmargin=20pt]
        \item
        $
        \Aut_{X'}^0(R) =    
        \left\{  \left( \begin{smallmatrix}
1 &  &  \\
 & e &  \\
 &  & i
\end{smallmatrix} \right) 
\in \PGL_3(R) \right\}
        $
        \item $(-2)$-curves: none
        \item $(-1)$-curves: $E_{1,0}^{},E_{2,0}^{}, E_{3,0}^{}, \ell_x^{(1)}, \ell_y^{(1)}, \ell_{z}^{(1)}$
        \item
        with configuration as in Figure \ref{Conf6A}.
    \end{itemize} 
    This is case \hyperref[Tab6A]{$6A$}.
        \end{minipage} \hspace{2mm} \begin{minipage}{0.3\textwidth} \begin{center} \includegraphics[width=0.72\textwidth]{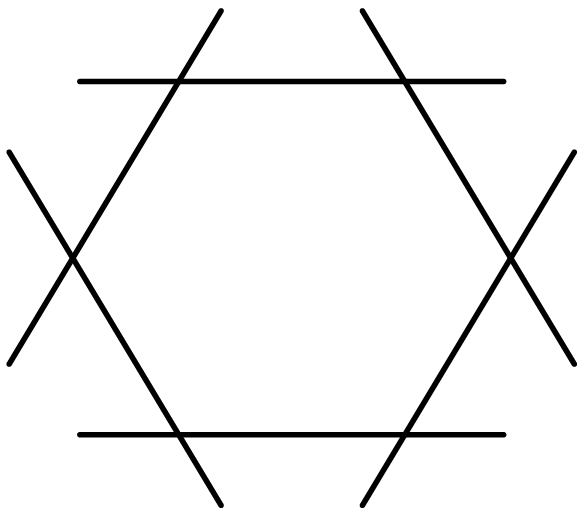} \vspace{-3.5mm}\captionof{figure}{} \label{Conf6A}  \end{center} \end{minipage} 
        
        \noindent     \begin{minipage}{0.65\textwidth}
        \item \vspace{2mm}
        $n = 2$: We may assume that $p_{1,0} = [1:0:0], p_{2,0} = [0:1:0]$.

    \begin{itemize}[leftmargin=20pt]
        \item
        $
        \Aut_{X'}^0(R) =    
        \left\{  \left( \begin{smallmatrix}
1 &  & c \\
 & e & f \\
 &  & i
\end{smallmatrix} \right) 
\in \PGL_3(R) \right\}
        $
        \item $(-2)$-curves: none
        \item $(-1)$-curves: $E_{1,0}^{}, E_{2,0}^{}, \ell_{z}^{(1)}$
        \item
        with configuration as in Figure \ref{Conf7A}.
    \end{itemize} 
    This is case \hyperref[Tab7A]{$7A$}.
        \end{minipage} \hspace{2mm} \begin{minipage}{0.3\textwidth} \begin{center} \includegraphics[width=0.85\textwidth]{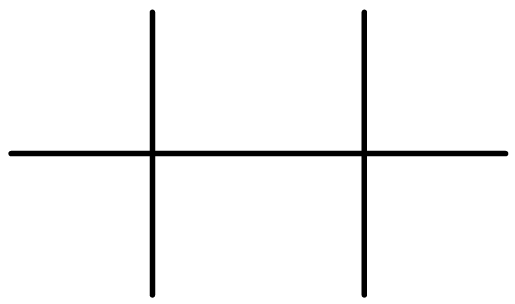} \vspace{-3.5mm}\captionof{figure}{}  \label{Conf7A} \end{center} \end{minipage} 
        
        \noindent     \begin{minipage}{0.65\textwidth}
        \item \vspace{2mm}
        $n = 1$: We may assume that $p_{1,0} = [1:0:0]$. 

    \begin{itemize}[leftmargin=20pt]
        \item
        $
        \Aut_{X'}^0(R) =    
        \left\{  \left( \begin{smallmatrix}
1 & b & c \\
 & e & f \\
 & h & i
\end{smallmatrix} \right) 
\in \PGL_3(R) \right\}
        $
        \item $(-2)$-curves: none
        \item $(-1)$-curves: $E_{1,0}^{}$
        \item
        with configuration as in Figure \ref{Conf8A}.
    \end{itemize} 
    This is case \hyperref[Tab8A]{$8A$}. 
        \end{minipage} \hspace{2mm} \begin{minipage}{0.3\textwidth} \begin{center} \includegraphics[width=0.9\textwidth]{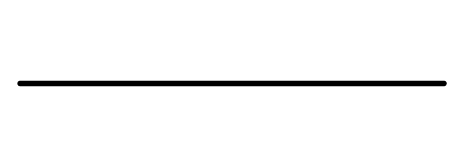} \vspace{-3.5mm}\captionof{figure}{}  \label{Conf8A} \end{center} \end{minipage} 
        
  \end{enumerate}

\vspace{2mm} 
\noindent
Summarizing, we obtain $\cal{L}_1 = \{X_{5A},X_{6C},X_{6A}, X_{7A},X_{8A}\}$.

\newpage

\subsection{Height 2}

 \subsubsection*{\underline{Case \hyperref[Tab5A]{$5A$}}}
We have $E= (\bigcup_{j=1}^{4} E_{j,0}^{})- \ell_{z}^{(1)}$.  Recall that the $R$-valued points of $\Aut_X^0$ are given by $\Aut_X^0(R) = \left\{  \left( \begin{smallmatrix}
1 &  &  \\
 & 1 &  \\
 &  & i
\end{smallmatrix} \right) 
\in \PGL_3(R) \right\}$.
We calculate the action of $\Aut_X^0$ on the $E_{j,0}^{}$ using adapted pencils:

\begin{itemize}[leftmargin=25pt]
    \item[-] $\lambda y + \mu z$ is $E_{1,0}^{}$-adapted and $\Aut_X^0(R)$ acts as
    $
    [\lambda:\mu] \mapsto [\lambda: i\mu]
    $
    \item[-] $\lambda x + \mu z$ is $E_{2,0}^{}$-adapted and $\Aut_X^0(R)$ acts as
    $
    [\lambda:\mu] \mapsto [\lambda: i\mu]
    $
    \item[-] $\lambda x + \mu y$ is $E_{3,0}^{}$-adapted and $\Aut_X^0(R)$ acts as
    $
    [\lambda:\mu] \mapsto [\lambda:\mu]
    $
    \item[-] $\lambda (x-y) + \mu z$ is $E_{4,0}^{}$-adapted and $\Aut_X^0(R)$ acts as
    $
    [\lambda:\mu] \mapsto [\lambda: i\mu]
    $
\end{itemize}

\noindent
In particular, there is one unique point with non-trivial stabilizer on $E \cap E_{1,0}^{}, E \cap E_{2,0}^{},$ and $E \cap E_{4,0}^{}$, respectively. Since $p_{1,0},p_{2,0}$ and $p_{4,0}$ can be interchanged by automorphisms of $\bbP^2$ preserving $p_{3,0}$, we have the following ten possibilities for $p_{1,1},\hdots,p_{n,1}$:

\begin{enumerate}[leftmargin=*]
\noindent     \begin{minipage}{0.65\textwidth}
    \item \vspace{2mm}

$ 
p_{1,1}= E_{1,0}^{} \cap \ell_{y}^{(1)}, 
p_{2,1}= E_{2,0}^{} \cap \ell_{x}^{(1)},
p_{3,1}= E_{3,0}^{} \cap \ell_{x+ \alpha y}^{(1)}$ \\with $\alpha \not\in \{0,-1\}$, $p_{4,1}= E_{4,0}^{} \cap \ell_{x-y}^{(1)}$
    \begin{itemize}[leftmargin=20pt]
        \item
        $
        \Aut_{X'}^0(R) =    
        \left\{  \left( \begin{smallmatrix}
1 &  &  \\
 & 1 &  \\
 &  & i
\end{smallmatrix} \right) 
\in \PGL_3(R) \right\}
        $
        \item $(-2)$-curves: $E_{1,0}^{(2)},E_{2,0}^{(2)},E_{3,0}^{(2)},E_{4,0}^{(2)}, \ell_{x}^{(2)}, \ell_{y}^{(2)}, \ell_{z}^{(2)}, \ell_{x-y}^{(2)}$
        \item $(-1)$-curves: $E_{1,1}^{}, E_{2,1}^{}, E_{3,1}^{}, E_{4,1}^{}, \ell_{x+\alpha y}^{(2)}$
        \item
        with configuration as in Figure \ref{Conf1A}.
    \end{itemize} 
    This is case \hyperref[Tab1A]{$1A$} and we see that we get a $1$-dimensional family of such surfaces $X_{1A,\alpha}$ depending on the parameter $\alpha$.     
        \end{minipage} \hspace{2mm} \begin{minipage}{0.3\textwidth} \begin{center} \includegraphics[width=0.95\textwidth]{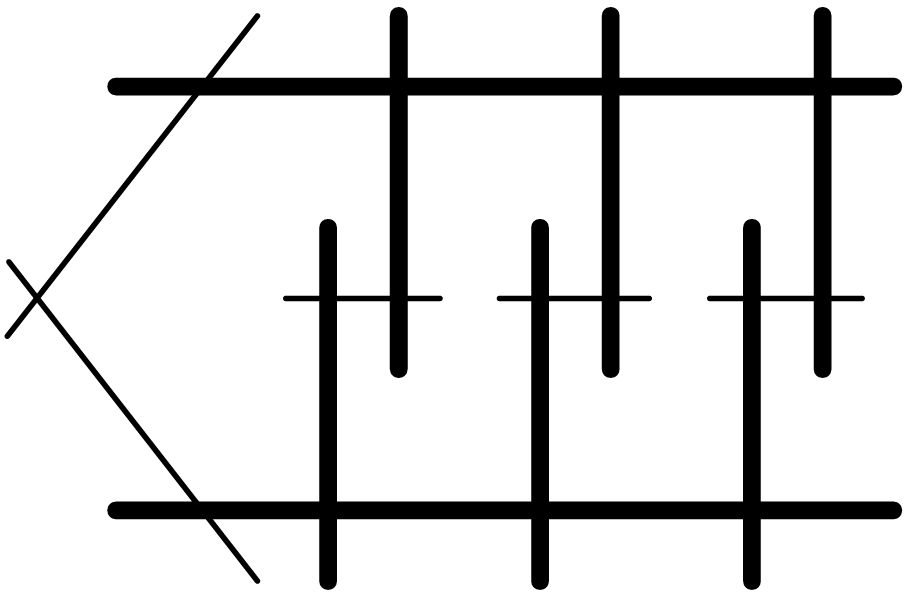} \vspace{-3.5mm}\captionof{figure}{} \label{Conf1A}  \end{center} \end{minipage}

    \noindent     \begin{minipage}{0.65\textwidth}
    \item \vspace{2mm}
    $
p_{1,1}= E_{1,0}^{} \cap \ell_{y}^{(1)}, 
p_{2,1}= E_{2,0}^{} \cap \ell_{x}^{(1)},
p_{3,1}= E_{3,0}^{} \cap \ell_{x+ \alpha y}^{(1)}$ \\with $\alpha \not\in \{0,-1\}$ 
    \begin{itemize}[leftmargin=20pt]
        \item
        $
        \Aut_{X'}^0(R) =    
        \left\{  \left( \begin{smallmatrix}
1 &  &  \\
 & 1 &  \\
 &  & i
\end{smallmatrix} \right) 
\in \PGL_3(R) \right\}
        $
        \item $(-2)$-curves: $E_{1,0}^{(2)},E_{2,0}^{(2)},E_{3,0}^{(2)}, \ell_{x}^{(2)}, \ell_{y}^{(2)}, \ell_{z}^{(2)}$
        \item $(-1)$-curves: $E_{1,1}^{}, E_{2,1}^{}, E_{3,1}^{}, E_{4,0}^{(2)},
        \ell_{x-y}^{(2)}, \ell_{x+ \alpha y}^{(2)}$
        \item
        with configuration as in Figure \ref{Conf2A}.
    \end{itemize} 
    This is case \hyperref[Tab2A]{$2A$} and we see that we get a $1$-dimensional family of such surfaces $X_{2A,\alpha}$ depending on the parameter $\alpha$.    
        \end{minipage} \hspace{2mm} \begin{minipage}{0.3\textwidth} \begin{center} \includegraphics[width=0.95\textwidth]{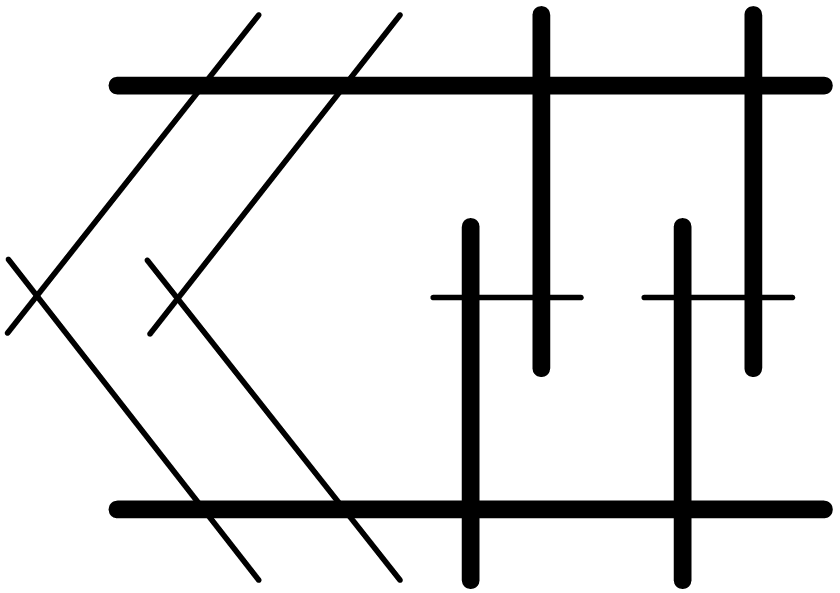} \vspace{-3.5mm}\captionof{figure}{}  \label{Conf2A} \end{center} \end{minipage} 
    
    \noindent     \begin{minipage}{0.65\textwidth}
    \item \vspace{2mm}
    $p_{1,1}= E_{1,0}^{} \cap \ell_{y}^{(1)}, 
        p_{2,1}= E_{2,0}^{} \cap \ell_{x}^{(1)}, p_{3,1}= E_{3,0}^{} \cap \ell_{x-y}^{(1)}$ 
    \begin{itemize}[leftmargin=20pt]
        \item
        $
        \Aut_{X'}^0(R) =    
        \left\{  \left( \begin{smallmatrix}
1 &  &  \\
 & 1 &  \\
 &  & i
\end{smallmatrix} \right) 
\in \PGL_3(R) \right\}
        $
        \item $(-2)$-curves: $E_{1,0}^{(2)},E_{2,0}^{(2)},E_{3,0}^{(2)}, \ell_{x}^{(2)}, \ell_{y}^{(2)}, \ell_{z}^{(2)},
        \ell_{x-y}^{(2)}$
        \item $(-1)$-curves: $E_{1,1}^{}, E_{2,1}^{}, E_{3,1}^{}, E_{4,0}^{(2)}$
        \item
        with configuration as in Figure \ref{Conf2D}.
    \end{itemize} 
    This is case \hyperref[Tab2D]{$2D$}.  
        \end{minipage} \hspace{2mm} \begin{minipage}{0.3\textwidth} \begin{center} \includegraphics[width=0.7\textwidth]{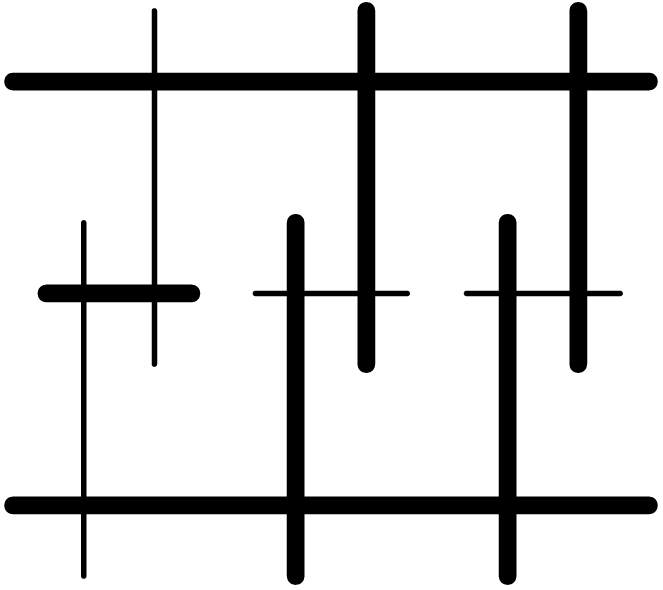} \vspace{-3.5mm}\captionof{figure}{} \label{Conf2D}  \end{center} \end{minipage} 
    
    \noindent     \begin{minipage}{0.65\textwidth}
    \item \vspace{2mm}
$p_{1,1}= E_{1,0}^{} \cap \ell_{y}^{(1)}, 
p_{2,1}= E_{2,0}^{} \cap \ell_{x}^{(1)}, 
p_{4,1}= E_{4,0}^{} \cap \ell_{x-y}^{(1)}$ 
    \begin{itemize}[leftmargin=20pt]
        \item
        $
        \Aut_{X'}^0(R) =    
        \left\{  \left( \begin{smallmatrix}
1 &  &  \\
 & 1 &  \\
 &  & i
\end{smallmatrix} \right) 
\in \PGL_3(R) \right\}
        $
        \item $(-2)$-curves: $E_{1,0}^{(2)},E_{2,0}^{(2)},E_{4,0}^{(2)}, \ell_{x}^{(2)}, \ell_{y}^{(2)}, \ell_{z}^{(2)},
        \ell_{x-y}^{(2)}$
        \item $(-1)$-curves: $E_{1,1}^{}, E_{2,1}^{}, E_{4,1}^{}, E_{3,0}^{(2)}$
        \item
        with configuration as in Figure \ref{Conf2E}.
    \end{itemize} 
    This is case \hyperref[Tab2E]{$2E$}.     
        \end{minipage} \hspace{2mm} \begin{minipage}{0.3\textwidth} \begin{center} \includegraphics[width=0.7\textwidth]{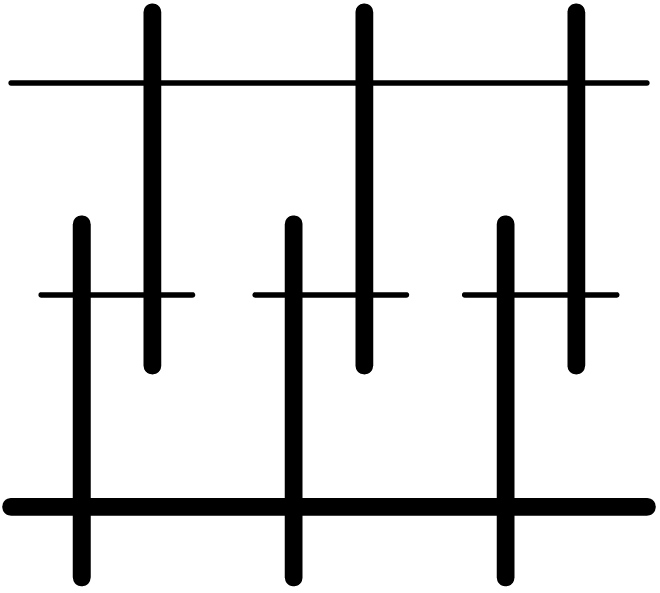} \vspace{-3.5mm}\captionof{figure}{}  \label{Conf2E} \end{center} \end{minipage}

    \noindent     \begin{minipage}{0.65\textwidth}
    \item \vspace{2mm}
    $p_{1,1}= E_{1,0}^{} \cap \ell_{y}^{(1)},  
    p_{3,1}= E_{3,0}^{} \cap \ell_{x+ \alpha y}^{(1)}$ with $\alpha \not\in \{0,-1\}$ 
    \begin{itemize}[leftmargin=20pt]
        \item
        $
        \Aut_{X'}^0(R) =    
        \left\{  \left( \begin{smallmatrix}
1 &  &  \\
 & 1 &  \\
 &  & i
\end{smallmatrix} \right) 
\in \PGL_3(R) \right\}
        $
        \item $(-2)$-curves: $E_{1,0}^{(2)},E_{3,0}^{(2)}, \ell_{y}^{(2)}, \ell_{z}^{(2)}$
        \item $(-1)$-curves: $E_{1,1}^{}, E_{3,1}^{}, E_{2,0}^{(2)}, E_{4,0}^{(2)}, \ell_{x}^{(2)}, \ell_{x-y}^{(2)}, \ell_{x+\alpha y}^{(2)}$
        \item
        with configuration as in Figure \ref{Conf3A}.
    \end{itemize} 
    This is case \hyperref[Tab3A]{$3A$} and we see that we get a $1$-dimensional family of such surfaces $X_{3A,\alpha}$ depending on the parameter $\alpha$.       
        \end{minipage} \hspace{2mm} \begin{minipage}{0.3\textwidth} \begin{center} \includegraphics[width=0.9\textwidth]{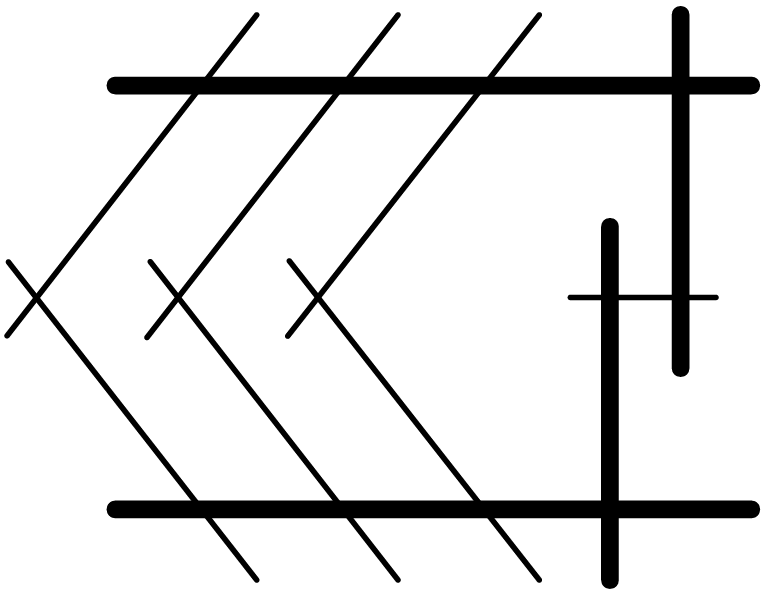} \vspace{-3.5mm}\captionof{figure}{} \label{Conf3A}  \end{center} \end{minipage}

    \noindent     \begin{minipage}{0.65\textwidth}
    \item \vspace{2mm}
    $p_{1,1}= E_{1,0}^{} \cap \ell_{y}^{(1)},  
    p_{3,1}= E_{3,0}^{} \cap \ell_{x}^{(1)}$
    \begin{itemize}[leftmargin=20pt]
        \item
        $
        \Aut_{X'}^0(R) =    
        \left\{  \left( \begin{smallmatrix}
1 &  &  \\
 & 1 &  \\
 &  & i
\end{smallmatrix} \right) 
\in \PGL_3(R) \right\}
        $
        \item $(-2)$-curves: $E_{1,0}^{(2)},E_{3,0}^{(2)}, \ell_{x}^{(2)},
        \ell_{y}^{(2)}, \ell_{z}^{(2)}$
        \item $(-1)$-curves: $E_{1,1}^{}, E_{3,1}^{}, E_{2,0}^{(2)}, E_{4,0}^{(2)},  \ell_{x-y}^{(2)}$
        \item
        with configuration as in Figure \ref{Conf3C}.
    \end{itemize} 
    This is case \hyperref[Tab3C]{$3C$}.       
        \end{minipage} \hspace{2mm} \begin{minipage}{0.3\textwidth} \begin{center} \includegraphics[width=0.73\textwidth]{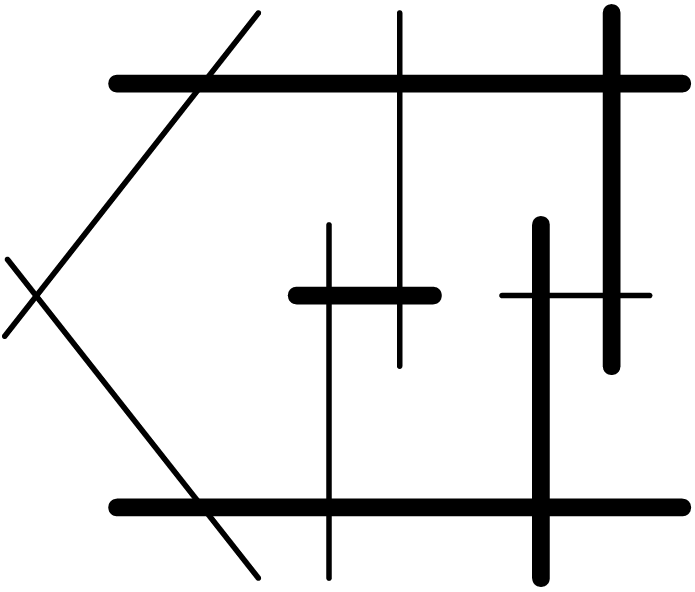} \vspace{-3.5mm}\captionof{figure}{}  \label{Conf3C} \end{center} \end{minipage}     
        
    \noindent     \begin{minipage}{0.65\textwidth}
    \item \vspace{2mm}
    $p_{1,1}= E_{1,0}^{} \cap \ell_{y}^{(1)}, 
    p_{2,1}= E_{2,0}^{} \cap \ell_{x}^{(1)}$
    \begin{itemize}[leftmargin=20pt]
        \item
        $
        \Aut_{X'}^0(R) =    
        \left\{  \left( \begin{smallmatrix}
1 &  &  \\
 & 1 &  \\
 &  & i
\end{smallmatrix} \right) 
\in \PGL_3(R) \right\}
        $
        \item $(-2)$-curves: $E_{1,0}^{(2)},E_{2,0}^{(2)}, \ell_{x}^{(2)},
        \ell_{y}^{(2)}, \ell_{z}^{(2)}$
        \item $(-1)$-curves: $E_{1,1}^{}, E_{2,1}^{}, E_{3,0}^{(2)}, E_{4,0}^{(2)},  \ell_{x-y}^{(2)}$
        \item
        with configuration as in Figure \ref{Conf3D}.
    \end{itemize} 
    This is case \hyperref[Tab3D]{$3D$}.
       \end{minipage} \hspace{2mm} \begin{minipage}{0.3\textwidth} \begin{center} \includegraphics[width=0.73\textwidth]{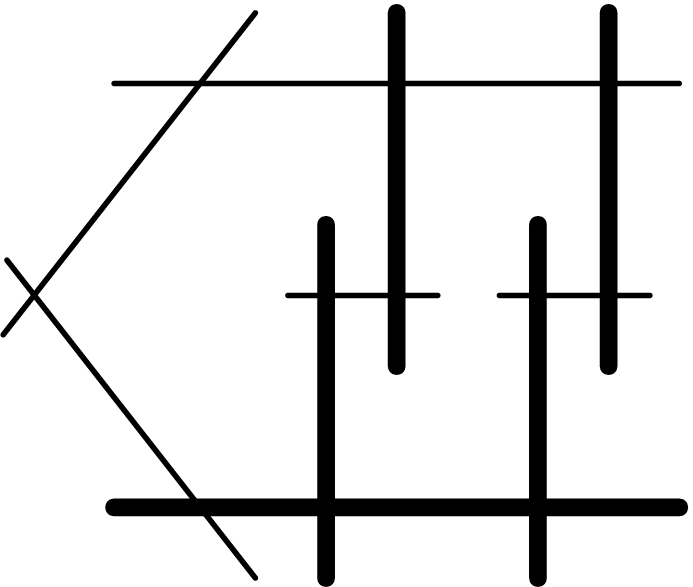} \vspace{-3.5mm}\captionof{figure}{}  \label{Conf3D} \end{center} \end{minipage}

    \noindent     \begin{minipage}{0.65\textwidth}
    \item \vspace{2mm}
    $p_{3,1}= E_{3,0}^{} \cap \ell_{x+\alpha y}^{(1)}$ with $\alpha \not\in \{0,-1\}$ 
    \begin{itemize}[leftmargin=20pt]
        \item
        $
        \Aut_{X'}^0(R) =    
        \left\{  \left( \begin{smallmatrix}
1 &  &  \\
 & 1 &  \\
 &  & i
\end{smallmatrix} \right) 
\in \PGL_3(R) \right\}
        $
        \item $(-2)$-curves: $E_{3,0}^{(2)}, \ell_{z}^{(2)}$
        \item $(-1)$-curves: $E_{3,1}^{}, E_{1,0}^{(2)},
        E_{2,0}^{(2)},
        E_{4,0}^{(2)},    
        \ell_{x}^{(2)},
        \ell_{y}^{(2)},
        \ell_{x-y}^{(2)},
        \ell_{x+ \alpha y}^{(2)}$
        \item
        with configuration as in Figure \ref{Conf4A}.
    \end{itemize} 
    This is case \hyperref[Tab4A]{$4A$} and we see that we get a $1$-dimensional family of such surfaces $X_{4A,\alpha}$ depending on the parameter $\alpha$.      
        \end{minipage} \hspace{2mm} \begin{minipage}{0.3\textwidth} \begin{center} \includegraphics[width=0.85\textwidth]{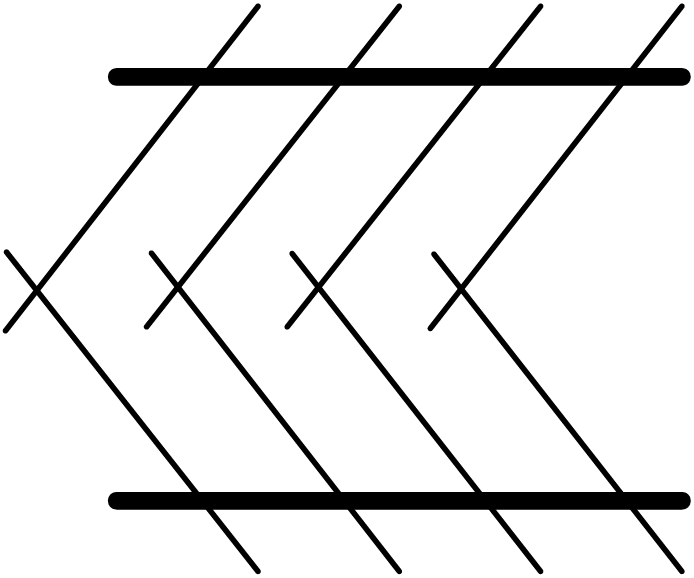} \vspace{-3.5mm}\captionof{figure}{}  \label{Conf4A} \end{center} \end{minipage}     
        
    \noindent     \begin{minipage}{0.65\textwidth}
    \item \vspace{2mm}
    $p_{3,1}= E_{3,0}^{} \cap \ell_{y}^{(1)}$ 
    \begin{itemize}[leftmargin=20pt]
        \item
        $
        \Aut_{X'}^0(R) =    
        \left\{  \left( \begin{smallmatrix}
1 &  &  \\
 & 1 &  \\
 &  & i
\end{smallmatrix} \right) 
\in \PGL_3(R) \right\}
        $
        \item $(-2)$-curves: $E_{3,0}^{(2)},
        \ell_{y}^{(2)},
        \ell_{z}^{(2)}$
        \item $(-1)$-curves: $E_{3,1}^{}, E_{1,0}^{(2)},
        E_{2,0}^{(2)},
        E_{4,0}^{(2)},    
        \ell_{x}^{(2)},
        \ell_{x-y}^{(2)}$
        \item
        with configuration as in Figure \ref{Conf4B}.
    \end{itemize} 
    This is case \hyperref[Tab4B]{$4B$}.
    \end{minipage} \hspace{2mm} \begin{minipage}{0.3\textwidth} \begin{center} \includegraphics[width=0.69\textwidth]{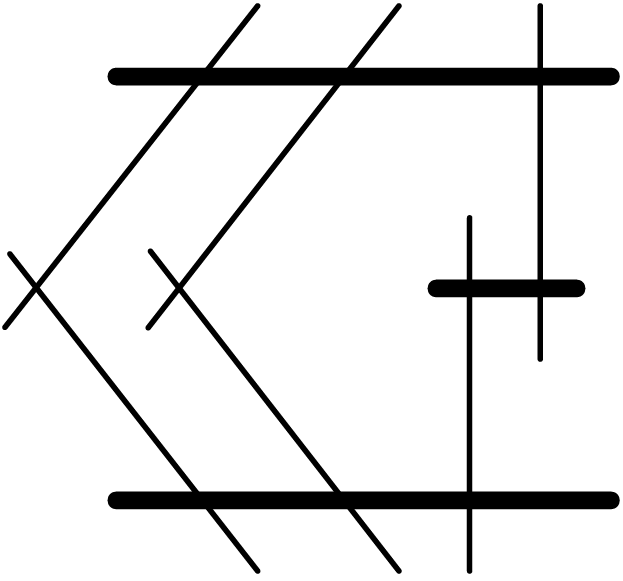} \vspace{-3.5mm}\captionof{figure}{} \label{Conf4B}  \end{center} \end{minipage} 

    \noindent     \begin{minipage}{0.65\textwidth}
    \item \vspace{2mm}
    $p_{1,1}= E_{1,0}^{} \cap \ell_{y}^{(1)}$ 
    \begin{itemize}[leftmargin=20pt]
        \item
        $
        \Aut_{X'}^0(R) =    
        \left\{  \left( \begin{smallmatrix}
1 &  &  \\
 & 1 &  \\
 &  & i
\end{smallmatrix} \right) 
\in \PGL_3(R) \right\}
        $
        \item $(-2)$-curves: $E_{1,0}^{(2)},
        \ell_{y}^{(2)},
        \ell_{z}^{(2)}$
        \item $(-1)$-curves: $E_{1,1}^{}, E_{2,0}^{(2)},
        E_{3,0}^{(2)},
        E_{4,0}^{(2)},    
        \ell_{x}^{(2)},
        \ell_{x-y}^{(2)}$
        \item
        with configuration as in Figure \ref{Conf4C}.
    \end{itemize} 
    This is case \hyperref[Tab4C]{$4C$}.
        \end{minipage} \hspace{2mm} \begin{minipage}{0.3\textwidth} \begin{center} \includegraphics[width=0.69\textwidth]{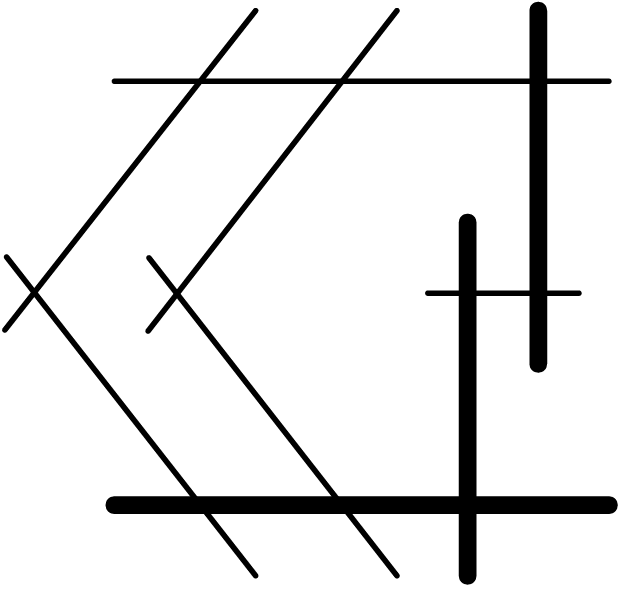} \vspace{-3.5mm}\captionof{figure}{} \label{Conf4C}  \end{center} \end{minipage} 

\end{enumerate}

 \subsubsection*{\underline{Case \hyperref[Tab6C]{$6C$}}}
We have $E= (\bigcup_{j=1}^{3} E_{j,0}^{})- \ell_{z}^{(1)}$ and $\Aut_X^0(R) = \left\{  \left( \begin{smallmatrix}
1 &  & c \\
 & 1 & f \\
 &  & i
\end{smallmatrix} \right) 
\in \PGL_3(R) \right\}$. 

\begin{itemize}[leftmargin=25pt]
    \item[-] $\lambda y + \mu z$ is $E_{1,0}^{}$-adapted and $\Aut_X^0(R)$ acts as 
    $[\lambda:\mu] \mapsto [\lambda: i\mu + f \lambda]$
    \item[-] $\lambda x + \mu z$ is $E_{2,0}^{}$-adapted and $\Aut_X^0(R)$ acts as $[\lambda:\mu] \mapsto [\lambda: i\mu + c \lambda]$
    \item[-] $\lambda (x-y) + \mu z$ is $E_{3,0}^{}$-adapted and $\Aut_X^0(R)$ acts as $[\lambda:\mu] \mapsto [\lambda: i\mu + (c-f) \lambda]$
\end{itemize}

\noindent
Since $p_{1,0},p_{2,0}$ and $p_{3,0}$ can be interchanged by automorphisms of $\bbP^2$ and the action of $\Aut_X^0$ is transitive on every $E \cap E_{i,0}$, we have the following three possibilities for $p_{1,1},\hdots,p_{n,1}$:

\begin{enumerate}[leftmargin=*]
\noindent     \begin{minipage}{0.65\textwidth}
    \item \vspace{2mm}
$p_{1,1}= E_{1,0}^{} \cap \ell_{y}^{(1)}, 
p_{2,1}= E_{2,0}^{} \cap \ell_{x}^{(1)},
p_{3,1}= E_{3,0}^{} \cap \ell_{x-y}^{(1)}$
    \begin{itemize}[leftmargin=20pt]
        \item
        $
        \Aut_{X'}^0(R) =    
        \left\{  \left( \begin{smallmatrix}
1 &  &  \\
 & 1 &  \\
 &  & i
\end{smallmatrix} \right) 
\in \PGL_3(R) \right\}
        $
        \item $(-2)$-curves: $E_{1,0}^{(2)},
        E_{2,0}^{(2)},
        E_{3,0}^{(2)}, 
        \ell_{z}^{(2)}$
        \item $(-1)$-curves: $E_{1,1}^{}, E_{2,1}^{},
        E_{3,1}^{},
        \ell_{x}^{(2)},
        \ell_{y}^{(2)},
        \ell_{x-y}^{(2)}$
        \item
        with configuration as in Figure \ref{Conf3B}.
    \end{itemize} 
    This is case \hyperref[Tab3B]{$3B$}.
    \end{minipage} \hspace{2mm} \begin{minipage}{0.3\textwidth} \begin{center} \includegraphics[width=0.9\textwidth]{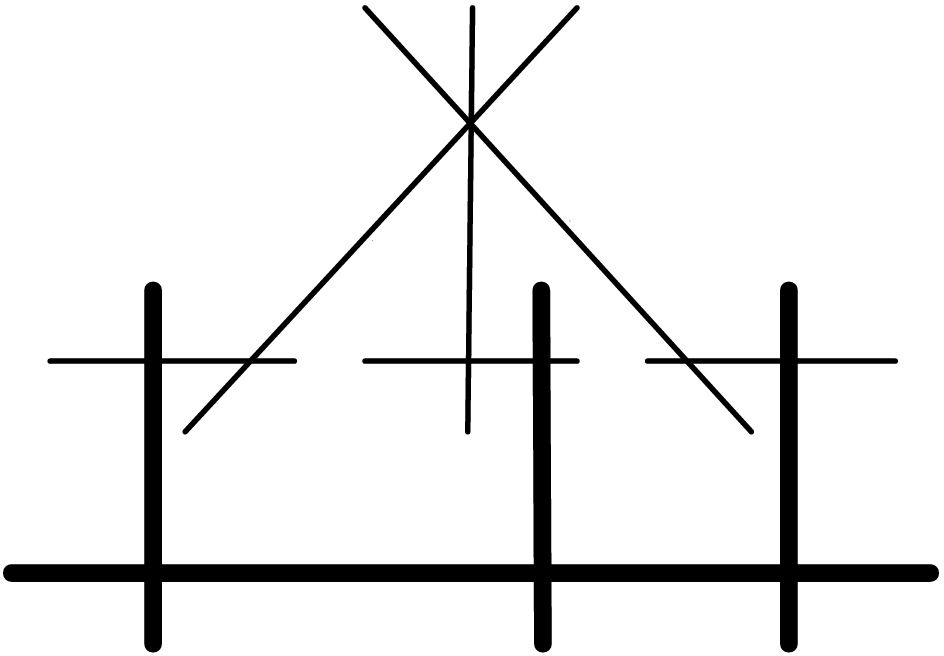} \vspace{-3.5mm}\captionof{figure}{} \label{Conf3B}  \end{center} \end{minipage}

    \noindent     \begin{minipage}{0.65\textwidth}
    \item \vspace{2mm}
$p_{1,1}= E_{1,0}^{} \cap \ell_{y}^{(1)}, 
p_{2,1}= E_{2,0}^{} \cap \ell_{x}^{(1)}$ 
    \begin{itemize}[leftmargin=20pt]
        \item
        $
        \Aut_{X'}^0(R) =    
        \left\{  \left( \begin{smallmatrix}
1 &  &  \\
 & 1 &  \\
 &  & i
\end{smallmatrix} \right) 
\in \PGL_3(R) \right\}
        $
        \item $(-2)$-curves: $E_{1,0}^{(2)},
        E_{2,0}^{(2)},
        \ell_{z}^{(2)}$
        \item $(-1)$-curves: $E_{1,1}^{}, E_{2,1}^{},
        E_{3,0}^{(2)},
        \ell_{x}^{(2)},
        \ell_{y}^{(2)}$
        \item
        with configuration as in Figure \ref{Conf4D}.
    \end{itemize} 
    This is case \hyperref[Tab4D]{$4D$}.    
        \end{minipage} \hspace{2mm} \begin{minipage}{0.3\textwidth} \begin{center} \includegraphics[width=0.76\textwidth]{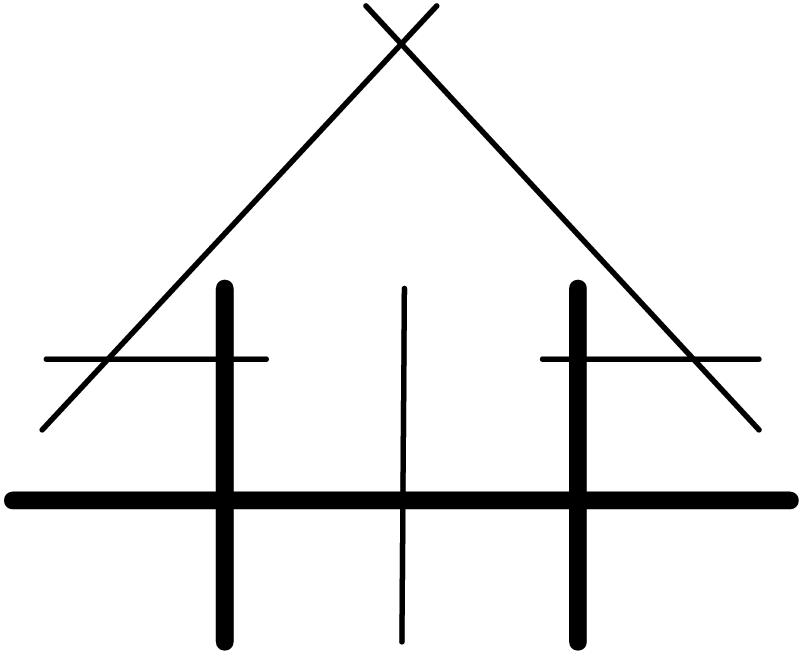} \vspace{-3.5mm}\captionof{figure}{}  \label{Conf4D}  \end{center} \end{minipage}

    \noindent     \begin{minipage}{0.65\textwidth}
    \item \vspace{2mm}
    $p_{1,1}= E_{1,0}^{} \cap \ell_{y}^{(1)}$
    \begin{itemize}[leftmargin=20pt]
        \item
        $
        \Aut_{X'}^0(R) =    
        \left\{  \left( \begin{smallmatrix}
1 &  & c \\
 & 1 &  \\
 &  & i
\end{smallmatrix} \right) 
\in \PGL_3(R) \right\}
        $
        \item $(-2)$-curves: $E_{1,0}^{(2)},
        \ell_{z}^{(2)}$
        \item $(-1)$-curves: $E_{1,1}^{}, E_{2,0}^{(2)},
        E_{3,0}^{(2)},
        \ell_{y}^{(2)}$
        \item
        with configuration as in Figure \ref{Conf5C}.
    \end{itemize} 
    This is case \hyperref[Tab5C]{$5C$}.

    \end{minipage} \hspace{2mm} \begin{minipage}{0.3\textwidth} \begin{center} \includegraphics[width=0.65\textwidth]{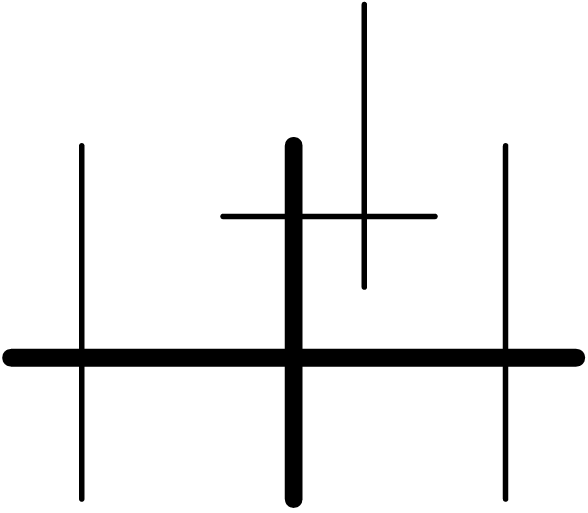} \vspace{-3.5mm}\captionof{figure}{}  \label{Conf5C} \end{center} \end{minipage} 
\end{enumerate}

 \subsubsection*{\underline{Case \hyperref[Tab6A]{$6A$}}}
We have $E= \bigcup_{j=0}^3 E_{j,0}$ and $\Aut_X^0(R) = \left\{  \left( \begin{smallmatrix}
1 &  &  \\
 & e &  \\
 &  & i
\end{smallmatrix} \right) 
\in \PGL_3(R) \right\}$.

\begin{itemize}[leftmargin=25pt]
    \item[-] $\lambda y + \mu z$ is $E_{1,0}^{}$-adapted and $\Aut_X^0(R)$ acts as $[\lambda:\mu] \mapsto [e \lambda: i \mu]$
    \item[-] $\lambda x + \mu z$ is $E_{2,0}^{}$-adapted and $\Aut_X^0(R)$ acts as $[\lambda:\mu] \mapsto [\lambda: i\mu]$
    \item[-] $\lambda x + \mu y$ is $E_{3,0}^{}$-adapted and $\Aut_X^0(R)$ acts as $[\lambda:\mu] \mapsto [\lambda: e\mu]$
\end{itemize}

\noindent
Since $p_{1,0},p_{2,0}$ and $p_{3,0}$ can be permuted arbitrarily by automorphisms of $\bbP^2$, we have the following nine possibilities for $p_{1,1},\hdots,p_{n,1}$:

\begin{enumerate}[leftmargin=*]
    \item \vspace{2mm}
    $p_{1,1}  =E_{1,0}^{} \cap \ell_{y-z}^{(1)}, 
    p_{2,1} = E_{2,0}^{} \cap \ell_{z}^{(1)}, 
    p_{3,1} = E_{3,0}^{} \cap \ell_{x}^{(1)}$
    \begin{itemize}[leftmargin=20pt]
    \noindent \begin{minipage}{0.5\textwidth}
        \item
        $
        \Aut_{X'}^0(R) =    
        \left\{  \left( \begin{smallmatrix}
1 &  &  \\
 & e &  \\
 &  & e
\end{smallmatrix} \right) 
\in \PGL_3(R) \right\}
        $
        \item $(-2)$-curves: $E_{1,0}^{(2)}, E_{2,0}^{(2)}, E_{3,0}^{(2)}, \ell_{x}^{(2)}, \ell_z^{(2)}$
        
        \end{minipage} \noindent \begin{minipage}{0.5\textwidth}
        
        \item $(-1)$-curves: $E_{1,1}^{}, E_{2,1}^{}, E_{3,1}^{}, \ell_y^{(2)}, \ell_{y-z}^{(2)}$
        \item
        with configuration as in Figure \ref{Conf3C}, that is, as in \\ case \hyperref[Tab3C]{$3C$}.
    \end{minipage}
    \end{itemize}
    \vspace{1mm} \noindent As explained in Remark \ref{R IsomorphismCheck}, one can check that $X' \cong X_{3C}$.

        \item \vspace{2mm}
     $p_{1,1}  =E_{1,0}^{} \cap \ell_{y-z}^{(1)}, 
    p_{2,1} = E_{2,0}^{} \cap \ell_{z}^{(1)}, 
    p_{3,1} = E_{3,0}^{} \cap \ell_{y}^{(1)}$
    \begin{itemize}[leftmargin=20pt]
    \noindent \begin{minipage}{0.5\textwidth}
        \item
        $
        \Aut_{X'}^0(R) =    
        \left\{  \left( \begin{smallmatrix}
1 &  &  \\
 & e &  \\
 &  & e
\end{smallmatrix} \right) 
\in \PGL_3(R) \right\}
        $
        \item $(-2)$-curves: $E_{1,0}^{(2)}, E_{2,0}^{(2)}, E_{3,0}^{(2)}, \ell_{y}^{(2)}, \ell_z^{(2)}$
        
        \end{minipage} \noindent \begin{minipage}{0.5\textwidth}
        
        \item $(-1)$-curves: $E_{1,1}^{}, E_{2,1}^{}, E_{3,1}^{}, \ell_x^{(2)}, \ell_{y-z}^{(1)}$
    \item
        with configuration as in Figure \ref{Conf3D}, that is, as in \\ case \hyperref[Tab3D]{$3D$}.
    \end{minipage}
    \end{itemize}
    \vspace{1mm} \noindent As explained in Remark \ref{R IsomorphismCheck}, one can check that $X' \cong X_{3D}$.

    \noindent     \begin{minipage}{0.65\textwidth}
        \item \vspace{2mm}
    $p_{1,1}  =E_{1,0}^{} \cap \ell_{z}^{(1)}, 
    p_{2,1} = E_{2,0}^{} \cap \ell_{x}^{(1)}, 
    p_{3,1} = E_{3,0}^{} \cap \ell_{y}^{(1)}$
    \begin{itemize}[leftmargin=20pt]
        \item
        $
        \Aut_{X'}^0(R) =    
        \left\{  \left( \begin{smallmatrix}
1 &  &  \\
 & e &  \\
 &  & i
\end{smallmatrix} \right) 
\in \PGL_3(R) \right\}
        $
        \item $(-2)$-curves: $E_{1,0}^{(2)}, E_{2,0}^{(2)}, E_{3,0}^{(2)}, \ell_{x}^{(2)}, \ell_{y}^{(2)}, \ell_z^{(2)}$
        \item $(-1)$-curves: $E_{1,1}^{}, E_{2,1}^{}, E_{3,1}^{}$
        \item
        with configuration as in Figure \ref{Conf3H}.
    \end{itemize}
    This is case \hyperref[Tab3H]{$3H$}.
        \end{minipage} \hspace{2mm} \begin{minipage}{0.3\textwidth} \begin{center} \includegraphics[width=0.68\textwidth]{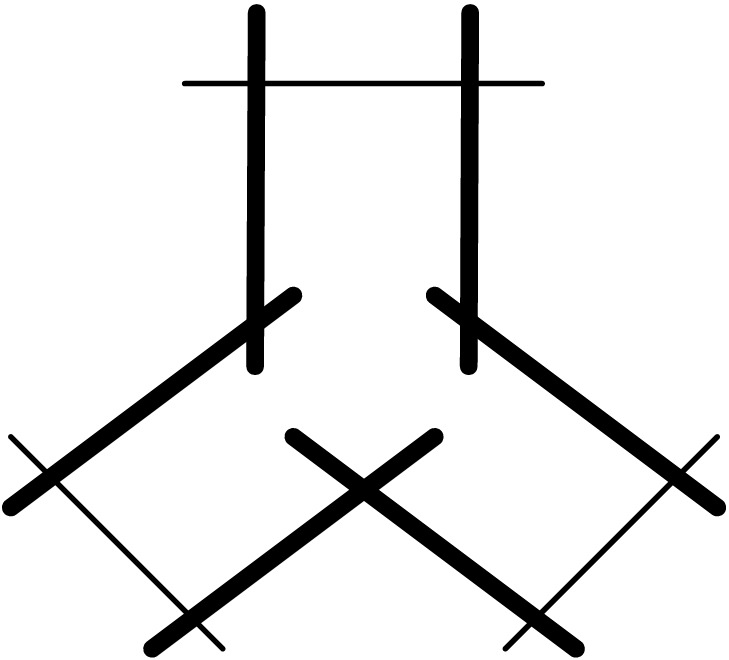} \vspace{-3.5mm}\captionof{figure}{}  \label{Conf3H} \end{center} \end{minipage}

        \item \vspace{2mm}
   $p_{1,1}  =E_{1,0}^{} \cap \ell_{y-z}^{(1)}, 
    p_{2,1} = E_{2,0}^{} \cap \ell_{z}^{(1)}$
    \begin{itemize}[leftmargin=20pt]
    \noindent \begin{minipage}{0.5\textwidth}
        \item
        $
        \Aut_{X'}^0(R) =    
        \left\{  \left( \begin{smallmatrix}
1 &  &  \\
 & e &  \\
 &  & e
\end{smallmatrix} \right) 
\in \PGL_3(R) \right\}
        $
        \item $(-2)$-curves: $E_{1,0}^{(2)}, E_{2,0}^{(2)}, \ell_z^{(2)}$
        
        \end{minipage} \noindent \begin{minipage}{0.5\textwidth}
        
        \item $(-1)$-curves: $E_{1,1}^{}, E_{2,1}^{}, E_{3,0}^{(2)}, \ell_x^{(2)}, \ell_y^{(2)}, \ell_{y-z}^{(2)}$
    \item
        with configuration as in Figure \ref{Conf4C}, that is, as in \\ case \hyperref[Tab4C]{$4C$}.
    \end{minipage}
    \end{itemize}
    \vspace{1mm} \noindent As explained in Remark \ref{R IsomorphismCheck}, one can check that $X' \cong X_{4C}$.

    \item \vspace{2mm}
    $p_{1,1}  =E_{1,0}^{} \cap \ell_{y-z}^{(1)}, 
    p_{2,1} = E_{2,0}^{} \cap \ell_{x}^{(1)}$
    \begin{itemize}[leftmargin=20pt]
    \noindent \begin{minipage}{0.5\textwidth}
        \item
        $
        \Aut_{X'}^0(R) =    
        \left\{  \left( \begin{smallmatrix}
1 &  &  \\
 & e &  \\
 &  & e
\end{smallmatrix} \right) 
\in \PGL_3(R) \right\}
        $
        \item $(-2)$-curves: $E_{1,0}^{(2)}, E_{2,0}^{(2)}, \ell_x^{(2)}$
        
        \end{minipage} \noindent \begin{minipage}{0.5\textwidth}
        
        \item $(-1)$-curves: $E_{1,1}^{}, E_{2,1}^{}, E_{3,0}^{(2)}, \ell_y^{(2)}, \ell_z^{(2)}, \ell_{y-z}^{(2)}$
    \item
        with configuration as in Figure \ref{Conf4B}, that is, as in \\ case \hyperref[Tab4B]{$4B$}.
    \end{minipage}
    \end{itemize}
    \vspace{1mm} \noindent As explained in Remark \ref{R IsomorphismCheck}, one can check that $X' \cong X_{4B}$.

    \noindent     \begin{minipage}{0.65\textwidth}
    \item \vspace{2mm}
     $p_{1,1}  =E_{1,0}^{} \cap \ell_{z}^{(1)}, 
    p_{2,1} = E_{2,0}^{} \cap \ell_{x}^{(1)}$
    \begin{itemize}[leftmargin=20pt]
        \item
        $
        \Aut_{X'}^0(R) =    
        \left\{  \left( \begin{smallmatrix}
1 &  &  \\
 & e &  \\
 &  & i
\end{smallmatrix} \right) 
\in \PGL_3(R) \right\}
        $
        \item $(-2)$-curves: $E_{1,0}^{(2)}, E_{2,0}^{(2)},  \ell_{x}^{(2)},  \ell_z^{(2)}$
        \item $(-1)$-curves: $E_{1,1}^{}, E_{2,1}^{}, E_{3,0}^{(2)}, \ell_y^{(2)}$
        \item
        with configuration as in Figure \ref{Conf4G}.
    \end{itemize} 
    This is case \hyperref[Tab4G]{$4G$}.
        \end{minipage} \hspace{2mm} \begin{minipage}{0.3\textwidth} \begin{center} \includegraphics[width=0.45\textwidth]{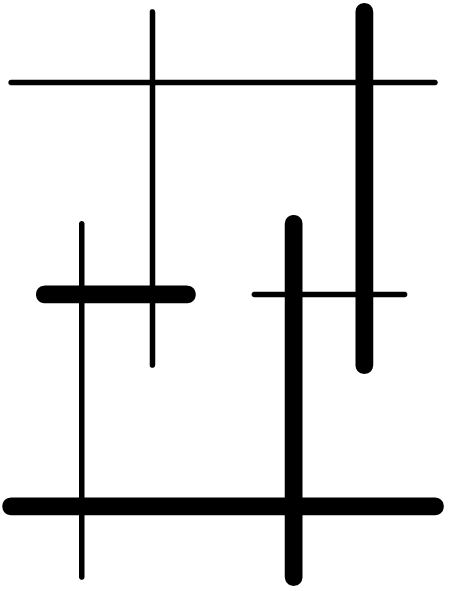} \vspace{-3.5mm}\captionof{figure}{}  \label{Conf4G} \end{center} \end{minipage}

      \noindent     \begin{minipage}{0.65\textwidth}
    \item \vspace{2mm}
     $p_{1,1}  =E_{1,0}^{} \cap \ell_{y}^{(1)}, 
    p_{2,1} = E_{2,0}^{} \cap \ell_{x}^{(1)}$
    \begin{itemize}[leftmargin=20pt]
        \item
        $
        \Aut_{X'}^0(R) =    
        \left\{  \left( \begin{smallmatrix}
1 &  &  \\
 & e &  \\
 &  & i
\end{smallmatrix} \right) 
\in \PGL_3(R) \right\}
        $
        \item $(-2)$-curves: $E_{1,0}^{(2)},E_{2,0}^{(2)}, \ell_x^{(2)}, \ell_y^{(2)}$
        \item $(-1)$-curves: $ E_{1,1}^{}, E_{2,1}^{}, E_{3,0}^{(2)}, \ell_{z}^{(2)}$
        \item
        with configuration as in Figure \ref{Conf4F}.
    \end{itemize} 
    This is case \hyperref[Tab4F]{$4F$}.
        \end{minipage} \hspace{2mm} \begin{minipage}{0.3\textwidth} \begin{center} \includegraphics[width=0.62\textwidth]{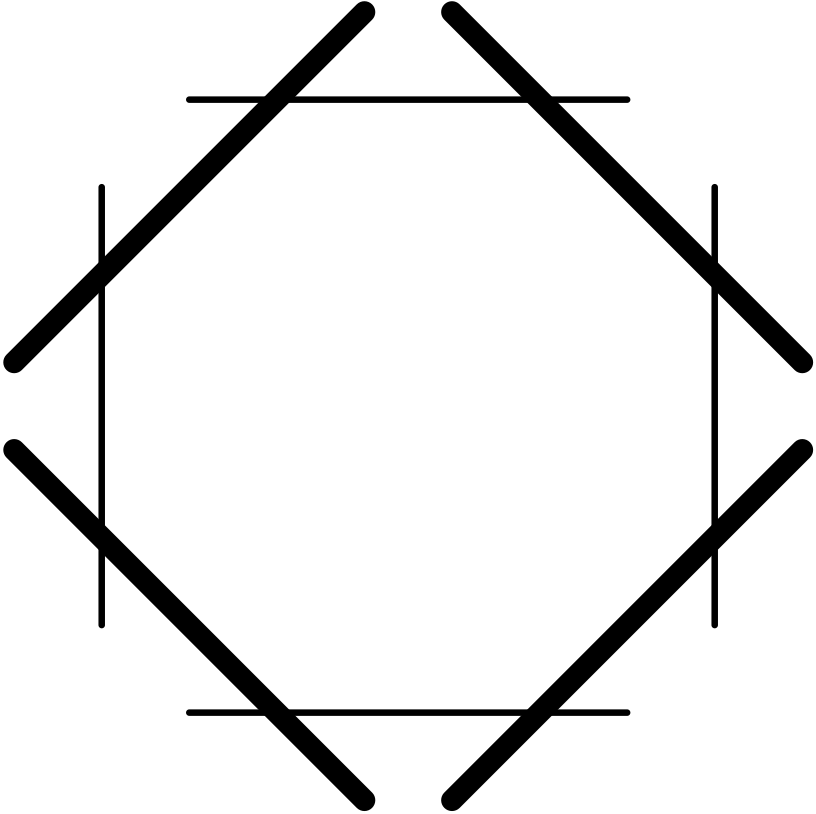} \vspace{-3.5mm}\captionof{figure}{}  \label{Conf4F} \end{center} \end{minipage}

    \item \vspace{2mm}
    $p_{1,1}  =E_{1,0}^{} \cap \ell_{y-z}^{(1)}$
    \begin{itemize}[leftmargin=20pt]
    \noindent \begin{minipage}{0.5\textwidth}
        \item
        $
        \Aut_{X'}^0(R) =    
        \left\{  \left( \begin{smallmatrix}
1 &  &  \\
 & e &  \\
 &  & e
\end{smallmatrix} \right) 
\in \PGL_3(R) \right\}
        $
        \item $(-2)$-curves: $E_{1,0}^{(2)}$
        \end{minipage} \noindent \begin{minipage}{0.5\textwidth}
        \item $(-1)$-curves: $E_{1,1}^{}, E_{2,0}^{(2)}, E_{3,0}^{(2)}, \ell_{x}^{(2)},  \ell_y^{(2)}, \ell_z^{(2)}, \ell_{y-z}^{(2)}$
     \item
        with configuration as in Figure \ref{Conf5A}, that is, as in \\ case \hyperref[Tab5A]{$5A$}.
    \end{minipage}
    \end{itemize}
    \vspace{1mm} \noindent As explained in Remark \ref{R IsomorphismCheck}, one can check that $X' \cong X_{5A}$.

    \noindent     \begin{minipage}{0.65\textwidth}
    \item \vspace{2mm}
    $p_{1,1}  =E_{1,0}^{} \cap \ell_{z}^{(1)}$
    \begin{itemize}[leftmargin=20pt]
        \item
        $
        \Aut_{X'}^0(R) =    
        \left\{  \left( \begin{smallmatrix}
1 &  &  \\
 & e &  \\
 &  & i
\end{smallmatrix} \right) 
\in \PGL_3(R) \right\}
        $
        \item $(-2)$-curves: $E_{1,0}^{(2)}, \ell_z^{(2)}$
        \item $(-1)$-curves: $E_{1,1}^{}, E_{2,0}^{(2)}, E_{3,0}^{(2)}, \ell_x^{(2)}, \ell_y^{(2)}$
        \item
        with configuration as in Figure \ref{Conf5B}.
    \end{itemize} 
    This is case \hyperref[Tab5B]{$5B$}.
    \end{minipage} \hspace{2mm} \begin{minipage}{0.3\textwidth} \begin{center} \includegraphics[width=0.45\textwidth]{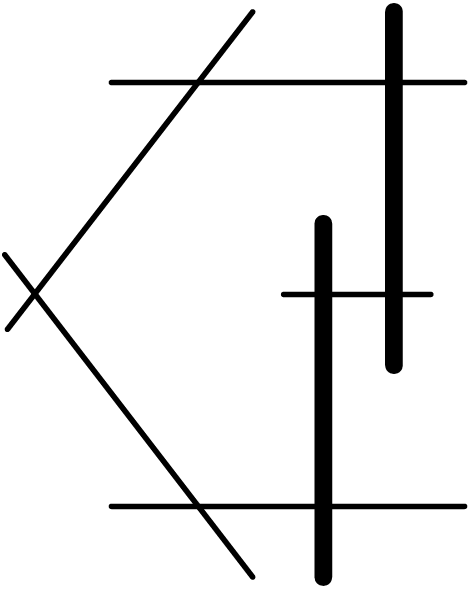} \vspace{-3.5mm}\captionof{figure}{} \label{Conf5B}  \end{center} \end{minipage} 
    
\end{enumerate}

\vspace{2mm} \subsubsection*{\underline{Case \hyperref[Tab7A]{$7A$}}}
We have $E = E_{1,0}^{} \cup E_{2,0}^{}$ and $\Aut_X^0(R) = \left\{  \left( \begin{smallmatrix}
1 &  & c \\
 & e & f \\
 &  & i
\end{smallmatrix} \right) 
\in \PGL_3(R) \right\}$.

\begin{itemize}[leftmargin=25pt]
    \item[-] $\lambda y + \mu z$ is $E_{1,0}^{}$-adapted and $\Aut_X^0(R)$ acts as $[\lambda:\mu] \mapsto [e \lambda: i\mu + f \lambda]$
    \item[-] $\lambda x + \mu z$ is $E_{2,0}^{}$-adapted and $\Aut_X^0(R)$ acts as $[\lambda:\mu] \mapsto [\lambda: i\mu + c \lambda]$
\end{itemize}

\noindent
Since $p_{1,0}$ and $p_{2,0}$ can be interchanged by an automorphism of $\bbP^2$, we have the following four possibilities for $p_{1,1},\hdots,p_{n,1}$:

\begin{enumerate}[leftmargin=*]

\item \vspace{2mm}
    $p_{1,1}  =E_{1,0}^{} \cap \ell_{y}^{(1)}, 
    p_{2,1} = E_{2,0}^{} \cap \ell_{x}^{(1)}$
    \begin{itemize}[leftmargin=20pt]
    \noindent \begin{minipage}{0.5\textwidth}
        \item
        $
        \Aut_{X'}^0(R) =    
        \left\{  \left( \begin{smallmatrix}
1 &  &  \\
 & e &  \\
 &  & i
\end{smallmatrix} \right) 
\in \PGL_3(R) \right\}
        $
        \item $(-2)$-curves: $E_{1,0}^{(2)}, E_{2,0}^{(2)}$
        
        \end{minipage} \noindent \begin{minipage}{0.5\textwidth}
        
        \item $(-1)$-curves: $E_{1,1}^{}, E_{2,1}^{},  \ell_x^{(2)}, \ell_y^{(2)}, \ell_z^{(2)}$
         \item
        with configuration as in Figure \ref{Conf5B}, that is, as in \\ case \hyperref[Tab5B]{$5B$}.
    \end{minipage}
    \end{itemize}
    \vspace{1mm} \noindent As explained in Remark \ref{R IsomorphismCheck}, one can check that $X' \cong X_{5B}$.

\noindent     \begin{minipage}{0.65\textwidth}
       \item \vspace{2mm}
    $p_{1,1}  =E_{1,0}^{} \cap \ell_{z}^{(1)}, 
    p_{2,1} = E_{2,0}^{} \cap \ell_{x}^{(1)}$
    \begin{itemize}[leftmargin=20pt]
        \item
        $
        \Aut_{X'}^0(R) =    
        \left\{  \left( \begin{smallmatrix}
1 &  &  \\
 & e & f \\
 &  & i
\end{smallmatrix} \right) 
\in \PGL_3(R) \right\}
        $
        \item $(-2)$-curves: $E_{1,0}^{(2)}, E_{2,0}^{(2)}, \ell_z^{(2)}$
        \item $(-1)$-curves: $E_{1,1}^{}, E_{2,1}^{},  \ell_x^{(2)}$
        \item
        with configuration as in Figure \ref{Conf5D}.
    \end{itemize} 
    This is case \hyperref[Tab5D]{$5D$}.
        \end{minipage} \hspace{2mm} \begin{minipage}{0.3\textwidth} \begin{center} \includegraphics[width=0.45\textwidth]{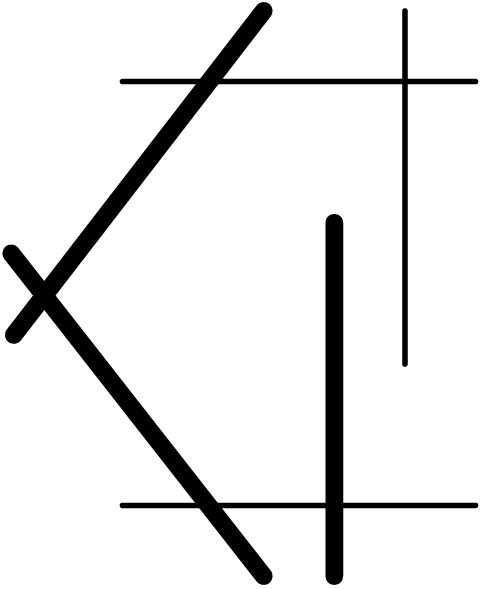} \vspace{-3.5mm}\captionof{figure}{} \label{Conf5D}  \end{center} \end{minipage} 
    
    \noindent     \begin{minipage}{0.65\textwidth}
    \item \vspace{2mm}
    $p_{1,1}  =E_{1,0}^{} \cap \ell_{y}^{(1)}$
    \begin{itemize}[leftmargin=20pt]
        \item
        $
        \Aut_{X'}^0(R) =    
        \left\{  \left( \begin{smallmatrix}
1 &  & c \\
 & e &  \\
 &  & i
\end{smallmatrix} \right) 
\in \PGL_3(R) \right\}
        $
        \item $(-2)$-curves: $E_{1,0}^{(2)}$
        \item $(-1)$-curves: $E_{1,1}^{}, E_{2,0}^{(2)},  \ell_y^{(2)}, \ell_z^{(2)}$
        \item
        with configuration as in Figure \ref{Conf6B}.
    \end{itemize} 
    This is case \hyperref[Tab6B]{$6B$}.
        \end{minipage} \hspace{2mm} \begin{minipage}{0.3\textwidth} \begin{center} \includegraphics[width=0.7\textwidth]{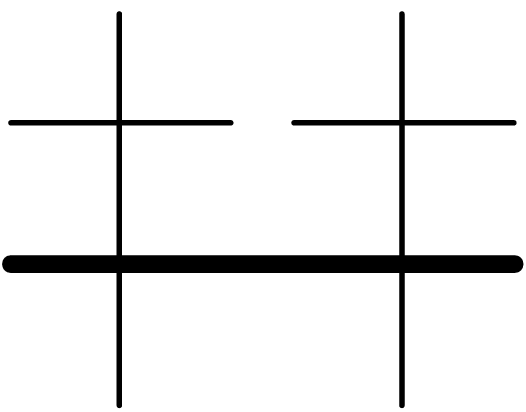} \vspace{-3.5mm}\captionof{figure}{} \label{Conf6B}  \end{center} \end{minipage}

    \noindent     \begin{minipage}{0.65\textwidth}
   \item \vspace{2mm}
    $p_{1,1}  =E_{1,0}^{} \cap \ell_{z}^{(1)}$
    \begin{itemize}[leftmargin=20pt]
        \item
        $
        \Aut_{X'}^0(R) =    
        \left\{  \left( \begin{smallmatrix}
1 &  & c \\
 & e & f \\
 &  & i
\end{smallmatrix} \right) 
\in \PGL_3(R) \right\}
        $
        \item $(-2)$-curves: $E_{1,0}^{(2)}, \ell_z^{(2)}$
        \item $(-1)$-curves: $E_{1,1}^{}, E_{2,0}^{(2)}$
        \item
        with configuration as in Figure \ref{Conf6D}.
    \end{itemize} 
    This is case \hyperref[Tab6D]{$6D$}.
        \end{minipage} \hspace{2mm} \begin{minipage}{0.3\textwidth} \begin{center} \includegraphics[width=0.9\textwidth]{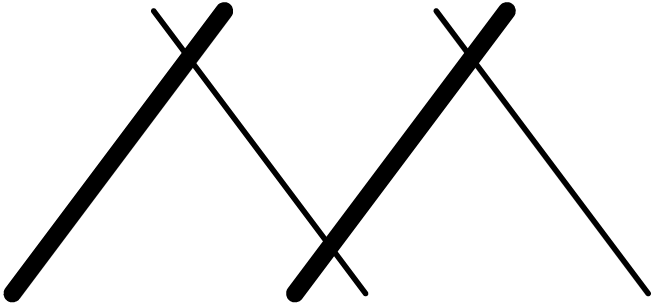} \vspace{-3.5mm}\captionof{figure}{}  \label{Conf6D} \end{center} \end{minipage} 

\end{enumerate}

\vspace{2mm} \subsubsection*{\underline{Case \hyperref[Tab8A]{$8A$}}} 
We have $E= E_{1,0}^{}$ and $\Aut_X^0(R) = \left\{  \left( \begin{smallmatrix}
1 & b & c \\
 & e & f \\
 & h & i
\end{smallmatrix} \right) 
\in \PGL_3(R) \right\}$.

\begin{itemize}[leftmargin=25pt]
    \item[-] $\lambda y + \mu z$ is $E_{1,0}^{}$-adapted and $\Aut_X^0(R)$ acts as $[\lambda:\mu] \mapsto [e \lambda + h \mu : i \mu + f \lambda]$
\end{itemize}

\noindent
Therefore, there is a unique possibility for $p_{1,1},\hdots,p_{n,1}$ up to isomorphism:

\begin{enumerate}[leftmargin=*]

\noindent     \begin{minipage}{0.65\textwidth}
    \item \vspace{2mm}
    $p_{1,1}  =E_{1,0}^{} \cap \ell_{z}^{(1)}$
    \begin{itemize}[leftmargin=20pt]
        \item
        $
        \Aut_{X'}^0(R) =    
        \left\{  \left( \begin{smallmatrix}
1 & b & c \\
 & e & f \\
 &  & i
\end{smallmatrix} \right) 
\in \PGL_3(R) \right\}
        $
        \item $(-2)$-curves: $E_{1,0}^{(2)}$
        \item $(-1)$-curves: $E_{1,1}^{}, \ell_z^{(2)}$
        \item
        with configuration as in Figure \ref{Conf7B}.
    \end{itemize} 
    This is case \hyperref[Tab7B]{$7B$}.
        \end{minipage} \hspace{2mm} \begin{minipage}{0.3\textwidth} \begin{center} \includegraphics[width=0.42\textwidth]{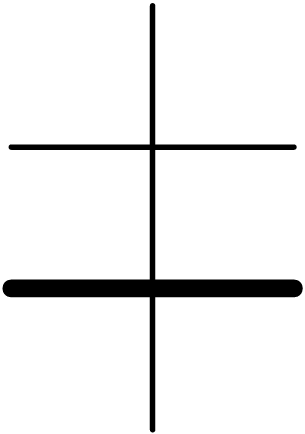} \vspace{-3.5mm}\captionof{figure}{}  \label{Conf7B} \end{center} \end{minipage} 
\end{enumerate}

\vspace{2mm} 
\noindent
Summarizing, we obtain
\begin{eqnarray*}
\cal{L}_2 &=& \{  X_{1A,\alpha}, X_{2A,\alpha}, X_{2D}, X_{2E}, X_{3A,\alpha},  X_{3C}, X_{3D}, X_{4A,\alpha}, X_{4B}, X_{4C},  \\ && X_{3B}, X_{4D}, X_{5C}, X_{3H}, X_{4G}, X_{4F}, X_{5B}, X_{5D}, X_{6B}, X_{6D}, X_{7B}\}.
\end{eqnarray*}

\subsection{Height 3}

\vspace{2mm} \subsubsection*{\underline{Case \hyperref[Tab2A]{$2A$}}}
We have $E = \bigcup_{j=1}^{3} E_{j,1}^{} - (\bigcup_{j=1}^{3} E_{j,0}^{(2)} \cup \ell_{x}^{(2)} \cup \ell_{y}^{(2)})$ and $\Aut_X^0(R) = \left\{  \left( \begin{smallmatrix}
1 &  &  \\
 & 1 &  \\
 &  & i
\end{smallmatrix} \right) 
\in \PGL_3(R) \right\}$.

\begin{itemize}[leftmargin=25pt]
    \item[-] $\lambda xy + \mu z^2$ is $E_{1,1}^{}$-adapted and $E_{2,1}^{}$-adapted and $\Aut_X^0(R)$ acts as $[\lambda:\mu] \mapsto [ \lambda: i^2 \mu]$
    \item[-] $\lambda y^2 + \mu (x + \alpha y)z$ is $E_{3,1}^{}$-adapted and $\Aut_X^0(R)$ acts as $[\lambda:\mu] \mapsto [\lambda: i\mu]$
\end{itemize}

\noindent 
Note that $X$ has degree $2$, therefore we are only allowed to blow up one more point $p_{j,2}$. Moreover, the involution $x \leftrightarrow \alpha y$ of $\bbP^2$ lifts to an involution of $X$ interchanging $E_{1,1}^{}$ and $E_{2,1}^{}$, thus we may assume without loss of generality that $j = 1$ or $j = 3$. Finally, if $j = 3$, then the stabilizer of $p_{3,2} \in E \cap E_{3,1}$ is trivial unless $p_{3,2}$ lies on the strict transform of $\ell_{x+\alpha y}$. Moreover, $\Aut_{X}^0$ acts transitively on $E \cap E_{1,1}$. Hence, we have the following two possibilities:

\begin{enumerate}[leftmargin=*]

        \item \vspace{2mm}
    $p_{3,2}  =E_{3,1}^{} \cap \ell_{x+ \alpha y}^{(2)}$ with $\alpha \not\in \{0, -1\}$
    \begin{itemize}[leftmargin=20pt]
    \noindent \begin{minipage}{0.5\textwidth}
        \item
        $
        \Aut_{X'}^0(R) =    
        \left\{  \left( \begin{smallmatrix}
1 &  &  \\
 & 1 &  \\
 &  & i
\end{smallmatrix} \right) 
\in \PGL_3(R) \right\}
        $
        \item $(-2)$-curves: $E_{1,0}^{(3)},E_{2,0}^{(3)},E_{3,0}^{(3)}, E_{3,1}^{(3)}, \ell_{x}^{(3)}, \ell_{y}^{(3)}, \\ 
        \ell_{z}^{(3)}, \ell_{x+ \alpha y}^{(3)}$
        
        \end{minipage} \noindent \begin{minipage}{0.5\textwidth}
        
        \item $(-1)$-curves: $E_{3,2}, E_{1,1}^{(3)}, E_{2,1}^{(3)},  E_{4,0}^{(3)},
        \ell_{x-y}^{(3)}$
        \item
        with configuration as in Figure \ref{Conf1A}, that is, as in \\ case \hyperref[Tab1A]{$1A$}.
    \end{minipage}
    \end{itemize} 
   \vspace{1mm} \noindent As explained in Remark \ref{R IsomorphismCheck}, one can check that $X' \cong X_{1A,\alpha'}$ for some $\alpha'$.
    
    \noindent     \begin{minipage}{0.65\textwidth}
         \item \vspace{2mm}
    $
p_{1,2}= E_{1,1}^{} \cap C_1^{(2)}$ with $C_1=\cal{V} (xy+z^2) $
    \begin{itemize}[leftmargin=20pt]
        \item
        $\Aut_{X'}^0(R) =
        \begin{cases}   
        \{\id\} & \text{ if } p \neq 2 \\
        \left\{  \left( \begin{smallmatrix}
1 &  &  \\
 & 1 &  \\
 &  & i
\end{smallmatrix} \right) 
\in \PGL_3(R) \bigg| i^2=1 \right\}
         & \text{ if } p=2
        \end{cases} $
        \\
        Hence, $X'$ has global vector fields only if $p=2$. Therefore, we assume $p=2$ when describing the configuration of negative curves.
        \item $(-2)$-curves: $E_{1,0}^{(3)},E_{2,0}^{(3)},E_{3,0}^{(3)}, E_{1,1}^{(3)}, \ell_{x}^{(3)}, \ell_{y}^{(3)}, \ell_{z}^{(3)}$
        \item $(-1)$-curves: $E_{1,2}^{}, E_{2,1}^{(3)}, E_{3,1}^{(3)}, E_{4,0}^{(3)},
        \ell_{x-y}^{(3)}, \ell_{x+ \alpha y}^{(3)}, C_1^{(3)}, C_2^{(3)}$ with $\alpha \not\in \{0, -1\}$ and $C_2=\cal{V}( x^3y+ xy^3+ x^2z^2+\alpha^2 y^2 z^2) $
        \item
        with configuration as in Figure \ref{Conf1L}.
    \end{itemize} 
    This is case \hyperref[Tab1L]{$1L$} and we see that we get a $1$-dimensional family of such surfaces $X_{1L,\alpha}$ depending on the parameter $\alpha$.  
        \end{minipage} \hspace{2mm} \begin{minipage}{0.3\textwidth} \begin{center} \includegraphics[width=0.98\textwidth]{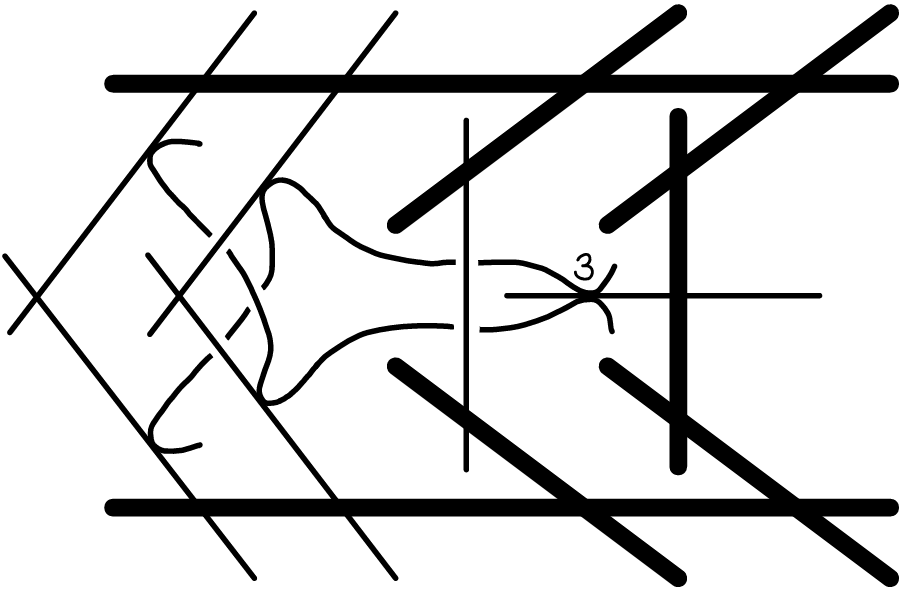} \vspace{-3.5mm}\captionof{figure}{}  \label{Conf1L} \end{center} \end{minipage}

\end{enumerate}

\vspace{2mm} \subsubsection*{\underline{Case \hyperref[Tab2D]{$2D$}}}
We have $E = \bigcup_{j=1}^{3} E_{j,1}^{} - (\bigcup_{j=1}^{3} E_{j,0}^{(2)} \cup \ell_{x}^{(2)} \cup \ell_{y}^{(2)} \cup \ell_{x-y}^{(2)})$ and $\Aut_X^0(R) = \left\{  \left( \begin{smallmatrix}
1 &  &  \\
 & 1 &  \\
 &  & i
\end{smallmatrix} \right) 
\in \PGL_3(R) \right\}$.

\begin{itemize}[leftmargin=25pt]
    \item[-] $\lambda xy + \mu z^2$ is $E_{1,1}^{}$-adapted and $E_{2,1}^{}$-adapted and $\Aut_X^0(R)$ acts as $[\lambda:\mu] \mapsto [\lambda: i^2 \mu]$
    \item[-] $\lambda y^2 + \mu (x - y)z$ is $E_{3,1}^{}$-adapted and $\Aut_X^0(R)$ acts as $[\lambda:\mu] \mapsto [\lambda: i \mu]$
\end{itemize}

\noindent
Note that $X$ has degree $2$, thus we are only allowed to blow up one more point $p_{j,2}$. Next, note that the stabilizer of every point on $E \cap E_{3,1}^{}$ is trivial, hence we may assume $j = 1$ or $j = 2$. Similar to Case $2A$, the involution $x \leftrightarrow y$ of $\bbP^2$ lifts to an involution of $X$ interchanging $E_{1,1}^{}$ and $E_{2,1}^{}$, thus we may assume without loss of generality that $j = 1$. Hence, there is the following unique choice for $p_{j,2}$ up to isomorphism:

\begin{enumerate}[leftmargin=*]
\noindent     \begin{minipage}{0.65\textwidth}
    \item \vspace{2mm}
    $p_{1,2}= E_{1,1}^{} \cap C^{(2)}$ with $C= \cal{V}(xy+z^2) $ 
    \begin{itemize}[leftmargin=20pt]
        \item
        $\Aut_{X'}^0(R) =
        \begin{cases}   
        \{\id\} & \text{ if } p \neq 2 \\
        \left\{  \left( \begin{smallmatrix}
1 &  &  \\
 & 1 &  \\
 &  & i
\end{smallmatrix} \right) 
\in \PGL_3(R) \bigg| i^2=1 \right\}
         & \text{ if } p=2
        \end{cases} $
        \\ Hence, $X'$ has global vector fields only if $p=2$. Therefore, we assume $p=2$ when describing the configuration of negative curves.
        \item $(-2)$-curves: $E_{1,0}^{(3)},E_{2,0}^{(3)},E_{3,0}^{(3)}, E_{1,1}^{(3)}, \ell_{x}^{(3)}, \ell_{y}^{(3)}, \ell_{z}^{(3)},
        \ell_{x-y}^{(3)}$
        \item $(-1)$-curves: $E_{1,2}^{},  E_{2,1}^{(3)}, E_{3,1}^{(3)}, E_{4,0}^{(3)}, C^{(3)}$
        \item
        with configuration as in Figure \ref{Conf1O}.
    \end{itemize} 
    This is case \hyperref[Tab1O]{$1O$}. 
        \end{minipage} \hspace{2mm} \begin{minipage}{0.3\textwidth} \begin{center} \includegraphics[width=0.98\textwidth]{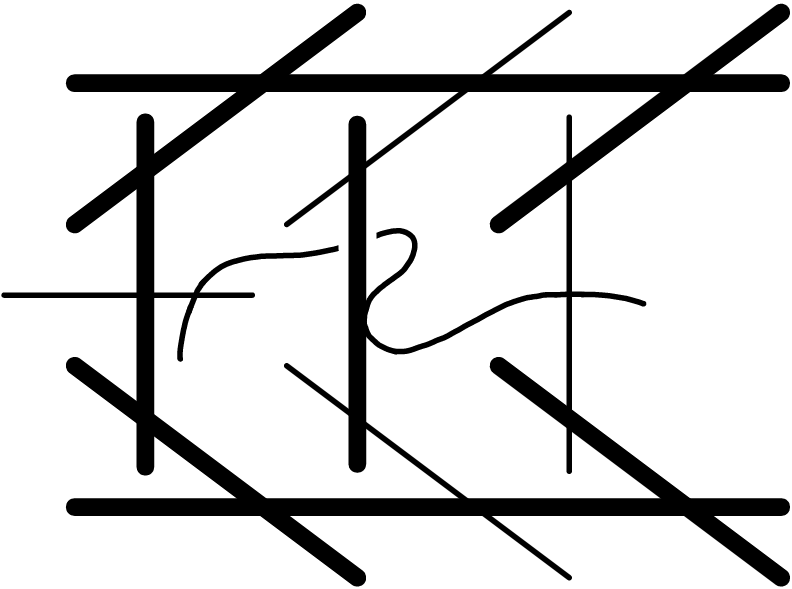} \vspace{-3.5mm}\captionof{figure}{}  \label{Conf1O} \end{center} \end{minipage} 
\end{enumerate}

\vspace{2mm} \subsubsection*{\underline{Case \hyperref[Tab2E]{$2E$}}}
We have 
$
E = (E_{1,1}^{} \cup E_{2,1}^{} \cup E_{4,1}^{}) - (E_{1,0}^{(2)} \cup E_{2,0}^{(2)} \cup E_{4,0}^{(2)} \cup \ell_x^{(2)} \cup \ell_y^{(2)} \cup \ell_{x-y}^{(2)}) 
$ and \\$\Aut_X^0(R) = \left\{  \left( \begin{smallmatrix}
1 &  &  \\
 & 1 &  \\
 &  & i
\end{smallmatrix} \right) 
\in \PGL_3(R) \right\}$.

\begin{itemize}[leftmargin=25pt]
    \item[-] $\lambda xy + \mu z^2$ is $E_{1,1}^{}$-adapted and $E_{2,1}^{}$-adapted and $\Aut_X^0(R)$ acts as $[\lambda:\mu] \mapsto [\lambda: i^2 \mu]$
    \item[-] $\lambda (x-y)x + \mu z^2$ is $E_{4,1}^{}$-adapted and $\Aut_X^0(R)$ acts as $[\lambda:\mu] \mapsto [\lambda: i^2 \mu]$
\end{itemize}

\noindent 
Note that $X$ has degree $2$, thus we are only allowed to blow up one more point $p_{j,2}$. Next, the automorphisms of $\bbP^2$ interchanging $p_{1,0},p_{2,0}$ and $p_{4,0}$ and preserving $p_{3,0}$ lift to $X$ and interchange $E_{1,1}^{},E_{2,1}^{}$ and $E_{4,1}^{}$, thus we may assume $j = 1$. Finally, $\Aut_X^0$ acts transitively on $E \cap E_{1,1}$, hence we have the following unique choice for $p_{j,2}$ up to isomorphism:

\begin{enumerate}[leftmargin=*]
\noindent     \begin{minipage}{0.65\textwidth}
    \item \vspace{2mm}
$p_{1,2}= E_{1,1}^{} \cap C_1^{(2)}$ with $C_1= \cal{V}(xy+z^2) $
    \begin{itemize}[leftmargin=20pt]
        \item
        $\Aut_{X'}^0(R) =
        \begin{cases}   
        \{\id\} & \text{ if } p \neq 2 \\
        \left\{  \left( \begin{smallmatrix}
1 &  &  \\
 & 1 &  \\
 &  & i
\end{smallmatrix} \right) 
\in \PGL_3(R) \bigg| i^2=1 \right\}
         & \text{ if } p=2
        \end{cases} $
        Hence, $X'$ has global vector fields only if $p=2$. Therefore, we assume $p=2$ when describing the configuration of negative curves.
        \item $(-2)$-curves: $E_{1,0}^{(3)},E_{2,0}^{(3)},E_{4,0}^{(3)}, E_{1,1}^{(3)}, \ell_{x}^{(3)}, \ell_{y}^{(3)}, \ell_{z}^{(3)},
        \ell_{x-y}^{(3)}$
        \item $(-1)$-curves: $E_{1,2}^{}, E_{2,1}^{(3)}, E_{4,1}^{(3)}, E_{3,0}^{(3)}, C_1^{(3)}, C_2^{(3)}$ with \\ $C_2= \cal{V} (xy+y^2+z^2)$
        \item
        with configuration as in Figure \ref{Conf1N}.
    \end{itemize} 
    This is case \hyperref[Tab1N]{$1N$}. 
        \end{minipage} \hspace{2mm} \begin{minipage}{0.3\textwidth} \begin{center} \includegraphics[width=0.97\textwidth]{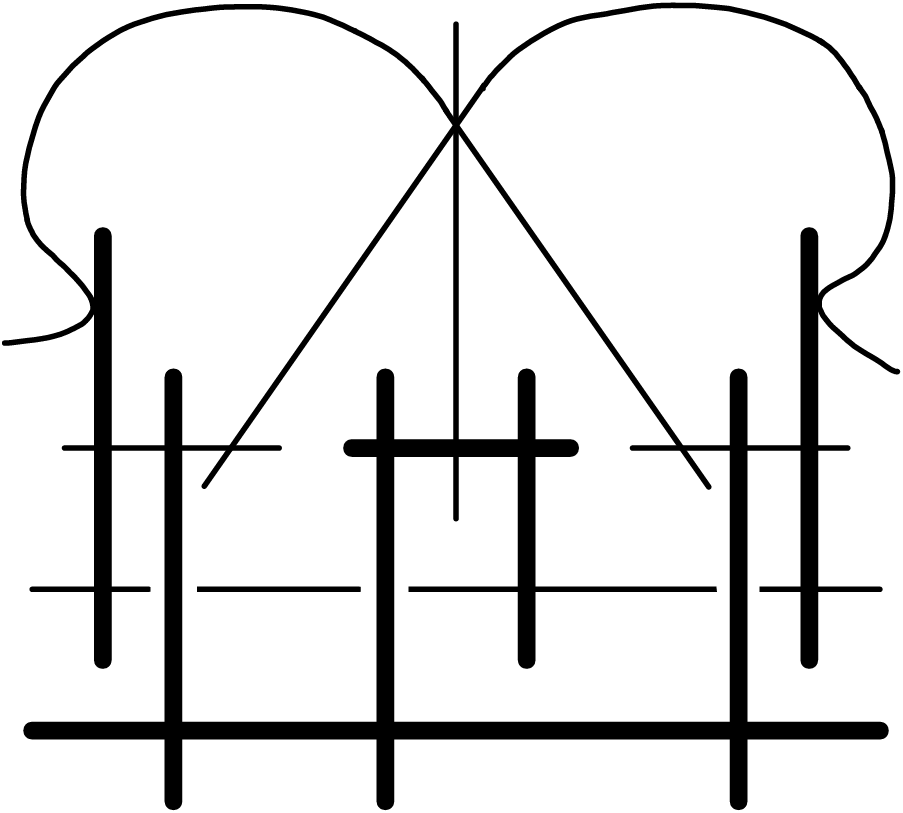} \vspace{-3.5mm}\captionof{figure}{}  \label{Conf1N} \end{center} \end{minipage} 
    
\end{enumerate}

\vspace{2mm} \subsubsection*{\underline{Case \hyperref[Tab3A]{$3A$}}}
We have 
$
E = (E_{1,1}^{} \cup E_{3,1}^{}) - (E_{1,0}^{(2)} \cup E_{3,0}^{(2)} \cup \ell_y^{(2)})
$ and $\Aut_X^0(R) = \left\{  \left( \begin{smallmatrix}
1 &  &  \\
 & 1 &  \\
 &  & i
\end{smallmatrix} \right) 
\in \PGL_3(R) \right\}$.

\begin{itemize}[leftmargin=25pt]
    \item[-] $\lambda xy + \mu z^2$ is $E_{1,1}^{}$-adapted and $\Aut_X^0(R)$ acts as $[\lambda:\mu] \mapsto [\lambda: i^2 \mu]$
    \item[-] $\lambda y^2 + \mu (x + \alpha y)z$ is $E_{3,1}^{}$-adapted and $\Aut_X^0(R)$ acts as $[\lambda:\mu] \mapsto [\lambda: i \mu]$
\end{itemize}

\noindent Note that there is one unique point with non-trivial stabilizer on $E \cap E_{3,1}^{}$, while $\Aut_X^0$ acts transitively on $E \cap E_{1,1}$. Hence, we have the following three choices up to isomorphism:

\begin{enumerate}[leftmargin=*]

 \item \vspace{2mm}
    $p_{1,2}= E_{1,1}^{} \cap C_1^{(2)}, 
    p_{3,2}= E_{3,1}^{} \cap \ell_{x+ \alpha y}^{(2)}$ with $C_1= \cal{V}(xy+z^2) $ and $\alpha \not\in \{0,-1\}$ 
    \begin{itemize}[leftmargin=20pt]
        \item
        $\Aut_{X'}^0(R) =
        \begin{cases}   
        \{\id\} & \text{ if } p \neq 2 \\
        \left\{  \left( \begin{smallmatrix}
1 &  &  \\
 & 1 &  \\
 &  & i
\end{smallmatrix} \right) 
\in \PGL_3(R) \bigg| i^2=1 \right\}
         & \text{ if } p=2
        \end{cases} $
        \\ Hence, $X'$ has global vector fields only if $p=2$. Therefore, we assume $p=2$ when describing the configuration of negative curves.
        \item $(-2)$-curves: $E_{1,0}^{(3)},E_{3,0}^{(3)}, E_{1,1}^{(3)}, E_{3,1}^{(3)}, \ell_{y}^{(3)}, \ell_{z}^{(3)}, \ell_{x+\alpha y}^{(3)}$
        \item $(-1)$-curves: $E_{1,2}^{}, E_{3,2}^{},   E_{2,0}^{(3)}, E_{4,0}^{(3)}, \ell_{x}^{(3)}, \ell_{x-y}^{(3)}, C_2^{(3)}, C_3^{(3)}$ with $C_2= \cal{V}(x^2y+xz^2+ \alpha yz^2), 
        \\C_3= \cal{V}(x^2y+ xz^2+ \alpha yz^2+y^3) $
        \item
        with configuration as in Figure \ref{Conf1L}, that is, as in case \hyperref[Tab1L]{$1L$}.
    \end{itemize} 
     \noindent As explained in Remark \ref{R IsomorphismCheck}, one can check that $X' \cong X_{1L,\alpha'}$ for some $\alpha'$.

 \item \vspace{2mm}
    $p_{3,2}= E_{3,1}^{} \cap \ell_{x+ \alpha y}^{(2)}$ with $\alpha \not\in \{0,-1\}$ 
    \begin{itemize}[leftmargin=20pt]
    \noindent \begin{minipage}{0.5\textwidth}
        \item
        $
        \Aut_{X'}^0(R) =    
        \left\{  \left( \begin{smallmatrix}
1 &  &  \\
 & 1 &  \\
 &  & i
\end{smallmatrix} \right) 
\in \PGL_3(R) \right\}
        $
        \item $(-2)$-curves: $E_{1,0}^{(3)},E_{3,0}^{(3)}, E_{3,1}^{(3)}, \ell_{y}^{(3)}, \ell_{z}^{(3)}, \ell_{x+\alpha y}^{(3)}$
        
        \end{minipage} \noindent \begin{minipage}{0.5\textwidth}
        
        \item $(-1)$-curves: $E_{3,2}^{}, E_{1,1}^{(3)},  E_{2,0}^{(3)}, E_{4,0}^{(3)}, \ell_{x}^{(3)}, \ell_{x-y}^{(3)}$
        \item
        with configuration as in Figure \ref{Conf2A}, that is, as in \\ case \hyperref[Tab2A]{$2A$}.
    \end{minipage}
    \end{itemize} 
  \vspace{1mm} \noindent As explained in Remark \ref{R IsomorphismCheck}, one can check that $X' \cong X_{2A,\alpha'}$ for some $\alpha'$.

\noindent     \begin{minipage}{0.65\textwidth}
 \item \vspace{2mm}
    $p_{1,2}= E_{1,1}^{} \cap C^{(2)}$ with $C= \cal{V}(xy+z^2) $
    \begin{itemize}[leftmargin=20pt]
        \item
        $\Aut_{X'}^0(R) =
        \begin{cases}   
        \{\id\} & \text{ if } p \neq 2 \\
        \left\{  \left( \begin{smallmatrix}
1 &  &  \\
 & 1 &  \\
 &  & i
\end{smallmatrix} \right) 
\in \PGL_3(R) \bigg| i^2=1 \right\}
         & \text{ if } p=2
        \end{cases} $
        \\Hence, $X'$ has global vector fields only if $p=2$. Therefore, we assume $p=2$ when describing the configuration of negative curves.
        \item $(-2)$-curves: $E_{1,0}^{(3)},E_{3,0}^{(3)}, E_{1,1}^{(3)}, \ell_{y}^{(3)}, \ell_{z}^{(3)}$
        \item $(-1)$-curves: $E_{1,2}^{}, E_{3,1}^{(3)}, E_{2,0}^{(3)}, E_{4,0}^{(3)}, \ell_{x}^{(3)}, \ell_{x-y}^{(3)}, \ell_{x+\alpha y}^{(3)}$ with 
        \\$\alpha \not\in \{0,-1\}$
        \item
        with configuration as in Figure \ref{Conf2N}.
    \end{itemize} 
    This is case \hyperref[Tab2N]{$2N$} and we see that we get a $1$-dimensional family of such surfaces $X_{2N,\alpha}$ depending on the parameter $\alpha$.     
        \end{minipage} \hspace{2mm} \begin{minipage}{0.3\textwidth} \begin{center} \includegraphics[width=0.98\textwidth]{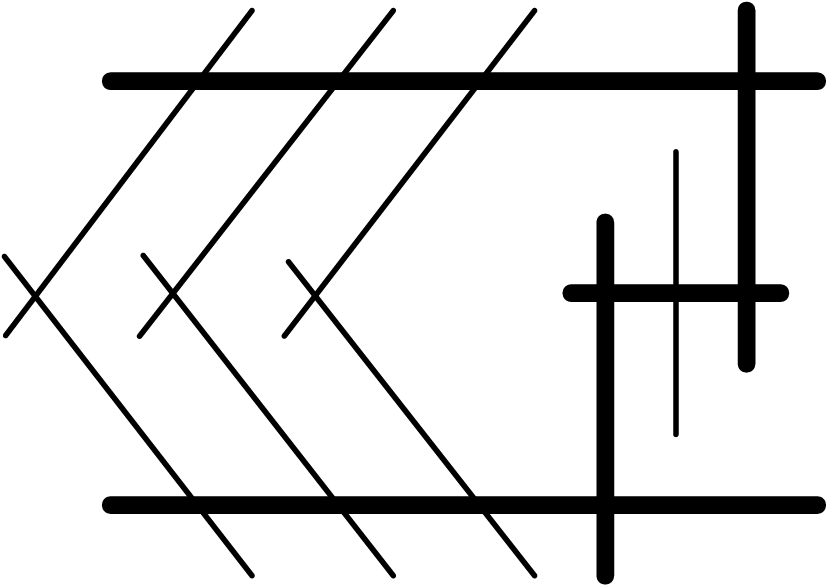} \vspace{-3.5mm}\captionof{figure}{}  \label{Conf2N} \end{center} \end{minipage} 
\end{enumerate}

\vspace{2mm} \subsubsection*{\underline{Case \hyperref[Tab3C]{$3C$}}}
We have 
$
E = (E_{1,1}^{} \cup E_{3,1}^{}) - (E_{1,0}^{(2)} \cup E_{3,0}^{(2)} \cup \ell_y^{(2)} \cup \ell_x^{(2)})
$ and $\Aut_X^0(R) = \left\{  \left( \begin{smallmatrix}
1 &  &  \\
 & 1 &  \\
 &  & i
\end{smallmatrix} \right) 
\in \PGL_3(R) \right\}$.

\begin{itemize}[leftmargin=25pt]
    \item[-] $\lambda xy + \mu z^2$ is $E_{1,1}^{}$-adapted and $\Aut_X^0(R)$ acts as $[\lambda:\mu] \mapsto [\lambda: i^2\mu]$
    \item[-] $\lambda xz + \mu y^2$ is $E_{3,1}^{}$-adapted and $\Aut_X^0(R)$ acts as $[\lambda:\mu] \mapsto [i \lambda: \mu]$
\end{itemize}

\noindent Note that the stabilizer of every point in $E \cap E_{3,1}^{}$ is trivial while $\Aut_X^0$ acts transitively on $E \cap E_{1,1}$. Hence, we have the following unique choice for $p_{1,2}$ up to isomorphism:

\begin{enumerate}[leftmargin=*]
\noindent     \begin{minipage}{0.65\textwidth}
       \item \vspace{2mm}
    $p_{1,2}= E_{1,1}^{} \cap C^{(2)}$ with $C=\cal{V}(xy+z^2) $  
    \begin{itemize}[leftmargin=20pt]
        \item
         $\Aut_{X'}^0(R) =
        \begin{cases}   
        \{\id\} & \text{ if } p \neq 2 \\
        \left\{  \left( \begin{smallmatrix}
1 &  &  \\
 & 1 &  \\
 &  & i
\end{smallmatrix} \right) 
\in \PGL_3(R) \bigg| i^2=1 \right\}
         & \text{ if } p=2
        \end{cases} $
        \\Hence, $X'$ has global vector fields only if $p=2$. Therefore, we assume $p=2$ when describing the configuration of negative curves.
        \item $(-2)$-curves: $E_{1,0}^{(3)},E_{3,0}^{(3)},E_{1,1}^{(3)}, \ell_{x}^{(3)},
        \ell_{y}^{(3)}, \ell_{z}^{(3)}$
        \item $(-1)$-curves: $E_{1,2}^{}, E_{3,1}^{(3)}, E_{2,0}^{(3)}, E_{4,0}^{(3)},  \ell_{x-y}^{(3)}$
        \item
        with configuration as in Figure \ref{Conf2Q}.
    \end{itemize} 
    This is case \hyperref[Tab2Q]{$2Q$}. 
        \end{minipage} \hspace{2mm} \begin{minipage}{0.3\textwidth} \begin{center} \includegraphics[width=0.98\textwidth]{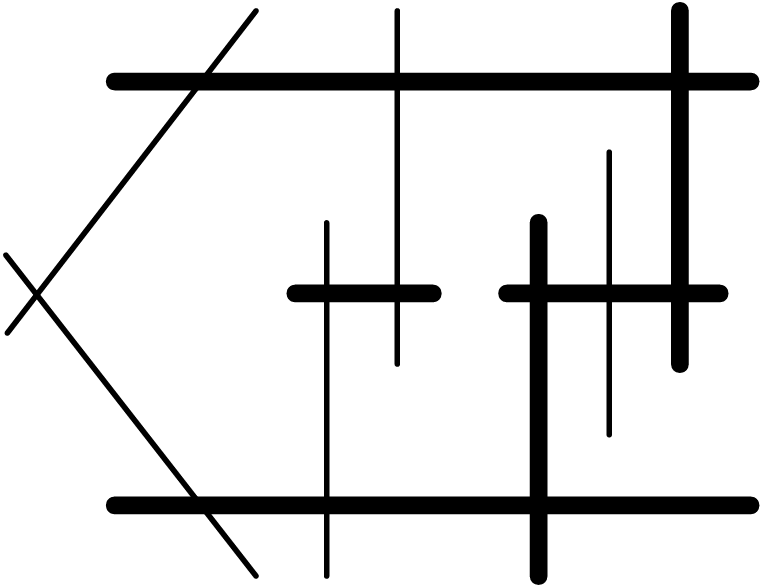} \vspace{-3.5mm}\captionof{figure}{}  \label{Conf2Q} \end{center} \end{minipage} 
\end{enumerate}

\vspace{2mm} \subsubsection*{\underline{Case \hyperref[Tab3D]{$3D$}}}
We have 
$
E = (E_{1,1}^{} \cup E_{2,1}^{}) - (E_{1,0}^{(2)} \cup E_{2,0}^{(2)} \cup \ell_y^{(2)} \cup \ell_x^{(2)})
$ and $\Aut_X^0(R) = \left\{  \left( \begin{smallmatrix}
1 &  &  \\
 & 1 &  \\
 &  & i
\end{smallmatrix} \right) 
\in \PGL_3(R) \right\}$.

\begin{itemize}[leftmargin=25pt]
    \item[-] $\lambda xy + \mu z^2$ is $E_{1,1}^{}$-adapted and $E_{2,1}^{}$-adapted and $\Aut_X^0(R)$ acts as $[\lambda:\mu] \mapsto [\lambda: i^2 \mu]$
\end{itemize}

\noindent Note that the involution $x \leftrightarrow y$ of $\bbP^2$ lifts to an involution of $X$ interchanging $E_{1,1}^{}$ and $E_{2,1}^{}$. Moreover, $\Aut_X^0$ acts transitively and with finite stabilizers on both $E \cap E_{1,1}$ and $E \cap E_{2,1}$. Hence, we have the following three possibilities for $p_{1,2},\hdots,p_{n,2}$ up to isomorphism:

\begin{enumerate}[leftmargin=*]
  \item \vspace{2mm}
    $p_{1,2}= E_{1,1}^{} \cap C^{(2)}, 
    p_{2,2}= E_{2,1}^{} \cap C^{(2)}$ with $C=\cal{V}(xy+z^2) $
    \begin{itemize}[leftmargin=20pt]
        \item
        $\Aut_{X'}^0(R) =
        \begin{cases}   
        \{\id\} & \text{ if } p \neq 2 \\
        \left\{  \left( \begin{smallmatrix}
1 &  &  \\
 & 1 &  \\
 &  & i
\end{smallmatrix} \right) 
\in \PGL_3(R) \bigg| i^2=1 \right\}
         & \text{ if } p=2
        \end{cases} $
        \\Hence, $X'$ has global vector fields only if $p=2$. Therefore, we assume $p=2$ when describing the configuration of negative curves.
        \item $(-2)$-curves: $E_{1,0}^{(3)},E_{2,0}^{(3)},E_{1,1}^{(3)}, E_{2,1}^{(3)}, \ell_{x}^{(3)},
        \ell_{y}^{(3)}, \ell_{z}^{(3)}, C^{(3)}$
        \item $(-1)$-curves: $E_{1,2}^{}, E_{2,2}^{}, E_{3,0}^{(3)}, E_{4,0}^{(3)},  \ell_{x-y}^{(3)}$
        \item
        with configuration as in Figure \ref{Conf1O}, that is, as in case \hyperref[Tab1O]{$1O$}.
    \end{itemize} 
        \noindent As explained in Remark \ref{R IsomorphismCheck}, one can check that $X' \cong X_{1O}$.

  \item \vspace{2mm}
    $p_{1,2}= E_{1,1}^{} \cap C_1^{(2)}, 
    p_{2,2}= E_{2,1}^{} \cap C_2^{(2)}$ with $C_1=\cal{V}(xy+z^2), C_2=\cal{V}(xy+ \alpha z^2) , \alpha \not\in \{0, 1\}$
    \begin{itemize}[leftmargin=20pt]
        \item
        $\Aut_{X'}^0(R) =
        \begin{cases}   
        \{\id\} & \text{ if } p \neq 2 \\
        \left\{  \left( \begin{smallmatrix}
1 &  &  \\
 & 1 &  \\
 &  & i
\end{smallmatrix} \right) 
\in \PGL_3(R) \bigg| i^2=1 \right\}
         & \text{ if } p=2
        \end{cases} $
        \\Hence, $X'$ has global vector fields only if $p=2$. Therefore, we assume $p=2$ when describing the configuration of negative curves.
        \item $(-2)$-curves: $E_{1,0}^{(3)},E_{2,0}^{(3)},E_{1,1}^{(3)}, E_{2,1}^{(3)}, \ell_{x}^{(3)},
        \ell_{y}^{(3)}, \ell_{z}^{(3)}$
        \item $(-1)$-curves: $E_{1,2}^{}, E_{2,2}^{}, E_{3,0}^{(3)}, E_{4,0}^{(3)},  \ell_{x-y}^{(3)}, C_1^{(3)}, C_2^{(3)}, C_3^{(3)}$ with 
        \\$C_{3}=\cal{V}(x^3y^2 + x^2y^3 + xz^4 + \alpha^2 yz^4) $ 
        \item
        with configuration as in Figure \ref{Conf1L}, that is, as in case \hyperref[Tab1L]{$1L$}.
    \end{itemize} 
   \noindent As explained in Remark \ref{R IsomorphismCheck}, one can check that $X' \cong X_{1L,\alpha'}$ for some $\alpha'$.
    
    \noindent     \begin{minipage}{0.65\textwidth}
     \item \vspace{2mm}
    $p_{1,2}= E_{1,1}^{} \cap C^{(2)}$ with $C= \cal{V}(xy+z^2)$
    \begin{itemize}[leftmargin=20pt]
        \item
        $\Aut_{X'}^0(R) =
        \begin{cases}   
        \{\id\} & \text{ if } p \neq 2 \\
        \left\{  \left( \begin{smallmatrix}
1 &  &  \\
 & 1 &  \\
 &  & i
\end{smallmatrix} \right) 
\in \PGL_3(R) \bigg| i^2=1 \right\}
         & \text{ if } p=2
        \end{cases} $
        \\Hence, $X'$ has global vector fields only if $p=2$. Therefore, we assume $p=2$ when describing the configuration of negative curves.
        \item $(-2)$-curves: $E_{1,0}^{(3)},E_{2,0}^{(3)},E_{1,1}^{(3)}, \ell_{x}^{(3)},
        \ell_{y}^{(3)}, \ell_{z}^{(3)}$
        \item $(-1)$-curves: $E_{1,2}^{}, E_{2,1}^{(3)}, E_{3,0}^{(3)}, E_{4,0}^{(3)},  \ell_{x-y}^{(3)}, C^{(3)}$
        \item
        with configuration as in Figure \ref{Conf2P}.
    \end{itemize} 
    This is case \hyperref[Tab2P]{$2P$}.
        \end{minipage} \hspace{2mm} \begin{minipage}{0.3\textwidth} \begin{center} \includegraphics[width=0.98\textwidth]{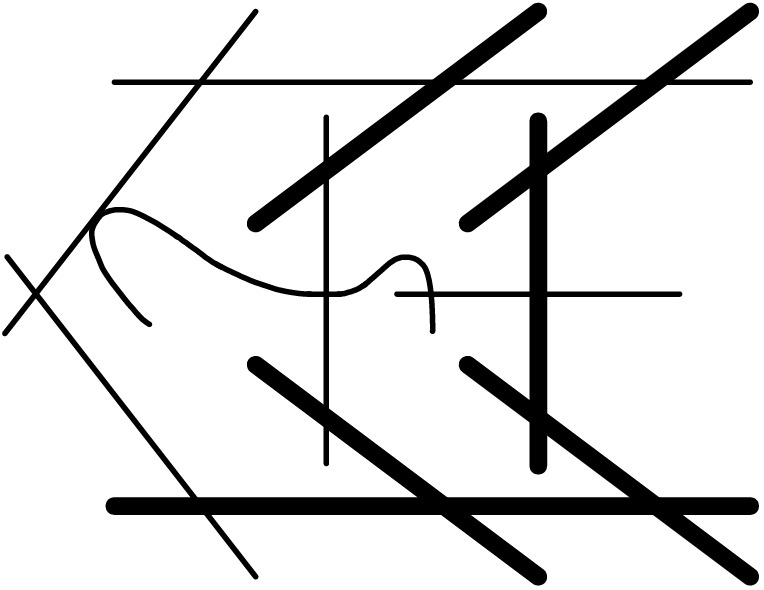} \vspace{-3.5mm}\captionof{figure}{}  \label{Conf2P} \end{center} \end{minipage}

\end{enumerate}

\vspace{2mm} \subsubsection*{\underline{Case \hyperref[Tab4A]{$4A$}}}
We have 
$
E = E_{3,1}^{} - E_{3,0}^{(2)}
$ and $\Aut_X^0(R) = 
 \left\{  \left( \begin{smallmatrix}
1 &  &  \\
 & 1 &  \\
 &  & i
\end{smallmatrix} \right) 
\in \PGL_3(R) \right\}$.

\begin{itemize}[leftmargin=25pt]
    \item[-] $\lambda y^2 + \mu (x+\alpha y)z$ is $E_{3,1}^{}$-adapted and $\Aut_X^0(R)$ acts as $[\lambda:\mu] \mapsto [\lambda: i \mu]$
\end{itemize}

\noindent Note that there is one unique point on $E \cap E_{3,1}^{}$ with non-trivial stabilizer, leading to the following unique choice for $p_{3,2}$:

\begin{enumerate}[leftmargin=*]

  \item \vspace{2mm}
    $p_{3,2}= E_{3,1}^{} \cap \ell_{x+\alpha y}^{(2)}$ with $\alpha \not\in \{0,-1\}$ 
    \begin{itemize}[leftmargin=20pt]
    \noindent \begin{minipage}{0.5\textwidth}
        \item
        $
        \Aut_{X'}^0(R) =    
        \left\{  \left( \begin{smallmatrix}
1 &  &  \\
 & 1 &  \\
 &  & i
\end{smallmatrix} \right) 
\in \PGL_3(R) \right\}
        $
        \item $(-2)$-curves: $E_{3,0}^{(3)},E_{3,1}^{(3)}, \ell_{z}^{(3)},
        \ell_{x+ \alpha y}^{(3)}$
        
        \end{minipage} \noindent \begin{minipage}{0.5\textwidth}
        
        \item $(-1)$-curves: $E_{3,2}^{}, E_{1,0}^{(3)},
        E_{2,0}^{(3)},
        E_{4,0}^{(3)},    
        \ell_{x}^{(3)},
        \ell_{y}^{(3)},
        \\ \ell_{x-y}^{(3)}$
        \item
        with configuration as in Figure \ref{Conf3A}, that is, as in \\ case \hyperref[Tab3A]{$3A$}.
    \end{minipage}
    \end{itemize} 
   \vspace{1mm} \noindent As explained in Remark \ref{R IsomorphismCheck}, one can check that $X' \cong X_{3A,\alpha'}$ for some $\alpha'$.

\end{enumerate}

\vspace{2mm} \subsubsection*{\underline{Case \hyperref[Tab4B]{$4B$}}}
We have 
$
E = E_{3,1}^{} - (E_{3,0}^{(2)} \cup \ell_y^{(2)})
$ and $\Aut_X^0(R) = \left\{  \left( \begin{smallmatrix}
1 &  &  \\
 & 1 &  \\
 &  & i
\end{smallmatrix} \right) 
\in \PGL_3(R) \right\}$.

\begin{itemize}[leftmargin=25pt]
    \item[-] $\lambda x^2 + \mu yz$ is $E_{3,1}^{}$-adapted and $\Aut_X^0(R)$ acts as $[\lambda:\mu] \mapsto [\lambda: i \mu]$
\end{itemize}

\noindent There is no point on $E \cap E_{3,1}^{}$ with non-trivial stabilizer, so we get no new cases by further blowing up $X$.

\vspace{2mm} \subsubsection*{\underline{Case \hyperref[Tab4C]{$4C$}}}
We have 
$
E = E_{1,1}^{} - (E_{1,0}^{(2)} \cup \ell_y^{(2)})
$ and $\Aut_X^0(R) = \left\{  \left( \begin{smallmatrix}
1 &  &  \\
 & 1 &  \\
 &  & i
\end{smallmatrix} \right) 
\in \PGL_3(R) \right\}$.

\begin{itemize}[leftmargin=25pt]
    \item[-] $\lambda xy + \mu z^2$ is $E_{1,1}^{}$-adapted and $\Aut_X^0(R)$ acts as $[\lambda:\mu] \mapsto [\lambda: i^2 \mu]$
\end{itemize}

\noindent In particular, $\Aut_X^0$ acts transitively on $E \cap E_{1,1}$. We get the following unique choice for $p_{1,2}$ up to isomorphism:

\begin{enumerate}[leftmargin=*]
\noindent     \begin{minipage}{0.65\textwidth}
    \item \vspace{2mm}
    $p_{1,2}= E_{1,1}^{} \cap C^{(2)}$ with $C= \cal{V}(xy+z^2)$ 
    \begin{itemize}[leftmargin=20pt]
        \item
        $\Aut_{X'}^0(R) =
        \begin{cases}   
        \{\id\} & \text{ if } p \neq 2 \\
        \left\{  \left( \begin{smallmatrix}
1 &  &  \\
 & 1 &  \\
 &  & i
\end{smallmatrix} \right) 
\in \PGL_3(R) \bigg| i^2=1 \right\}
         & \text{ if } p=2
        \end{cases} $
        \\ Hence, $X'$ has global vector fields only if $p=2$. Therefore, we assume $p=2$ when describing the configuration of negative curves.
        \item $(-2)$-curves: $E_{1,0}^{(3)},E_{1,1}^{(3)},
        \ell_{y}^{(3)},
        \ell_{z}^{(3)}$
        \item $(-1)$-curves: $E_{1,2}^{}, E_{2,0}^{(3)},
        E_{3,0}^{(3)},
        E_{4,0}^{(3)},    
        \ell_{x}^{(3)},
        \ell_{x-y}^{(3)}$
        \item
        with configuration as in Figure \ref{Conf3N}.
    \end{itemize} 
    This is case \hyperref[Tab3N]{$3N$}.
        \end{minipage} \hspace{2mm} \begin{minipage}{0.3\textwidth} \begin{center} \includegraphics[width=0.98\textwidth]{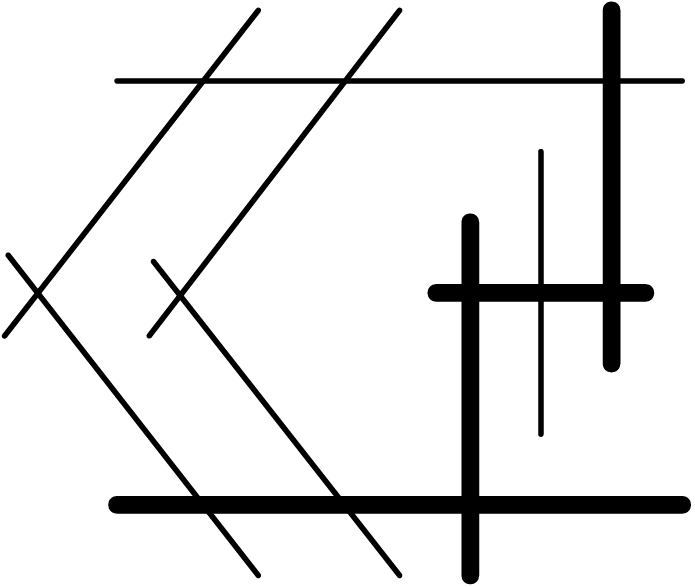} \vspace{-3.5mm}\captionof{figure}{}  \label{Conf3N} \end{center} \end{minipage} 
\end{enumerate}

\vspace{2mm} \subsubsection*{\underline{Case \hyperref[Tab3B]{$3B$}}}
We have 
$
E = \bigcup_{j=1}^{3} E_{j,1}^{} - \bigcup_{j=1}^{3} E_{j,0}^{(2)}
$ and $\Aut_X^0(R) = \left\{  \left( \begin{smallmatrix}
1 &  &  \\
 & 1 &  \\
 &  & i
\end{smallmatrix} \right) 
\in \PGL_3(R) \right\}$.

\begin{itemize}[leftmargin=25pt]
    \item[-] $\lambda xy + \mu z^2$ is $E_{1,1}^{}$-adapted and $E_{2,1}^{}$-adapted and  $\Aut_X^0(R)$ acts as $[\lambda:\mu] \mapsto [\lambda: i^2 \mu]$
    \item[-] $\lambda (x-y)x + \mu z^2$ is $E_{3,1}^{}$-adapted and  $\Aut_X^0(R)$ acts as $[\lambda:\mu] \mapsto [\lambda: i^2 \mu]$
\end{itemize}

\noindent 
Note that automorphisms of $\bbP^2$ fixing $[0:0:1]$ and interchanging the $p_{j,0}$ lift to automorphisms of $X$ interchanging the $E_{j,1}^{}$. Moreover, since $X$ has degree $3$, we are only allowed to blow up two more points. Finally, on every $E \cap E_{j,1}$, the action of $\Aut_X^0$ has two orbits and one of them is a fixed point. Hence, we get the following six possibilities for $p_{1,2},\hdots,p_{3,2}$ up to isomorphism:

\begin{enumerate}[leftmargin=*]

\noindent     \begin{minipage}{0.65\textwidth}
    \item \vspace{2mm}
$p_{1,2}= E_{1,1}^{} \cap C_1^{(2)}, 
p_{2,2}= E_{2,1}^{} \cap C_1^{(2)}$ with $C_1=\cal{V}(xy+z^2)$
    \begin{itemize}[leftmargin=20pt]
        \item
        $\Aut_{X'}^0(R) =
        \begin{cases}   
        \{\id\} & \text{ if } p \neq 2 \\
        \left\{  \left( \begin{smallmatrix}
1 &  &  \\
 & 1 &  \\
 &  & i
\end{smallmatrix} \right) 
\in \PGL_3(R) \bigg| i^2=1 \right\}
         & \text{ if } p=2
        \end{cases} $
        \\Hence, $X'$ has global vector fields only if $p=2$. Therefore, we assume $p=2$ when describing the configuration of negative curves.
        \item $(-2)$-curves: $E_{1,0}^{(3)},
        E_{2,0}^{(3)},
        E_{3,0}^{(3)}, E_{1,1}^{(3)}, E_{2,1}^{(3)},
        \ell_{z}^{(3)}, C_1^{(3)}$
        \item $(-1)$-curves: $ E_{1,2}^{}, E_{2,2}^{},
        E_{3,1}^{(3)},
        \ell_{x}^{(3)},
        \ell_{y}^{(3)},
        \ell_{x-y}^{(3)}, C_2^{(3)}, C_3^{(3)}$ with 
        \\$C_2= \cal{V}(xy+y^2+z^2), C_3=\cal{V}(xy+x^2+z^2)$
        \item
        with configuration as in Figure \ref{Conf1K}.
    \end{itemize} 
    This is case \hyperref[Tab1K]{$1K$}.
        \end{minipage} \hspace{2mm} \begin{minipage}{0.3\textwidth} \begin{center} \includegraphics[width=0.98\textwidth]{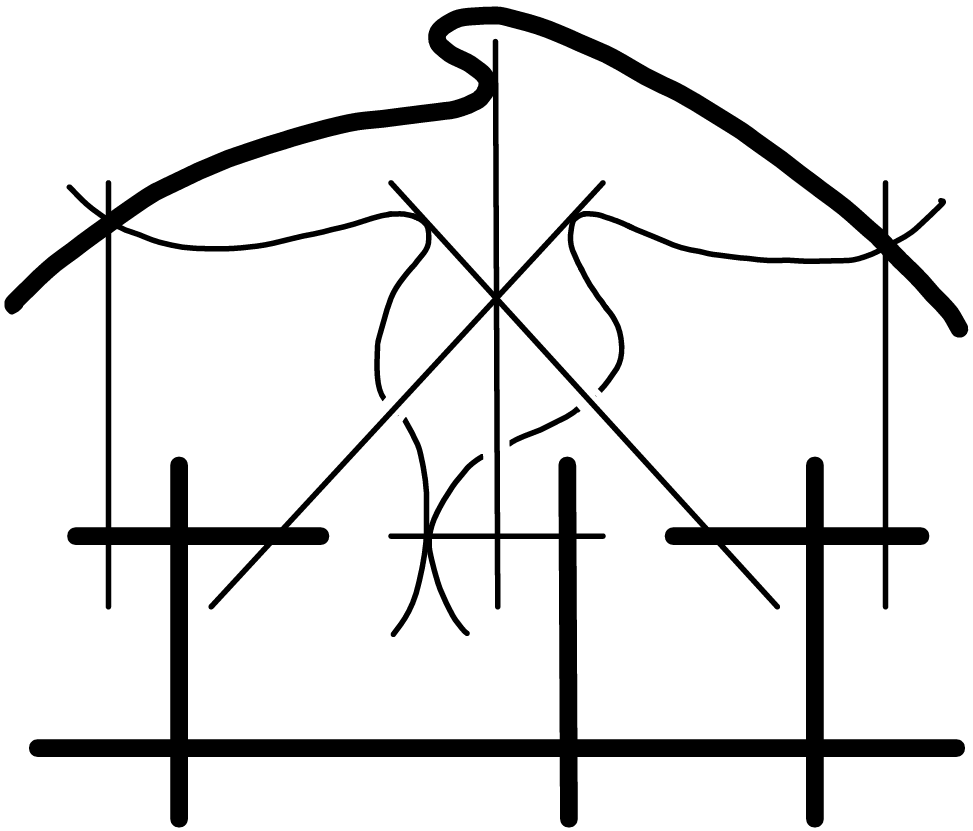} \vspace{-3.5mm}\captionof{figure}{}  \label{Conf1K} \end{center} \end{minipage}

    \noindent     \begin{minipage}{0.65\textwidth}
    \item \vspace{2mm}
$p_{1,2}= E_{1,1}^{} \cap C_1^{(2)}, 
p_{2,2}= E_{2,1}^{} \cap C_2^{(2)}$ with $C_1=\cal{V}(xy+z^2), C_2=\cal{V}(xy+ \alpha z^2), \alpha \not\in \{0, 1 \}$
    \begin{itemize}[leftmargin=20pt]
        \item
        $\Aut_{X'}^0(R) =
        \begin{cases}   
        \{\id\} & \text{ if } p \neq 2 \\
        \left\{  \left( \begin{smallmatrix}
1 &  &  \\
 & 1 &  \\
 &  & i
\end{smallmatrix} \right) 
\in \PGL_3(R) \bigg| i^2=1 \right\}
         & \text{ if } p=2
        \end{cases} $
        \\Hence, $X'$ has global vector fields only if $p=2$. Therefore, we assume $p=2$ when describing the configuration of negative curves.
        \item $(-2)$-curves: $E_{1,0}^{(3)},
        E_{2,0}^{(3)},
        E_{3,0}^{(3)}, E_{1,1}^{(3)}, E_{2,1}^{(3)},
        \ell_{z}^{(3)}$
        \item $(-1)$-curves: $E_{1,2}^{}, E_{2,2}^{},
        E_{3,1}^{(3)},
        \ell_{x}^{(3)},
        \ell_{y}^{(3)},
        \ell_{x-y}^{(3)}, C_1^{(3)}, C_2^{(3)}, C_3^{(3)}, 
        \\C_4^{(3)}, C_5^{(3)}, C_6^{(3)}, C_7^{(3)}$ with $C_3= \cal{V}(xy+y^2+z^2), 
        \\C_4= \cal{V}(xy+x^2+ \alpha z^2), C_5= \cal{V}(x^2y^2+ xy^3+ \alpha y^2z^2+ z^4), 
        \\C_6= \cal{V}(x^2y^2 + x^3y + x^2z^2+ \alpha^2 z^4), 
        \\C_7= \cal{V}(x^3y^2+x^2y^3+xz^4+ \alpha^2 yz^4)$
        \item
        with configuration as in Figure \ref{Conf1J}.
    \end{itemize} 
    This is case \hyperref[Tab1J]{$1J$} and we see that we get a $1$-dimensional family of such surfaces $X_{1J,\alpha}$ depending on the parameter $\alpha$. 
        \end{minipage} \hspace{2mm} \begin{minipage}{0.3\textwidth} \begin{center} \includegraphics[width=0.98\textwidth]{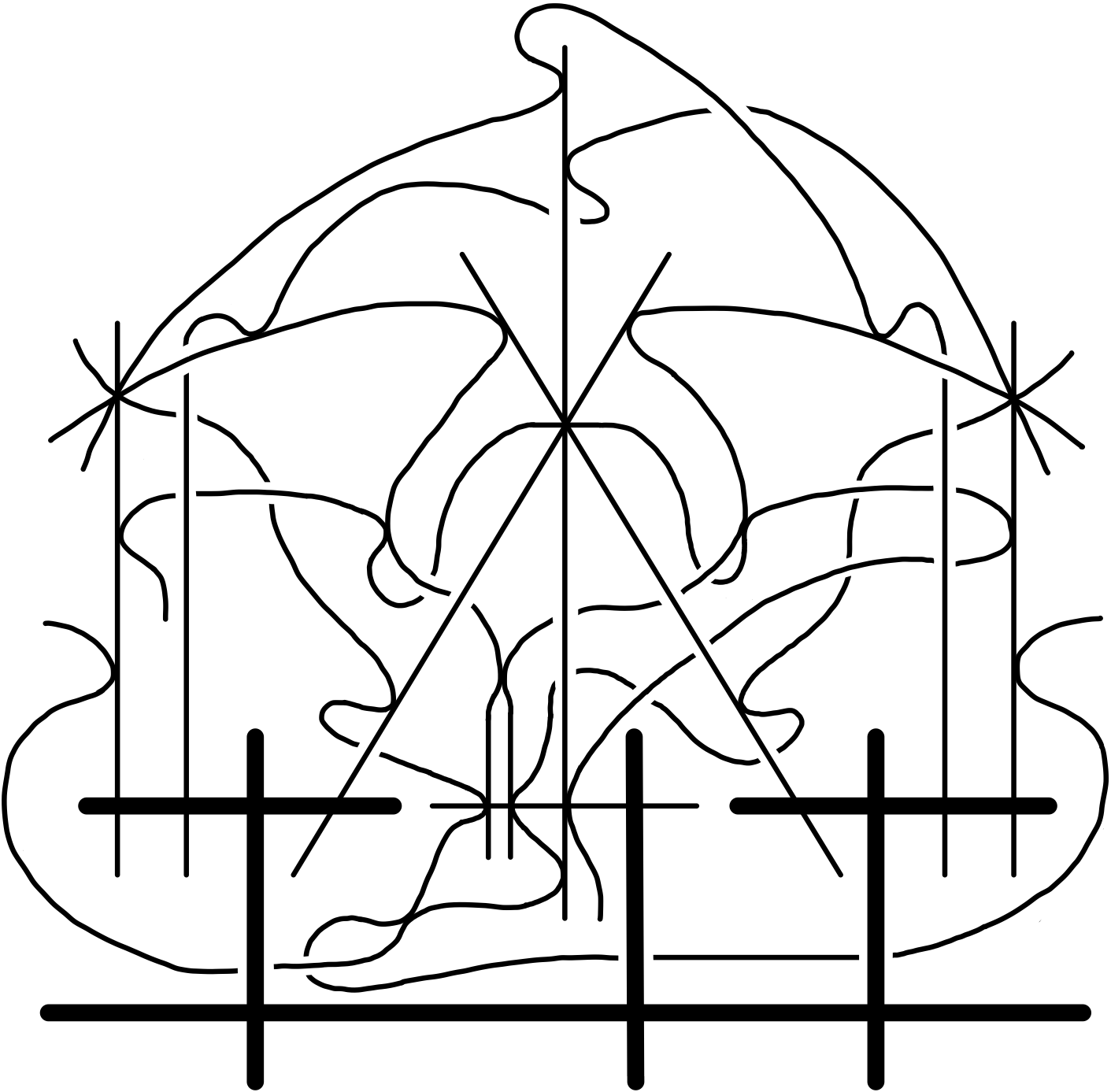} \vspace{-3.5mm}\captionof{figure}{}  \label{Conf1J} \end{center} \end{minipage}

    \item \vspace{2mm}
$p_{1,2}= E_{1,1}^{} \cap C_1^{(2)}, 
p_{2,2}= E_{2,1}^{} \cap \ell_{x}^{(2)}$ with $C_1= \cal{V}(xy+z^2)$
    \begin{itemize}[leftmargin=20pt]
        \item
        $\Aut_{X'}^0(R) =
        \begin{cases}   
        \{\id\} & \text{ if } p \neq 2 \\
        \left\{  \left( \begin{smallmatrix}
1 &  &  \\
 & 1 &  \\
 &  & i
\end{smallmatrix} \right) 
\in \PGL_3(R) \bigg| i^2=1 \right\}
         & \text{ if } p=2
        \end{cases} $
        \\Hence, $X'$ has global vector fields only if $p=2$. Therefore, we assume $p=2$ when describing the configuration of negative curves.
        \item $(-2)$-curves: $E_{1,0}^{(3)},
        E_{2,0}^{(3)},
        E_{3,0}^{(3)}, E_{1,1}^{(3)}, E_{2,1}^{(3)},\ell_{x}^{(3)},
        \ell_{z}^{(3)}$
        \item $(-1)$-curves: $E_{1,2}^{}, E_{2,2}^{},
        E_{3,1}^{(3)},
        \ell_{y}^{(3)},
        \ell_{x-y}^{(3)}, C_1^{(3)}, C_2^{(3)}, C_3^{(3)}$ with $C_2= \cal{V}(xy+y^2+z^2), 
        \\C_3^{(3)}=\cal{V}(x^2y^2+xy^3+z^4)$
        \item
        with configuration as in Figure \ref{Conf1K}, that is, as in case \hyperref[Tab1K]{$1K$}.
    \end{itemize} 
      \noindent As explained in Remark \ref{R IsomorphismCheck}, one can check that $X' \cong X_{1K}$.

    \noindent     \begin{minipage}{0.65\textwidth}
   \item \vspace{2mm}
$p_{1,2}= E_{1,1}^{} \cap \ell_{y}^{(2)}, 
p_{2,2}= E_{2,1}^{} \cap \ell_{x}^{(2)}$
    \begin{itemize}[leftmargin=20pt]
        \item
        $
        \Aut_{X'}^0(R) =    
        \left\{  \left( \begin{smallmatrix}
1 &  &  \\
 & 1 &  \\
 &  & i
\end{smallmatrix} \right) 
\in \PGL_3(R) \right\}
        $
        \item $(-2)$-curves: $E_{1,0}^{(3)},
        E_{2,0}^{(3)},
        E_{3,0}^{(3)}, E_{1,1}^{(3)}, E_{2,1}^{(3)},\ell_{x}^{(3)},
        \ell_{y}^{(3)},
        \ell_{z}^{(3)}$
        \item $(-1)$-curves: $E_{1,2}^{}, E_{2,2}^{},
        E_{3,1}^{(3)},
        \ell_{x-y}^{(3)}$
        \item
        with configuration as in Figure \ref{Conf1B}.
    \end{itemize} 
    This is case \hyperref[Tab1B]{$1B$}.
        \end{minipage} \hspace{2mm} \begin{minipage}{0.3\textwidth} \begin{center} \includegraphics[width=0.6\textwidth]{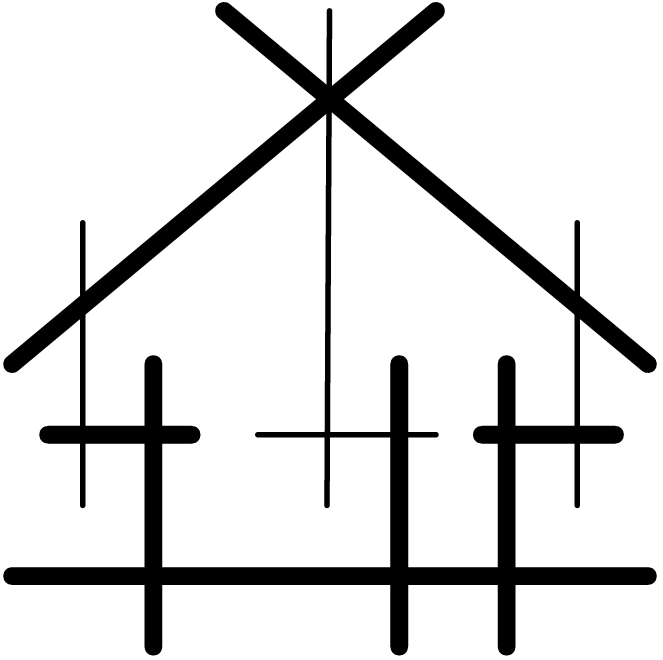} \vspace{-3.5mm}\captionof{figure}{} \label{Conf1B}  \end{center} \end{minipage}

    \noindent     \begin{minipage}{0.65\textwidth}
    \item \vspace{2mm}
$p_{1,2}= E_{1,1}^{} \cap C_1^{(2)}$ with $C_1= \cal{V}(xy+z^2)$
    \begin{itemize}[leftmargin=20pt]
        \item
        $\Aut_{X'}^0(R) =
        \begin{cases}   
        \{\id\} & \text{ if } p \neq 2 \\
        \left\{  \left( \begin{smallmatrix}
1 &  &  \\
 & 1 &  \\
 &  & i
\end{smallmatrix} \right) 
\in \PGL_3(R) \bigg| i^2=1 \right\}
         & \text{ if } p=2
        \end{cases} $
        \\Hence, $X'$ has global vector fields only if $p=2$. Therefore, we assume $p=2$ when describing the configuration of negative curves.
        \item $(-2)$-curves: $E_{1,0}^{(3)},
        E_{2,0}^{(3)},
        E_{3,0}^{(3)}, E_{1,1}^{(3)},
        \ell_{z}^{(3)}$
        \item $(-1)$-curves: $E_{1,2}^{}, E_{2,1}^{(3)},
        E_{3,1}^{(3)},
        \ell_{x}^{(3)},
        \ell_{y}^{(3)},
        \ell_{x-y}^{(3)}, C_1^{(3)}, C_2^{(3)}$ with \\$C_2= \cal{V}(xy+y^2+z^2)$
        \item
        with configuration as in Figure \ref{Conf2O}.
    \end{itemize} 
    This is case \hyperref[Tab2O]{$2O$}.
        \end{minipage} \hspace{2mm} \begin{minipage}{0.3\textwidth} \begin{center} \includegraphics[width=0.96\textwidth]{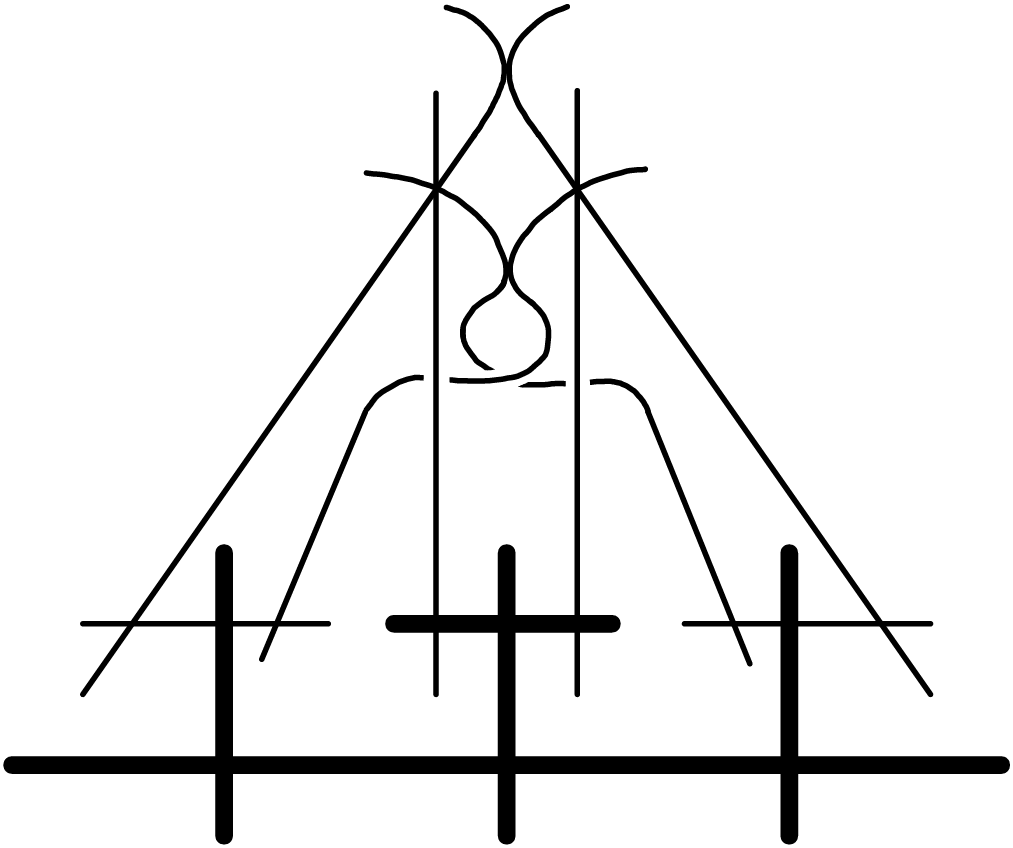} \vspace{-3.5mm}\captionof{figure}{}  \label{Conf2O} \end{center} \end{minipage}

    \noindent     \begin{minipage}{0.65\textwidth}
    \item \vspace{2mm}
$p_{1,2}= E_{1,1}^{} \cap \ell_{y}^{(2)}$
    \begin{itemize}[leftmargin=20pt]
        \item
        $
        \Aut_{X'}^0(R) =    
        \left\{  \left( \begin{smallmatrix}
1 &  &  \\
 & 1 &  \\
 &  & i
\end{smallmatrix} \right) 
\in \PGL_3(R) \right\}
        $
        \item $(-2)$-curves: $E_{1,0}^{(3)},
        E_{2,0}^{(3)},
        E_{3,0}^{(3)}, E_{1,1}^{(3)},
        \ell_{y}^{(3)},
        \ell_{z}^{(3)}$
        \item $(-1)$-curves: $E_{1,2}^{}, E_{2,1}^{(3)},
        E_{3,1}^{(3)},
        \ell_{x}^{(3)},
        \ell_{x-y}^{(3)}$
        \item
        with configuration as in Figure \ref{Conf2B}.
    \end{itemize} 
    This is case \hyperref[Tab2B]{$2B$}.
        \end{minipage} \hspace{2mm} \begin{minipage}{0.3\textwidth} \begin{center} \includegraphics[width=0.8\textwidth]{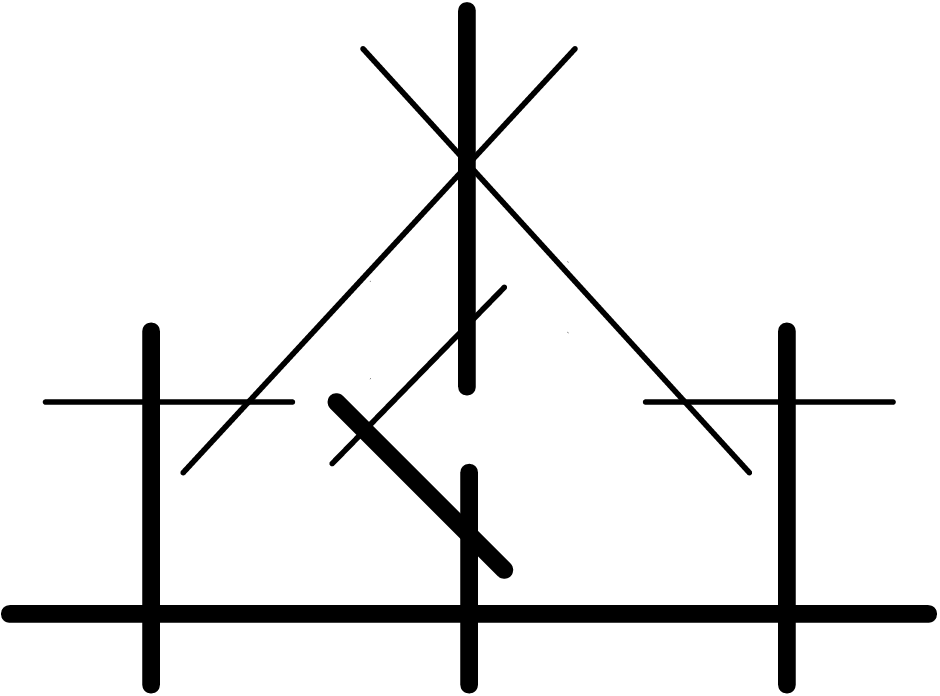} \vspace{-3.5mm}\captionof{figure}{}  \label{Conf2B} \end{center} \end{minipage}

\end{enumerate}

\vspace{2mm} \subsubsection*{\underline{Case \hyperref[Tab4D]{$4D$}}}
We have 
$
E = \bigcup_{j=1}^{2} E_{j,1}^{} - \bigcup_{j=1}^{2} E_{j,0}^{(2)}
$ and $\Aut_X^0(R) = \left\{  \left( \begin{smallmatrix}
1 &  &  \\
 & 1 &  \\
 &  & i
\end{smallmatrix} \right) 
\in \PGL_3(R) \right\}$.

\begin{itemize}[leftmargin=25pt]
    \item[-] $\lambda xy + \mu z^2$ is $E_{1,1}^{}$-adapted and $E_{2,1}^{}$-adapted and  $\Aut_X^0(R)$ acts as $[\lambda:\mu] \mapsto [\lambda: i^2 \mu]$
\end{itemize}

\noindent 
Note that automorphisms of $\bbP^2$ fixing $[0:0:1]$ and interchanging $p_{1,0}$ and $p_{2,0}$ lift to automorphisms of $X$ interchanging $E_{1,1}^{}$ and $E_{2,1}^{}$. Moreover, $\Aut_X^0$ has two orbits on each $E \cap E_{j,1}$, one of which is a fixed point. Hence, we get the following six possibilities for $p_{1,2},p_{2,2}$ up to isomorphism:

\begin{enumerate}[leftmargin=*]

 \item \vspace{2mm}
$p_{1,2}= E_{1,1}^{} \cap C^{(2)}, 
p_{2,2}= E_{2,1}^{} \cap C^{(2)}$ with $C= \cal{V}(xy+z^2)$ 
    \begin{itemize}[leftmargin=20pt]
        \item
        $\Aut_{X'}^0(R) =
        \begin{cases}   
        \{\id\} & \text{ if } p \neq 2 \\
        \left\{  \left( \begin{smallmatrix}
1 &  &  \\
 & 1 &  \\
 &  & i
\end{smallmatrix} \right) 
\in \PGL_3(R) \bigg| i^2=1 \right\}
         & \text{ if } p=2
        \end{cases} $
        \\Hence, $X'$ has global vector fields only if $p=2$. Therefore, we assume $p=2$ when describing the configuration of negative curves.
        \item $(-2)$-curves: $E_{1,0}^{(3)},
        E_{2,0}^{(3)},E_{1,1}^{(3)}, E_{2,1}^{(3)},
        \ell_{z}^{(3)}, C^{(3)}$
        \item $(-1)$-curves: $E_{1,2}^{}, E_{2,2}^{},
        E_{3,0}^{(3)},
        \ell_{x}^{(3)},
        \ell_{y}^{(3)}$
        \item
        with configuration as in Figure \ref{Conf2Q}, that is, as in case \hyperref[Tab2Q]{$2Q$}.
    \end{itemize} 
        \noindent As explained in Remark \ref{R IsomorphismCheck}, one can check that $X' \cong X_{2Q}$.
    
     \item \vspace{2mm}
$p_{1,2}= E_{1,1}^{} \cap C_1^{(2)}, 
p_{2,2}= E_{2,1}^{} \cap C_2^{(2)}$ with $C_1= \cal{V}(xy+z^2), C_2=\cal{V}(xy+ \alpha z^2), \alpha \not\in \{0, 1\} $ 
    \begin{itemize}[leftmargin=20pt]
        \item
        $\Aut_{X'}^0(R) =
        \begin{cases}   
        \{\id\} & \text{ if } p \neq 2 \\
        \left\{  \left( \begin{smallmatrix}
1 &  &  \\
 & 1 &  \\
 &  & i
\end{smallmatrix} \right) 
\in \PGL_3(R) \bigg| i^2=1 \right\}
        & \text{ if } p=2
        \end{cases} $
        \\Hence, $X'$ has global vector fields only if $p=2$. Therefore, we assume $p=2$ when describing the configuration of negative curves.
        \item $(-2)$-curves: $E_{1,0}^{(3)},
        E_{2,0}^{(3)},E_{1,1}^{(3)}, E_{2,1}^{(3)},
        \ell_{z}^{(3)}$
        \item $(-1)$-curves: $E_{1,2}^{}, E_{2,2}^{},
        E_{3,0}^{(3)},
        \ell_{x}^{(3)},
        \ell_{y}^{(3)}, C_1^{(3)}, C_2^{(3)}$
        \item
        with configuration as in Figure \ref{Conf2N}, that is, as in case \hyperref[Tab2N]{$2N$}.
    \end{itemize} 
         \noindent As explained in Remark \ref{R IsomorphismCheck}, one can check that $X' \cong X_{2N,\alpha'}$ for some $\alpha'$.

     \item \vspace{2mm}
$p_{1,2}= E_{1,1}^{} \cap C^{(2)}, 
p_{2,2}= E_{2,1}^{} \cap \ell_{x}^{(2)}$ with $C= \cal{V}(xy+z^2)$ 
    \begin{itemize}[leftmargin=20pt]
        \item
        $\Aut_{X'}^0(R) =
        \begin{cases}   
        \{\id\} & \text{ if } p \neq 2 \\
        \left\{  \left( \begin{smallmatrix}
1 &  &  \\
 & 1 &  \\
 &  & i
\end{smallmatrix} \right) 
\in \PGL_3(R) \bigg| i^2=1 \right\}
         & \text{ if } p=2
        \end{cases} $
        \\Hence, $X'$ has global vector fields only if $p=2$. Therefore, we assume $p=2$ when describing the configuration of negative curves.
        \item $(-2)$-curves: $E_{1,0}^{(3)},
        E_{2,0}^{(3)},E_{1,1}^{(3)}, E_{2,1}^{(3)},
        \ell_{x}^{(3)},
        \ell_{z}^{(3)}$
        \item $(-1)$-curves: $E_{1,2}^{}, E_{2,2}^{},
        E_{3,0}^{(3)},
        \ell_{y}^{(3)}, C^{(3)}$
        \item
        with configuration as in Figure \ref{Conf2Q}, that is, as in case \hyperref[Tab2Q]{$2Q$}.
    \end{itemize} 
       \noindent As explained in Remark \ref{R IsomorphismCheck}, one can check that $X' \cong X_{2Q}$.
    
    \noindent     \begin{minipage}{0.65\textwidth}
     \item \vspace{2mm}
$p_{1,2}= E_{1,1}^{} \cap \ell_{y}^{(2)}, 
p_{2,2}= E_{2,1}^{} \cap \ell_{x}^{(2)}$ 
    \begin{itemize}[leftmargin=20pt]
        \item
        $
        \Aut_{X'}^0(R) =    
        \left\{  \left( \begin{smallmatrix}
1 &  &  \\
 & 1 &  \\
 &  & i
\end{smallmatrix} \right) 
\in \PGL_3(R) \right\}
        $
        \item $(-2)$-curves: $E_{1,0}^{(3)},
        E_{2,0}^{(3)},E_{1,1}^{(3)}, E_{2,1}^{(3)},
        \ell_{x}^{(3)},
        \ell_{y}^{(3)},
        \ell_{z}^{(3)}$
        \item $(-1)$-curves: $E_{1,2}^{}, E_{2,2}^{},
        E_{3,0}^{(3)}$
        \item
        with configuration as in Figure \ref{Conf2F}.
    \end{itemize} 
    This is case \hyperref[Tab2F]{$2F$}.     \end{minipage} \hspace{2mm} \begin{minipage}{0.3\textwidth} \begin{center} \includegraphics[width=0.55\textwidth]{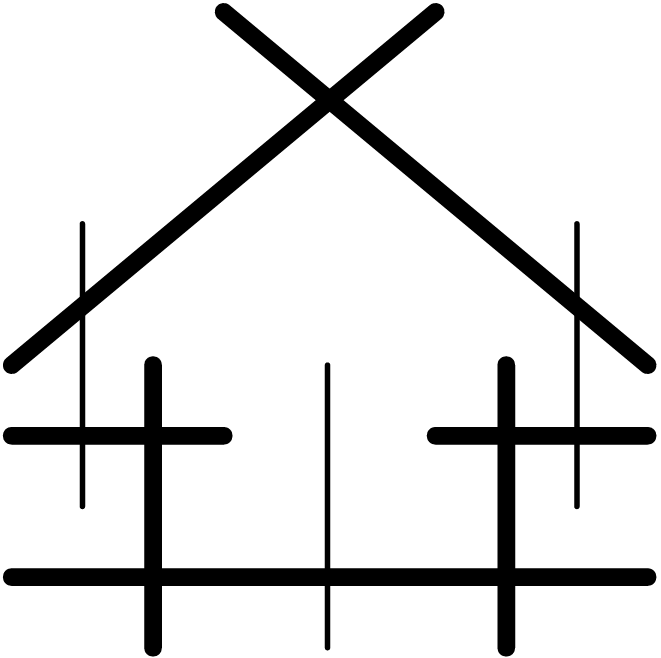} \vspace{-3.5mm}\captionof{figure}{}  \label{Conf2F} \end{center} \end{minipage}

     \item \vspace{2mm}
$p_{1,2}= E_{1,1}^{} \cap C^{(2)}$ with $C=\cal{V}(xy+z^2)$
    \begin{itemize}[leftmargin=20pt]
        \item
        $\Aut_{X'}^0(R) =
        \begin{cases}   
        \{\id\} & \text{ if } p \neq 2 \\
        \left\{  \left( \begin{smallmatrix}
1 &  &  \\
 & 1 &  \\
 &  & i
\end{smallmatrix} \right) 
\in \PGL_3(R) \bigg| i^2=1 \right\}
         & \text{ if } p=2
        \end{cases} $
        \\Hence, $X'$ has global vector fields only if $p=2$. Therefore, we assume $p=2$ when describing the configuration of negative curves.
        \item $(-2)$-curves: $E_{1,0}^{(3)},
        E_{2,0}^{(3)},E_{1,1}^{(3)},
        \ell_{z}^{(3)}$
        \item $(-1)$-curves: $E_{1,2}^{}, E_{2,1}^{(3)},
        E_{3,0}^{(3)},
        \ell_{x}^{(3)},
        \ell_{y}^{(3)}, C^{(3)}$
        \item
        with configuration as in Figure \ref{Conf3N}, that is, as in case \hyperref[Tab3N]{$3N$}.
    \end{itemize} 
     \noindent As explained in Remark \ref{R IsomorphismCheck}, one can check that $X' \cong X_{3N}$.
    
    \noindent     \begin{minipage}{0.65\textwidth}
    \item \vspace{2mm}
$p_{1,2}= E_{1,1}^{} \cap \ell_{y}^{(2)}$ 
    \begin{itemize}[leftmargin=20pt]
        \item
        $
        \Aut_{X'}^0(R) =    
        \left\{  \left( \begin{smallmatrix}
1 &  &  \\
 & 1 &  \\
 &  & i
\end{smallmatrix} \right) 
\in \PGL_3(R) \right\}
        $
        \item $(-2)$-curves: $E_{1,0}^{(3)},
        E_{2,0}^{(3)},E_{1,1}^{(3)},
        \ell_{y}^{(3)},
        \ell_{z}^{(3)}$
        \item $(-1)$-curves: $E_{1,2}^{}, E_{2,1}^{(3)},
        E_{3,0}^{(3)},
        \ell_{x}^{(3)}$
        \item
        with configuration as in Figure \ref{Conf3E}.
    \end{itemize} 
    This is case \hyperref[Tab3E]{$3E$}. 
        \end{minipage} \hspace{2mm} \begin{minipage}{0.3\textwidth} \begin{center} \includegraphics[width=0.8\textwidth]{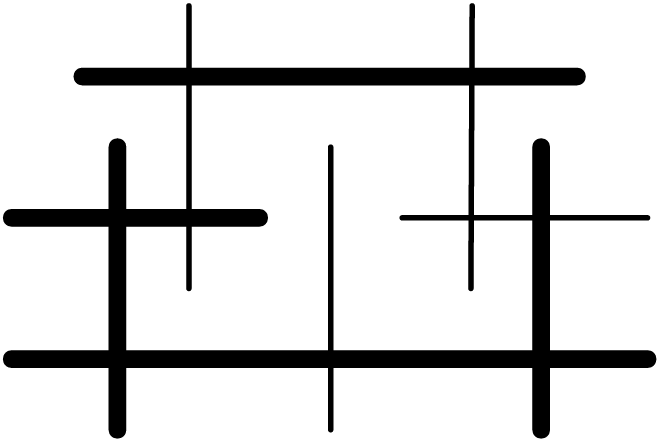} \vspace{-3.5mm}\captionof{figure}{} \label{Conf3E}  \end{center} \end{minipage}

\end{enumerate}

\vspace{2mm} \subsubsection*{\underline{Case \hyperref[Tab5C]{$5C$}}}
We have 
$
E = E_{1,1}^{} - E_{1,0}^{(2)}
$ and $\Aut_X^0(R) = \left\{  \left( \begin{smallmatrix}
1 &  & c \\
 & 1 &  \\
 &  & i
\end{smallmatrix} \right) 
\in \PGL_3(R) \right\}$.

\begin{itemize}[leftmargin=25pt]
    \item[-] $\lambda xy + \mu z^2$ is $E_{1,1}^{}$-adapted and  $\Aut_X^0(R)$ acts as $[\lambda:\mu] \mapsto [\lambda: i^2 \mu]$
\end{itemize}

\noindent 
Note that this is the first case in which we could not choose an $\Aut_X^0$-stable $E_{1,1}$-adapted pencil of conics. Therefore, to obtain the above description of the $\Aut_X^0$-action on $E_{1,1}$, Remark \ref{R obwohlsnichtadaptedist} must be applied. This will happen increasingly often from this point on, and we will no longer mention that we are applying Remark \ref{R obwohlsnichtadaptedist}. 

Since $\Aut_X^0$ has two orbits on $E \cap E_{1,1}$, we get the following two possibilities for $p_{1,2}$ up to isomorphism:

\begin{enumerate}[leftmargin=*]

\noindent     \begin{minipage}{0.65\textwidth}
    \item \vspace{2mm}
    $p_{1,2}= E_{1,1}^{} \cap C^{(2)}$ with $C=\cal{V}(xy+z^2)$
    \begin{itemize}[leftmargin=20pt]
        \item
        $\Aut_{X'}^0(R) =
        \begin{cases}   
        \left\{  \left( \begin{smallmatrix}
1 &  & c \\
 & 1 &  \\
 &  & 1
\end{smallmatrix} \right) 
\in \PGL_3(R) \right\} & \text{ if } p \neq 2 \\
        \left\{  \left( \begin{smallmatrix}
1 &  & c \\
 & 1 &  \\
 &  & i
\end{smallmatrix} \right) 
\in \PGL_3(R) \bigg| i^2=1 \right\}
         & \text{ if } p=2
        \end{cases} $ 
        \\We describe the configurations of negative curves on $X'$ for $p \neq 2$ and $p=2$ simultaneously:
        \item $(-2)$-curves: $E_{1,0}^{(3)},E_{1,1}^{(3)},
        \ell_{z}^{(3)}$
        \item $(-1)$-curves: $E_{1,2}^{}, E_{2,0}^{(3)},
        E_{3,0}^{(3)},
        \ell_{y}^{(3)}$
        \item
        with configuration as in Figure \ref{Conf4E4M}.
    \end{itemize} 
    This is case \hyperref[Tab4E]{$4E$} if $p \neq 2$, and case \hyperref[Tab4M]{$4M$} if $p=2$.
        \end{minipage} \hspace{2mm} \begin{minipage}{0.3\textwidth} \begin{center} \includegraphics[width=0.8\textwidth]{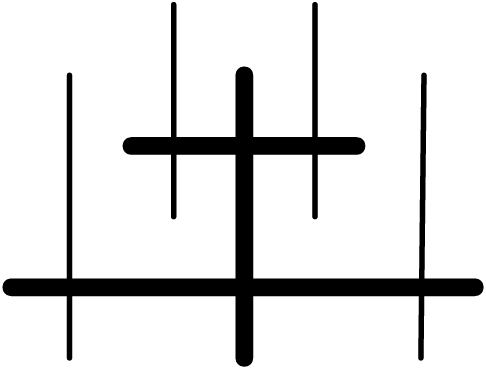} \vspace{-3.5mm}\captionof{figure}{} \label{Conf4E4M}  \end{center} \end{minipage}

    \noindent     \begin{minipage}{0.65\textwidth}
       \item \vspace{2mm}
    $p_{1,2}= E_{1,1}^{} \cap \ell_{y}^{(2)}$
    \begin{itemize}[leftmargin=20pt]
        \item
        $
        \Aut_{X'}^0(R) =    
        \left\{  \left( \begin{smallmatrix}
1 &  & c \\
 & 1 &  \\
 &  & i
\end{smallmatrix} \right) 
\in \PGL_3(R) \right\}
        $
        \item $(-2)$-curves: $E_{1,0}^{(3)},E_{1,1}^{(3)},
        \ell_{y}^{(3)},
        \ell_{z}^{(3)}$
        \item $(-1)$-curves: $E_{1,2}^{}, E_{2,0}^{(3)},
        E_{3,0}^{(3)}$
        \item
        with configuration as in Figure \ref{Conf4H}.
    \end{itemize} 
    This is case \hyperref[Tab4H]{$4H$}.
        \end{minipage} \hspace{2mm} \begin{minipage}{0.3\textwidth} \begin{center} \includegraphics[width=0.7\textwidth]{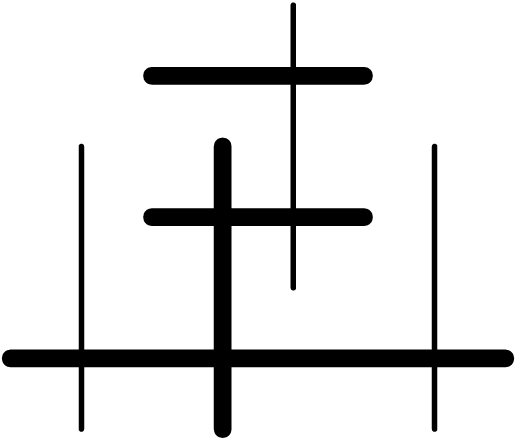} \vspace{-3.5mm}\captionof{figure}{}  \label{Conf4H} \end{center} \end{minipage} 
    
\end{enumerate}

\vspace{2mm} \subsubsection*{\underline{Case \hyperref[Tab3H]{$3H$}}}
We have 
$
E = \bigcup_{j=1}^{3} E_{j,1}^{}- (\bigcup_{j=1}^{3} E_{j,0}^{(2)} \cup \ell_x^{(2)} \cup \ell_y^{(2)} \cup \ell_z^{(2)})
$ and $\Aut_X^0(R) = \left\{  \left( \begin{smallmatrix}
1 &  &  \\
 & e &  \\
 &  & i
\end{smallmatrix} \right) 
\in \PGL_3(R) \right\}$.

\begin{itemize}[leftmargin=25pt]
    \item[-] $\lambda xz + \mu y^2$ is $E_{1,1}^{}$-adapted and  $\Aut_X^0(R)$ acts as $[\lambda:\mu] \mapsto [i \lambda: e^2 \mu]$
    \item[-] $\lambda xy + \mu z^2$ is $E_{2,1}^{}$-adapted and  $\Aut_X^0(R)$ acts as $[\lambda:\mu] \mapsto [e \lambda: i^2 \mu]$
    \item[-] $\lambda yz + \mu x^2$ is $E_{3,1}^{}$-adapted and  $\Aut_X^0(R)$ acts as $[\lambda:\mu] \mapsto [ei \lambda: \mu]$
\end{itemize}

\noindent 
Note that all automorphisms of $\bbP^2$ inducing cyclic permutations of $p_{1,0},p_{2,0},$ and $p_{3,0}$ lift to automorphisms of $X$ and since $X$ has degree $3$, we can only blow up two additional points. Moreover, $\Aut_X^0$ acts transitively on every $E \cap E_{j,1}$.
Hence, we get the following two possibilities for $p_{1,2},\hdots,p_{3,2}$ up to isomorphism:

\begin{enumerate}[leftmargin=*]
\noindent     \begin{minipage}{0.65\textwidth}
    \item \vspace{2mm}
    $p_{1,2}  =E_{1,1}^{} \cap C_1^{(2)}, 
    p_{2,2} = E_{2,1}^{} \cap C_2^{(2)}$ with $C_1= \cal{V}(xz+y^2),
    \\C_2=\cal{V}(xy+z^2)$
    \begin{itemize}[leftmargin=20pt]
        \item
        $\Aut_{X'}^0(R) =
        \begin{cases}   
        \{\id\} & \text{ if } p \neq 3 \\
        \left\{  \left( \begin{smallmatrix}
1 &  &  \\
 & e &  \\
 &  & e^2
\end{smallmatrix} \right) 
\in \PGL_3(R) \bigg| e^3=1 \right\}
         & \text{ if } p=3
        \end{cases} $
        \\Hence, $X'$ has global vector fields only if $p=3$. Therefore, we assume $p=3$ when describing the configuration of negative curves.
        \item $(-2)$-curves: $E_{1,0}^{(3)}, E_{2,0}^{(3)}, E_{3,0}^{(3)},E_{1,1}^{(3)}, E_{2,1}^{(3)}, \ell_{x}^{(3)}, \ell_{y}^{(3)}, \ell_z^{(3)}$
        \item $(-1)$-curves: $E_{1,2}^{}, E_{2,2}^{}, E_{3,1}^{(3)}$
        \item
        with configuration as in Figure \ref{Conf1G}.
    \end{itemize}
    This is case \hyperref[Tab1G]{$1G$}.
        \end{minipage} \hspace{2mm} \begin{minipage}{0.3\textwidth} \begin{center} \includegraphics[width=0.98\textwidth]{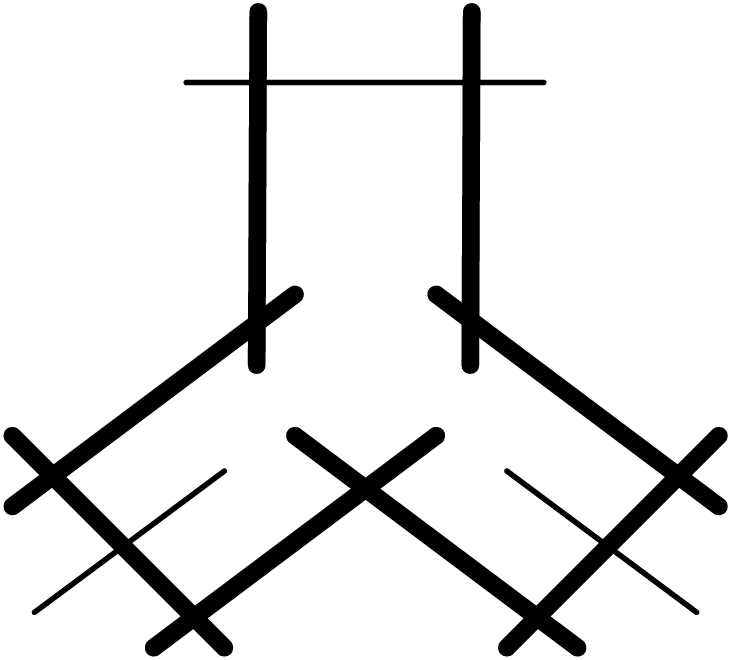} \vspace{-3.5mm}\captionof{figure}{}  \label{Conf1G} \end{center} \end{minipage}

   \item \vspace{2mm}
    $p_{1,2}  =E_{1,1}^{} \cap C^{(2)}$ with $C=\cal{V}(xz+y^2)$
    \begin{itemize}[leftmargin=20pt]
    \noindent \begin{minipage}{0.5\textwidth}
        \item
        $
        \Aut_{X'}^0(R) =    
        \left\{  \left( \begin{smallmatrix}
1 &  &  \\
 & e &  \\
 &  & e^2
\end{smallmatrix} \right) 
\in \PGL_3(R) \right\}
        $
        \item $(-2)$-curves: $E_{1,0}^{(3)}, E_{2,0}^{(3)}, E_{3,0}^{(3)},E_{1,1}^{(3)}, \ell_{x}^{(3)}, \ell_{y}^{(3)}, \ell_z^{(3)}$
        
        \end{minipage} \noindent \begin{minipage}{0.5\textwidth}
        
        \item $(-1)$-curves: $E_{1,2}^{}, E_{2,1}^{(3)}, E_{3,1}^{(3)}$
        \item
        with configuration as in Figure \ref{Conf2F}, that is, as in \\ case \hyperref[Tab2F]{$2F$}.
    \end{minipage}
    \end{itemize}
   \vspace{1mm} \noindent As explained in Remark \ref{R IsomorphismCheck}, one can check that $X' \cong X_{2F}$.
    
\end{enumerate}

\newpage

\vspace{2mm} \subsubsection*{\underline{Case \hyperref[Tab4G]{$4G$}}}
We have 
$
E = \bigcup_{j=1}^{2} E_{j,1}^{}- (\bigcup_{j=1}^{2} E_{j,0}^{(2)} \cup \ell_z^{(2)} \cup \ell_x^{(2)})
$ and $\Aut_X^0(R) = \left\{  \left( \begin{smallmatrix}
1 &  &  \\
 & e &  \\
 &  & i
\end{smallmatrix} \right) 
\in \PGL_3(R) \right\}$.

\begin{itemize}[leftmargin=25pt]
    \item[-] $\lambda xz + \mu y^2$ is $E_{1,1}^{}$-adapted and  $\Aut_X^0(R)$ acts as $[\lambda:\mu] \mapsto [i \lambda: e^2 \mu]$
    \item[-] $\lambda xy + \mu z^2$ is $E_{2,1}^{}$-adapted and  $\Aut_X^0(R)$ acts as $[\lambda:\mu] \mapsto [e \lambda: i^2 \mu]$
\end{itemize}

\noindent 
Since $\Aut_X^0$ acts transitively on every $E \cap E_{j,1}$, we get the following three possibilities for $p_{1,2},p_{2,2}$ up to isomorphism:

\begin{enumerate}[leftmargin=*]

\noindent     \begin{minipage}{0.65\textwidth}
    \item \vspace{2mm}
     $p_{1,2}  =E_{1,1}^{} \cap C_1^{(2)}, 
    p_{2,2} = E_{2,1}^{} \cap C_2^{(2)}$ with $C_1= \cal{V}(xz+y^2), 
    \\C_2= \cal{V}(xy+z^2)$
    \begin{itemize}[leftmargin=20pt]
        \item
        $\Aut_{X'}^0(R) =
        \begin{cases}   
        \{\id\} & \text{ if } p \neq 3 \\
        \left\{  \left( \begin{smallmatrix}
1 &  &  \\
 & e &  \\
 &  & e^2
\end{smallmatrix} \right) 
\in \PGL_3(R) \bigg| e^3=1 \right\}
         & \text{ if } p=3
        \end{cases} $
        \\Hence, $X'$ has global vector fields only if $p=3$. Therefore, we assume $p=3$ when describing the configuration of negative curves.
        \item $(-2)$-curves: $E_{1,0}^{(3)}, E_{2,0}^{(3)}, E_{1,1}^{(3)}, E_{2,1}^{(3)}, \ell_{x}^{(3)},  \ell_z^{(3)}$
        \item $(-1)$-curves: $E_{1,2}^{}, E_{2,2}^{}, E_{3,0}^{(3)}, \ell_y^{(3)}$
        \item
        with configuration as in Figure \ref{Conf2J}.
    \end{itemize} 
    This is case \hyperref[Tab2J]{$2J$}.
        \end{minipage} \hspace{2mm} \begin{minipage}{0.3\textwidth} \begin{center} \includegraphics[width=0.65\textwidth]{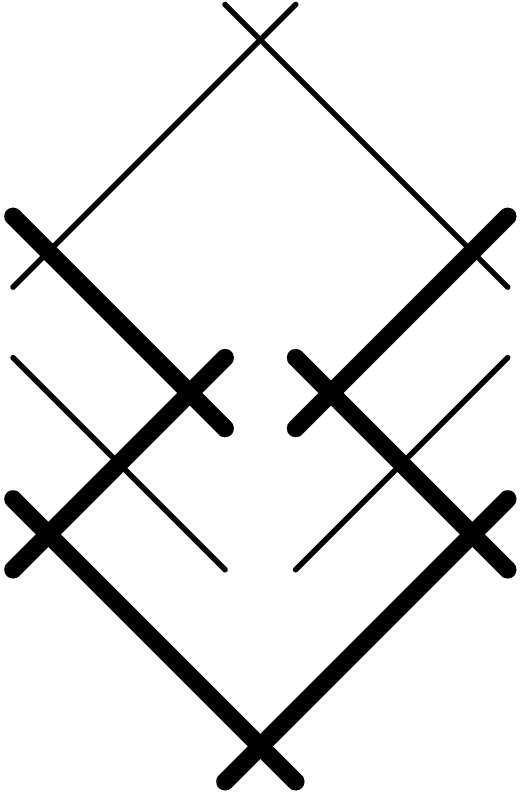} \vspace{-3.5mm}\captionof{figure}{}  \label{Conf2J} \end{center} \end{minipage}

   \item \vspace{2mm}
     $p_{2,2} = E_{2,1}^{} \cap C^{(2)}$ with $C= \cal{V}(xy+z^2)$
    \begin{itemize}[leftmargin=20pt]
    \noindent \begin{minipage}{0.5\textwidth}
        \item
        $
        \Aut_{X'}^0(R) =    
        \left\{  \left( \begin{smallmatrix}
1 &  &  \\
 & i^2 &  \\
 &  & i
\end{smallmatrix} \right) 
\in \PGL_3(R) \right\}
        $
        \item $(-2)$-curves: $E_{1,0}^{(3)}, E_{2,0}^{(3)}, E_{2,1}^{(3)}, \ell_{x}^{(3)},  \ell_z^{(3)}$
        
        \end{minipage} \noindent \begin{minipage}{0.5\textwidth}
        
        \item $(-1)$-curves: $E_{2,2}^{}, E_{1,1}^{(3)}, E_{3,0}^{(3)}, \ell_y^{(3)}$
        \item
        with configuration as in Figure \ref{Conf3E}, that is, as in \\ case \hyperref[Tab3E]{$3E$}.
    \end{minipage}
    \end{itemize} 
    \vspace{1mm} \noindent As explained in Remark \ref{R IsomorphismCheck}, one can check that $X' \cong X_{3E}$.
    
   \item \vspace{2mm}
     $p_{1,2}  =E_{1,1}^{} \cap C^{(2)}$ with $C= \cal{V}(xz+y^2)$
    \begin{itemize}[leftmargin=20pt]
    \noindent \begin{minipage}{0.5\textwidth}
        \item
        $
        \Aut_{X'}^0(R) =    
        \left\{  \left( \begin{smallmatrix}
1 &  &  \\
 & e &  \\
 &  & e^2
\end{smallmatrix} \right) 
\in \PGL_3(R) \right\}
        $
        \item $(-2)$-curves: $E_{1,0}^{(3)}, E_{2,0}^{(3)}, E_{1,1}^{(3)}, \ell_{x}^{(3)},  \ell_z^{(3)}$
        
        \end{minipage} \noindent \begin{minipage}{0.5\textwidth}
        
        \item $(-1)$-curves: $E_{1,2}^{}, E_{2,1}^{(3)}, E_{3,0}^{(3)}, \ell_y^{(3)}$
        \item
        with configuration as in Figure \ref{Conf3E}, that is, as in \\ case \hyperref[Tab3E]{$3E$}.
    \end{minipage}
    \end{itemize} 
    \vspace{1mm} \noindent As explained in Remark \ref{R IsomorphismCheck}, one can check that $X' \cong X_{3E}$.
     
\end{enumerate}

\vspace{2mm} \subsubsection*{\underline{Case \hyperref[Tab4F]{$4F$}}}
We have 
$
E = (E_{1,1}^{} \cup E_{2,1})- (E_{1,0}^{(2)} \cup E_{2,0}^{(2)} \cup \ell_x^{(2)} \cup \ell_y^{(2)})
$ and $\Aut_X^0(R) = \left\{  \left( \begin{smallmatrix}
1 &  &  \\
 & e &  \\
 &  & i
\end{smallmatrix} \right) 
\in \PGL_3(R) \right\}$.

\begin{itemize}[leftmargin=25pt]
    \item[-] $\lambda xy + \mu z^2$ is $E_{1,1}^{}$-adapted and $E_{2,1}^{}$-adapted and $\Aut_X^0(R)$ acts as $[\lambda:\mu] \mapsto [e \lambda: i^2 \mu]$
\end{itemize}

Note that the involution $x \leftrightarrow y$ of $\bbP^2$ lifts to an involution of $X$ interchanging $E_{1,1}$ and $E_{2,1}$. Moreover, $\Aut_X^0$ acts transitively on both $E \cap E_{1,1}$ and $E \cap E_{2,1}$, but the stabilizer of every point on $E \cap E_{1,1}$ acts trivially on $E \cap E_{2,1}$. Hence, we have the following three possibilities up to isomorphism:

\begin{enumerate}[leftmargin=*]
   \item \vspace{2mm}
    $p_{1,2}  =E_{1,1}^{} \cap C^{(2)}, p_{2,2}=E_{2,1}^{} \cap C^{(2)}$ with $C= \cal{V}(xy+z^2)$
    \begin{itemize}[leftmargin=20pt]
    \noindent \begin{minipage}{0.5\textwidth}
        \item
        $
        \Aut_{X'}^0(R) =    
        \left\{  \left( \begin{smallmatrix}
1 &  &  \\
 & i^2 &  \\
 &  & i
\end{smallmatrix} \right) 
\in \PGL_3(R) \right\}
        $
       \item $(-2)$-curves: $E_{1,0}^{(3)},E_{2,0}^{(3)}, E_{1,1}^{(3)}, E_{2,1}^{(3)}, \ell_x^{(3)}, \ell_y^{(3)}, C^{(3)}$
        
        \end{minipage} \noindent \begin{minipage}{0.5\textwidth}
        
        \item $(-1)$-curves: $E_{1,2}^{}, E_{2,2}^{},  E_{3,0}^{(3)}, \ell_{z}^{(3)}$
        \item
        with configuration as in Figure \ref{Conf2D}, that is, as in \\ case \hyperref[Tab2D]{$2D$}.
    \end{minipage}
    \end{itemize} 
   \vspace{1mm} \noindent As explained in Remark \ref{R IsomorphismCheck}, one can check that $X' \cong X_{2D}$.

   \item \vspace{2mm}
    $p_{1,2}  =E_{1,1}^{} \cap C_1^{(2)}, p_{2,2}=E_{2,1}^{} \cap C_2^{(2)}$ with $C_1= \cal{V}(xy+z^2), C_2= \cal{V}(xy+ \alpha z^2), \alpha \not\in \{0, 1\}$
    \begin{itemize}[leftmargin=20pt]
    \noindent \begin{minipage}{0.5\textwidth}
        \item
        $
        \Aut_{X'}^0(R) =    
        \left\{  \left( \begin{smallmatrix}
1 &  &  \\
 & i^2 &  \\
 &  & i
\end{smallmatrix} \right) 
\in \PGL_3(R) \right\}
        $
       \item $(-2)$-curves: $E_{1,0}^{(3)},E_{2,0}^{(3)}, E_{1,1}^{(3)}, E_{2,1}^{(3)}, \ell_x^{(3)}, \ell_y^{(3)}$
        
        \end{minipage} \noindent \begin{minipage}{0.5\textwidth}
        
        \item $(-1)$-curves: $E_{1,2}^{}, E_{2,2}^{}, E_{3,0}^{(3)}, \ell_{z}^{(3)}, C_1^{(3)}, C_2^{(3)}$
        \item
        with configuration as in Figure \ref{Conf2A}, that is, as in \\ case \hyperref[Tab2A]{$2A$}.
    \end{minipage}
    \end{itemize} 
   \vspace{1mm} \noindent As explained in Remark \ref{R IsomorphismCheck}, one can check that $X' \cong X_{2A, \alpha'}$ for some $\alpha'$.
   
   \newpage
   
   \item \vspace{2mm}
    $p_{1,2}  =E_{1,1}^{} \cap C^{(2)}$ with $C= \cal{V}(xy+z^2)$
    \begin{itemize}[leftmargin=20pt]
    \noindent \begin{minipage}{0.5\textwidth}
        \item
        $
        \Aut_{X'}^0(R) =    
        \left\{  \left( \begin{smallmatrix}
1 &  &  \\
 & i^2 &  \\
 &  & i
\end{smallmatrix} \right) 
\in \PGL_3(R) \right\}
        $
        \item $(-2)$-curves: $E_{1,0}^{(3)},E_{2,0}^{(3)}, E_{1,1}^{(3)}, \ell_x^{(3)}, \ell_y^{(3)}$

        \end{minipage} \noindent \begin{minipage}{0.5\textwidth}
        
        \item $(-1)$-curves: $E_{1,2}^{}, E_{2,1}^{(3)}, E_{3,0}^{(3)}, \ell_{z}^{(3)}, C^{(3)}$
        \item
        with configuration as in Figure \ref{Conf3D}, that is, as in \\ case \hyperref[Tab3D]{$3D$}.
    \end{minipage}
    \end{itemize} 
   \vspace{1mm} \noindent As explained in Remark \ref{R IsomorphismCheck}, one can check that $X' \cong X_{3D}$.
   
\end{enumerate}

\vspace{2mm} \subsubsection*{\underline{Case \hyperref[Tab5B]{$5B$}}}
We have 
$
E = E_{1,1}^{}- (E_{1,0}^{(2)} \cup \ell_z^{(2)})
$ and $\Aut_X^0(R) = \left\{  \left( \begin{smallmatrix}
1 &  &  \\
 & e &  \\
 &  & i
\end{smallmatrix} \right) 
\in \PGL_3(R) \right\}$.

\begin{itemize}[leftmargin=25pt]
    \item[-] $\lambda xz + \mu y^2$ is $E_{1,1}^{}$-adapted and  $\Aut_X^0(R)$ acts as $[\lambda:\mu] \mapsto [i \lambda: e^2 \mu]$
\end{itemize}

\noindent 
Since $\Aut_X^0$ acts transitively on $E \cap E_{1,1}$, we have the following unique choice for $p_{1,2}$ up to isomorphism:

\begin{enumerate}[leftmargin=*]
   \item \vspace{2mm}
    $p_{1,2}  =E_{1,1}^{} \cap C^{(2)}$ with $C= \cal{V}(xz+y^2)$
    \begin{itemize}[leftmargin=20pt]
    \noindent \begin{minipage}{0.5\textwidth}
        \item
        $
        \Aut_{X'}^0(R) =    
        \left\{  \left( \begin{smallmatrix}
1 &  &  \\
 & e &  \\
 &  & e^2
\end{smallmatrix} \right) 
\in \PGL_3(R) \right\}
        $
        \item $(-2)$-curves: $E_{1,0}^{(3)}, E_{1,1}^{(3)},\ell_z^{(3)}$
        
        \end{minipage} \noindent \begin{minipage}{0.5\textwidth}
        
        \item $(-1)$-curves: $E_{1,2}^{}, E_{2,0}^{(3)}, E_{3,0}^{(3)}, \ell_x^{(3)}, \ell_y^{(3)}$
        \item
        with configuration as in Figure \ref{Conf4D}, that is, as in \\ case \hyperref[Tab4D]{$4D$}.
    \end{minipage}
    \end{itemize} 
   \vspace{1mm} \noindent As explained in Remark \ref{R IsomorphismCheck}, one can check that $X' \cong X_{4D}$.
\end{enumerate}

\vspace{2mm} \subsubsection*{\underline{Case \hyperref[Tab5D]{$5D$}}}
We have 
$
E = \bigcup_{j=1}^{2} E_{j,1}^{} - (\bigcup_{j=1}^{2} E_{j,0}^{(2)} \cup \ell_z^{(2)})
$ and $\Aut_X^0(R) = \left\{  \left( \begin{smallmatrix}
1 &  &  \\
 & e & f \\
 &  & i
\end{smallmatrix} \right) 
\in \PGL_3(R) \right\}$.

\begin{itemize}[leftmargin=25pt]
    \item[-] $\lambda xz + \mu y^2$ is $E_{1,1}^{}$-adapted and  $\Aut_X^0(R)$ acts as $[\lambda:\mu] \mapsto [i \lambda: e^2 \mu]$
    \item[-] $\lambda xy + \mu z^2$ is $E_{2,1}^{}$-adapted and  $\Aut_X^0(R)$ acts as $[\lambda:\mu] \mapsto [e \lambda: i^2 \mu]$
\end{itemize}

\noindent 
Note that $\Aut_X^0$ acts transitively on $E \cap E_{1,1}$, and with two orbits, one of which is a fixed point, on $E \cap E_{2,1}$.
Hence, we have the following five choices for $p_{1,2},p_{2,2}$ up to isomorphism:

\begin{enumerate}[leftmargin=*]

\noindent     \begin{minipage}{0.65\textwidth}
    \item \vspace{2mm}
    $p_{1,2}  =E_{1,1}^{} \cap C_1^{(2)}, 
    p_{2,2} = E_{2,1}^{} \cap C_2^{(2)}$ with $C_1= \cal{V}(xz+y^2), 
    \\C_2= \cal{V}(xy+z^2)$
    \begin{itemize}[leftmargin=20pt]
        \item
        $\Aut_{X'}^0(R) =
        \begin{cases}   
        \left\{  \left( \begin{smallmatrix}
1 &  &  \\
 & 1 & f \\
 &  & 1
\end{smallmatrix} \right) 
\in \PGL_3(R) \right\} & \text{ if } p \neq 3 \\
        \left\{  \left( \begin{smallmatrix}
1 &  &  \\
 & e & f \\
 &  & e^2
\end{smallmatrix} \right) 
\in \PGL_3(R) \bigg| e^3=1 \right\}
         & \text{ if } p=3
        \end{cases} $ 
        \\We describe the configurations of negative curves on $X'$ for $p \neq 3$ and $p=3$ simultaneously:
        \item $(-2)$-curves: $E_{1,0}^{(3)}, E_{2,0}^{(3)},E_{1,1}^{(3)}, E_{2,1}^{(3)}, \ell_z^{(3)}$
        \item $(-1)$-curves: $E_{1,2}^{}, E_{2,2}^{},  \ell_x^{(3)}$
        \item
        with configuration as in Figure \ref{Conf3F3K}.
    \end{itemize} 
    This is case \hyperref[Tab3F]{$3F$} if $p \neq 3$, and case \hyperref[Tab3K]{$3K$} if $p=3$.
        \end{minipage} \hspace{2mm} \begin{minipage}{0.3\textwidth} \begin{center} \includegraphics[width=0.98\textwidth]{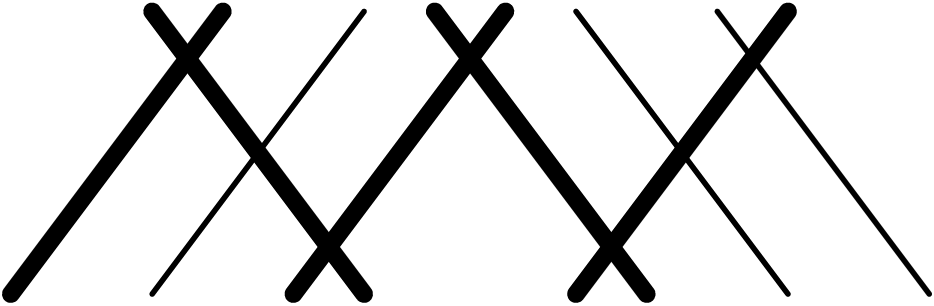} \vspace{-3.5mm}\captionof{figure}{}  \label{Conf3F3K} \end{center} \end{minipage}

    \noindent     \begin{minipage}{0.65\textwidth}
    \item \vspace{2mm}
    $p_{1,2}  =E_{1,1}^{} \cap C^{(2)}, 
    p_{2,2} = E_{2,1}^{} \cap \ell_{x}^{(2)}$ with $C= \cal{V}(xz+y^2)$
    \begin{itemize}[leftmargin=20pt]
        \item
        $
        \Aut_{X'}^0(R) =    
        \left\{  \left( \begin{smallmatrix}
1 &  &  \\
 & e & f \\
 &  & e^2
\end{smallmatrix} \right) 
\in \PGL_3(R) \right\}
        $
        \item $(-2)$-curves: $E_{1,0}^{(3)}, E_{2,0}^{(3)}, E_{1,1}^{(3)}, E_{2,1}^{(3)},  \ell_x^{(3)}, \ell_z^{(3)}$
        \item $(-1)$-curves: $E_{1,2}^{}, E_{2,2}^{}$
        \item
        with configuration as in Figure \ref{Conf3I}.
    \end{itemize} 
    This is case \hyperref[Tab3I]{$3I$}.
        \end{minipage} \hspace{2mm} \begin{minipage}{0.3\textwidth} \begin{center} \includegraphics[width=0.98\textwidth]{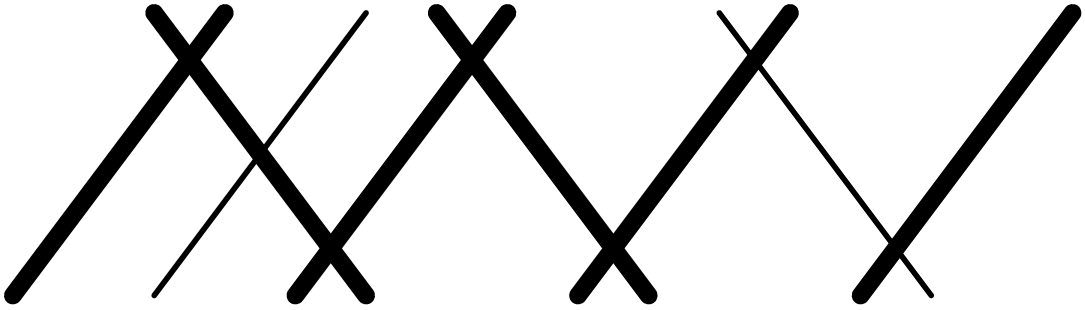} \vspace{-3.5mm}\captionof{figure}{}  \label{Conf3I} \end{center} \end{minipage}

   \item \vspace{2mm}
    $p_{2,2} = E_{2,1}^{} \cap C^{(2)}$ with $C=\cal{V}(xy+z^2)$
    \begin{itemize}[leftmargin=20pt]
    \noindent \begin{minipage}{0.5\textwidth}
        \item
        $
        \Aut_{X'}^0(R) =    
        \left\{  \left( \begin{smallmatrix}
1 &  &  \\
 & i^2 & f \\
 &  & i
\end{smallmatrix} \right) 
\in \PGL_3(R) \right\}
        $
        \item $(-2)$-curves: $E_{1,0}^{(3)}, E_{2,0}^{(3)},E_{2,1}^{(3)}, \ell_z^{(3)}$
        
        \end{minipage} \noindent \begin{minipage}{0.5\textwidth}
        
        \item $(-1)$-curves: $ E_{2,2}^{}, E_{1,1}^{(3)}, \ell_x^{(3)}$
        \item
        with configuration as in Figure \ref{Conf4H}, that is, as in \\ case \hyperref[Tab4H]{$4H$}.
    \end{minipage}
    \end{itemize} 
   \vspace{1mm} \noindent As explained in Remark \ref{R IsomorphismCheck}, one can check that $X' \cong X_{4H}$.

    \noindent     \begin{minipage}{0.65\textwidth}
    \item \vspace{2mm}
    $p_{2,2} = E_{2,1}^{} \cap \ell_{x}^{(2)}$
    \begin{itemize}[leftmargin=20pt]
        \item
        $
        \Aut_{X'}^0(R) =    
        \left\{  \left( \begin{smallmatrix}
1 &  &  \\
 & e & f \\
 &  & i
\end{smallmatrix} \right) 
\in \PGL_3(R) \right\}
        $
        \item $(-2)$-curves: $E_{1,0}^{(3)}, E_{2,0}^{(3)},E_{2,1}^{(3)},\ell_x^{(3)}, \ell_z^{(3)}$
        \item $(-1)$-curves: $E_{2,2}^{}, E_{1,1}^{(3)}$
        \item
        with configuration as in Figure \ref{Conf4K}.
    \end{itemize} 
    This is case \hyperref[Tab4K]{$4K$}.
        \end{minipage} \hspace{2mm} \begin{minipage}{0.3\textwidth} \begin{center} \includegraphics[width=0.8\textwidth]{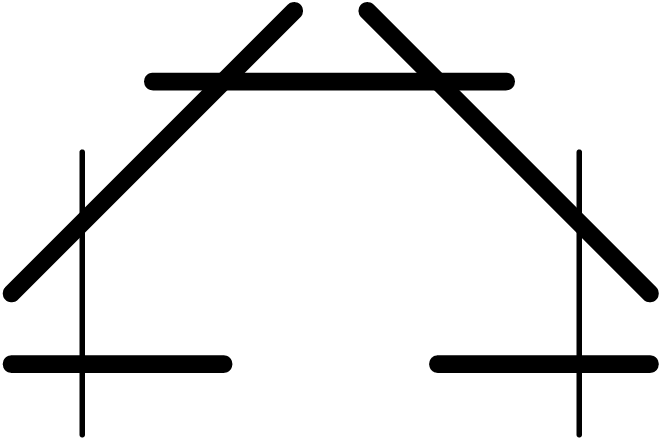} \vspace{-3.5mm}\captionof{figure}{}  \label{Conf4K} \end{center} \end{minipage}

    \noindent     \begin{minipage}{0.65\textwidth}
    \item \vspace{2mm}
    $p_{1,2}  =E_{1,1}^{} \cap C^{(2)}$ with $C= \cal{V}(xz+y^2)$
    \begin{itemize}[leftmargin=20pt]
        \item
        $
        \Aut_{X'}^0(R) =    
        \left\{  \left( \begin{smallmatrix}
1 &  &  \\
 & e & f \\
 &  & e^2
\end{smallmatrix} \right) 
\in \PGL_3(R) \right\}
        $
        \item $(-2)$-curves: $E_{1,0}^{(3)}, E_{2,0}^{(3)},E_{1,1}^{(3)}, \ell_z^{(3)}$
        \item $(-1)$-curves: $E_{1,2}^{}, E_{2,1}^{(3)},  \ell_x^{(3)}$
        \item
        with configuration as in Figure \ref{Conf4I}.
    \end{itemize} 
    This is case \hyperref[Tab4I]{$4I$}.
    \end{minipage} \hspace{2mm} \begin{minipage}{0.3\textwidth} \begin{center} \includegraphics[width=0.6\textwidth]{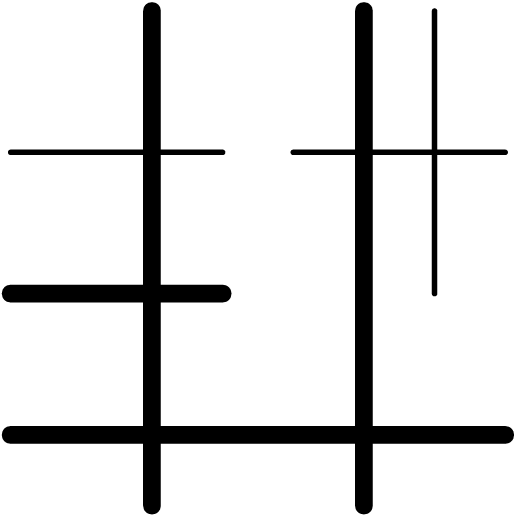} \vspace{-3.5mm}\captionof{figure}{}  \label{Conf4I} \end{center} \end{minipage}

\end{enumerate}

\vspace{2mm} \subsubsection*{\underline{Case \hyperref[Tab6B]{$6B$}}}
We have 
$
E = E_{1,1}^{} - E_{1,0}^{(2)}
$ and $\Aut_X^0(R) = \left\{  \left( \begin{smallmatrix}
1 &  & c \\
 & e &  \\
 &  & i
\end{smallmatrix} \right) 
\in \PGL_3(R) \right\}$.

\begin{itemize}[leftmargin=25pt]
    \item[-] $\lambda xy + \mu z^2$ is $E_{1,1}^{}$-adapted and  $\Aut_X^0(R)$ acts as $[\lambda:\mu] \mapsto [e \lambda: i^2 \mu]$
\end{itemize}

\noindent 
 Since $\Aut_X^0$ has two orbits on $E \cap E_{1,1}$, we have the following two choices for $p_{1,2}$ up to isomorphism:

\begin{enumerate}[leftmargin=*]

    \item \vspace{2mm}
    $p_{1,2}  =E_{1,1}^{} \cap C^{(2)}$ with $C=\cal{V}(xy+z^2)$
    \begin{itemize}[leftmargin=20pt]
    \noindent \begin{minipage}{0.5\textwidth}
        \item
        $
        \Aut_{X'}^0(R) =    
        \left\{  \left( \begin{smallmatrix}
1 &  & c \\
 & i^2 &  \\
 &  & i
\end{smallmatrix} \right) 
\in \PGL_3(R) \right\}
        $
        \item $(-2)$-curves: $E_{1,0}^{(3)},E_{1,1}^{(3)}$
        
        \end{minipage} \noindent \begin{minipage}{0.5\textwidth}
        
        \item $(-1)$-curves: $E_{1,2}^{}, E_{2,0}^{(3)},  \ell_y^{(3)}, \ell_z^{(3)}$
        \item
        with configuration as in Figure \ref{Conf5C}, that is, as in \\ case \hyperref[Tab5C]{$5C$}.
    \end{minipage}
    \end{itemize} 
    \vspace{1mm} \noindent As explained in Remark \ref{R IsomorphismCheck}, one can check that $X' \cong X_{5C}$.
    
    \item \vspace{2mm}
    $p_{1,2}  =E_{1,1}^{} \cap \ell_{y}^{(2)}$
    \begin{itemize}[leftmargin=20pt]
    \noindent \begin{minipage}{0.5\textwidth}
        \item
        $
        \Aut_{X'}^0(R) =    
        \left\{  \left( \begin{smallmatrix}
1 &  & c \\
 & e &  \\
 &  & i
\end{smallmatrix} \right) 
\in \PGL_3(R) \right\}
        $
        \item $(-2)$-curves: $E_{1,0}^{(3)}, E_{1,1}^{(3)},  \ell_y^{(3)}$
        
        \end{minipage} \noindent \begin{minipage}{0.5\textwidth}
        
        \item $(-1)$-curves: $E_{1,2}^{}, E_{2,0}^{(3)}, \ell_z^{(3)}$
        \item
        with configuration as in Figure \ref{Conf5D}, that is, as in \\ case \hyperref[Tab5D]{$5D$}.
    \end{minipage}
    \end{itemize} 
    \vspace{1mm} \noindent As explained in Remark \ref{R IsomorphismCheck}, one can check that $X' \cong X_{5D}$.
\end{enumerate}

\vspace{2mm} \subsubsection*{\underline{Case \hyperref[Tab6D]{$6D$}}}
We have 
$
E = E_{1,1}^{} - (E_{1,0}^{(2)} \cup \ell_z^{(2)})
$ and $\Aut_X^0(R) = \left\{  \left( \begin{smallmatrix}
1 &  & c \\
 & e & f \\
 &  & i
\end{smallmatrix} \right) 
\in \PGL_3(R) \right\}$.

\begin{itemize}[leftmargin=25pt]
    \item[-] $\lambda xz + \mu y^2$ is $E_{1,1}^{}$-adapted and  $\Aut_X^0(R)$ acts as $[\lambda:\mu] \mapsto [i \lambda: e^2 \mu]$
\end{itemize}

\noindent 
 Since $\Aut_X^0$ acts transitively on $E \cap E_{1,1}$, there is only one choice for $p_{1,2}$ up to isomorphism:

\begin{enumerate}[leftmargin=*]
\noindent     \begin{minipage}{0.65\textwidth}
    \item \vspace{2mm}
    $p_{1,2}  =E_{1,1}^{} \cap C^{(2)}$ with $C= \cal{V}(xz+y^2)$
    \begin{itemize}[leftmargin=20pt]
        \item
        $
        \Aut_{X'}^0(R) =    
        \left\{  \left( \begin{smallmatrix}
1 &  & c \\
 & e & f \\
 &  & e^2
\end{smallmatrix} \right) 
\in \PGL_3(R) \right\}
        $
        \item $(-2)$-curves: $E_{1,0}^{(3)}, E_{1,1}^{(3)},\ell_z^{(3)}$
        \item $(-1)$-curves: $E_{1,2}^{}, E_{2,0}^{(3)}$
        \item
        with configuration as in Figure \ref{Conf5E}.
    \end{itemize} 
    This is case \hyperref[Tab5E]{$5E$}.
        \end{minipage} \hspace{2mm} \begin{minipage}{0.3\textwidth} \begin{center} \includegraphics[width=0.8\textwidth]{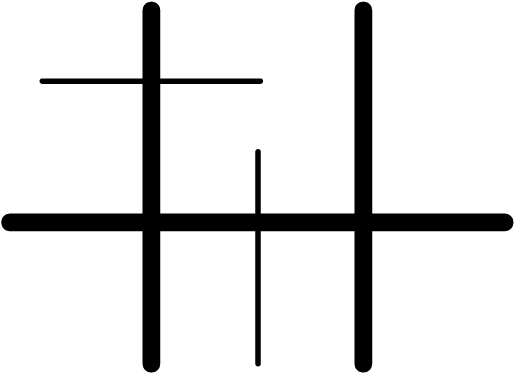} \vspace{-3.5mm}\captionof{figure}{}  \label{Conf5E}  \end{center} \end{minipage} 
\end{enumerate}

\vspace{2mm} \subsubsection*{\underline{Case \hyperref[Tab7B]{$7B$}}}
We have 
$
E = E_{1,1}^{} - E_{1,0}^{(2)}
$ and $\Aut_X^0(R) = \left\{  \left( \begin{smallmatrix}
1 & b & c \\
 & e & f \\
 &  & i
\end{smallmatrix} \right) 
\in \PGL_3(R) \right\}$.

\begin{itemize}[leftmargin=25pt]
    \item[-] $\lambda xz + \mu y^2$ is $E_{1,1}^{}$-adapted and  $\Aut_X^0(R)$ acts as $[\lambda:\mu] \mapsto [i \lambda: e^2 \mu]$
\end{itemize}

\noindent 
 Since $\Aut_X^0$ has two orbits on $E \cap E_{1,1}$, there are the following two choices for $p_{1,2}$ up to isomorphism:

\begin{enumerate}[leftmargin=*]

\noindent     \begin{minipage}{0.65\textwidth}
\item \vspace{2mm}
    $p_{1,2}  =E_{1,1}^{} \cap C^{(2)}$ with $C=\cal{V}(xz+y^2)$
    \begin{itemize}[leftmargin=20pt]
        \item
        $
        \Aut_{X'}^0(R) =    
        \left\{  \left( \begin{smallmatrix}
1 & b & c \\
 & e & f \\
 &  & e^2
\end{smallmatrix} \right) 
\in \PGL_3(R) \right\}
        $
        \item $(-2)$-curves: $E_{1,0}^{(3)}, E_{1,1}^{(3)}$
        \item $(-1)$-curves: $E_{1,2}^{}, \ell_z^{(3)}$
        \item
        with configuration as in Figure \ref{Conf6E}.
    \end{itemize} 
    This is case \hyperref[Tab6E]{$6E$}.
        \end{minipage} \hspace{2mm} \begin{minipage}{0.3\textwidth} \begin{center} \includegraphics[width=0.95\textwidth]{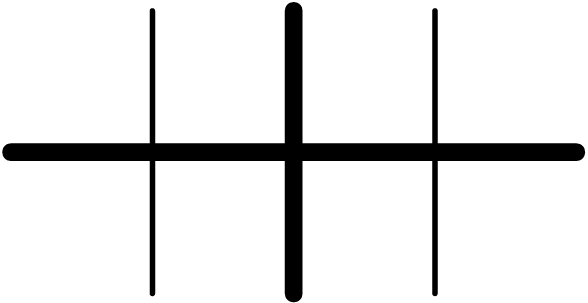} \vspace{-3.5mm}\captionof{figure}{} \label{Conf6E}  \end{center} \end{minipage} 
    
    \noindent     \begin{minipage}{0.65\textwidth}
    \item \vspace{2mm}
    $p_{1,2}  =E_{1,1}^{} \cap \ell_{z}^{(2)}$
    \begin{itemize}[leftmargin=20pt]
        \item
        $
        \Aut_{X'}^0(R) =    
        \left\{  \left( \begin{smallmatrix}
1 & b & c \\
 & e & f \\
 &  & i
\end{smallmatrix} \right) 
\in \PGL_3(R) \right\}
        $
        \item $(-2)$-curves: $E_{1,0}^{(3)}, E_{1,1}^{(3)}, \ell_z^{(3)}$
        \item $(-1)$-curves: $E_{1,2}^{}$
        \item
        with configuration as in Figure \ref{Conf6F}.
    \end{itemize} 
    This is case \hyperref[Tab6F]{$6F$}.    \end{minipage} \hspace{2mm} \begin{minipage}{0.3\textwidth} \begin{center} \includegraphics[width=0.95\textwidth]{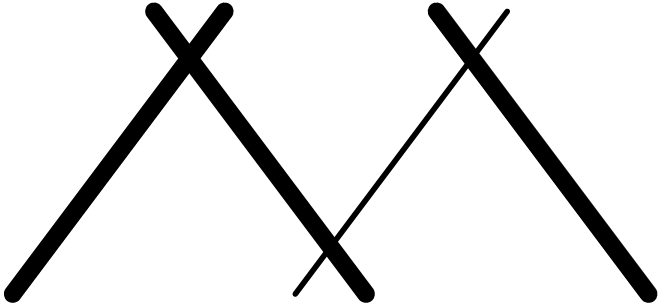} \vspace{-3.5mm}\captionof{figure}{}  \label{Conf6F} \end{center} \end{minipage}

\end{enumerate}

\vspace{2mm} 
\noindent
Summarizing, we obtain
\begin{eqnarray*}
\cal{L}_3 &=& \{  X_{1L,\alpha}, X_{1O}, X_{1N}, X_{2N,\alpha}, X_{2Q}, X_{2P}, X_{3N}, X_{1K}, X_{1J,\alpha}, X_{1B}, X_{2O}, X_{2B}, X_{2F}, \\ &&  X_{3E}, X_{4E}, X_{4M}, X_{4H}, X_{1G}, X_{2J}, X_{3F}, X_{3K}, X_{3I}, X_{4K}, X_{4I}, X_{5E}, X_{6E}, X_{6F}  \}.
\end{eqnarray*}

\subsection{Height 4}

\vspace{2mm} \subsubsection*{\underline{Case \hyperref[Tab2N]{$2N$}}}
This case exists only if $p = 2$. 
\\ \noindent We have $E =
E_{1,2}^{} - E_{1,1}^{(3)}$ and $\Aut_X^0(R) = \left\{  \left( \begin{smallmatrix}
1 &  &  \\
 & 1 &  \\
 &  & i
\end{smallmatrix} \right) 
\in \PGL_3(R) \bigg| i^2=1 \right\}$.

\begin{itemize}[leftmargin=25pt]
    \item[-] $\lambda (x^2y + xz^2)+ \mu z^3$ is $E_{1,2}^{}$-adapted and $\Aut_X^0(R)$ acts as $[\lambda:\mu] \mapsto [\lambda: i \mu]$.
\end{itemize}

\noindent
Note that there is only one point on $E \cap E_{1,2}^{}$ with non-trivial stabilizer, hence we have the following unique choice for $p_{1,3}$:

\begin{enumerate}[leftmargin=*]
   \item \vspace{2mm}
    $p_{1,3}= E_{1,2}^{} \cap C_1^{(3)}$ with $C_1= \cal{V}(xy+z^2) $
    \begin{itemize}[leftmargin=20pt]
        \item
        $\Aut_{X'}^0(R) =
        \left\{  \left( \begin{smallmatrix}
1 &  &  \\
 & 1 &  \\
 &  & i
\end{smallmatrix} \right) 
\in \PGL_3(R) \bigg| i^2=1 \right\}$
        \item $(-2)$-curves: $E_{1,0}^{(4)},E_{3,0}^{(4)}, E_{1,1}^{(4)},E_{1,2}^{(4)}, \ell_{y}^{(4)}, \ell_{z}^{(4)}$
        \item $(-1)$-curves: $E_{1,3}^{}, E_{3,1}^{(4)}, E_{2,0}^{(4)}, E_{4,0}^{(4)}, \ell_{x}^{(4)}, \ell_{x-y}^{(4)}, \ell_{x+\alpha y}^{(4)}, C_1^{(4)}, C_2^{(4)}, C_3^{(4)}, C_4^{(4)}, C_5^{(4)}, C_6^{(4)}$ with
        \\$C_2= \cal{V}(xy+y^2+z^2), C_3= \cal{V}(x^2y+xz^2+ \alpha y z^2), C_4= \cal{V}(x^2y+xz^2+y^3+ \alpha y z^2), 
        \\C_5= \cal{V}(x^2y^2+ x^2z^2+ x^3y+ \alpha^2 y^2z^2), C_6= \cal{V}(xy^3+ x^2z^2+ x^3y+ \alpha^2 y^2z^2), \alpha \not\in \{0,-1\}$
        \item
        with configuration as in Figure \ref{Conf1J}, that is, as in case \hyperref[Tab1J]{$1J$}.
    \end{itemize} 
  \noindent As explained in Remark \ref{R IsomorphismCheck}, one can check that $X' \cong X_{1J,\alpha'}$ for some $\alpha'$.
\end{enumerate}

\vspace{2mm} \subsubsection*{\underline{Case \hyperref[Tab2Q]{$2Q$}}}
This case exists only if $p = 2$. 
\\ \noindent We have $E =
E_{1,2}^{} - E_{1,1}^{(3)}$ and $\Aut_X^0(R) = \left\{  \left( \begin{smallmatrix}
1 &  &  \\
 & 1 &  \\
 &  & i
\end{smallmatrix} \right) 
\in \PGL_3(R) \bigg| i^2=1 \right\}$.

\begin{itemize}[leftmargin=25pt]
    \item[-] $\lambda (x^2y + xz^2) + \mu z^3$ is $E_{1,2}^{}$-adapted and $\Aut_X^0(R)$ acts as $[\lambda:\mu] \mapsto [\lambda: i \mu]$.
\end{itemize}

\noindent
Note that there is only one point on $E \cap E_{1,2}^{}$ with non-trivial stabilizer, hence we have the following unique choice for $p_{1,3}$:

\begin{enumerate}[leftmargin=*]
    \item \vspace{2mm}
    $p_{1,3}= E_{1,2}^{} \cap C_1^{(3)}$ with $C_1= \cal{V}(xy+z^2) $
    \begin{itemize}[leftmargin=20pt]
    \noindent \begin{minipage}{0.5\textwidth}
        \item
        $
        \Aut_{X'}^0(R) =    
        \left\{  \left( \begin{smallmatrix}
1 &  &  \\
 & 1 &  \\
 &  & i
\end{smallmatrix} \right) 
\in \PGL_3(R) \bigg| i^2=1 \right\}
        $
        \item $(-2)$-curves: $E_{1,0}^{(4)},E_{3,0}^{(4)},E_{1,1}^{(4)},E_{1,2}^{(4)}, \ell_{x}^{(4)},
        \ell_{y}^{(4)}, \ell_{z}^{(4)}$
        
        \end{minipage} \noindent \begin{minipage}{0.5\textwidth}
        
        \item $(-1)$-curves: $E_{1,3}^{}, E_{3,1}^{(4)}, E_{2,0}^{(4)}, E_{4,0}^{(4)},  \ell_{x-y}^{(4)}, C_1^{(4)}, 
        \\ C_2^{(4)}, C_3^{(4)}$ with $C_2= \cal{V}(xy+y^2+z^2), 
        \\ C_3= \cal{V}(xz^2+x^2y+y^3)$
        \item
        with configuration as in Figure \ref{Conf1K}, that is, as in \\ case \hyperref[Tab1K]{$1K$}.
    \end{minipage}
    \end{itemize} 
  \vspace{1mm} \noindent As explained in Remark \ref{R IsomorphismCheck}, one can check that $X' \cong X_{1K}$.
\end{enumerate}

\vspace{2mm} \subsubsection*{\underline{Case \hyperref[Tab2P]{$2P$}}}
This case exists only if $p = 2$.  
\\ \noindent We have $E =
E_{1,2}^{}  - E_{1,1}^{(3)}$ and $\Aut_X^0(R) = \left\{  \left( \begin{smallmatrix}
1 &  &  \\
 & 1 &  \\
 &  & i
\end{smallmatrix} \right) 
\in \PGL_3(R) \bigg| i^2=1 \right\}$.

\begin{itemize}[leftmargin=25pt]
    \item[-] $\lambda (x^2y + xz^2) + \mu z^3$ is $E_{1,2}^{}$-adapted and $\Aut_X^0(R)$ acts as $[\lambda:\mu] \mapsto [i \lambda: \mu]$.
\end{itemize}
\noindent Note that there is only one point on $E \cap E_{1,2}^{}$ with non-trivial stabilizer, hence we have the following unique choice for $p_{1,3}$:

\begin{enumerate}[leftmargin=*]
   \item \vspace{2mm}
    $p_{1,3}= E_{1,2}^{} \cap C_1^{(3)}$ with $C_1= \cal{V}(xy+z^2) $
    \begin{itemize}[leftmargin=20pt]
    \noindent \begin{minipage}{0.5\textwidth}
        \item
        $\Aut_{X'}^0(R) =
        \left\{  \left( \begin{smallmatrix}
1 &  &  \\
 & 1 &  \\
 &  & i
\end{smallmatrix} \right) 
\in \PGL_3(R) \bigg| i^2=1 \right\}$
        \item $(-2)$-curves: $E_{1,0}^{(4)},E_{2,0}^{(4)},E_{1,1}^{(4)},E_{1,2}^{(4)}, \ell_{x}^{(4)},
        \ell_{y}^{(4)}, 
        \\ \ell_{z}^{(4)}, C_1^{(4)}$
        
        \end{minipage} \noindent \begin{minipage}{0.5\textwidth}
        
        \item $(-1)$-curves: $E_{1,3}^{}, E_{2,1}^{(4)}, E_{3,0}^{(4)}, E_{4,0}^{(4)},  \ell_{x-y}^{(4)},  C_2^{(4)}$ 
        \\with $C_2=\cal{V}(xy+y^2+z^2)$
        \item
        with configuration as in Figure \ref{Conf1N}, that is, as in \\ case \hyperref[Tab1N]{$1N$}.
    \end{minipage}
    \end{itemize} 
  \vspace{1mm} \noindent As explained in Remark \ref{R IsomorphismCheck}, one can check that $X' \cong X_{1N}$.
\end{enumerate}

\vspace{2mm} \subsubsection*{\underline{Case \hyperref[Tab3N]{$3N$}}}
This case exists only if $p = 2$.  
\\ \noindent We have $E =
E_{1,2}^{} - E_{1,1}^{(3)}$ and $\Aut_X^0(R) = \left\{  \left( \begin{smallmatrix}
1 &  &  \\
 & 1 &  \\
 &  & i
\end{smallmatrix} \right) 
\in \PGL_3(R) \bigg| i^2=1 \right\}$.

\begin{itemize}[leftmargin=25pt]
    \item[-] $\lambda (x^2y + xz^2) + \mu z^3$ is $E_{1,2}^{}$-adapted and $\Aut_X^0(R)$ acts as $[\lambda:\mu] \mapsto [\lambda: i \mu]$.
\end{itemize}
\noindent Note that there is only one point on $E \cap E_{1,2}^{}$ with non-trivial stabilizer, hence we have the following unique choice for $p_{1,3}$:

\begin{enumerate}[leftmargin=*]
   \item \vspace{2mm}
    $p_{1,3}= E_{1,2}^{} \cap C_1^{(3)}$ with $C_1= \cal{V}(xy+z^2) $
    \begin{itemize}[leftmargin=20pt]
    \noindent \begin{minipage}{0.5\textwidth}
        \item
        $\Aut_{X'}^0(R) =
        \left\{  \left( \begin{smallmatrix}
1 &  &  \\
 & 1 &  \\
 &  & i
\end{smallmatrix} \right) 
\in \PGL_3(R) \bigg| i^2=1 \right\}$
        \item $(-2)$-curves: $E_{1,0}^{(4)},E_{1,1}^{(4)},E_{1,2}^{(4)},
        \ell_{y}^{(4)},
        \ell_{z}^{(4)}$
        
        \end{minipage} \noindent \begin{minipage}{0.5\textwidth}
        
        \item $(-1)$-curves: $E_{1,3}^{}, E_{2,0}^{(4)},
        E_{3,0}^{(4)},
        E_{4,0}^{(4)},    
        \ell_{x}^{(4)},
        \ell_{x-y}^{(4)}, 
        \\ C_1^{(4)}, C_2^{(4)}$ with $C_2= \cal{V}(xy+y^2+z^2)$
        \item
        with configuration as in Figure \ref{Conf2O}, that is, as in \\ case \hyperref[Tab2O]{$2O$}.
    \end{minipage}
    \end{itemize} 
 \vspace{1mm} \noindent As explained in Remark \ref{R IsomorphismCheck}, one can check that $X' \cong X_{2O}$.
\end{enumerate}

\vspace{2mm} \subsubsection*{\underline{Case \hyperref[Tab2O]{$2O$}}}
This case exists only if $p = 2$.  
\\ \noindent We have $E =
E_{1,2}^{} - E_{1,1}^{(3)}$ and $\Aut_X^0(R) = \left\{  \left( \begin{smallmatrix}
1 &  &  \\
 & 1 &  \\
 &  & i
\end{smallmatrix} \right) 
\in \PGL_3(R) \bigg| i^2=1 \right\}$.

\begin{itemize}[leftmargin=25pt]
    \item[-] $\lambda (x^2y + xz^2) + \mu z^3$ is $E_{1,2}^{}$-adapted and $\Aut_X^0(R)$ acts as $[\lambda:\mu] \mapsto [\lambda: i \mu]$.
\end{itemize}
\noindent Note that there is only one point on $E \cap E_{1,2}^{}$ with non-trivial stabilizer, hence we have the following unique choice for $p_{1,3}$:

\begin{enumerate}[leftmargin=*]

    \item \vspace{2mm}
$p_{1,3}= E_{1,2}^{} \cap C_1^{(3)}$ with $C_1= \cal{V}(xy+z^2)$
    \begin{itemize}[leftmargin=20pt]
    \noindent \begin{minipage}{0.5\textwidth}
        \item
        $\Aut_{X'}^0(R) =
        \left\{  \left( \begin{smallmatrix}
1 &  &  \\
 & 1 &  \\
 &  & i
\end{smallmatrix} \right) 
\in \PGL_3(R) \bigg| i^2=1 \right\}$
        \item $(-2)$-curves: $E_{1,0}^{(4)},
        E_{2,0}^{(4)},
        E_{3,0}^{(4)}, E_{1,1}^{(4)},E_{1,2}^{(4)},
        \ell_{z}^{(4)}, 
        \\C_1^{(4)}, C_2^{(4)}$ with $C_2= \cal{V}(xy+y^2+z^2)$
        
        \end{minipage} \noindent \begin{minipage}{0.5\textwidth}
        
        \item $(-1)$-curves: $E_{1,3}^{}, E_{2,1}^{(4)},
        E_{3,1}^{(4)},
        \ell_{x}^{(4)},
        \ell_{y}^{(4)},
        \ell_{x-y}^{(4)}$
        \item
        with configuration as in Figure \ref{Conf1N}, that is, as in \\ case \hyperref[Tab1N]{$1N$}.
    \end{minipage}
    \end{itemize} 
  \vspace{1mm} \noindent As explained in Remark \ref{R IsomorphismCheck}, one can check that $X' \cong X_{1N}$.
\end{enumerate}

\vspace{2mm} \subsubsection*{\underline{Case \hyperref[Tab2B]{$2B$}}}
We have $E =
E_{1,2}^{} - (E_{1,1}^{(3)} \cup \ell_y^{(3)})$ and $\Aut_X^0(R) = \left\{  \left( \begin{smallmatrix}
1 &  &  \\
 & 1 &  \\
 &  & i
\end{smallmatrix} \right) 
\in \PGL_3(R) \right\}$.

\begin{itemize}[leftmargin=25pt]
    \item[-] $\lambda x^2y + \mu z^3$ is $E_{1,2}^{}$-adapted and $\Aut_X^0(R)$ acts as $[\lambda:\mu] \mapsto [\lambda: i^3 \mu]$.
\end{itemize}
\noindent Hence, we have the following unique choice for $p_{1,3}$ up to isomorphism:

\begin{enumerate}[leftmargin=*]
\noindent     \begin{minipage}{0.65\textwidth}
   \item \vspace{2mm}
$p_{1,3}= E_{1,2}^{} \cap C^{(3)}$ with $C=\cal{V}(x^2y+z^3)$
    \begin{itemize}[leftmargin=20pt]
        \item
        $\Aut_{X'}^0(R) =
        \begin{cases}   
        \{\id\} & \text{ if } p \neq 3 \\
        \left\{  \left( \begin{smallmatrix}
1 &  &  \\
 & 1 &  \\
 &  & i
\end{smallmatrix} \right) 
\in \PGL_3(R) \bigg| i^3=1 \right\}
         & \text{ if } p=3
        \end{cases} $
        \\Hence, $X'$ has global vector fields only if $p=3$. Therefore, we assume $p=3$ when describing the configuration of negative curves.
        \item $(-2)$-curves: $E_{1,0}^{(4)},
        E_{2,0}^{(4)},
        E_{3,0}^{(4)}, E_{1,1}^{(4)},E_{1,2}^{(4)},
        \ell_{y}^{(4)},
        \ell_{z}^{(4)}$
        \item $(-1)$-curves: $E_{1,3}^{}, E_{2,1}^{(4)},
        E_{3,1}^{(4)},
        \ell_{x}^{(4)},
        \ell_{x-y}^{(4)}$
        \item
        with configuration as in Figure \ref{Conf1E}.
    \end{itemize} 
    This is case \hyperref[Tab1E]{$1E$}.
        \end{minipage} \hspace{2mm} \begin{minipage}{0.3\textwidth} \begin{center} \includegraphics[width=0.98\textwidth]{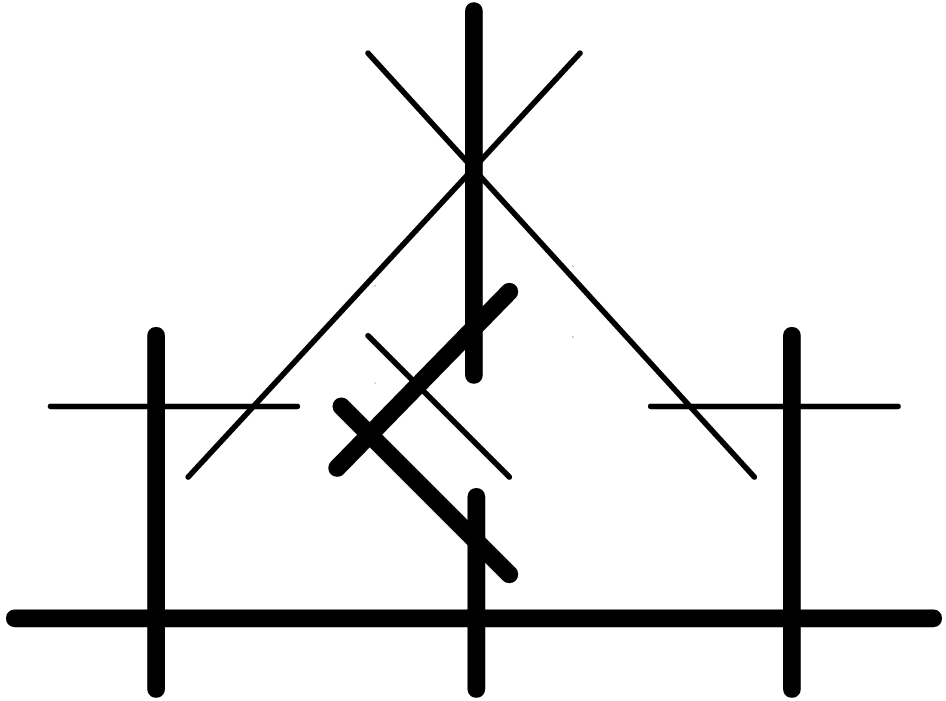} \vspace{-3.5mm}\captionof{figure}{}  \label{Conf1E} \end{center} \end{minipage} 
\end{enumerate}

\vspace{2mm} \subsubsection*{\underline{Case \hyperref[Tab2F]{$2F$}}}
We have $E =
(E_{1,2}^{} \cup E_{2,2}^{}) - (E_{1,1}^{(3)} \cup E_{2,1}^{(3)} \cup \ell_{x}^{(3)} \cup \ell_y^{(3)})$ and $\Aut_X^0(R) = \left\{  \left( \begin{smallmatrix}
1 &  &  \\
 & 1 &  \\
 &  & i
\end{smallmatrix} \right) 
\in \PGL_3(R) \right\}$.

\begin{itemize}[leftmargin=25pt]
    \item[-] $\lambda x^2y + \mu z^3$ is $E_{1,2}^{}$-adapted and $\Aut_X^0(R)$ acts as $[\lambda:\mu] \mapsto [\lambda: i^3 \mu]$.
    \item[-] $\lambda xy^2 + \mu z^3$ is $E_{2,2}^{}$-adapted and $\Aut_X^0(R)$ acts as $[\lambda:\mu] \mapsto [\lambda: i^3 \mu]$.
\end{itemize}
\noindent Note that the involution $x \leftrightarrow y$ of $\bbP^2$ lifts to an automorphism of $X$ interchanging $E_{1,2}^{}$ and $E_{2,2}^{}$. Moreover, since $X$ has degree $2$, we are only allowed to blow up one more point. Hence, we have the following unique choice for $p_{1,3},p_{2,3}$ up to isomorphism:

\begin{enumerate}[leftmargin=*]
    \item \vspace{2mm}
$p_{1,3}= E_{1,2}^{} \cap C^{(3)}$ with $C=\cal{V}(x^2y+z^3)$
    \begin{itemize}[leftmargin=20pt]
        \item
        $\Aut_{X'}^0(R) =
        \begin{cases}   
        \{\id\} & \text{ if } p \neq 3 \\
        \left\{  \left( \begin{smallmatrix}
1 &  &  \\
 & 1 &  \\
 &  & i
\end{smallmatrix} \right) 
\in \PGL_3(R) \bigg| i^3=1 \right\}
         & \text{ if } p=3
        \end{cases} $
        \\Hence, $X'$ has global vector fields only if $p=3$. Therefore, we assume $p=3$ when describing the configuration of negative curves.
        \item $(-2)$-curves: $E_{1,0}^{(4)},
        E_{2,0}^{(4)},E_{1,1}^{(4)}, E_{2,1}^{(4)},E_{1,2}^{(4)},
        \ell_{x}^{(4)},
        \ell_{y}^{(4)},
        \ell_{z}^{(4)}$
        \item $(-1)$-curves: $E_{1,3}^{}, E_{2,2}^{(4)},
        E_{3,0}^{(4)}$
        \item
        with configuration as in Figure \ref{Conf1G}, that is, as in case \hyperref[Tab1G]{$1G$}.
    \end{itemize} 
\noindent As explained in Remark \ref{R IsomorphismCheck}, one can check that $X' \cong X_{1G}$.
\end{enumerate}

\vspace{2mm} \subsubsection*{\underline{Case \hyperref[Tab3E]{$3E$}}}
We have $E =
E_{1,2}^{} - (E_{1,1}^{(3)} \cup \ell_y^{(3)})$ and $\Aut_X^0(R) = \left\{  \left( \begin{smallmatrix}
1 &  &  \\
 & 1 &  \\
 &  & i
\end{smallmatrix} \right) 
\in \PGL_3(R) \right\}$.

\begin{itemize}[leftmargin=25pt]
    \item[-] $\lambda x^2y + \mu z^3$ is $E_{1,2}^{}$-adapted and $\Aut_X^0(R)$ acts as $[\lambda:\mu] \mapsto [\lambda: i^3 \mu]$.
\end{itemize}
Hence, we have the following unique choice for $p_{1,3}$ up to isomorphism:

\newpage

\begin{enumerate}[leftmargin=*]
   \item \vspace{2mm}
$p_{1,3}= E_{1,2}^{} \cap C^{(3)}$ with $C=\cal{V}(x^2y+z^3)$
    \begin{itemize}[leftmargin=20pt]
        \item
        $\Aut_{X'}^0(R) =
        \begin{cases}   
        \{\id\} & \text{ if } p \neq 3 \\
        \left\{  \left( \begin{smallmatrix}
1 &  &  \\
 & 1 &  \\
 &  & i
\end{smallmatrix} \right) 
\in \PGL_3(R) \bigg| i^3=1 \right\}
         & \text{ if } p=3
        \end{cases} $
        \\Hence, $X'$ has global vector fields only if $p=3$. Therefore, we assume $p=3$ when describing the configuration of negative curves.
        \item $(-2)$-curves: $E_{1,0}^{(4)},
        E_{2,0}^{(4)},E_{1,1}^{(4)},E_{1,2}^{(4)},
        \ell_{y}^{(4)},
        \ell_{z}^{(4)}$
        \item $(-1)$-curves: $E_{1,3}^{}, E_{2,1}^{(4)},
        E_{3,0}^{(4)},
        \ell_{x}^{(4)}$
        \item
        with configuration as in Figure \ref{Conf2J}, that is, as in case \hyperref[Tab2J]{$2J$}.
    \end{itemize} 
\noindent As explained in Remark \ref{R IsomorphismCheck}, one can check that $X' \cong X_{2J}$.
\end{enumerate}

\vspace{2mm} \subsubsection*{\underline{Case \hyperref[Tab4E]{$4E$}}}
This case exists only if $p \neq 2$.  
\\ \noindent  
We have $E =
E_{1,2}^{} - E_{1,1}^{(3)}$ and $\Aut_X^0(R) = \left\{  \left( \begin{smallmatrix}
1 &  & c \\
 & 1 &  \\
 &  & 1
\end{smallmatrix} \right) 
\in \PGL_3(R) \right\}$.

\begin{itemize}[leftmargin=25pt]
    \item[-] $\lambda (x^2y + xz^2) + \mu z^3$ is $E_{1,2}^{}$-adapted and $\Aut_X^0(R)$ acts as 
    $[\lambda:\mu] \mapsto [\lambda: \mu - c \lambda ]$
\end{itemize}
In particular, the stabilizer of every point on $E \cap E_{1,2}^{}$ is trivial, hence this case does not lead to additional weak del Pezzo surfaces with global vector fields.

\vspace{2mm} \subsubsection*{\underline{Case \hyperref[Tab4M]{$4M$}}}
This case exists only if $p = 2$.  
\\ \noindent We have $E =
E_{1,2}^{} - E_{1,1}^{(3)}$ and $\Aut_X^0(R) = \left\{  \left( \begin{smallmatrix}
1 &  & c \\
 & 1 &  \\
 &  & i
\end{smallmatrix} \right) 
\in \PGL_3(R) \bigg| i^2=1 \right\}$.

\begin{itemize}[leftmargin=25pt]
    \item[-] $\lambda (x^2y + xz^2) + \mu z^3$ is $E_{1,2}^{}$-adapted and $\Aut_X^0(R)$ acts as
    $[\lambda:\mu] \mapsto [\lambda: i\mu + c \lambda ]$    
\end{itemize}
In particular, $\Aut_X^0$ acts transitively on $E \cap E_{1,2}$, so there is the following unique possibility for $p_{1,3}$ up to isomorphism:

\begin{enumerate}[leftmargin=*]
   \item \vspace{2mm}
    $p_{1,3}= E_{1,2}^{} \cap C_1^{(3)}$ with $C_1=\cal{V}(x y +  z^2)$
    \begin{itemize}[leftmargin=20pt]
    \noindent \begin{minipage}{0.5\textwidth}
        \item
        $\Aut_{X'}^0(R) =
        \left\{  \left( \begin{smallmatrix}
1 &  &  \\
 & 1 &  \\
 &  & i
\end{smallmatrix} \right) 
\in \PGL_3(R) \bigg| i^2=1 \right\}
       $ 
        \item $(-2)$-curves: $E_{1,0}^{(4)},E_{1,1}^{(4)},E_{1,2}^{(4)},
        \ell_{z}^{(4)}$
        
        \end{minipage} \noindent \begin{minipage}{0.5\textwidth}

        \item $(-1)$-curves: $E_{1,3}^{}, E_{2,0}^{(4)},
        E_{3,0}^{(4)},
        \ell_{y}^{(4)}, C_1^{(4)}, C_2^{(4)}$ 
        \\with $C_2= \cal{V}(xy+y^2+z^2)$
        \item
        with configuration as in Figure \ref{Conf3N}, that is, as in \\ case \hyperref[Tab3N]{$3N$}.
    \end{minipage}
    \end{itemize} 
\vspace{1mm} \noindent As explained in Remark \ref{R IsomorphismCheck}, one can check that $X' \cong X_{3N}$.
\end{enumerate}

\vspace{2mm} \subsubsection*{\underline{Case \hyperref[Tab4H]{$4H$}}}
We have $E =
E_{1,2}^{} - (E_{1,1}^{(3)} \cup \ell_y^{(3)})$ and $\Aut_X^0(R) = \left\{  \left( \begin{smallmatrix}
1 &  & c \\
 & 1 &  \\
 &  & i
\end{smallmatrix} \right) 
\in \PGL_3(R) \right\}$.

\begin{itemize}[leftmargin=25pt]
    \item[-] $\lambda x^2y + \mu z^3$ is $E_{1,2}^{}$-adapted 
    and $\Aut_X^0(R)$ acts as
    $[\lambda:\mu] \mapsto [\lambda: i^3\mu]$  
\end{itemize}
\noindent Since $\Aut_X^0$ acts transitively on $E \cap E_{1,2}$, there is the following unique possibility for $p_{1,3}$ up to isomorphism:

\begin{enumerate}[leftmargin=*]
    \item \vspace{2mm}
    $p_{1,3}= E_{1,2}^{} \cap C^{(3)}$ with $C= \cal{V}(x^2y+z^3)$
    \begin{itemize}[leftmargin=20pt]
        \item
        $\Aut_{X'}^0(R) =
        \begin{cases}   
        \left\{  \left( \begin{smallmatrix}
1 &  & c \\
 & 1 &  \\
 &  & 1
\end{smallmatrix} \right) 
\in \PGL_3(R) \right\} & \text{ if } p \neq 3 \\
        \left\{  \left( \begin{smallmatrix}
1 &  & c \\
 & 1 &  \\
 &  & i
\end{smallmatrix} \right) 
\in \PGL_3(R) \bigg| i^3=1 \right\}
         & \text{ if } p=3
        \end{cases} $ 
        \\We describe the configurations of negative curves on $X'$ for $p \neq 3$ and $p=3$ simultaneously:
        \item $(-2)$-curves: $E_{1,0}^{(4)},E_{1,1}^{(4)},E_{1,2}^{(4)},
        \ell_{y}^{(4)},
        \ell_{z}^{(4)}$
        \item $(-1)$-curves: $E_{1,3}^{}, E_{2,0}^{(4)},
        E_{3,0}^{(4)}$
        \item
        with configuration as in Figure \ref{Conf3F3K}, that is, as in case \hyperref[Tab3F]{$3F$} or \hyperref[Tab3K]{$3K$}.
    \end{itemize} 
    \noindent As explained in Remark \ref{R IsomorphismCheck}, one can check that $X' \cong X_{3F}$ if $p \neq 3$, and $X' \cong X_{3K}$ if $p = 3$.
\end{enumerate}

\vspace{2mm} \subsubsection*{\underline{Case \hyperref[Tab2J]{$2J$}}}
This case exists only if $p=3$.  
\\ \noindent We have $E =
(E_{1,2}^{} \cup E_{2,2}^{}) - (E_{1,1}^{(3)} \cup E_{2,1}^{(3)})$ and $\Aut_X^0(R) = \left\{  \left( \begin{smallmatrix}
1 &  &  \\
 & e &  \\
 &  & e^2
\end{smallmatrix} \right) 
\in \PGL_3(R) \bigg| e^3 = 1 \right\}$.

\begin{itemize}[leftmargin=25pt]
    \item[-] $\lambda (x^2z + xy^2) + \mu y^3$ is $E_{1,2}^{}$-adapted and $\Aut_X^0(R)$ acts as $[\lambda:\mu] \mapsto [e^2 \lambda: \mu]$.
    \item[-] $\lambda (xy^2 + yz^2) + \mu z^3$ is $E_{2,2}^{}$-adapted and $\Aut_X^0(R)$ acts as $[\lambda:\mu] \mapsto [e^2 \lambda: \mu]$.
    
\end{itemize}
\noindent Note that $X$ has degree $2$, hence we are only allowed to blow up one more point. Moreover, there is a unique point on $E \cap E_{1,2}^{}$ and on $E \cap E_{2,2}^{}$ with non-trivial stabilizer. Therefore, we have the following two possibilities for $p_{1,3}$ and $p_{2,3}$:

\begin{enumerate}[leftmargin=*]

   \item \vspace{2mm}
     $ p_{2,3} = E_{2,2}^{} \cap C^{(3)}$ with $C= \cal{V}(xy+z^2)$
    \begin{itemize}[leftmargin=20pt]
    \noindent \begin{minipage}{0.5\textwidth}
        \item
        $\Aut_{X'}^0(R) =
        \left\{  \left( \begin{smallmatrix}
1 &  &  \\
 & e &  \\
 &  & e^2
\end{smallmatrix} \right) 
\in \PGL_3(R) \bigg| e^3=1 \right\}
        $
        \item $(-2)$-curves: $E_{1,0}^{(4)}, E_{2,0}^{(4)}, E_{1,1}^{(4)}, E_{2,1}^{(4)},E_{2,2}^{(4)}, \ell_{x}^{(4)},  \ell_z^{(4)}$
        
        \end{minipage} \noindent \begin{minipage}{0.5\textwidth}
        
        \item $(-1)$-curves: $ E_{2,3}^{}, E_{1,2}^{(4)}, E_{3,0}^{(4)}, \ell_y^{(4)}, C^{(4)}$
        \item
        with configuration as in Figure \ref{Conf1E}, that is, as in \\ case \hyperref[Tab1E]{$1E$}.
    \end{minipage}
    \end{itemize} 
\vspace{1mm} \noindent As explained in Remark \ref{R IsomorphismCheck}, one can check that $X' \cong X_{1E}$.

    \item \vspace{2mm}
     $p_{1,3}  =E_{1,2}^{} \cap C^{(3)}$ with $C= \cal{V}(xz+y^2)$
    \begin{itemize}[leftmargin=20pt]
    \noindent \begin{minipage}{0.5\textwidth}
        \item
        $\Aut_{X'}^0(R) =
        \left\{  \left( \begin{smallmatrix}
1 &  &  \\
 & e &  \\
 &  & e^2
\end{smallmatrix} \right) 
\in \PGL_3(R) \bigg| e^3=1 \right\}
        $
        \item $(-2)$-curves: $E_{1,0}^{(4)}, E_{2,0}^{(4)}, E_{1,1}^{(4)}, E_{2,1}^{(4)}, E_{1,2}^{(4)},\ell_{x}^{(4)},  \ell_z^{(4)}$
        
        \end{minipage} \noindent \begin{minipage}{0.5\textwidth}
        
        \item $(-1)$-curves: $E_{1,3}^{}, E_{2,2}^{(4)}, E_{3,0}^{(4)}, \ell_y^{(4)}, C^{(4)}$
       \item
        with configuration as in Figure \ref{Conf1E}, that is, as in \\ case \hyperref[Tab1E]{$1E$}.
    \end{minipage}
    \end{itemize} 
\vspace{1mm} \noindent As explained in Remark \ref{R IsomorphismCheck}, one can check that $X' \cong X_{1E}$.

\end{enumerate}

\vspace{2mm} \subsubsection*{\underline{Case \hyperref[Tab3F]{$3F$}}}
This case exists only if $p \neq 3$.  
\\ \noindent We have $E =
(E_{1,2}^{} \cup E_{2,2}^{}) - (E_{1,1}^{(3)} \cup E_{2,1}^{(3)})$ and $\Aut_X^0(R) = \left\{  \left( \begin{smallmatrix}
1 &  &  \\
 & 1 & f \\
 &  & 1
\end{smallmatrix} \right) 
\in \PGL_3(R) \right\}$.

\begin{itemize}[leftmargin=25pt]
    \item[-] $\lambda (x^2z + xy^2) + \mu y^3$ is $E_{1,2}^{}$-adapted and $\Aut_X^0(R)$ acts as
    $[\lambda:\mu] \mapsto [\lambda: \mu - 2 f \lambda ]$   
    \item[-] $\lambda (xy^2 + yz^2) + \mu z^3$ is $E_{2,2}^{}$-adapted and $\Aut_X^0(R)$ acts as
    $[\lambda:\mu] \mapsto [\lambda: \mu - f \lambda ]$    
\end{itemize}
If $p \neq 2$, then $\Aut_X^0$ acts simply transitively on both $E \cap E_{1,2}$ and $E \cap E_{2,2}$, hence we cannot blow up $X$ any further and still obtain a weak del Pezzo surface with global vector fields. If $p = 2$, then $\Aut_X^0$ still acts transitively on $E \cap E_{2,2}$, but now it acts trivially on $E \cap E_{1,2}$. This leads to the following possibilities for $p_{1,3}$:

\begin{enumerate}[leftmargin=*]
\noindent     \begin{minipage}{0.65\textwidth}
\item \vspace{2mm}
    $p_{1,3}  =E_{1,2}^{} \cap C^{(3)}$ with $C= \cal{V}(x^2z+xy^2+ \alpha y^3)$
    \begin{itemize}[leftmargin=20pt]
        \item
        $\Aut_{X'}^0(R) =
        \begin{cases}   
        \{ \id \} & \text{ if } p \neq 2, 3 \\
        \left\{  \left( \begin{smallmatrix}
1 &  &  \\
 & 1 & f \\
 &  & 1
\end{smallmatrix} \right) 
\in \PGL_3(R)  \right\}
         & \text{ if } p=2
        \end{cases} $ 
        \\Hence, $X'$ has global vector fields only if $p=2$. Therefore, we assume $p=2$ when describing the configuration of negative curves.
        \item $(-2)$-curves: $E_{1,0}^{(4)}, E_{2,0}^{(4)},E_{1,1}^{(4)}, E_{2,1}^{(4)}, E_{1,2}^{(4)},\ell_z^{(4)}$
        \item $(-1)$-curves: $E_{1,3}^{}, E_{2,2}^{(4)},  \ell_x^{(4)}$
        \item
        with configuration as in Figure \ref{Conf2K2R}.
    \end{itemize} 
    This is case \hyperref[Tab2R]{$2R$} and we see that we get a $1$-dimensional family of such surfaces $X_{2R,\alpha}$ depending on the parameter $\alpha$.
        \end{minipage} \hspace{2mm} \begin{minipage}{0.3\textwidth} \begin{center} \includegraphics[width=0.97\textwidth]{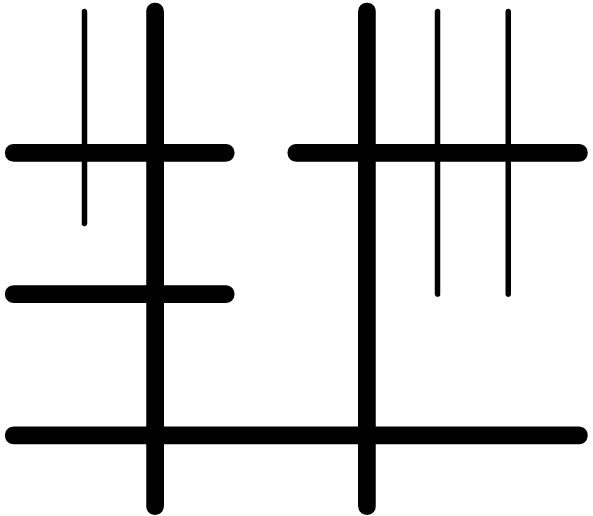} \vspace{-3.5mm}\captionof{figure}{}  \label{Conf2K2R} \end{center} \end{minipage} 
\end{enumerate}

\vspace{2mm} \subsubsection*{\underline{Case \hyperref[Tab3K]{$3K$}}}
This case exists only if $p = 3$.  
\\ \noindent We have $E =
(E_{1,2}^{} \cup E_{2,2}^{}) - (E_{1,1}^{(3)} \cup E_{2,1}^{(3)})$ and $\Aut_X^0(R) = \left\{  \left( \begin{smallmatrix}
1 &  &  \\
 & e & f \\
 &  & e^2
\end{smallmatrix} \right) 
\in \PGL_3(R) \bigg| e^3 = 1 \right\}$.

\begin{itemize}[leftmargin=25pt]
    \item[-] $\lambda (x^2z + xy^2) + \mu y^3$ is $E_{1,2}^{}$-adapted and $\Aut_X^0(R)$ acts as
    $[\lambda:\mu] \mapsto [e^2\lambda: \mu - 2ef \lambda ]$
    \item[-] $\lambda (xy^2 + yz^2) + \mu z^3$ is $E_{2,2}^{}$-adapted and $\Aut_X^0(R)$ acts as
    $[\lambda:\mu] \mapsto [e^2\lambda: \mu - ef \lambda ]$
\end{itemize}
Note that $\Aut_X^0$ acts transitively on both $E \cap E_{1,2}$ and $E \cap E_{2,2}$. The stabilizer of every point on $E \cap E_{1,2}$ is isomorphic to $\mu_3$ and this $\mu_3$ has a unique fixed point on $E \cap E_{2,2}$. This leads to the following three possibilities for $p_{1,3},p_{2,3}$ up to isomorphism:

\begin{enumerate}[leftmargin=*]

\item \vspace{2mm}
    $p_{1,3}  =E_{1,2}^{} \cap C_1^{(3)}, 
    p_{2,3} = E_{2,2}^{} \cap C_2^{(2)}$ with $C_1= \cal{V}(xz+y^2), C_2= \cal{V}(xy+z^2)$
    \begin{itemize}[leftmargin=20pt]
    \noindent \begin{minipage}{0.5\textwidth}
        \item
        $\Aut_{X'}^0(R) =
        \left\{  \left( \begin{smallmatrix}
1 &  &  \\
 & e &  \\
 &  & e^2
\end{smallmatrix} \right) 
\in \PGL_3(R) \bigg| e^3=1 \right\}
        $ 
        \item $(-2)$-curves: $E_{1,0}^{(4)}, E_{2,0}^{(4)},E_{1,1}^{(4)}, E_{2,1}^{(4)}, E_{1,2}^{(4)}, 
        \\E_{2,2}^{(4)},\ell_z^{(4)}$
        
        \end{minipage} \noindent \begin{minipage}{0.5\textwidth}
        
        \item $(-1)$-curves: $E_{1,3}^{}, E_{2,3}^{},  \ell_x^{(4)}, C_2^{(4)}, C_3^{(4)}$ with 
        \\$C_3= \cal{V}(x^2y^2+x^3z+z^4)$
        \item
        with configuration as in Figure \ref{Conf1E}, that is, as in \\ case \hyperref[Tab1E]{$1E$}.
    \end{minipage}
    \end{itemize} 
\vspace{1mm} \noindent As explained in Remark \ref{R IsomorphismCheck}, one can check that $X' \cong X_{1E}$.
    
    \item \vspace{2mm}
    $p_{2,3} = E_{2,2}^{} \cap C^{(3)}$ with $C= \cal{V}(xy+z^2)$
    \begin{itemize}[leftmargin=20pt]
    \noindent \begin{minipage}{0.5\textwidth}
        \item
        $\Aut_{X'}^0(R) =
        \left\{  \left( \begin{smallmatrix}
1 &  &  \\
 & e &  \\
 &  & e^2
\end{smallmatrix} \right) 
\in \PGL_3(R) \bigg| e^3=1 \right\}
        $ 
        \item $(-2)$-curves: $E_{1,0}^{(4)}, E_{2,0}^{(4)},E_{1,1}^{(4)}, E_{2,1}^{(4)}, E_{2,2}^{(4)},\ell_z^{(4)}$
        
        \end{minipage} \noindent \begin{minipage}{0.5\textwidth}
        
        \item $(-1)$-curves: $ E_{2,3}^{}, E_{1,2}^{(4)}, \ell_x^{(4)}, C^{(4)}$
        \item
        with configuration as in Figure \ref{Conf2J}, that is, as in \\ case \hyperref[Tab2J]{$2J$}.
    \end{minipage}
    \end{itemize} 
\vspace{1mm} \noindent As explained in Remark \ref{R IsomorphismCheck}, one can check that $X' \cong X_{2J}$.

    \item \vspace{2mm}
    $p_{1,3}  =E_{1,2}^{} \cap C^{(3)}$ with $C= \cal{V}(xz+y^2)$.

    \begin{itemize}[leftmargin=20pt]
    \noindent \begin{minipage}{0.5\textwidth}
        \item
        $\Aut_{X'}^0(R) =
        \left\{  \left( \begin{smallmatrix}
1 &  &  \\
 & e &  \\
 &  & e^2
\end{smallmatrix} \right) 
\in \PGL_3(R) \bigg| e^3=1 \right\}
        $ 
        \item $(-2)$-curves: $E_{1,0}^{(4)}, E_{2,0}^{(4)},E_{1,1}^{(4)}, E_{2,1}^{(4)},E_{1,2}^{(4)},  \ell_z^{(4)}$
        
        \end{minipage} \noindent \begin{minipage}{0.5\textwidth}
        
        \item $(-1)$-curves: $E_{1,3}^{}, E_{2,2}^{(4)},  \ell_x^{(4)}$
        \item
        with configuration as in Figure \ref{Conf2K2R}.
    \end{minipage}
    \end{itemize} 
    \vspace{1mm} \noindent This is case \hyperref[Tab2K]{$2K$}.

\end{enumerate}

\vspace{2mm} \subsubsection*{\underline{Case \hyperref[Tab3I]{$3I$}}}
We have $E =
(E_{1,2}^{} \cup E_{2,2}^{}) - (E_{1,1}^{(3)} \cup E_{2,1}^{(3)} \cup \ell_x^{(3)})$ and $\Aut_X^0(R) = \left\{  \left( \begin{smallmatrix}
1 &  &  \\
 & e & f \\
 &  & e^2
\end{smallmatrix} \right) 
\in \PGL_3(R) \right\}$.

\begin{itemize}[leftmargin=25pt]
    \item[-] $\lambda (x^2z + xy^2) + \mu y^3$ is $E_{1,2}^{}$-adapted and $\Aut_X^0(R)$ acts as
    $[\lambda:\mu] \mapsto [e^2\lambda: e^3\mu - 2ef \lambda ]$ 
    \item[-] $\lambda xy^2  + \mu z^3$ is $E_{2,2}^{}$-adapted 
    and $\Aut_X^0(R)$ acts as
    $[\lambda:\mu] \mapsto [e^2\lambda: e^6\mu]$  
\end{itemize}
\noindent 
Note that $\Aut_X^0$ acts transitively on $E \cap E_{2,2}$. If $p \neq 2$ (resp. $p = 2$), then $\Aut_X^0$ acts transitively (resp. with two orbits) on $E \cap E_{1,2}$.
We have the following five possibilities for $p_{1,3},p_{2,3}$ up to isomorphism:

\begin{enumerate}[leftmargin=*]
    
    \noindent     \begin{minipage}{0.65\textwidth}
    \item \vspace{2mm}
    $p_{1,3}  =E_{1,2}^{} \cap C_1^{(3)}, 
    p_{2,3} = E_{2,2}^{} \cap C_2^{(3)}$ with 
    \\$C_1= \cal{V}(x^2z+xy^2+y^3), C_2= \cal{V}(xy^2+ \alpha z^3), \alpha \neq 0$
    \begin{itemize}[leftmargin=20pt]
        \item
        $\Aut_{X'}^0(R) =
        \begin{cases}   
        \{ \id \} & \text{ if } p \neq 2 \\
        \left\{  \left( \begin{smallmatrix}
1 &  &  \\
 & 1 & f \\
 &  & 1
\end{smallmatrix} \right) 
\in \PGL_3(R)  \right\}
         & \text{ if } p=2
        \end{cases} $ 
        \\Hence, $X'$ has global vector fields only if $p=2$. Therefore, we assume $p=2$ when describing the configuration of negative curves.
        \item $(-2)$-curves: $E_{1,0}^{(4)}, E_{2,0}^{(4)}, E_{1,1}^{(4)}, E_{2,1}^{(4)},E_{1,2}^{(4)}, E_{2,2}^{(4)},  \ell_x^{(4)}, \ell_z^{(4)}$
        \item $(-1)$-curves: $E_{1,3}^{}, E_{2,3}^{}$
        \item
        with configuration as in Figure \ref{Conf1Q1R}.
    \end{itemize} 
    This is case \hyperref[Tab1Q]{$1Q$} and we see that we get a $1$-dimensional family of such surfaces $X_{1Q,\alpha}$ depending on the parameter $\alpha$. 
        \end{minipage} \hspace{2mm} \begin{minipage}{0.3\textwidth} \begin{center} \includegraphics[width=0.98\textwidth]{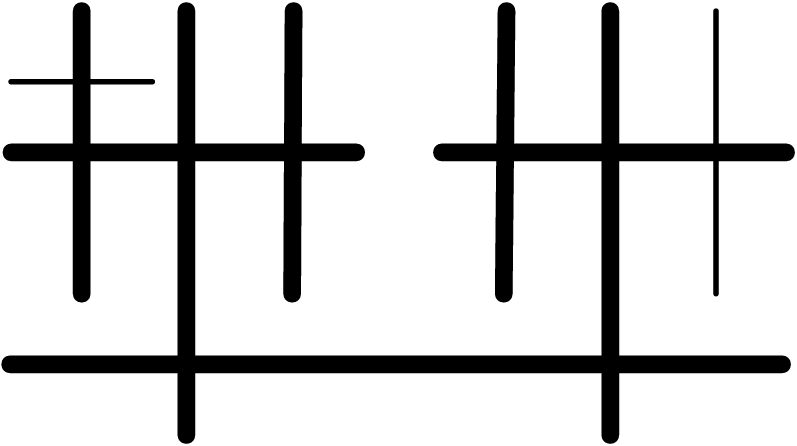} \vspace{-3.5mm}\captionof{figure}{}  \label{Conf1Q1R} \end{center} \end{minipage} 
    
    \newpage
    
    \item \vspace{2mm}
    $p_{1,3}  =E_{1,2}^{} \cap C_1^{(3)}, 
    p_{2,3} = E_{2,2}^{} \cap C_2^{(3)}$ with $C_1= \cal{V}(xz+y^2), C_2= \cal{V}(xy^2+z^3)$
    \begin{itemize}[leftmargin=20pt]
        \item
        $\Aut_{X'}^0(R) =
        \begin{cases}   
        \{ \id \} & \text{ if } p \neq 2 \\
        \left\{  \left( \begin{smallmatrix}
1 &  &  \\
 & e & f \\
 &  & e^2
\end{smallmatrix} \right) 
\in \PGL_3(R) \bigg| e^4=1 \right\}
         & \text{ if } p=2
        \end{cases} $ 
        \\Hence, $X'$ has global vector fields only if $p=2$. Therefore, we assume $p=2$ when describing the configuration of negative curves.
        \item $(-2)$-curves: $E_{1,0}^{(4)}, E_{2,0}^{(4)}, E_{1,1}^{(4)}, E_{2,1}^{(4)},E_{1,2}^{(4)}, E_{2,2}^{(4)},  \ell_x^{(4)}, \ell_z^{(4)}$
        \item $(-1)$-curves: $E_{1,3}^{}, E_{2,3}^{}$
        \item
        with configuration as in Figure \ref{Conf1Q1R}.
    \end{itemize} 
    This is case \hyperref[Tab1R]{$1R$}.

    \noindent     \begin{minipage}{0.65\textwidth}
    \item \vspace{2mm}
    $ p_{2,3} = E_{2,2}^{} \cap C^{(3)}$ with $C= \cal{V}(xy^2+z^3)$
    \begin{itemize}[leftmargin=20pt]
        \item
        $\Aut_{X'}^0(R) =
        \begin{cases}   
        \left\{  \left( \begin{smallmatrix}
1 &  &  \\
 & 1 & f \\
 &  & 1
\end{smallmatrix} \right) 
\in \PGL_3(R) \right\} & \text{ if } p \neq 2 \\
        \left\{  \left( \begin{smallmatrix}
1 &  &  \\
 & e & f \\
 &  & e^2
\end{smallmatrix} \right) 
\in \PGL_3(R) \bigg| e^4=1 \right\}
         & \text{ if } p=2
        \end{cases} $ 
        \\We describe the configurations of negative curves on $X'$ for $p \neq 2$ and $p=2$ simultaneously:
        \item $(-2)$-curves: $E_{1,0}^{(4)}, E_{2,0}^{(4)}, E_{1,1}^{(4)}, E_{2,1}^{(4)}, E_{2,2}^{(4)}, \ell_x^{(4)}, \ell_z^{(4)}$
        \item $(-1)$-curves: $E_{2,3}^{}, E_{1,2}^{(4)}$
        \item
        with configuration as in Figure \ref{Conf2H2V}.
    \end{itemize} 
    This is case \hyperref[Tab2H]{$2H$} if $p \neq 2$, and case \hyperref[Tab2V]{$2V$} if $p=2$.
        \end{minipage} \hspace{2mm} \begin{minipage}{0.3\textwidth} \begin{center} \includegraphics[width=0.98\textwidth]{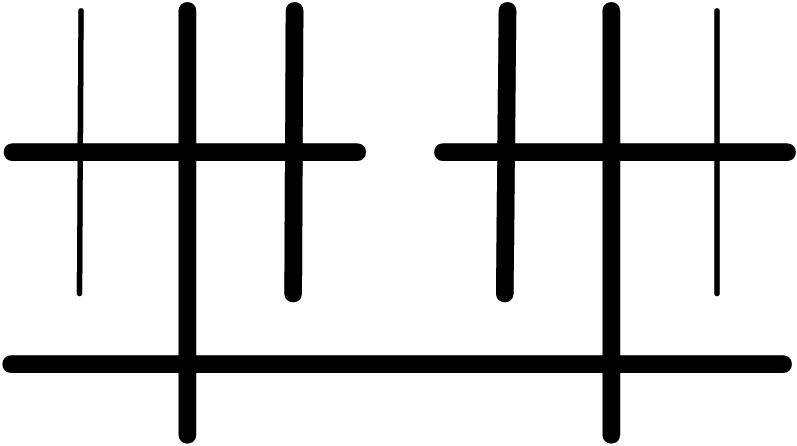} \vspace{-3.5mm}\captionof{figure}{}  \label{Conf2H2V} \end{center} \end{minipage}

    \noindent     \begin{minipage}{0.65\textwidth}
    \item \vspace{2mm}
    $p_{1,3}  =E_{1,2}^{} \cap C^{(3)}$ with $C= \cal{V}(xz+y^2)$
    \begin{itemize}[leftmargin=20pt]
        \item
        $\Aut_{X'}^0(R) =
        \begin{cases}   
        \left\{  \left( \begin{smallmatrix}
1 &  &  \\
 & e &  \\
 &  & e^2
\end{smallmatrix} \right) 
\in \PGL_3(R) \right\} & \text{ if } p \neq 2 \\
        \left\{  \left( \begin{smallmatrix}
1 &  &  \\
 & e & f \\
 &  & e^2
\end{smallmatrix} \right) 
\in \PGL_3(R)  \right\}
         & \text{ if } p=2
        \end{cases} $ 
        \\We describe the configurations of negative curves on $X'$ for $p \neq 2$ and $p=2$ simultaneously:
        \item $(-2)$-curves: $E_{1,0}^{(4)}, E_{2,0}^{(4)}, E_{1,1}^{(4)}, E_{2,1}^{(4)}, E_{1,2}^{(4)}, \ell_x^{(4)}, \ell_z^{(4)}$
        \item $(-1)$-curves: $E_{1,3}^{}, E_{2,2}^{(4)}$
        \item
        with configuration as in Figure \ref{Conf2G2T2U}.
    \end{itemize} 
    This is case \hyperref[Tab2G]{$2G$} if $p \neq 2$, and case \hyperref[Tab2U]{$2U$} if $p=2$.
        \end{minipage} \hspace{2mm} \begin{minipage}{0.3\textwidth} \begin{center} \includegraphics[width=0.98\textwidth]{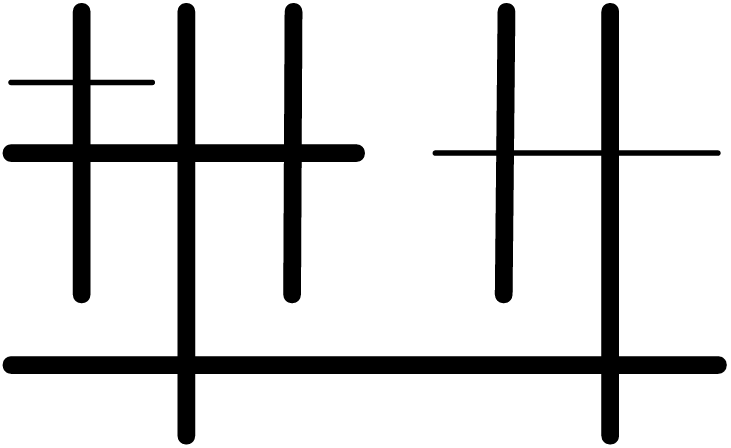} \vspace{-3.5mm}\captionof{figure}{}  \label{Conf2G2T2U} \end{center} \end{minipage} 
   
    \item \vspace{2mm}
    Let $p=2$ and
    $p_{1,3}  =E_{1,2}^{} \cap C^{(3)}$ with $C= \cal{V}(x^2z+xy^2+y^3)$.
    \begin{itemize}[leftmargin=20pt]
    \noindent \begin{minipage}{0.5\textwidth}
        \item
        $
        \Aut_{X'}^0(R) =    
        \left\{  \left( \begin{smallmatrix}
1 &  &  \\
 & 1 & f \\
 &  & 1
\end{smallmatrix} \right) 
\in \PGL_3(R) \right\}
        $
        \item $(-2)$-curves: $E_{1,0}^{(4)}, E_{2,0}^{(4)}, E_{1,1}^{(4)}, E_{2,1}^{(4)}, E_{1,2}^{(4)}, \ell_x^{(4)}, \ell_z^{(4)}$
        
        \end{minipage} \noindent \begin{minipage}{0.5\textwidth}
        
        \item $(-1)$-curves: $E_{1,3}^{}, E_{2,2}^{(4)}$
        \item
        with configuration as in Figure \ref{Conf2G2T2U}. 
    \end{minipage}
    \end{itemize} 
    \vspace{1mm} \noindent This is case \hyperref[Tab2T]{$2T$}.

\end{enumerate}

\vspace{2mm} \subsubsection*{\underline{Case \hyperref[Tab4K]{$4K$}}}
We have $E =
E_{2,2}^{} - (E_{2,1}^{(3)} \cup \ell_x^{(3)})$ and $\Aut_X^0(R) = \left\{  \left( \begin{smallmatrix}
1 &  &  \\
 & e & f \\
 &  & i
\end{smallmatrix} \right) 
\in \PGL_3(R) \right\}$.

\begin{itemize}[leftmargin=25pt]
    \item[-] $\lambda xy^2 + \mu z^3$ is $E_{2,2}^{}$-adapted 
    and $\Aut_X^0(R)$ acts as
    $[\lambda:\mu] \mapsto [e^2\lambda: i^3\mu]$
\end{itemize}
\noindent Since $\Aut_X^0$ acts transitively on $E \cap E_{2,2}$, there is a unique possibility for $p_{2,3}$ up to isomorphism:

\begin{enumerate}[leftmargin=*]
\item \vspace{2mm}
    $p_{2,3} = E_{2,2}^{} \cap C^{(3)}$ with $C= \cal{V}(xy^2+z^3)$ 
    \begin{itemize}[leftmargin=20pt]
    \noindent \begin{minipage}{0.5\textwidth}
        \item
        $
        \Aut_{X'}^0(R) =    
        \left\{  \left( \begin{smallmatrix}
1 &  &  \\
 & e & f \\
 &  & i
\end{smallmatrix} \right) 
\in \PGL_3(R) \bigg| e^2= i^3\right\}
        $
        \item $(-2)$-curves: $E_{1,0}^{(4)}, E_{2,0}^{(4)},E_{2,1}^{(4)}, E_{2,2}^{(4)},\ell_x^{(4)}, \ell_z^{(4)}$
        
        \end{minipage} \noindent \begin{minipage}{0.5\textwidth}
        
        \item $(-1)$-curves: $E_{2,3}^{}, E_{1,1}^{(4)}$
        \item
        with configuration as in Figure \ref{Conf3I}, that is, as in \\ case \hyperref[Tab3I]{$3I$}.
    \end{minipage}
    \end{itemize} 
    \vspace{1mm} \noindent As explained in Remark \ref{R IsomorphismCheck}, one can check that $X' \cong X_{3I}$.
\end{enumerate}

\vspace{2mm} \subsubsection*{\underline{Case \hyperref[Tab4I]{$4I$}}}
We have $E =
E_{1,2}^{} - E_{1,1}^{(3)}$ and $\Aut_X^0(R) = 
 \left\{  \left( \begin{smallmatrix}
1 &  &  \\
 & e & f \\
 &  & e^2
\end{smallmatrix} \right) 
\in \PGL_3(R) \right\}$.

\begin{itemize}[leftmargin=25pt]
    \item[-] $\lambda (x^2z + xy^2) + \mu y^3$ is
    $E_{1,2}^{}$-adapted and $\Aut_X^0(R)$ acts as
    $[\lambda:\mu] \mapsto [e^2\lambda: e^3\mu - 2ef \lambda ]$  
\end{itemize}
If $p \neq 2$, then $\Aut_X^0$ acts transitively on $E \cap E_{1,2}$, while if $p = 2$, then the $\Aut_X^0$ has two orbits on $E \cap E_{1,2}$. Hence, if $p = 2$, there is only one possibility for $p_{1,3}$ and if $p = 2$, there are two possibilities up to isomorphism:

\begin{enumerate}[leftmargin=*]
    
    \noindent     \begin{minipage}{0.65\textwidth}
   \item \vspace{2mm}
    $p_{1,3}  =E_{1,2}^{} \cap C^{(3)}$ with $C=\cal{V}(xz+y^2)$
    \begin{itemize}[leftmargin=20pt]
        \item
        $\Aut_{X'}^0(R) =
        \begin{cases}   
        \left\{  \left( \begin{smallmatrix}
1 &  &  \\
 & e &  \\
 &  & e^2
\end{smallmatrix} \right) 
\in \PGL_3(R) \right\} & \text{ if } p \neq 2 \\
        \left\{  \left( \begin{smallmatrix}
1 &  &  \\
 & e & f \\
 &  & e^2
\end{smallmatrix} \right) 
\in \PGL_3(R)  \right\}
         & \text{ if } p=2
        \end{cases} $ 
        \\We describe the configurations of negative curves on $X'$ for $p \neq 2$ and $p=2$ simultaneously:
        \item $(-2)$-curves: $E_{1,0}^{(4)}, E_{2,0}^{(4)},E_{1,1}^{(4)}, E_{1,2}^{(4)},\ell_z^{(4)}$
        \item $(-1)$-curves: $E_{1,3}^{}, E_{2,1}^{(4)},  \ell_x^{(4)}$
        \item
        with configuration as in Figure \ref{Conf3G3O3P}.
    \end{itemize} 
    This is case \hyperref[Tab3G]{$3G$} if $p \neq 2$, and case \hyperref[Tab3P]{$3P$} if $p=2$.
    \end{minipage} \hspace{2mm} \begin{minipage}{0.3\textwidth} \begin{center} \includegraphics[width=0.85\textwidth]{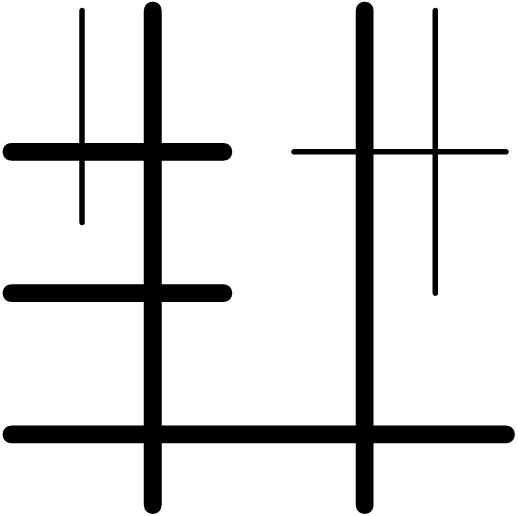} \vspace{-3.5mm}\captionof{figure}{}  \label{Conf3G3O3P} \end{center} \end{minipage}

    \item \vspace{2mm}
    Let $p=2$ and
    $p_{1,3}  =E_{1,2}^{} \cap C^{(3)}$ with $C= \cal{V}(x^2z+xy^2+y^3)$.
    \begin{itemize}[leftmargin=20pt]
    \noindent \begin{minipage}{0.5\textwidth}
        \item
        $
        \Aut_{X'}^0(R) =    
        \left\{  \left( \begin{smallmatrix}
1 &  &  \\
 & 1 & f \\
 &  & 1
\end{smallmatrix} \right) 
\in \PGL_3(R) \right\}
        $
        \item $(-2)$-curves: $E_{1,0}^{(4)}, E_{2,0}^{(4)},E_{1,1}^{(4)},E_{1,2}^{(4)}, \ell_z^{(4)}$
        
        \end{minipage} \noindent \begin{minipage}{0.5\textwidth}
        
        \item $(-1)$-curves: $E_{1,3}^{}, E_{2,1}^{(4)},  \ell_x^{(4)}$
        \item
        with configuration as in Figure \ref{Conf3G3O3P}.
    \end{minipage}
    \end{itemize} 
    \vspace{1mm} \noindent This is case \hyperref[Tab3O]{$3O$}.

\end{enumerate}

\vspace{2mm} \subsubsection*{\underline{Case \hyperref[Tab5E]{$5E$}}}
We have $E =
E_{1,2}^{} - E_{1,1}^{(3)}$ and $\Aut_X^0(R) =
 \left\{  \left( \begin{smallmatrix}
1 &  & c \\
 & e & f \\
 &  & e^2
\end{smallmatrix} \right) 
\in \PGL_3(R) \right\}$.

\begin{itemize}[leftmargin=25pt]
    \item[-] $\lambda (x^2z + xy^2) + \mu y^3$ is $E_{1,2}^{}$-adapted and $\Aut_X^0(R)$ acts as
    $[\lambda:\mu] \mapsto [e^2\lambda: e^3\mu - 2ef\lambda ]$
\end{itemize}

\noindent
As in the previous case, if $p \neq 2$, there is only one possibility for $p_{1,3}$ up to isomorphism, and if $p = 2$, there are two possibilities up to isomorphism:

\begin{enumerate}[leftmargin=*]

\noindent     \begin{minipage}{0.65\textwidth}
\item \vspace{2mm}
    $p_{1,3}  =E_{1,2}^{} \cap C^{(3)}$ with $C= \cal{V}(xz+y^2)$
    \begin{itemize}[leftmargin=20pt]
        \item
        $\Aut_{X'}^0(R) =
        \begin{cases}   
        \left\{  \left( \begin{smallmatrix}
1 &  & c \\
 & e &  \\
 &  & e^2
\end{smallmatrix} \right) 
\in \PGL_3(R) \right\} & \text{ if } p \neq 2 \\
        \left\{  \left( \begin{smallmatrix}
1 &  & c \\
 & e & f \\
 &  & e^2
\end{smallmatrix} \right) 
\in \PGL_3(R)  \right\}
         & \text{ if } p=2
        \end{cases} $ 
        \\We describe the configurations of negative curves on $X'$ for $p \neq 2$ and $p=2$ simultaneously:
        \item $(-2)$-curves: $E_{1,0}^{(4)}, E_{1,1}^{(4)},E_{1,2}^{(4)},\ell_z^{(4)}$
        \item $(-1)$-curves: $E_{1,3}^{}, E_{2,0}^{(4)}$
        \item
        with configuration as in Figure \ref{Conf4N4J4O}.
    \end{itemize} 
    This is case \hyperref[Tab4J]{$4J$} if $p \neq 2$, and case \hyperref[Tab4O]{$4O$} if $p=2$.
    \end{minipage} \hspace{2mm} \begin{minipage}{0.3\textwidth} \begin{center} \includegraphics[width=0.85\textwidth]{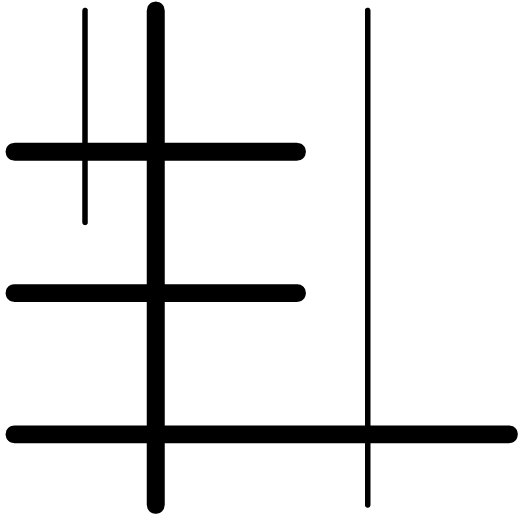} \vspace{-3.5mm}\captionof{figure}{}  \label{Conf4N4J4O} \end{center} \end{minipage}

    \item \vspace{2mm}
    Let $p=2$ and
    $p_{1,3}  =E_{1,2}^{} \cap C^{(3)}$ with $C= \cal{V}(x^2z+xy^2+y^3)$.
    \begin{itemize}[leftmargin=20pt]
    \noindent \begin{minipage}{0.5\textwidth}
        \item
        $
        \Aut_{X'}^0(R) =    
        \left\{  \left( \begin{smallmatrix}
1 &  & c \\
 & 1 & f \\
 &  & 1
\end{smallmatrix} \right) 
\in \PGL_3(R) \right\}
        $
        \item $(-2)$-curves: $E_{1,0}^{(4)}, E_{1,1}^{(4)},E_{1,2}^{(4)},\ell_z^{(4)}$
        
        \end{minipage} \noindent \begin{minipage}{0.5\textwidth}
        
        \item $(-1)$-curves: $E_{1,3}^{}, E_{2,0}^{(4)}$
        \item
        with configuration as in Figure \ref{Conf4N4J4O}.
    \end{minipage}
    \end{itemize} 
    \vspace{1mm} \noindent This is case \hyperref[Tab4N]{$4N$}.
\end{enumerate}

\vspace{2mm} \subsubsection*{\underline{Case \hyperref[Tab6E]{$6E$}}}
We have $E =
E_{1,2}^{} - E_{1,1}^{(3)}$ and $\Aut_X^0(R) = 
 \left\{  \left( \begin{smallmatrix}
1 & b & c \\
 & e & f \\
 &  & e^2
\end{smallmatrix} \right) 
\in \PGL_3(R) \right\}$.

\begin{itemize}[leftmargin=25pt]
    \item[-] $\lambda (x^2z + xy^2) + \mu y^3$ is $E_{1,2}^{}$-adapted and $\Aut_X^0(R)$ acts as $[\lambda:\mu] \mapsto [e^2 \lambda: e^3 \mu - be^2 \lambda - 2ef \lambda]$.
\end{itemize}

Since $\Aut_X^0$ acts transitively on $E \cap E_{1,2}$, there is a unique possibility for $p_{1,3}$ up to isomorphism:

\begin{enumerate}[leftmargin=*]

\item \vspace{2mm}
    $p_{1,3}  =E_{1,2}^{} \cap C^{(3)}$ with $C=\cal{V}(xz+y^2)$
    \begin{itemize}[leftmargin=20pt]
    \noindent \begin{minipage}{0.5\textwidth}
        \item
        $
        \Aut_{X'}^0(R) =    
        \left\{  \left( \begin{smallmatrix}
1 & -2fe^{-1} & c \\
 & e & f \\
 &  & e^2
\end{smallmatrix} \right) 
\in \PGL_3(R) \right\}
        $
        \item $(-2)$-curves: $E_{1,0}^{(4)}, E_{1,1}^{(4)}, E_{1,2}^{(4)}$
        
        \end{minipage} \noindent \begin{minipage}{0.5\textwidth}
        
        \item $(-1)$-curves: $E_{1,3}^{}, \ell_z^{(4)}$
        \item
        with configuration as in Figure \ref{Conf5E}, that is, as in \\ case \hyperref[Tab5E]{$5E$}.
    \end{minipage}
    \end{itemize} 
\vspace{1mm} \noindent As explained in Remark \ref{R IsomorphismCheck}, one can check that $X' \cong X_{5E}$.
    
\end{enumerate}

\vspace{2mm} \subsubsection*{\underline{Case \hyperref[Tab6F]{$6F$}}}
We have $E =
E_{1,2}^{} - (E_{1,1}^{(3)} \cup \ell_{z}^{(3)})$ and $\Aut_X^0(R) = \left\{  \left( \begin{smallmatrix}
1 & b & c \\
 & e & f \\
 &  & i
\end{smallmatrix} \right) 
\in \PGL_3(R) \right\}$.

\begin{itemize}[leftmargin=25pt]
    \item[-] $\lambda x^2z + \mu y^3$ is $E_{1,2}^{}$-adapted and $\Aut_X^0(R)$ acts as $[\lambda:\mu] \mapsto [i \lambda: e^3 \mu]$.
\end{itemize}

Since $\Aut_X^0$ acts transitively on $E \cap E_{1,2}$, there is a unique possibility for $p_{1,3}$ up to isomorphism:

\begin{enumerate}[leftmargin=*]
\noindent     \begin{minipage}{0.65\textwidth}
    \item \vspace{2mm}
    $p_{1,3}  =E_{1,2}^{} \cap C^{(3)}$ with $C= \cal{V}(x^2z+y^3)$
    \begin{itemize}[leftmargin=20pt]
        \item
        $
        \Aut_{X'}^0(R) =    
        \left\{  \left( \begin{smallmatrix}
1 & b & c \\
 & e & f \\
 &  & e^3
\end{smallmatrix} \right) 
\in \PGL_3(R) \right\}
        $
        \item $(-2)$-curves: $E_{1,0}^{(4)}, E_{1,1}^{(4)},E_{1,2}^{(4)}, \ell_z^{(4)}$
        \item $(-1)$-curves: $E_{1,3}^{}$
        \item
        with configuration as in Figure \ref{Conf5F}.
    \end{itemize} 
    This is case \hyperref[Tab5F]{$5F$}.
        \end{minipage} \hspace{2mm} \begin{minipage}{0.3\textwidth} \begin{center} \includegraphics[width=0.98\textwidth]{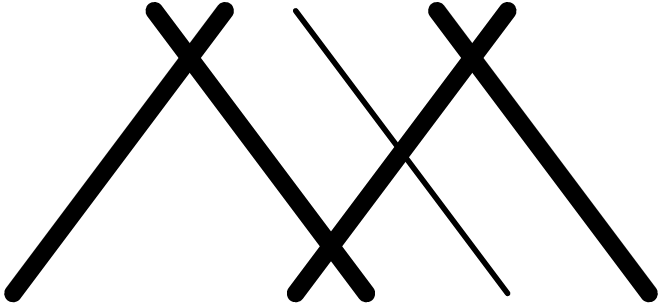} \vspace{-3.5mm}\captionof{figure}{}  \label{Conf5F} \end{center} \end{minipage} 
\end{enumerate}

\vspace{2mm} 
\noindent
Summarizing, we obtain
\begin{eqnarray*}
\cal{L}_4 &=& \{  X_{1E}, X_{2R,\alpha}, X_{2K}, X_{1Q,\alpha}, X_{1R}, X_{2H},X_{2V}, X_{2G},X_{2U}, X_{2T}, X_{3G},X_{3P}, X_{3O}, X_{4J},X_{4O}, X_{4N}, X_{5F}\}.
\end{eqnarray*}

\subsection{Height 5}

\vspace{2mm} \subsubsection*{\underline{Case \hyperref[Tab2R]{$2R$}}}
This case exists only if $p = 2$.  
\\ \noindent We have $E =
E_{1,3}^{} - E_{1,2}^{(4)}$ and $\Aut_X^0(R) = \left\{  \left( \begin{smallmatrix}
1 &  &  \\
 & 1 & f \\
 &  & 1
\end{smallmatrix} \right) 
\in \PGL_3(R) \right\}$.

\begin{itemize}[leftmargin=25pt]
    \item[-] $\lambda (x+\alpha y)^2(xz + y^2 + \alpha yz) + \mu y^4$ is $E_{1,3}^{}$-adapted and $\Aut_X^0(R)$ acts as $[\lambda:\mu] \mapsto [\lambda: \mu + (\alpha f + f^2)\lambda]$.
\end{itemize}

\noindent
Therefore, if $\alpha \neq 0$, then the identity component of the stabilizer of every point on $E \cap E_{1,3}^{}$ is trivial, hence there is no way of further blowing up $X$ and still obtaining a weak del Pezzo surface with global vector fields. If $\alpha = 0$, then there is the following unique possibility for $p_{1,4}$ up to isomorphism:

\begin{enumerate}[leftmargin=*]
\noindent     \begin{minipage}{0.65\textwidth}
\item \vspace{2mm}
    $p_{1,4}  =E_{1,3}^{} \cap C_1^{(4)}$ with $C_1= \cal{V}(xz+y^2)$
    \begin{itemize}[leftmargin=20pt]
        \item
        $\Aut_{X'}^0(R) =
        \left\{  \left( \begin{smallmatrix}
1 &  &  \\
 & 1 & f \\
 &  & 1
\end{smallmatrix} \right) 
\in \PGL_3(R) \bigg| f^2=0  \right\}
       $ 
        \item $(-2)$-curves: $E_{1,0}^{(5)}, E_{2,0}^{(5)},E_{1,1}^{(5)}, E_{2,1}^{(5)}, E_{1,2}^{(5)},E_{1,3}^{(5)},\ell_z^{(5)}$
        \item $(-1)$-curves: $E_{1,4}^{}, E_{2,2}^{(5)},  \ell_x^{(5)}, C_1^{(5)}, C_2^{(5)}$ with 
        \\$C_2= \cal{V}(x^2y^2+x^3z+z^4)$
        \item
        with configuration as in Figure \ref{Conf1M}.
    \end{itemize} 
    This is case \hyperref[Tab1M]{$1M$}.
        \end{minipage} \hspace{2mm} \begin{minipage}{0.3\textwidth} \begin{center} \includegraphics[width=0.98\textwidth]{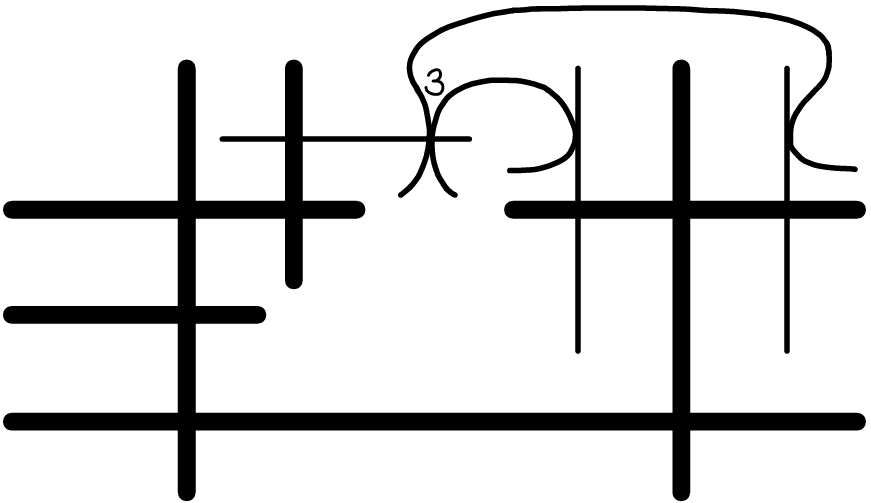} \vspace{-3.5mm}\captionof{figure}{}  \label{Conf1M} \end{center} \end{minipage} 
\end{enumerate}

\vspace{2mm} \subsubsection*{\underline{Case \hyperref[Tab2K]{$2K$}}}
This case exists only if $p = 3$.  
\\ \noindent We have $E =
E_{1,3}^{} - E_{1,2}^{(4)}$ and $\Aut_X^0(R) = \left\{  \left( \begin{smallmatrix}
1 &  &  \\
 & e &  \\
 &  & e^2
\end{smallmatrix} \right) 
\in \PGL_3(R) \bigg| e^3 = 1 \right\}$.

\begin{itemize}[leftmargin=25pt]
    \item[-] $\lambda x^2(xz + y^2) + \mu y^4$ is $E_{1,3}^{}$-adapted and $\Aut_X^0(R)$ acts as $[\lambda:\mu] \mapsto [e^2\lambda: e\mu ]$.
\end{itemize}
Note that there is a unique point on $E \cap E_{1,3}^{}$ with non-trivial stabilizer. This leads to the following unique possibility for $p_{1,4}$:

\begin{enumerate}[leftmargin=*]
\noindent     \begin{minipage}{0.65\textwidth}
\item \vspace{2mm}
    $p_{1,4}  =E_{1,3}^{} \cap C_1^{(4)}$ with $C_1= \cal{V}(xz+y^2)$ 
    \begin{itemize}[leftmargin=20pt]
        \item
        $\Aut_{X'}^0(R) =
        \left\{  \left( \begin{smallmatrix}
1 &  &  \\
 & e &  \\
 &  & e^2
\end{smallmatrix} \right) 
\in \PGL_3(R) \bigg| e^3=1 \right\}
        $ 
        \item $(-2)$-curves: $E_{1,0}^{(5)}, E_{2,0}^{(5)},E_{1,1}^{(5)}, E_{2,1}^{(5)},E_{1,2}^{(5)},E_{1,3}^{(5)},  \ell_z^{(5)}$
        \item $(-1)$-curves: $E_{1,4}^{}, E_{2,2}^{(5)},  \ell_x^{(5)}, C_1^{(5)}, C_2^{(5)}$ with 
        \\$C_2= \cal{V}(x^2y^2+x^3z+z^4+2xyz^2)$
        \item
        with configuration as in Figure \ref{Conf1F}.
    \end{itemize} 
    This is case \hyperref[Tab1F]{$1F$}. 
        \end{minipage} \hspace{2mm} \begin{minipage}{0.3\textwidth} \begin{center} \includegraphics[width=0.98\textwidth]{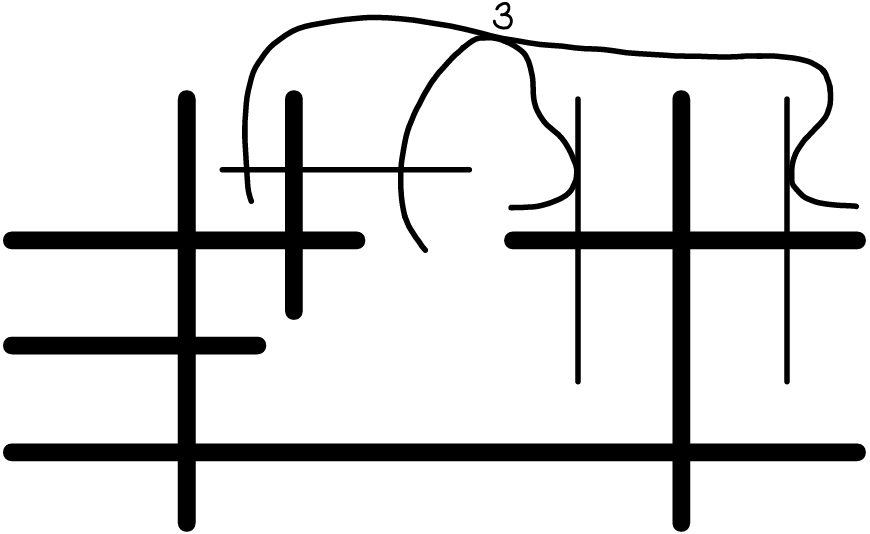} \vspace{-3.5mm}\captionof{figure}{}  \label{Conf1F} \end{center} \end{minipage} 
\end{enumerate}

\vspace{2mm} \subsubsection*{\underline{Case \hyperref[Tab2H]{$2H$}}}
This case exists only if $p \neq 2$.  
\\ \noindent We have $E =
E_{2,3}^{} - E_{2,2}^{(4)}$ and $\Aut_X^0(R) = \left\{  \left( \begin{smallmatrix}
1 &  &  \\
 & 1 & f\\
 &  & 1
\end{smallmatrix} \right) 
\in \PGL_3(R) \right\}$.

\begin{itemize}[leftmargin=25pt]
    \item[-] $\lambda y(xy^2 + z^3) + \mu z^4$ is $E_{2,3}^{}$-adapted and $\Aut_X^0(R)$ acts as $[\lambda:\mu] \mapsto [\lambda : \mu  - 2  f \lambda]$.
\end{itemize}

\noindent In particular, since $p \neq 2$ the stabilizer of every point on $E \cap E_{2,3}^{}$ is trivial, hence there is no way of further blowing up $X$ and obtaining a weak del Pezzo surface with global vector fields.

\vspace{2mm} \subsubsection*{\underline{Case \hyperref[Tab2V]{$2V$}}}
This case exists only if $p = 2$.  
\\ \noindent We have $E =
E_{2,3}^{} - E_{2,2}^{(4)}$ and $\Aut_X^0(R) = \left\{  \left( \begin{smallmatrix}
1 &  &  \\
 & e & f\\
 &  & e^2
\end{smallmatrix} \right) 
\in \PGL_3(R) \bigg| e^4 = 1 \right\}$.

\begin{itemize}[leftmargin=25pt]
    \item[-] $\lambda y(xy^2 + z^3) + \mu z^4$ is $E_{2,3}^{}$-adapted and $\Aut_X^0(R)$ acts as $[\lambda:\mu] \mapsto [e^3 \lambda : \mu ]$.
\end{itemize}
This leads to the following possibilities for $p_{1,4}$:

\begin{enumerate}[leftmargin=*]
\item \vspace{2mm}
    $ p_{2,4} = E_{2,3}^{} \cap C^{(4)}$ with $C= \cal{V}(xy^3+yz^3 + \alpha z^4), \alpha \neq 0$
    \begin{itemize}[leftmargin=20pt]
    \noindent \begin{minipage}{0.5\textwidth}
        \item
        $\Aut_{X'}^0(R) = \left\{ \left( \begin{smallmatrix}
1 &  &  \\
 & 1 & f \\
 &  & 1
\end{smallmatrix} \right) 
\in \PGL_3(R) \right\} $ 
        \item $(-2)$-curves: $E_{1,0}^{(5)}, E_{2,0}^{(5)}, E_{1,1}^{(5)}, E_{2,1}^{(5)}, E_{2,2}^{(5)},E_{2,3}^{(5)}, 
        \\ \ell_x^{(5)}, \ell_z^{(5)}$
        
        \end{minipage} \noindent \begin{minipage}{0.5\textwidth}
        
        \item $(-1)$-curves: $E_{2,4}^{}, E_{1,2}^{(5)}$
        \item
        with configuration as in Figure \ref{Conf1Q1R}, that is, as in \\ case \hyperref[Tab1Q]{$1Q$}.
    \end{minipage}
    \end{itemize} 
    \vspace{1mm} \noindent As explained in Remark \ref{R IsomorphismCheck}, one can check that $X' \cong X_{1Q,\alpha'}$ for some $\alpha'$.

    \item \vspace{2mm}
    $ p_{2,4} = E_{2,3}^{} \cap C^{(4)}$ with $C= \cal{V}(xy^2+z^3)$
    \begin{itemize}[leftmargin=20pt]
    \noindent \begin{minipage}{0.5\textwidth}
        \item
        $\Aut_{X'}^0(R) =
        \left\{  \left( \begin{smallmatrix}
1 &  &  \\
 & e & f \\
 &  & e^2
\end{smallmatrix} \right) 
\in \PGL_3(R) \bigg| e^4=1 \right\} $ 
        \item $(-2)$-curves: $E_{1,0}^{(5)}, E_{2,0}^{(5)}, E_{1,1}^{(5)}, E_{2,1}^{(5)}, E_{2,2}^{(5)},E_{2,3}^{(5)}, \\ \ell_x^{(5)}, \ell_z^{(5)}$
        
        \end{minipage} \noindent \begin{minipage}{0.5\textwidth}
        
        \item $(-1)$-curves: $E_{2,4}^{}, E_{1,2}^{(5)}$
         \item
        with configuration as in Figure \ref{Conf1Q1R}, that is, as in \\ case \hyperref[Tab1R]{$1R$}.
    \end{minipage}
    \end{itemize} 
\vspace{1mm} \noindent As explained in Remark \ref{R IsomorphismCheck}, one can check that $X' \cong X_{1R}$.
    
    \end{enumerate}

\vspace{2mm} \subsubsection*{\underline{Case \hyperref[Tab2G]{$2G$}}}
This case exists only if $p \neq 2$.  
\\ \noindent We have $E =
E_{1,3}^{} - E_{1,2}^{(4)}$ and $\Aut_X^0(R) = \left\{  \left( \begin{smallmatrix}
1 &  &  \\
 & e & \\
 &  & e^2
\end{smallmatrix} \right) 
\in \PGL_3(R) \right\}$.

\begin{itemize}[leftmargin=25pt]
    \item[-] $\lambda x^2(xz + y^2) + \mu y^4$ is $E_{1,3}^{}$-adapted and $\Aut_X^0(R)$ acts as $[\lambda:\mu] \mapsto [\lambda : e^2 \mu ]$.
\end{itemize}
Since $p \neq 2$, there is a unique point on $E \cap E_{1,3}^{}$ such that the identity component of its stabilizer is non-trivial. This leads to the following unique possibility for $p_{1,4}$:

\begin{enumerate}[leftmargin=*]
\noindent     \begin{minipage}{0.65\textwidth}
   \item \vspace{2mm}
    $p_{1,4}  =E_{1,3}^{} \cap C^{(4)}$ with $C= \cal{V}(xz+y^2)$
    \begin{itemize}[leftmargin=20pt]
        \item
        $\Aut_{X'}^0(R) =
        \left\{  \left( \begin{smallmatrix}
1 &  &  \\
 & e &  \\
 &  & e^2
\end{smallmatrix} \right) 
\in \PGL_3(R) \right\} $ 
        \item $(-2)$-curves: $E_{1,0}^{(5)}, E_{2,0}^{(5)}, E_{1,1}^{(5)}, E_{2,1}^{(5)}, E_{1,2}^{(5)}, E_{1,3}^{(5)},\ell_x^{(5)}, \ell_z^{(5)}$
        \item $(-1)$-curves: $E_{1,4}^{}, E_{2,2}^{(5)}, C^{(5)}$
        \item
        with configuration as in Figure \ref{Conf1C1P}.
    \end{itemize} 
    This is case \hyperref[Tab1C]{$1C$}.
        \end{minipage} \hspace{2mm} \begin{minipage}{0.3\textwidth} \begin{center} \includegraphics[width=0.65\textwidth]{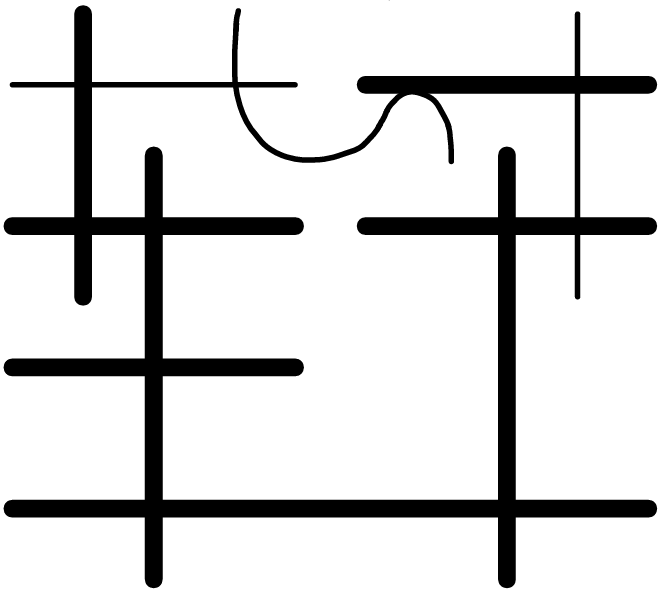} \vspace{-3.5mm}\captionof{figure}{} \label{Conf1C1P}  \end{center} \end{minipage} 
\end{enumerate}

\vspace{2mm} \subsubsection*{\underline{Case \hyperref[Tab2U]{$2U$}}}
This case exists only if $p = 2$.  
\\ \noindent We have $E =
E_{1,3}^{} - E_{1,2}^{(4)}$ and $\Aut_X^0(R) = \left\{  \left( \begin{smallmatrix}
1 &  &  \\
 & e & f\\
 &  & e^2
\end{smallmatrix} \right) 
\in \PGL_3(R) \right\}$.

\begin{itemize}[leftmargin=25pt]
    \item[-] $\lambda x^2(xz + y^2) + \mu y^4$ is $E_{1,3}^{}$-adapted and $\Aut_X^0(R)$ acts as $[\lambda:\mu] \mapsto [e^2 \lambda : e^4 \mu +  f^2\lambda]$.
\end{itemize}
Since $\Aut_X^0$ acts transitively on $E \cap E_{1,3}$, there is a unique possibility for $p_{1,4}$ up to isomorphism:

\begin{enumerate}[leftmargin=*]
    \item \vspace{2mm}
    $p_{1,4}  =E_{1,3}^{} \cap C^{(4)}$ with $C= \cal{V}(xz+y^2)$
    \begin{itemize}[leftmargin=20pt]
    \noindent \begin{minipage}{0.5\textwidth}
        \item
        $\Aut_{X'}^0(R) =
        \left\{  \left( \begin{smallmatrix}
1 &  &  \\
 & e & f \\
 &  & e^2
\end{smallmatrix} \right) 
\in \PGL_3(R) \bigg| f^2=0 \right\}
        $ 
        \item $(-2)$-curves: $E_{1,0}^{(5)}, E_{2,0}^{(5)}, E_{1,1}^{(5)}, E_{2,1}^{(5)}, E_{1,2}^{(5)},E_{1,3}^{(5)}, 
        \\ \ell_x^{(5)}, \ell_z^{(5)}$
        
        \end{minipage} \noindent \begin{minipage}{0.5\textwidth}
        
        \item $(-1)$-curves: $E_{1,4}^{}, E_{2,2}^{(5)}, C^{(5)}$
        \item
        with configuration as in Figure \ref{Conf1C1P}. 
    \end{minipage}
    \end{itemize} 
    \vspace{1mm} \noindent This is case \hyperref[Tab1P]{$1P$}.
     
\end{enumerate}

\vspace{2mm} \subsubsection*{\underline{Case \hyperref[Tab2T]{$2T$}}}
This case exists only if $p = 2$.  
\\ \noindent We have $E =
E_{1,3}^{} - E_{1,2}^{(4)}$ and $\Aut_X^0(R) = \left\{  \left( \begin{smallmatrix}
1 &  &  \\
 & 1 & f\\
 &  & 1
\end{smallmatrix} \right) 
\in \PGL_3(R) \right\}$.

\begin{itemize}[leftmargin=25pt]
    \item[-] $\lambda (x+y)(x^2z + xy^2 + y^3 + y^2z) + \mu y^4$ is $E_{1,3}^{}$-adapted and $\Aut_X^0(R)$ acts as $[\lambda:\mu] \mapsto [\lambda : \mu + (f + f^2)\lambda]$.
\end{itemize}
Note that the identity component of the stabilizer of every point on $E \cap E_{1,3}^{}$ is trivial, hence we cannot blow up further and still obtain a weak del Pezzo surface with global vector fields.

\vspace{2mm} \subsubsection*{\underline{Case \hyperref[Tab3G]{$3G$}}}
This case exists only if $p \neq 2$.  
\\ \noindent We have $E =
E_{1,3}^{} - E_{1,2}^{(4)}$ and $\Aut_X^0(R) = \left\{  \left( \begin{smallmatrix}
1 &  &  \\
 & e & \\
 &  & e^2
\end{smallmatrix} \right) 
\in \PGL_3(R) \right\}$.

\begin{itemize}[leftmargin=25pt]
    \item[-] $\lambda x^2(xz + y^2) + \mu y^4$ is $E_{1,3}^{}$-adapted and $\Aut_X^0(R)$ acts as $[\lambda:\mu] \mapsto [\lambda : e^2\mu]$.
\end{itemize}
Since $p \neq 2$, there is a unique point on $E \cap E_{1,3}^{}$ for which the identity component of the stabilizer is non-trivial. This leads to the following unique possibility for $p_{1,4}$:

\begin{enumerate}[leftmargin=*]
\noindent     \begin{minipage}{0.65\textwidth}
    \item \vspace{2mm}
    $p_{1,4}  =E_{1,3}^{} \cap C^{(4)}$ with $C=\cal{V}(xz+y^2)$
    \begin{itemize}[leftmargin=20pt]
        \item
        $\Aut_{X'}^0(R) =
        \left\{  \left( \begin{smallmatrix}
1 &  &  \\
 & e &  \\
 &  & e^2
\end{smallmatrix} \right) 
\in \PGL_3(R) \right\} $ 
        \item $(-2)$-curves: $E_{1,0}^{(5)}, E_{2,0}^{(5)},E_{1,1}^{(5)}, E_{1,2}^{(5)},E_{1,3}^{(5)},\ell_z^{(5)}$
        \item $(-1)$-curves: $E_{1,4}^{}, E_{2,1}^{(5)},  \ell_x^{(5)}, C^{(5)}$
        \item
        with configuration as in Figure \ref{Conf2C2S}.
    \end{itemize} 
    This is case \hyperref[Tab2C]{$2C$}.
        \end{minipage} \hspace{2mm} \begin{minipage}{0.3\textwidth} \begin{center} \includegraphics[width=0.65\textwidth]{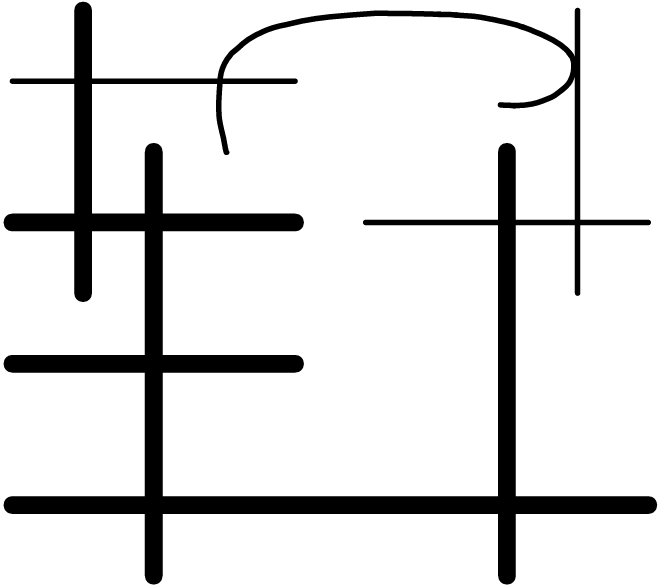} \vspace{-3.5mm}\captionof{figure}{} \label{Conf2C2S}  \end{center} \end{minipage} 
\end{enumerate}

\vspace{2mm} \subsubsection*{\underline{Case \hyperref[Tab3P]{$3P$}}}
This case exists only if $p = 2$.  
\\ \noindent We have $E =
E_{1,3}^{} - E_{1,2}^{(4)}$ and $\Aut_X^0(R) = \left\{  \left( \begin{smallmatrix}
1 &  &  \\
 & e & f \\
 &  & e^2
\end{smallmatrix} \right) 
\in \PGL_3(R) \right\}$.

\begin{itemize}[leftmargin=25pt]
    \item[-] $\lambda x^2(xz + y^2) + \mu y^4$ is $E_{1,3}^{}$-adapted and $\Aut_X^0(R)$ acts as $[\lambda:\mu] \mapsto [e^2\lambda : e^4 \mu + f^2 \lambda]$.
\end{itemize}
Since $\Aut_X^0$ acts transitively on $E \cap E_{1,3}$, there is a unique possibility for $p_{1,4}$ up to isomorphism:

\begin{enumerate}[leftmargin=*]

   \item \vspace{2mm}
    $p_{1,4}  =E_{1,3}^{} \cap C^{(4)}$ with $C=\cal{V}(xz+y^2)$
    \begin{itemize}[leftmargin=20pt]
    \noindent \begin{minipage}{0.5\textwidth}
        \item
        $\Aut_{X'}^0(R) =
        \left\{  \left( \begin{smallmatrix}
1 &  &  \\
 & e & f \\
 &  & e^2
\end{smallmatrix} \right) 
\in \PGL_3(R)  \bigg| f^2=0 \right\}
        $ 
        \item $(-2)$-curves: $E_{1,0}^{(5)}, E_{2,0}^{(5)},E_{1,1}^{(5)}, E_{1,2}^{(5)},E_{1,3}^{(5)},\ell_z^{(5)}$
        
        \end{minipage} \noindent \begin{minipage}{0.5\textwidth}
        
        \item $(-1)$-curves: $E_{1,4}^{}, E_{2,1}^{(5)},  \ell_x^{(5)}, C^{(5)}$
        \item
        with configuration as in Figure \ref{Conf2C2S}. 
    \end{minipage}
    \end{itemize} 
    \vspace{1mm} \noindent This is case \hyperref[Tab2S]{$2S$}.
\end{enumerate}

\vspace{2mm} \subsubsection*{\underline{Case \hyperref[Tab3O]{$3O$}}}
This case exists only if $p = 2$.  
\\ \noindent We have $E =
E_{1,3}^{} - E_{1,2}^{(4)}$ and $\Aut_X^0(R) = \left\{  \left( \begin{smallmatrix}
1 &  &  \\
 & 1 & f \\
 &  & 1
\end{smallmatrix} \right) 
\in \PGL_3(R) \right\}$.

\begin{itemize}[leftmargin=25pt]
    \item[-] $\lambda (x+y)(x^2z +xy^2 + y^3 + y^2z) + \mu y^4$ is $E_{1,3}^{}$-adapted and $\Aut_X^0(R)$ acts as $[\lambda:\mu] \mapsto [\lambda : \mu + (f + f^2)\lambda]$.
\end{itemize}
In particular, the identity component of the stabilizer of every point on $E \cap E_{1,3}^{}$ is trivial, hence we cannot blow up further.

\vspace{2mm} \subsubsection*{\underline{Case \hyperref[Tab4J]{$4J$}}}
This case exists only if $p \neq 2$.  
\\ \noindent We have $E =
E_{1,3}^{} - E_{1,2}^{(4)}$ and $\Aut_X^0(R) = \left\{  \left( \begin{smallmatrix}
1 &  & c \\
 & e&  \\
 &  & e^2
\end{smallmatrix} \right) 
\in \PGL_3(R) \right\}$.

\begin{itemize}[leftmargin=25pt]
    \item[-] $\lambda x^2(xz +y^2) + \mu y^4$ is $E_{1,3}^{}$-adapted and $\Aut_X^0(R)$ acts as $[\lambda:\mu] \mapsto [\lambda : e^2\mu + c \lambda]$.
\end{itemize}
Since $\Aut_X^0$ acts transitively on $E \cap E_{1,3}$, we have the following unique possibility for $p_{1,4}$ up to isomorphism:

\begin{enumerate}[leftmargin=*]
   \item \vspace{2mm}
    $p_{1,4}  =E_{1,3}^{} \cap C^{(4)}$ with $C= \cal{V}(xz+y^2)$
    \begin{itemize}[leftmargin=20pt]
    \noindent \begin{minipage}{0.5\textwidth}
        \item
        $\Aut_{X'}^0(R) =
        \left\{  \left( \begin{smallmatrix}
1 &  &  \\
 & e &  \\
 &  & e^2
\end{smallmatrix} \right) 
\in \PGL_3(R) \right\} $ 
        \item $(-2)$-curves: $E_{1,0}^{(5)}, E_{1,1}^{(5)},E_{1,2}^{(5)},E_{1,3}^{(5)},\ell_z^{(5)}$
        
        \end{minipage} \noindent \begin{minipage}{0.5\textwidth}
        
        \item $(-1)$-curves: $E_{1,4}^{}, E_{2,0}^{(5)}, C^{(5)}$
        \item
        with configuration as in Figure \ref{Conf3G3O3P}, that is, as in \\ case \hyperref[Tab3G]{$3G$}.
    \end{minipage}
    \end{itemize} 
   \vspace{1mm} \noindent As explained in Remark \ref{R IsomorphismCheck}, one can check that $X' \cong X_{3G}$.
\end{enumerate}

\vspace{2mm} \subsubsection*{\underline{Case \hyperref[Tab4O]{$4O$}}}
This case exists only if $p = 2$.  
\\ \noindent We have $E =
E_{1,3}^{} - E_{1,2}^{(4)}$ and $\Aut_X^0(R) = \left\{  \left( \begin{smallmatrix}
1 &  & c \\
 & e& f \\
 &  & e^2
\end{smallmatrix} \right) 
\in \PGL_3(R) \right\}$.

\begin{itemize}[leftmargin=25pt]
    \item[-] $\lambda x^2(xz +y^2) + \mu y^4$ is $E_{1,3}^{}$-adapted and $\Aut_X^0(R)$ acts as $[\lambda:\mu] \mapsto [e^2\lambda : e^4\mu + (ce^2 + f^2)\lambda ]$.
\end{itemize}
Since $\Aut_X^0$ acts transitively on $E \cap E_{1,3}$, we have the following unique possibility for $p_{1,4}$ up to isomorphism:

\begin{enumerate}[leftmargin=*]
    \item \vspace{2mm}
    $p_{1,4}  =E_{1,3}^{} \cap C^{(4)}$ with $C= \cal{V}(xz+y^2)$
    \begin{itemize}[leftmargin=20pt]
    \noindent \begin{minipage}{0.5\textwidth}
        \item
        $\Aut_{X'}^0(R) =
        \left\{  \left( \begin{smallmatrix}
1 &  & f^2 e^{-2} \\
 & e & f \\
 &  & e^2
\end{smallmatrix} \right) 
\in \PGL_3(R)   \right\}
        $ 
        \item $(-2)$-curves: $E_{1,0}^{(5)}, E_{1,1}^{(5)},E_{1,2}^{(5)},E_{1,3}^{(5)},\ell_z^{(5)}$
        
        \end{minipage} \noindent \begin{minipage}{0.5\textwidth}
        
        \item $(-1)$-curves: $E_{1,4}^{}, E_{2,0}^{(5)}, C^{(5)}$
        \item
        with configuration as in Figure \ref{Conf3G3O3P}, that is, as in \\ case \hyperref[Tab3P]{$3P$}.
    \end{minipage}
    \end{itemize} 
  \vspace{1mm} \noindent As explained in Remark \ref{R IsomorphismCheck}, one can check that $X' \cong X_{3P}$.
\end{enumerate}

\vspace{2mm} \subsubsection*{\underline{Case \hyperref[Tab4N]{$4N$}}}
This case exists only if $p = 2$.  
\\ \noindent We have $E =
E_{1,3}^{} - E_{1,2}^{(4)}$ and $\Aut_X^0(R) = \left\{  \left( \begin{smallmatrix}
1 &  & c \\
 & 1& f \\
 &  & 1
\end{smallmatrix} \right) 
\in \PGL_3(R) \right\}$.

\begin{itemize}[leftmargin=25pt]
    \item[-] $\lambda (x+y)(x^2z +xy^2 + y^3 + y^2z) + \mu y^4$ is $E_{1,3}^{}$-adapted and $\Aut_X^0(R)$ acts as 
    $[\lambda:\mu]  \mapsto  [\lambda : \mu + (c + f + f^2)\lambda]$.
\end{itemize}
Since $\Aut_X^0$ acts transitively on $E \cap E_{1,3}$, we have the following unique possibility for $p_{1,4}$ up to isomorphism:

\begin{enumerate}[leftmargin=*]
    \item \vspace{2mm}
    $p_{1,4}  =E_{1,3}^{} \cap C_1^{(4)}$ with $C_1= \cal{V}(x^2z+xy^2+y^3)$
    \begin{itemize}[leftmargin=20pt]
    \noindent \begin{minipage}{0.5\textwidth}
        \item
        $
        \Aut_{X'}^0(R) =    
        \left\{  \left( \begin{smallmatrix}
1 &  & f+f^2 \\
 & 1 & f \\
 &  & 1
\end{smallmatrix} \right) 
\in \PGL_3(R) \right\}
        $
        \item $(-2)$-curves: $E_{1,0}^{(5)}, E_{1,1}^{(5)},E_{1,2}^{(5)},E_{1,3}^{(5)},\ell_z^{(5)}$
        
        \end{minipage} \noindent \begin{minipage}{0.5\textwidth}
        
        \item $(-1)$-curves: $E_{1,4}^{}, E_{2,0}^{(5)}, C_2^{(5)}$ with 
        \\$C_2= \cal{V}(xz+yz+y^2)$
        \item
        with configuration as in Figure \ref{Conf3G3O3P}, that is, as in \\ case \hyperref[Tab3O]{$3O$}.
    \end{minipage}
    \end{itemize} 
    \vspace{1mm} \noindent As explained in Remark \ref{R IsomorphismCheck}, one can check that $X' \cong X_{3O}$.
\end{enumerate}

\vspace{2mm} \subsubsection*{\underline{Case \hyperref[Tab5F]{$5F$}}}
We have $E =
E_{1,3}^{} - E_{1,2}^{(4)}$ and $\Aut_X^0(R) = 
 \left\{  \left( \begin{smallmatrix}
1 & b & c \\
 & e & f \\
 &  & e^3
\end{smallmatrix} \right) 
\in \PGL_3(R) \right\}$. 

\begin{itemize}[leftmargin=25pt]
    \item[-] $\lambda x(x^2z + y^3) + \mu y^4$ is $E_{1,3}^{}$-adapted and $\Aut_X^0(R)$ acts as $[\lambda:\mu] \mapsto [\lambda : e \mu - 2  b \lambda]$.
\end{itemize}
Therefore, if $p \neq 2$, we have one unique possibility for $p_{1,4} \in E \cap E_{1,3}$, while if $p = 2$, there are two possibilities:

\begin{enumerate}[leftmargin=*]

\noindent     \begin{minipage}{0.65\textwidth}
    \item \vspace{2mm}
    $p_{1,4}  =E_{1,3}^{} \cap C^{(4)}$ with $C= \cal{V}(x^2z+y^3)$
    \begin{itemize}[leftmargin=20pt]
        \item
        $\Aut_{X'}^0(R) =
        \begin{cases}   
        \left\{  \left( \begin{smallmatrix}
1 &  & c \\
 & e & f \\
 &  & e^3
\end{smallmatrix} \right) 
\in \PGL_3(R)  \right\} & \text{ if } p \neq 2 \\
        \left\{  \left( \begin{smallmatrix}
1 & b & c \\
 & e & f \\
 &  & e^3
\end{smallmatrix} \right) \in \PGL_3(R)  \right\}
         & \text{ if } p=2
        \end{cases} $ 
        \\We describe the configurations of negative curves on $X'$ for $p \neq 2$ and $p=2$ simultaneously:
        \item $(-2)$-curves: $E_{1,0}^{(5)}, E_{1,1}^{(5)},E_{1,2}^{(5)}, E_{1,3}^{(5)}, \ell_z^{(5)}$
        \item $(-1)$-curves: $E_{1,4}^{}$
        \item
        with configuration as in Figure \ref{Conf4P4Q4L}.
    \end{itemize} 
    This is case \hyperref[Tab4L]{$4L$} if $p \neq 2$, and case \hyperref[Tab4Q]{$4Q$} if $p=2$.
    \end{minipage} \hspace{2mm} \begin{minipage}{0.3\textwidth} \begin{center} \includegraphics[width=0.98\textwidth]{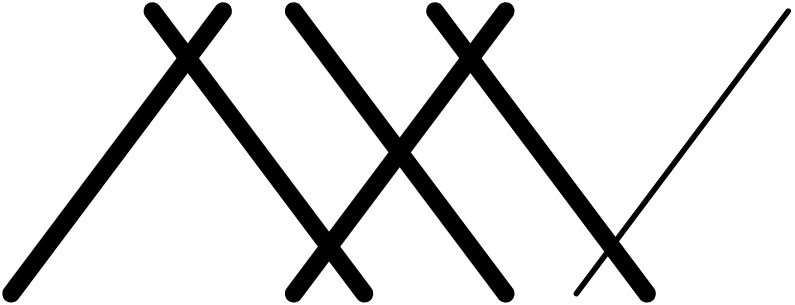} \vspace{-3.5mm}\captionof{figure}{}  \label{Conf4P4Q4L} \end{center} \end{minipage}

    \item \vspace{2mm}
    Let $p=2$ and
    $p_{1,4}  =E_{1,3}^{} \cap C^{(4)}$ with $C= \cal{V}(x^3z + xy^3 + y^4)$.
    \begin{itemize}[leftmargin=20pt]
    \noindent \begin{minipage}{0.5\textwidth}
        \item
        $
        \Aut_{X'}^0(R) =    
        \left\{  \left( \begin{smallmatrix}
1 & b & c \\
 & 1 & f \\
 &  & 1
\end{smallmatrix} \right) 
\in \PGL_3(R) \right\}
        $
        \item $(-2)$-curves: $E_{1,0}^{(5)}, E_{1,1}^{(5)},E_{1,2}^{(5)}, E_{1,3}^{(5)}, \ell_z^{(5)}$
        
        \end{minipage} \noindent \begin{minipage}{0.5\textwidth}
        
        \item $(-1)$-curves: $E_{1,4}^{}$
        \item
        with configuration as in Figure \ref{Conf4P4Q4L}.
    \end{minipage}
    \end{itemize} 
    \vspace{1mm} \noindent This is case \hyperref[Tab4P]{$4P$}.

\end{enumerate}

\vspace{2mm} 
\noindent
Summarizing, we obtain
\begin{eqnarray*}
\cal{L}_5 &=& \{  X_{1M}, X_{1F}, X_{1C}, X_{1P}, X_{2C}, X_{2S}, X_{4L}, X_{4Q}, X_{4P}\}.
\end{eqnarray*}

\subsection{Height 6}

\vspace{2mm} \subsubsection*{\underline{Case \hyperref[Tab2C]{$2C$}}}
This case exists only if $p \neq 2$.  
\\ \noindent We have $E =
E_{1,4}^{} - E_{1,3}^{(5)}$ and $\Aut_X^0(R) = \left\{  \left( \begin{smallmatrix}
1 &  &  \\
 & e&  \\
 &  & e^2
\end{smallmatrix} \right) 
\in \PGL_3(R) \right\}$.

\begin{itemize}[leftmargin=25pt]
    \item[-] $\lambda x^3(xz + y^2) + \mu y^5$ is $E_{1,4}^{}$-adapted and $\Aut_X^0(R)$ acts as $[\lambda:\mu] \mapsto [\lambda : e^3 \mu]$.
\end{itemize}
Note that if $p \neq 3$, then there is a unique point on $E \cap E_{1,4}^{}$ such that the identity component of its stabilizer is non-trivial. If $p = 3$, this identity component is non-trivial for every point. In all characteristics, the action of $\Aut_X^0$ on $E \cap E_{1,4}$ has two orbits. Hence, we have the following two possibilities for $p_{1,5}$ up to isomorphism:

\begin{enumerate}[leftmargin=*]
    \item \vspace{2mm}
    $p_{1,5}  =E_{1,4}^{} \cap C_1^{(5)}$ with $C_1=\cal{V}(x^4z+x^3y^2+y^5)$
    \begin{itemize}[leftmargin=20pt]
        \item
        $\Aut_{X'}^0(R) =
        \begin{cases}   
        \{ \id \} & \text{ if } p \neq 2,3 \\
        \left\{  \left( \begin{smallmatrix}
1 &  &  \\
 & e &  \\
 &  & e^2
\end{smallmatrix} \right) \in \PGL_3(R) \bigg| e^3=1 \right\}
         & \text{ if } p=3
        \end{cases} $ 
        \\Hence, $X'$ has global vector fields only if $p=3$. Therefore, we assume $p=3$ when describing the configuration of negative curves.
        \item $(-2)$-curves: $E_{1,0}^{(6)}, E_{2,0}^{(6)},E_{1,1}^{(6)}, E_{1,2}^{(6)},E_{1,3}^{(6)},E_{1,4}^{(6)},\ell_z^{(6)}$
        \item $(-1)$-curves: $E_{1,5}^{}, E_{2,1}^{(6)},  \ell_x^{(6)}, C_2^{(6)}, C_3^{(6)}$ with $C_2= \cal{V}(xz+y^2), 
        \\C_3= \cal{V}(xy^4-xyz^3-x^2y^2z+x^3z^2-y^3z^2-z^5)$
        \item
        with configuration as in Figure \ref{Conf1F}, that is, as in case \hyperref[Tab1F]{$1F$}.
    \end{itemize} 
     \noindent As explained in Remark \ref{R IsomorphismCheck}, one can check that $X' \cong X_{1F}$.
    
    \item \vspace{2mm}
    $p_{1,5}  =E_{1,4}^{} \cap C^{(5)}$ with $C=\cal{V}(xz+y^2)$
    \begin{itemize}[leftmargin=20pt]
    \noindent \begin{minipage}{0.5\textwidth}
        \item
        $\Aut_{X'}^0(R) =
        \left\{  \left( \begin{smallmatrix}
1 &  &  \\
 & e &  \\
 &  & e^2
\end{smallmatrix} \right) 
\in \PGL_3(R) \right\} $ 
        \item $(-2)$-curves: $E_{1,0}^{(6)}, E_{2,0}^{(6)},E_{1,1}^{(6)}, E_{1,2}^{(6)},E_{1,3}^{(6)},E_{1,4}^{(6)},
        \\ \ell_z^{(6)}, C^{(6)}$
        
        \end{minipage} \noindent \begin{minipage}{0.5\textwidth}
        
        \item $(-1)$-curves: $E_{1,5}^{}, E_{2,1}^{(6)},  \ell_x^{(6)}$
        \item
        with configuration as in Figure \ref{Conf1C1P}, that is, as in \\ case \hyperref[Tab1C]{$1C$}.
    \end{minipage}
    \end{itemize} 
    \vspace{1mm} \noindent As explained in Remark \ref{R IsomorphismCheck}, one can check that $X' \cong X_{1C}$.
    
\end{enumerate}

\vspace{2mm} \subsubsection*{\underline{Case \hyperref[Tab2S]{$2S$}}}
This case exists only if $p = 2$.  
\\ \noindent We have $E =
E_{1,4}^{} - E_{1,3}^{(5)}$ and $\Aut_X^0(R) = \left\{  \left( \begin{smallmatrix}
1 &  &  \\
 & e& f \\
 &  & e^2
\end{smallmatrix} \right) 
\in \PGL_3(R) \bigg| f^2 = 0\right\}$.

\begin{itemize}[leftmargin=25pt]
    \item[-] $\lambda x^3(xz + y^2) + \mu y^5$ is $E_{1,4}^{}$-adapted and $\Aut_X^0(R)$ acts as $[\lambda:\mu] \mapsto [\lambda : e^3 \mu]$.
\end{itemize}
Since $\Aut_X^0$ acts on $E \cap E_{1,4}$ with two orbits, we have the following two possibilities for $p_{1,5}$ up to isomorphism:

\begin{enumerate}[leftmargin=*]

    \item \vspace{2mm}
    $p_{1,5}  =E_{1,4}^{} \cap C_1^{(5)}$ with $C_1=\cal{V}(x^4z+x^3y^2+y^5)$
    \begin{itemize}[leftmargin=20pt]
    \noindent \begin{minipage}{0.5\textwidth}
        \item
        $\Aut_{X'}^0(R) =
        \left\{  \left( \begin{smallmatrix}
1 &  &  \\
 & 1 & f \\
 &  & 1
\end{smallmatrix} \right) 
\in \PGL_3(R)  \bigg| f^2=0 \right\}
        $ 
        \item $(-2)$-curves: $E_{1,0}^{(6)}, E_{2,0}^{(6)},E_{1,1}^{(6)}, E_{1,2}^{(6)},E_{1,3}^{(6)},
        \\ E_{1,4}^{(6)},\ell_z^{(6)}$
        
        \end{minipage} \noindent \begin{minipage}{0.5\textwidth}
        
        \item $(-1)$-curves: $E_{1,5}^{}, E_{2,1}^{(6)},  \ell_x^{(6)}, C_2^{(6)}, C_3^{(6)}$ with 
        \\$C_2= \cal{V}(xz+y^2), C_3= \cal{V}(xy^4+x^3z^2+z^5)$
        \item
        with configuration as in Figure \ref{Conf1M}, that is, as in \\ case \hyperref[Tab1M]{$1M$}.
    \end{minipage}
    \end{itemize} 
    \vspace{1mm} \noindent As explained in Remark \ref{R IsomorphismCheck}, one can check that $X' \cong X_{1M}$.

    \item \vspace{2mm}
    $p_{1,5}  =E_{1,4}^{} \cap C^{(5)}$ with $C=\cal{V}(xz+y^2)$
    \begin{itemize}[leftmargin=20pt]
    \noindent \begin{minipage}{0.5\textwidth}
        \item
        $\Aut_{X'}^0(R) =
        \left\{  \left( \begin{smallmatrix}
1 &  &  \\
 & e & f \\
 &  & e^2
\end{smallmatrix} \right) 
\in \PGL_3(R)  \bigg| f^2=0 \right\}
        $ 
        \item $(-2)$-curves: $E_{1,0}^{(6)}, E_{2,0}^{(6)},E_{1,1}^{(6)}, E_{1,2}^{(6)},E_{1,3}^{(6)},E_{1,4}^{(6)},
        \\ \ell_z^{(6)}, C^{(6)}$
        
        \end{minipage} \noindent \begin{minipage}{0.5\textwidth}
        
        \item $(-1)$-curves: $E_{1,5}^{}, E_{2,1}^{(6)},  \ell_x^{(6)}$
        \item
        with configuration as in Figure \ref{Conf1C1P}, that is, as in \\ case \hyperref[Tab1P]{$1P$}.
    \end{minipage}
    \end{itemize} 
   \vspace{1mm} \noindent As explained in Remark \ref{R IsomorphismCheck}, one can check that $X' \cong X_{1P}$.
    
\end{enumerate}

\newpage

\vspace{2mm} \subsubsection*{\underline{Case \hyperref[Tab4L]{$4L$}}}
This case exists only if $p \neq 2$.  
\\ \noindent We have $E =
E_{1,4}^{} - E_{1,3}^{(5)}$ and $\Aut_X^0(R) = \left\{  \left( \begin{smallmatrix}
1 &  & c \\
 & e& f \\
 &  & e^3
\end{smallmatrix} \right) 
\in \PGL_3(R) \right\}$. 

\begin{itemize}[leftmargin=25pt]
    \item[-] $\lambda x^2(x^2z + y^3) + \mu y^5$ is $E_{1,4}^{}$-adapted and $\Aut_X^0(R)$ acts as $[\lambda:\mu] \mapsto [e\lambda : e^3 \mu - 3 f \lambda]$. 
\end{itemize}
In particular, if $p \neq 3$, then $\Aut_X^0$ acts transitively on $E \cap E_{1,4}$ and we have only one choice for $p_{1,5}$ up to isomorphism, and if $p = 3$, then $\Aut_X^0$ acts with two orbits on $E \cap E_{1,4}$, hence we have two choices up to isomorphism:

\begin{enumerate}[leftmargin=*]

\noindent     \begin{minipage}{0.65\textwidth}
    \item \vspace{2mm}
    $p_{1,5}  =E_{1,4}^{} \cap C^{(5)}$ with $C= \cal{V}(x^2z+y^3)$
    \begin{itemize}[leftmargin=20pt]
        \item
        $\Aut_{X'}^0(R) =
        \begin{cases}   
        \left\{  \left( \begin{smallmatrix}
1 &  & c \\
 & e &  \\
 &  & e^3
\end{smallmatrix} \right) 
\in \PGL_3(R)  \right\} & \text{ if } p \neq 2,3 \\
        \left\{  \left( \begin{smallmatrix}
1 &  & c \\
 & e & f \\
 &  & e^3
\end{smallmatrix} \right) \in \PGL_3(R)  \right\}
         & \text{ if } p=3
        \end{cases} $ 
        \\We describe the configurations of negative curves on $X'$ for \\
        $p \neq 2,3$ and $p=3$ simultaneously:
        \item $(-2)$-curves: $E_{1,0}^{(6)}, E_{1,1}^{(6)},E_{1,2}^{(6)}, E_{1,3}^{(6)}, E_{1,4}^{(6)}, \ell_z^{(6)}$
        \item $(-1)$-curves: $E_{1,5}^{}$
        \item
        with configuration as in Figure \ref{Conf3J3L3M3R3Q}.
    \end{itemize} 
    This is case \hyperref[Tab3J]{$3J$} if $p \neq 2,3$, and case \hyperref[Tab3M]{$3M$} if $p=3$. 
    \end{minipage} \hspace{2mm} \begin{minipage}{0.3\textwidth} \begin{center} \includegraphics[width=0.98\textwidth]{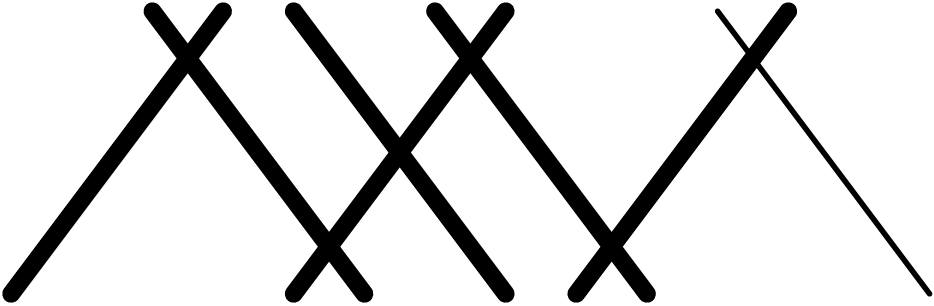} \vspace{-3.5mm}\captionof{figure}{}  \label{Conf3J3L3M3R3Q} \end{center} \end{minipage}

    \item \vspace{2mm}
    Let $p=3$ and
    $p_{1,5}  =E_{1,4}^{} \cap C^{(5)}$ with $C= \cal{V}(x^4z+x^2y^3+y^5)$.
    \begin{itemize}[leftmargin=20pt]
    \noindent \begin{minipage}{0.5\textwidth}
        \item
        $\Aut_{X'}^0(R) =
        \left\{  \left( \begin{smallmatrix}
1 &  & c \\
 & 1 & f \\
 &  & 1
\end{smallmatrix} \right) 
\in \PGL_3(R)  \right\} $ 
        \item $(-2)$-curves: $E_{1,0}^{(6)}, E_{1,1}^{(6)},E_{1,2}^{(6)}, E_{1,3}^{(6)},E_{1,4}^{(6)}, \ell_z^{(6)}$
        
        \end{minipage} \noindent \begin{minipage}{0.5\textwidth}
        
        \item $(-1)$-curves: $E_{1,5}^{}$
        \item
        with configuration as in Figure \ref{Conf3J3L3M3R3Q}.
    \end{minipage}
    \end{itemize} 
    \vspace{1mm} \noindent This is case \hyperref[Tab3L]{$3L$}.

\end{enumerate}

\vspace{2mm} \subsubsection*{\underline{Case \hyperref[Tab4Q]{$4Q$}}}
This case exists only if $p = 2$.  
\\ \noindent We have $E =
E_{1,4}^{} - E_{1,3}^{(5)}$ and $\Aut_X^0(R) = \left\{  \left( \begin{smallmatrix}
1 & b & c \\
 & e& f \\
 &  & e^3
\end{smallmatrix} \right) 
\in \PGL_3(R) \right\}$. 

\begin{itemize}[leftmargin=25pt]
    \item[-] $\lambda x^2(x^2z + y^3) + \mu y^5$ is $E_{1,4}^{}$-adapted and $\Aut_X^0(R)$ acts as $[\lambda:\mu] \mapsto [e \lambda : e^3 \mu + (b^2e + f) \lambda]$.
\end{itemize}
Since $\Aut_X^0$ acts transitively on $E \cap E_{1,4}$, there is a unique choice for $p_{1,5}$ up to isomorphism:

\begin{enumerate}[leftmargin=*]

   \item \vspace{2mm}
    $p_{1,5}  =E_{1,4}^{} \cap C^{(5)}$ with $C= \cal{V}(x^2z+y^3)$
    \begin{itemize}[leftmargin=20pt]
    \noindent \begin{minipage}{0.5\textwidth}
        \item
        $\Aut_{X'}^0(R) =
        \left\{  \left( \begin{smallmatrix}
1 & b & c \\
 & e & b^2e \\
 &  & e^3
\end{smallmatrix} \right) 
\in \PGL_3(R)  \right\}$ 
        \item $(-2)$-curves: $E_{1,0}^{(6)}, E_{1,1}^{(6)},E_{1,2}^{(6)}, E_{1,3}^{(6)},E_{1,4}^{(6)}, \ell_z^{(6)}$
        
        \end{minipage} \noindent \begin{minipage}{0.5\textwidth}
        
        \item $(-1)$-curves: $E_{1,5}^{}$
        \item
        with configuration as in Figure \ref{Conf3J3L3M3R3Q}.
    \end{minipage}
    \end{itemize} 
    \vspace{1mm} \noindent This is case \hyperref[Tab3R]{$3R$}. 

\end{enumerate}

\vspace{2mm} \subsubsection*{\underline{Case \hyperref[Tab4P]{$4P$}}}
This case exists only if $p = 2$.  
\\ \noindent We have $E =
E_{1,4}^{} - E_{1,3}^{(5)}$ and $\Aut_X^0(R) = \left\{  \left( \begin{smallmatrix}
1 & b & c \\
 & 1& f \\
 &  & 1
\end{smallmatrix} \right) 
\in \PGL_3(R) \right\}$.

\begin{itemize}[leftmargin=25pt]
    \item[-] $\lambda x(x^3z + xy^3 + y^4) + \mu y^5$ is $E_{1,4}^{}$-adapted and $\Aut_X^0(R)$ acts as $[\lambda:\mu] \mapsto [\lambda : \mu +  (b + b^2 + f)\lambda]$.
\end{itemize}
Since $\Aut_X^0$ acts transitively on $E \cap E_{1,4}$, we have the following unique choice for $p_{1,5}$ up to isomorphism:

\begin{enumerate}[leftmargin=*]

    \item \vspace{2mm}
    $p_{1,5}  =E_{1,4}^{} \cap C^{(5)}$ with $C= \cal{V}(x^3z + xy^3 + y^4)$
    \begin{itemize}[leftmargin=20pt]
    \noindent \begin{minipage}{0.5\textwidth}
        \item
        $
        \Aut_{X'}^0(R) =    
        \left\{  \left( \begin{smallmatrix}
1 & b & c \\
 & 1 & b^2+b \\
 &  & 1
\end{smallmatrix} \right) 
\in \PGL_3(R) \right\}
        $
        \item $(-2)$-curves: $E_{1,0}^{(6)}, E_{1,1}^{(6)},E_{1,2}^{(6)}, E_{1,3}^{(6)},E_{1,4}^{(6)}, \ell_z^{(6)}$
        
        \end{minipage} \noindent \begin{minipage}{0.5\textwidth}
        
        \item $(-1)$-curves: $E_{1,5}^{}$
        \item
        with configuration as in Figure \ref{Conf3J3L3M3R3Q}. 
    \end{minipage}
    \end{itemize} 
    \vspace{1mm} \noindent This is case \hyperref[Tab3Q]{$3Q$}.

\end{enumerate}

\vspace{2mm} 
\noindent
Summarizing, we obtain
\begin{eqnarray*}
\cal{L}_6 &=& \{  X_{3J},X_{3M},X_{3L}, X_{3R}, X_{3Q}\}.
\end{eqnarray*}

\subsection{Height 7}

\vspace{2mm} \subsubsection*{\underline{Case \hyperref[Tab3J]{$3J$}}}
This case exists only if $p \neq 2,3$.  
\\ \noindent We have $E =
E_{1,5}^{} - E_{1,4}^{(6)}$ and $\Aut_X^0(R) = \left\{  \left( \begin{smallmatrix}
1 &  & c \\
 & e&  \\
 &  & e^3
\end{smallmatrix} \right) 
\in \PGL_3(R) \right\}$.

\begin{itemize}[leftmargin=25pt]
    \item[-] $\lambda x^3(x^2z + y^3) + \mu y^6$ is $E_{1,5}^{}$-adapted and $\Aut_X^0(R)$ acts as $[\lambda:\mu] \mapsto [\lambda : e^3 \mu + 2 c \lambda$].
\end{itemize}
Since $p \neq 2$, $\Aut_X^0$ acts transitively on $E \cap E_{1,5}$, so there is a unique choice for $p_{1,6}$ up to isomorphism:

\begin{enumerate}[leftmargin=*]
\noindent     \begin{minipage}{0.65\textwidth}
    \item \vspace{2mm}
    $p_{1,6}  =E_{1,5}^{} \cap C^{(6)}$ with $C= \cal{V}(x^2z+y^3)$ 
    \begin{itemize}[leftmargin=20pt]
        \item
        $\Aut_{X'}^0(R) =
        \left\{  \left( \begin{smallmatrix}
1 &  &  \\
 & e &  \\
 &  & e^3
\end{smallmatrix} \right) 
\in \PGL_3(R)  \right\}$ 
        \item $(-2)$-curves: $E_{1,0}^{(7)}, E_{1,1}^{(7)},E_{1,2}^{(7)}, E_{1,3}^{(7)}, E_{1,4}^{(7)},E_{1,5}^{(7)},  \ell_z^{(7)}$
        \item $(-1)$-curves: $E_{1,6}^{}$
        \item
        with configuration as in Figure \ref{Conf2I2L2M2W2X2Y}.
    \end{itemize} 
    This is case \hyperref[Tab2I]{$2I$}.
        \end{minipage} \hspace{2mm} \begin{minipage}{0.3\textwidth} \begin{center} \includegraphics[width=0.98\textwidth]{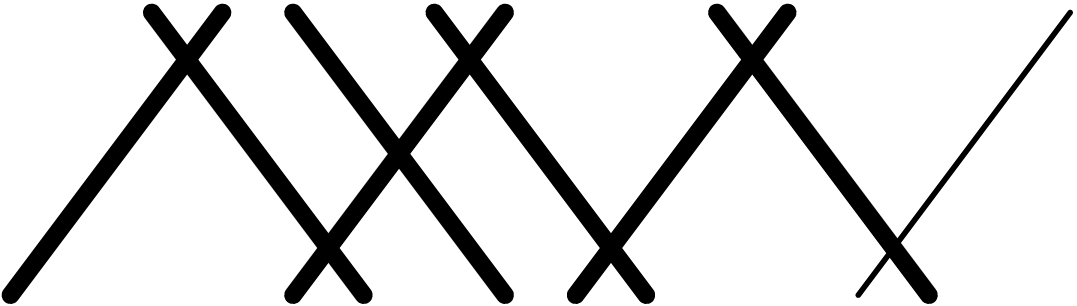} \vspace{-3.5mm}\captionof{figure}{}  \label{Conf2I2L2M2W2X2Y} \end{center} \end{minipage} 
\end{enumerate}

\vspace{2mm} \subsubsection*{\underline{Case \hyperref[Tab3M]{$3M$}}}
This case exists only if $p = 3$.  
\\ \noindent We have $E =
E_{1,5}^{} - E_{1,4}^{(6)}$ and $\Aut_X^0(R) = \left\{  \left( \begin{smallmatrix}
1 &  & c \\
 & e& f \\
 &  & e^3
\end{smallmatrix} \right) 
\in \PGL_3(R) \right\}$.

\begin{itemize}[leftmargin=25pt]
    \item[-] $\lambda x^3(x^2z + y^3) + \mu y^6$ is $E_{1,5}^{}$-adapted and $\Aut_X^0(R)$ acts as $[\lambda:\mu] \mapsto [\lambda : e^3 \mu + 2 c \lambda]$.
\end{itemize}
As in the previous case, there is a unique choice for $p_{1,6}$ up to isomorphism:

\begin{enumerate}[leftmargin=*]

   \item \vspace{2mm}
    $p_{1,6}  =E_{1,5}^{} \cap C^{(6)}$ with $C= \cal{V}(x^2z+y^3)$
    \begin{itemize}[leftmargin=20pt]
    \noindent \begin{minipage}{0.5\textwidth}
        \item
        $\Aut_{X'}^0(R) =
        \left\{  \left( \begin{smallmatrix}
1 &  &  \\
 & e & f \\
 &  & e^3
\end{smallmatrix} \right) 
\in \PGL_3(R)  \right\} $
        \item $(-2)$-curves: $E_{1,0}^{(7)}, E_{1,1}^{(7)},E_{1,2}^{(7)}, E_{1,3}^{(7)}, E_{1,4}^{(7)},
        E_{1,5}^{(7)}, \ell_z^{(7)}$
        
        \end{minipage} \noindent \begin{minipage}{0.5\textwidth}
        
        \item $(-1)$-curves: $E_{1,6}^{}$
        \item
        with configuration as in Figure \ref{Conf2I2L2M2W2X2Y}.
    \end{minipage}
    \end{itemize} 
    \vspace{1mm} \noindent This is case \hyperref[Tab2M]{$2M$}.

\end{enumerate}

\vspace{2mm} \subsubsection*{\underline{Case \hyperref[Tab3L]{$3L$}}}
This case exists only if $p = 3$.  
\\ \noindent We have $E =
E_{1,5}^{} - E_{1,4}^{(6)}$ and $\Aut_X^0(R) = \left\{  \left( \begin{smallmatrix}
1 &  & c \\
 & 1& f \\
 &  & 1
\end{smallmatrix} \right) 
\in \PGL_3(R) \right\}$.

\begin{itemize}[leftmargin=25pt]
    \item[-] $\lambda x(x^4z + x^2y^3 + y^5) + \mu y^6$ is $E_{1,5}^{}$-adapted and $\Aut_X^0(R)$ acts as $[\lambda:\mu] \mapsto [\lambda : \mu + 2 c \lambda]$.
\end{itemize}
As in the previous case, there is a unique choice for $p_{1,6}$ up to isomorphism:

\begin{enumerate}[leftmargin=*]

    \item \vspace{2mm}
    $p_{1,6}  =E_{1,5}^{} \cap C^{(6)}$ with $C= \cal{V}(x^4z+x^2y^3+y^5)$
    \begin{itemize}[leftmargin=20pt]
    \noindent \begin{minipage}{0.5\textwidth}
        \item
        $\Aut_{X'}^0(R) =
        \left\{  \left( \begin{smallmatrix}
1 &  &  \\
 & 1 & f \\
 &  & 1
\end{smallmatrix} \right) 
\in \PGL_3(R)  \right\} $ 
        \item $(-2)$-curves: $E_{1,0}^{(7)}, E_{1,1}^{(7)},E_{1,2}^{(7)}, E_{1,3}^{(7)},E_{1,4}^{(7)},
        E_{1,5}^{(7)}, \ell_z^{(7)}$
        
        \end{minipage} \noindent \begin{minipage}{0.5\textwidth}

        \item $(-1)$-curves: $E_{1,6}^{}$
        \item
        with configuration as in Figure \ref{Conf2I2L2M2W2X2Y}.
    \end{minipage}
    \end{itemize} 
    \vspace{1mm} \noindent This is case \hyperref[Tab2L]{$2L$}.

\end{enumerate}

\vspace{2mm} \subsubsection*{\underline{Case \hyperref[Tab3R]{$3R$}}}
This case exists only if $p = 2$.  
\\ \noindent We have $E =
E_{1,5}^{} - E_{1,4}^{(6)}$ and $\Aut_X^0(R) = \left\{  \left( \begin{smallmatrix}
1 & b & c \\
 & e& b^2e\\
 &  & e^3
\end{smallmatrix} \right) 
\in \PGL_3(R) \right\}$.

\begin{itemize}[leftmargin=25pt]
    \item[-] $\lambda x^3(x^2z + y^3) + \mu y^6$ is $E_{1,5}^{}$-adapted and $\Aut_X^0(R)$ acts as $[\lambda:\mu] \mapsto [\lambda : e^3 \mu]$.
\end{itemize}
Since $\Aut_X^0$ has two orbits on $E \cap E_{1,5}$, we have the following two choices for $p_{1,6}$ up to isomorphism:

\begin{enumerate}[leftmargin=*]
 
   \item \vspace{2mm}
    $p_{1,6}  =E_{1,5}^{} \cap C^{(6)}$ with $C= \cal{V}(x^5z+x^3y^3+y^6)$
    \begin{itemize}[leftmargin=20pt]
    \noindent \begin{minipage}{0.5\textwidth}
        \item
        $\Aut_{X'}^0(R) =
        \left\{  \left( \begin{smallmatrix}
1 & b & c \\
 & 1 & b^2 \\
 &  & 1
\end{smallmatrix} \right) 
\in \PGL_3(R)  \right\}$ 
        \item $(-2)$-curves: $E_{1,0}^{(7)}, E_{1,1}^{(7)},E_{1,2}^{(7)}, E_{1,3}^{(7)},E_{1,4}^{(7)}, 
        E_{1,5}^{(7)}, \ell_z^{(7)}$
        
        \end{minipage} \noindent \begin{minipage}{0.5\textwidth}
        
        \item $(-1)$-curves: $E_{1,6}^{}$
        \item
        with configuration as in Figure \ref{Conf2I2L2M2W2X2Y}.
    \end{minipage}
    \end{itemize} 
    \vspace{1mm} \noindent This is case \hyperref[Tab2X]{$2X$}.

    \item \vspace{2mm}
    $p_{1,6}  =E_{1,5}^{} \cap C^{(6)}$ with $C= \cal{V}(x^2z+y^3)$
    \begin{itemize}[leftmargin=20pt]
    \noindent \begin{minipage}{0.5\textwidth}
        \item
        $\Aut_{X'}^0(R) =
        \left\{  \left( \begin{smallmatrix}
1 & b & c \\
 & e & b^2e \\
 &  & e^3
\end{smallmatrix} \right) 
\in \PGL_3(R)  \right\}$ 
        \item $(-2)$-curves: $E_{1,0}^{(7)}, E_{1,1}^{(7)},E_{1,2}^{(7)}, E_{1,3}^{(7)},E_{1,4}^{(7)}, 
        E_{1,5}^{(7)}, \ell_z^{(7)}$
        
        \end{minipage} \noindent \begin{minipage}{0.5\textwidth}
        
        \item $(-1)$-curves: $E_{1,6}^{}$
        \item
        with configuration as in Figure \ref{Conf2I2L2M2W2X2Y}.
        \end{minipage}
        \end{itemize}
    \vspace{1mm} \noindent This is case \hyperref[Tab2Y]{$2Y$}. 

\end{enumerate}

\vspace{2mm} \subsubsection*{\underline{Case \hyperref[Tab3Q]{$3Q$}}}
This case exists only if $p = 2$.  
\\ \noindent We have $E =
E_{1,5}^{} - E_{1,4}^{(6)}$ and $\Aut_X^0(R) = \left\{  \left( \begin{smallmatrix}
1 & b & c \\
 & 1& b^2 + b\\
 &  & 1
\end{smallmatrix} \right) 
\in \PGL_3(R) \right\}$.

\begin{itemize}[leftmargin=25pt]
    \item[-] $\lambda x^2(x^3z + xy^3 + y^4) + \mu y^6$ is $E_{1,5}^{}$-adapted and $\Aut_X^0(R)$ acts as $[\lambda:\mu] \mapsto [\lambda : \mu + (b^2 + b)\lambda]$.
\end{itemize}
Since $\Aut_X^0$ acts transitively on $E \cap E_{1,5}$, we have the following unique choice for $p_{1,6}$ up to isomorphism:

\begin{enumerate}[leftmargin=*]

    \item \vspace{2mm}
    $p_{1,6}  =E_{1,5}^{} \cap C^{(6)}$ with $C= \cal{V}(x^3z + xy^3 + y^4)$ 
    \begin{itemize}[leftmargin=20pt]
    \noindent \begin{minipage}{0.5\textwidth}
        \item
        $
        \Aut_{X'}^0(R) =    
        \left\{  \left( \begin{smallmatrix}
1 &  & c \\
 & 1 &  \\
 &  & 1
\end{smallmatrix} \right) 
\in \PGL_3(R) \right\}
        $ 
        \item $(-2)$-curves: $E_{1,0}^{(7)}, E_{1,1}^{(7)},E_{1,2}^{(7)}, E_{1,3}^{(7)},E_{1,4}^{(7)},
        E_{1,5}^{(7)}, \ell_z^{(7)}$
        
        \end{minipage} \noindent \begin{minipage}{0.5\textwidth}
        
        \item $(-1)$-curves: $E_{1,6}^{}$
        \item
        with configuration as in Figure \ref{Conf2I2L2M2W2X2Y}.
    \end{minipage}
    \end{itemize} 
    \vspace{1mm} \noindent This is case \hyperref[Tab2W]{$2W$}. 

\end{enumerate}

\vspace{2mm} 
\noindent
Summarizing, we obtain
\begin{eqnarray*}
\cal{L}_7 &=& \{ X_{2I}, X_{2M}, X_{2L}, X_{2X}, X_{2Y}, X_{2W} \}.
\end{eqnarray*}

\subsection{Height 8}

\vspace{2mm} \subsubsection*{\underline{Case \hyperref[Tab2I]{$2I$}}}
This case exists only if $p \neq 2,3$.  
\\ \noindent We have $E =
E_{1,6}^{} - E_{1,5}^{(7)}$ and $\Aut_X^0(R) = \left\{  \left( \begin{smallmatrix}
1 &  &  \\
 & e&  \\
 &  & e^3
\end{smallmatrix} \right) 
\in \PGL_3(R) \right\}$.

\begin{itemize}[leftmargin=25pt]
    \item[-] $\lambda x^4(x^2z + y^3) + \mu y^7$ is $E_{1,6}^{}$-adapted and $\Aut_X^0(R)$ acts as $[\lambda:\mu] \mapsto [\lambda : e^4 \mu$].
\end{itemize}
Since $p \neq 2$, there is a unique point on $E \cap E_{1,6}$ whose stabilizer has non-trivial identity component. This leads to the following unique choice for $p_{1,7}$ up to isomorphism:

\begin{enumerate}[leftmargin=*]
\noindent     \begin{minipage}{0.65\textwidth}
    \item \vspace{2mm}
    $p_{1,7}  =E_{1,6}^{} \cap C^{(7)}$ with $C= \cal{V}(x^2z+y^3)$
    \begin{itemize}[leftmargin=20pt]
        \item
        $\Aut_{X'}^0(R) =
        \left\{  \left( \begin{smallmatrix}
1 &  &  \\
 & e &  \\
 &  & e^3
\end{smallmatrix} \right) 
\in \PGL_3(R)  \right\}$ 
        \item $(-2)$-curves: $E_{1,0}^{(8)}, E_{1,1}^{(8)},E_{1,2}^{(8)}, E_{1,3}^{(8)}, E_{1,4}^{(8)},E_{1,5}^{(8)}, E_{1,6}^{(8)}, \ell_z^{(8)}$
        \item $(-1)$-curves: $E_{1,7}^{}$
        \item
        with configuration as in Figure \ref{Conf1D1H1I1S1T}.
    \end{itemize} 
    This is case \hyperref[Tab1D]{$1D$}.
        \end{minipage} \hspace{2mm} \begin{minipage}{0.3\textwidth} \begin{center} \includegraphics[width=0.98\textwidth]{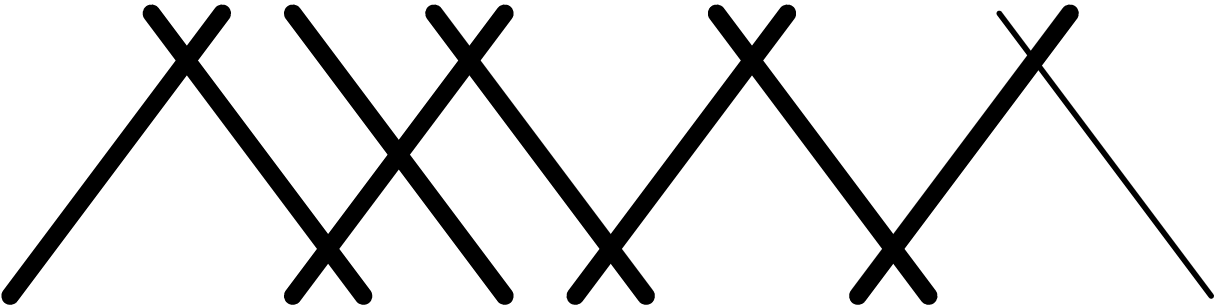} \vspace{-3.5mm}\captionof{figure}{}  \label{Conf1D1H1I1S1T} \end{center} \end{minipage} 
\end{enumerate}

\vspace{2mm} \subsubsection*{\underline{Case \hyperref[Tab2M]{$2M$}}}
This case exists only if $p = 3$.  
\\ \noindent We have $E =
E_{1,6}^{} - E_{1,5}^{(7)}$ and $\Aut_X^0(R) = \left\{  \left( \begin{smallmatrix}
1 &  &  \\
 & e&  f \\
 &  & e^3
\end{smallmatrix} \right) 
\in \PGL_3(R) \right\}$.

\begin{itemize}[leftmargin=25pt]
    \item[-] $\lambda x^4(x^2z + y^3) + \mu y^7$ is $E_{1,6}^{}$-adapted and $\Aut_X^0(R)$ acts as $[\lambda:\mu] \mapsto [\lambda : e^4 \mu$].
\end{itemize}
Since $\Aut_X^0$ acts with two orbits on $E \cap E_{1,6}$, we have the following two choices for $p_{1,7}$ up to isomorphism:

\begin{enumerate}[leftmargin=*]

    \item \vspace{2mm}
    $p_{1,7}  =E_{1,6}^{} \cap C^{(7)}$ with $C= \cal{V}(x^6z+x^4y^3+y^7)$
    \begin{itemize}[leftmargin=20pt]
    \noindent \begin{minipage}{0.5\textwidth}
        \item
        $\Aut_{X'}^0(R) =
        \left\{  \left( \begin{smallmatrix}
1 &  &  \\
 & 1 & f \\
 &  & 1
\end{smallmatrix} \right) 
\in \PGL_3(R)  \right\} $  
        \item $(-2)$-curves: $E_{1,0}^{(8)}, E_{1,1}^{(8)},E_{1,2}^{(8)}, E_{1,3}^{(8)}, E_{1,4}^{(8)},E_{1,5}^{(8)}, 
        \\E_{1,6}^{(8)}, \ell_z^{(8)}$
        
        \end{minipage} \noindent \begin{minipage}{0.5\textwidth}
        
        \item $(-1)$-curves: $E_{1,7}^{}$
        \item
        with configuration as in Figure \ref{Conf1D1H1I1S1T}.
    \end{minipage}
    \end{itemize} 
    \vspace{1mm} \noindent This is case \hyperref[Tab1H]{$1H$}.

    \item \vspace{2mm}
    $p_{1,7}  =E_{1,6}^{} \cap C^{(7)}$ with $C= \cal{V}(x^2z+y^3)$
    \begin{itemize}[leftmargin=20pt]
    \noindent \begin{minipage}{0.5\textwidth}
        \item
        $\Aut_{X'}^0(R) =
        \left\{  \left( \begin{smallmatrix}
1 &  &  \\
 & e & f \\
 &  & e^3
\end{smallmatrix} \right) 
\in \PGL_3(R)  \right\} $  
        \item $(-2)$-curves: $E_{1,0}^{(8)}, E_{1,1}^{(8)},E_{1,2}^{(8)}, E_{1,3}^{(8)}, E_{1,4}^{(8)},E_{1,5}^{(8)}, 
        \\E_{1,6}^{(8)}, \ell_z^{(8)}$
        
        \end{minipage} \noindent \begin{minipage}{0.5\textwidth}
        
        \item $(-1)$-curves: $E_{1,7}^{}$
        \item
        with configuration as in Figure \ref{Conf1D1H1I1S1T}.
    \end{minipage}
    \end{itemize} 
    \vspace{1mm} \noindent This is case \hyperref[Tab1I]{$1I$}.
        
\end{enumerate}

\vspace{2mm} \subsubsection*{\underline{Case \hyperref[Tab2L]{$2L$}}}
This case exists only if $p = 3$.  
\\ \noindent We have $E =
E_{1,6}^{} - E_{1,5}^{(7)}$ and $\Aut_X^0(R) = \left\{  \left( \begin{smallmatrix}
1 &  &  \\
 & 1& f \\
 &  & 1
\end{smallmatrix} \right) 
\in \PGL_3(R) \right\}$.

\begin{itemize}[leftmargin=25pt]
    \item[-] $\lambda x^2(x^4z + x^2y^3 + y^5) + \mu y^7$ is $E_{1,6}^{}$-adapted and $\Aut_X^0(R)$ acts as $[\lambda:\mu] \mapsto [\lambda : \mu +  f\lambda$].
\end{itemize}
Hence, the stabilizer of every point on $E \cap E_{1,6}^{}$ is trivial, therefore we cannot blow up $X$ further and still obtain a weak del Pezzo surface with global vector fields.

\vspace{2mm} \subsubsection*{\underline{Case \hyperref[Tab2X]{$2X$}}}
This case exists only if $p = 2$.  
\\ \noindent We have $E =
E_{1,6}^{} - E_{1,5}^{(7)}$ and $\Aut_X^0(R) = \left\{  \left( \begin{smallmatrix}
1 & b & c \\
 & 1& b^2 \\
 &  & 1
\end{smallmatrix} \right) 
\in \PGL_3(R) \right\}$.

\begin{itemize}[leftmargin=25pt]
    \item[-] $\lambda x(x^5z + x^3y^3 + y^6) + \mu y^7$ is $E_{1,6}^{}$-adapted and $\Aut_X^0(R)$ acts as $[\lambda:\mu] \mapsto [\lambda : \mu + (b + b^4)\lambda$].
\end{itemize}
Since $\Aut_X^0$ acts transitively on $E \cap E_{1,6}$, there is a unique choice for $p_{1,7}$ up to isomorphism:

\begin{enumerate}[leftmargin=*]

   \item \vspace{2mm}
    $p_{1,7}  =E_{1,6}^{} \cap C^{(7)}$ with $C= \cal{V}(x^5z+x^3y^3+y^6)$
    \begin{itemize}[leftmargin=20pt]
    \noindent \begin{minipage}{0.5\textwidth}
        \item
        $\Aut_{X'}^0(R) =
        \left\{  \left( \begin{smallmatrix}
1 &  & c \\
 & 1 &  \\
 &  & 1
\end{smallmatrix} \right) 
\in \PGL_3(R)  \right\}$ 
        \item $(-2)$-curves: $E_{1,0}^{(8)}, E_{1,1}^{(8)},E_{1,2}^{(8)}, E_{1,3}^{(8)},E_{1,4}^{(8)}, E_{1,5}^{(8)},
        \\E_{1,6}^{(8)}, \ell_z^{(8)}$
        
        \end{minipage} \noindent \begin{minipage}{0.5\textwidth}
        
        \item $(-1)$-curves: $E_{1,7}^{}$
        \item
        with configuration as in Figure \ref{Conf1D1H1I1S1T}.
    \end{minipage}
    \end{itemize} 
    \vspace{1mm} \noindent This is case \hyperref[Tab1S]{$1S$}. 
        
\end{enumerate}

\vspace{2mm} \subsubsection*{\underline{Case \hyperref[Tab2Y]{$2Y$}}}
This case exists only if $p = 2$.  
\\ \noindent We have $E =
E_{1,6}^{} - E_{1,5}^{(7)}$ and $\Aut_X^0(R) = \left\{  \left( \begin{smallmatrix}
1 & b & c \\
 & e& b^2 e\\
 &  & e^3
\end{smallmatrix} \right) 
\in \PGL_3(R) \right\}$.

\begin{itemize}[leftmargin=25pt]
    \item[-] $\lambda x^4(x^2z + y^3) + \mu y^7$ is $E_{1,6}^{}$-adapted and $\Aut_X^0(R)$ acts as $[\lambda:\mu] \mapsto [\lambda : e^4\mu +  b^4\lambda$].
\end{itemize}
Since $\Aut_X^0$ acts transitively on $E \cap E_{1,6}$, there is a unique choice for $p_{1,7}$ up to isomorphism:

\begin{enumerate}[leftmargin=*]

    \item \vspace{2mm}
    $p_{1,7}  =E_{1,6}^{} \cap C^{(7)}$ with $C= \cal{V}(x^2z+y^3)$
    \begin{itemize}[leftmargin=20pt]
    \noindent \begin{minipage}{0.5\textwidth}
        \item
        $\Aut_{X'}^0(R) =
        \left\{  \left( \begin{smallmatrix}
1 & b & c \\
 & e & b^2e \\
 &  & e^3
\end{smallmatrix} \right) 
\in \PGL_3(R) \bigg| b^4=0 \right\}$ 
        \item $(-2)$-curves: $E_{1,0}^{(8)}, E_{1,1}^{(8)},E_{1,2}^{(8)}, E_{1,3}^{(8)},E_{1,4}^{(8)}, E_{1,5}^{(8)},
        \\E_{1,6}^{(8)}, \ell_z^{(8)}$
        
        \end{minipage} \noindent \begin{minipage}{0.5\textwidth}
        
        \item $(-1)$-curves: $E_{1,7}^{}$
        \item
        with configuration as in Figure \ref{Conf1D1H1I1S1T}.
    \end{minipage}
    \end{itemize} 
    \vspace{1mm} \noindent This is case \hyperref[Tab1T]{$1T$}.
        
\end{enumerate}

\vspace{2mm} \subsubsection*{\underline{Case \hyperref[Tab2W]{$2W$}}}
This case exists only if $p = 2$.  
\\ \noindent We have $E =
E_{1,6}^{} - E_{1,5}^{(7)}$ and $\Aut_X^0(R) = \left\{  \left( \begin{smallmatrix}
1 &  & c \\
 & 1& \\
 &  & 1
\end{smallmatrix} \right) 
\in \PGL_3(R) \right\}$.

\begin{itemize}[leftmargin=25pt]
    \item[-] $\lambda x^3(x^3z + xy^3 + y^4) + \mu y^7$ is $E_{1,6}^{}$-adapted and $\Aut_X^0(R)$ acts as $[\lambda:\mu] \mapsto [\lambda : \mu +  c\lambda$].
\end{itemize}
In particular, the identity component of the stabilizer of every point on $E \cap E_{1,6}^{}$ is trivial, hence we cannot blow up further and still obtain a weak del Pezzo surface with global vector fields.

\vspace{2mm} 
\noindent
Summarizing, we obtain
\begin{eqnarray*}
\cal{L}_8 &=& \{ X_{1D}, X_{1H}, X_{1I}, X_{1S}, X_{1T} \}.
\end{eqnarray*}

\qed

\vspace{1.5cm}

\bibliographystyle{alpha} 
\bibliography{delPezzo}

\end{document}